\documentclass[11pt]{book}
\usepackage{amsthm,amsfonts,amsmath,amscd,amssymb,latexsym,epsfig}
\usepackage{mathrsfs}
\usepackage{mathdots}
\usepackage[all]{xy}
\usepackage{hyperref}
\title{\Huge The Moment-Weight Inequality \\ and the Hilbert--Mumford Criterion}

\author{
Valentina~Georgoulas\\ETH Z\"urich
\and Joel~W.~Robbin\\ UW Madison
\and Dietmar~A.~Salamon\\ETH Z\"urich
}

\date{30 May 2021}

\numberwithin{equation}{chapter} 
\newtheorem{PARA}{}[chapter]
\newtheorem{thm}[PARA]{Theorem}
\newtheorem{cor}[PARA]{Corollary}
\newtheorem{lem}[PARA]{Lemma}

\newtheorem{defn}[PARA]{Definition}
\newtheorem{remk}[PARA]{Remark}
\newtheorem{exple}[PARA]{Example}
\newcommand{\cA}{\mathcal{A}}

\newcommand{\cC}{\mathcal{C}}

\newcommand{\cF}{\mathcal{F}}

\newcommand{\cG}{\mathcal{G}}
\newcommand{\cH}{\mathcal{H}}

\newcommand{\cJ}{\mathcal{J}}

\newcommand{\cL}{\mathcal{L}}
\newcommand{\cM}{\mathcal{M}}

\newcommand{\cO}{\mathcal{O}}
\newcommand{\cP}{\mathcal{P}}

\newcommand{\tg}{{\widetilde{g}}}
\newcommand{\tu}{{\widetilde{u}}}
\newcommand{\tx}{{\widetilde{x}}}

\newcommand{\tga}{{\widetilde{\gamma}}}
\newcommand{\tmu}{{\widetilde{\mu}}}
\newcommand{\txi}{{\widetilde{\xi}}}
\newcommand{\teta}{{\widetilde{\eta}}}
\newcommand{\tzeta}{{\widetilde{\zeta}}}

\newcommand{\tom}{{\widetilde{\om}}}

\newcommand{\tPhi}{{\widetilde{\Phi}}}
\newcommand{\tPsi}{{\widetilde{\Psi}}}

\newcommand{\trG}{{\widetilde{\rG}}}
\newcommand{\trT}{{\widetilde{\rT}}}
\newcommand{\tcg}{{\widetilde{\cg}}}
\newcommand{\tct}{{\widetilde{\ct}}}

\newcommand{\heta}{{\widehat{\eta}}}

\newcommand{\one}{{{\mathchoice \mathrm{ 1\mskip-4mu l} \mathrm{ 1\mskip-4mu l}
\mathrm{ 1\mskip-4.5mu l} \mathrm{ 1\mskip-5mu l}}}}
\newcommand{\dslash}{/\mskip-6mu/}

\newcommand{\C}{{\mathbb{C}}}

\newcommand{\bbD}{{\mathbb{D}}}

\newcommand{\bbH}{{\mathbb{H}}}

\newcommand{\N}{{\mathbb{N}}}
\newcommand{\bbP}{{\mathbb{P}}}
\newcommand{\Q}{{\mathbb{Q}}}
\newcommand{\R}{{\mathbb{R}}}
\newcommand{\bbS}{{\mathbb{S}}}

\newcommand{\Z}{{\mathbb{Z}}}

\newcommand{\sL}{{\mathscr{L}}}

\newcommand{\sT}{{\mathscr{T}}}

\newcommand{\sZ}{{\mathscr{Z}}}

\renewcommand{\ss}{{\mathsf{ss}}}
\newcommand{\ps}{{\mathsf{ps}}}
\newcommand{\us}{{\mathsf{us}}}
\newcommand{\s}{{\mathsf{s}}}

\newcommand{\im}{\mathrm{im}}

\newcommand{\Crit}{\mathrm{Crit}}

\newcommand{\trace}{\mathrm{trace}}

\newcommand{\id}{\mathrm{id}} 

\newcommand{\dvol}{\mathrm{dvol}}    
\newcommand{\Vol}{\mathrm{Vol}}

\renewcommand{\Re}{\mathrm{Re}} 
\renewcommand{\Im}{\mathrm{Im}} 
\newcommand{\ad}{\mathrm{ad}}    

\newcommand{\rB}{\mathrm{B}}
\newcommand{\rC}{\mathrm{C}}
\newcommand{\rG}{\mathrm{G}}
\newcommand{\rH}{\mathrm{H}}
\newcommand{\rK}{\mathrm{K}}

\newcommand{\SO}{\mathrm{SO}}
\newcommand{\rP}{\mathrm{P}}
\newcommand{\rT}{\mathrm{T}}
\newcommand{\rU}{\mathrm{U}}
\newcommand{\rZ}{\mathrm{Z}}
\newcommand{\SU}{\mathrm{SU}}
\newcommand{\GL}{\mathrm{GL}}
\newcommand{\SL}{\mathrm{SL}}
\newcommand{\PSL}{\mathrm{PSL}}
\newcommand{\Lie}{\mathrm{Lie}}   
\newcommand{\Aut}{\mathrm{Aut}}   
\newcommand{\Der}{\mathrm{Der}} 
\newcommand{\Diff}{\mathrm{Diff}} 
\newcommand{\Vect}{\mathrm{Vect}} 
   
\newcommand{\Hom}{\mathrm{Hom}}

\newcommand{\cg}{\mathfrak{g}}
\newcommand{\ch}{\mathfrak{h}}

\newcommand{\cso}{\mathfrak{so}}

\newcommand{\cp}{\mathfrak{p}}
\newcommand{\ct}{\mathfrak{t}}
\newcommand{\cu}{\mathfrak{u}}
\newcommand{\csu}{\mathfrak{su}}
\newcommand{\cgl}{\mathfrak{gl}}

\newcommand{\eps}{{\varepsilon}}
\newcommand{\om}{{\omega}}
\newcommand{\Om}{{\Omega}}
\newcommand{\Ahat}{{\widehat{A}}}

\newcommand{\Phat}{{\widehat{P}}}

\newcommand{\ghat}{{\widehat{g}}}
\newcommand{\phat}{{\widehat{p}}}
\newcommand{\uhat}{{\widehat{u}}}
\newcommand{\vhat}{{\widehat{v}}}
\newcommand{\xhat}{{\widehat{x}}}
\newcommand{\yhat}{{\widehat{y}}}

\newcommand{\etahat}{{\widehat{\eta}}}
\newcommand{\oa}{{\overline{a}}}
\newcommand{\ob}{{\overline{b}}}
\newcommand{\ox}{{\overline{x}}}
\newcommand{\oz}{{\overline{z}}}
\newcommand{\ola}{{\overline{\lambda}}}
\newcommand{\ocM}{{\overline{\cM}}}
\newcommand{\bi}{{\mathbf{i}}}

\renewcommand{\i}{{\mathbf{i}}}

\newcommand{\inner}[2]{\langle #1, #2\rangle}
\newcommand{\Inner}[2]{\left\langle #1, #2\right\rangle}

\def\NABLA#1{{\mathop{\nabla\kern-.5ex\lower1ex\hbox{$#1$}}}}
\def\Nabla#1{\nabla\kern-.5ex{}_{#1}}
\def\Tabla#1{\Tilde\nabla\kern-.5ex{}_{#1}}
\def\abs#1{\mathopen|#1\mathclose|}
\def\Abs#1{\left|#1\right|}

\renewcommand{\Tilde}{\widetilde}

\newcommand{\p}{{\partial}}

\makeindex             

\begin{document}

\frontmatter

\maketitle

\chapter*{Preface}

This book gives an essentially self contained exposition 
(except for an appeal to the Lojasiewicz gradient inequality)
of geometric invariant theory from a differential
geometric viewpoint. Central ingredients are the 
moment-weight inequality (relating the Mumford 
numerical invariants to the norm of the moment map), 
the negative gradient flow of the moment map squared,
and the Kempf--Ness function.

The last author DAS owes a lot to lectures by and 
conversations with Simon Donaldson, with
G\'abor Sz\'ekelyhidi, and with Xiuxiong Chen and Song Sun.  
The second author JWR learned much from a course 
given by and conversations with Sean Paul
at the University of Wisconsin. 
Most of this book was written when JWR visited the 
Forschungsinstitut f\"ur Mathematik at ETH Z\"urich
and he thanks them for their hospitality.
The first version of this book was completed while DAS 
visited the IAS, Princeton, and the SCGP, Stony Brook; 
he thanks both institutes for their hospitality.
Thanks to Samuel Trautwein for helpful discussions.
Thanks to Amanda Jenny for pointing out errors
in earlier versions of the manuscript.
 
\vspace{\baselineskip}
\begin{flushright}\noindent
\hfill {\it Valentina Georgoulas}\\
Z\"urich,\hfill {\it Joel W.~Robbin}\\
November 2019\hfill {\it Dietmar A.~Salamon}\\
\end{flushright}

\tableofcontents


\mainmatter



\chapter{Introduction}\label{ch:INTRO}

Many important problems in geometry can be reduced to a partial
differential equation of the form 
$$
\mu(x)=0, 
$$
where $x$ ranges over a complexifed group orbit in an 
infinite-dimensional symplectic manifold $X$
and~${\mu:X\to\cg}$ is an associated moment map 
(see Calabi~\cite{CALABI0,CALABI1,CALABI}, Yau~\cite{YAU0,YAU1,YAU2}, 
Tian~\cite{TIAN}, Chen--Donaldson--Sun~\cite{CDS1,CDS2,CDS3},
Atiyah--Bott~\cite{AB}, Uhlenbeck--Yau~\cite{UY}, 
Donaldson~\cite{DON-NS,DON-asd,DON-momdiff,DON-K1,DON-K2}).  
Problems like this are extremely difficult.  The purpose of this book
is to explain the analogous finite-dimensional situation, which is the 
subject of Geometric Invariant Theory. 

GIT was originally developed to study actions of a complex reductive
Lie group $\rG^c$ on a projective algebraic variety~${X\subset\bbP(V)}$.
Here $\rG^c$ is the complexification of a compact Lie group $\rG$
and the Hermitian structure on~$V$ can be chosen so that~$\rG$ 
acts by unitary automorphisms. In the smooth case~$X$ inherits
the structure of a K\"ahler manifold from the standard K\"ahler 
structure on~$\bbP(V)$.  The $\rG$-action is generated by
the standard moment map 
$$
\mu:X\to\cg
$$ 
(with values in the Lie algebra of $\rG$).  
In the original treatment of Mumford~\cite{MUMFORD} 
the symplectic form and the moment map were not used. 
Subsequently several authors discovered the connection 
between the theory of the moment map and GIT 
(see Kirwan~\cite{KIRWAN} and Ness~\cite{NESS2}).
In the article of Lerman~\cite{LERMAN} it is noted that
both Guilleman--Sternberg~\cite{GS} and Ness~\cite{NESS2}
credit Mumford for the relation between the complex quotient 
and the Marsden--Weinstein quotient 
$$
X\dslash\rG := \mu^{-1}(0)/\rG.
$$

In our exposition we assume that~$X$ is a closed K\"ahler manifold
but do not assume that it is a projective variety.
In the latter case the aforementioned standard 
moment map satisfies certain rationality conditions 
(Chapter~\ref{ch:RAT}) which we do not use in our treatment.
As a result the Mumford numerical invariants 
$$
w_\mu(x,\xi) 
:= \lim_{t\to\infty}\inner{\mu(\exp(\i t\xi)x)}{\xi}
$$
associated to a point~${x\in X}$ and an
element
$
\xi\in\cg\setminus\{0\}
$ 
(Chapter~\ref{ch:WEIGHTS}) need not be integers as they are 
in traditional GIT. In the classical theory the Lie algebra 
element belongs to the set 
$$
\Lambda:=\left\{\xi\in\cg\setminus\{0\}\,|\,\exp(\xi)=\one\right\}
$$
and thus generates a one-parameter subgroup of $\rG^c$.
In our treatment $\xi$ can be a general nonzero element of $\cg$.

A central ingredient in our treatment is the
{\bf moment-weight inequality}
\begin{equation}\label{eq:MW}
\sup_{\xi\in\cg\setminus\{0\}}
\frac{-w_\mu(x,\xi)}{\abs{\xi}}
\le \inf_{g\in\rG^c}\abs{\mu(gx)}.
\end{equation}
We give two proofs of this inequality in Chapter~\ref{ch:MW1},
one due to Mumford~\cite{MUMFORD} and Ness~\cite[Lemma~3.1]{NESS2} 
and one due to Xiuxiong Chen~\cite{CIII}.  
(For an in depth discussion of this 
inequality see Atiyah--Bott~\cite{AB}, in the setting of bundles over 
Riemann surfaces, and Donaldson~\cite{DON-K4}, Sz\'ekelyhidi~\cite{S}, 
Chen~\cite{CIII,CIV}, in the setting of K\"ahler--Einstein geometry.)  
Following an argument of Chen--Sun~\cite{CS} we also prove 
that equality holds in~\eqref{eq:MW} whenever the right hand 
side is positive (Theorem~\ref{thm:KEMPF}).
We also prove that the supremum on the left
is always attained (Theorem~\ref{thm:SUPMAX})
and that the supremum over all $\xi\in\cg\setminus\{0\}$ agrees 
with the supremum over all $\xi\in\Lambda$ (Theorem~\ref{thm:M}). 
In the projective case the supremum is attained at an 
element~${\xi\in\Lambda}$ by a theorem of Kempf~\cite{KEMPF}, 
however, that need not be the case in our more general setting.
The {\bf Hilbert--Mumford numerical criterion} 
for $\mu$-semistability is an immediate consequence 
of the aforementioned results (Theorem~\ref{thm:HMss}). 
It asserts that
$$
\overline{\rG^c(x)}\cap\mu^{-1}(0)\ne\emptyset
\qquad\iff\qquad
w_\mu(x,\xi)\ge 0\;\forall\,\xi\in\Lambda.
$$
Further consequences of the moment-weight inequality 
include the {\bf Kirwan--Ness Inequality} which asserts that
if $x$ is a critical point of the moment map squared 
then~${\abs{\mu(x)}=\inf_{g\in\rG^c}\abs{\mu(gx)}}$ 
(Corollary~\ref{cor:KIRNESS1}), 
the {\bf Moment Limit Theorem} which asserts that each 
negative gradient flow line of the moment map squared 
converges to a minimum of the moment map squared on the 
complexified group orbit (Theorem~\ref{thm:MLT}),
and the {\bf Ness Uniqueness Theorem} which asserts that
any two critical points of the moment map squared 
in the same $\rG^c$-orbit in fact belong to the same 
$\rG$-orbit (Theorem~\ref{thm:NESS1}) and, moreover, 
that the minimum of the moment map squared on the 
closure of a $\rG^c$-orbit is taken on at a 
unique $\rG$-orbit (Theorem~\ref{thm:NESS2}).

A central ingredient in the proofs of these theorems
is the negative gradient flow of the moment map 
squared~${f:=\tfrac{1}{2}\abs{\mu}^2:X\to\R}$.
The gradient flow equation takes the form
\begin{equation}\label{eq:GRADFLOW}
\dot x = -JL_x\mu(x)
\end{equation}
where~${L_x:\cg\to T_xX}$ denotes the infinitesimal 
action of the Lie algebra on~$X$. 
Each $\rG^c$-orbit~${\rG^c(x)\subset X}$ 
is invariant under this flow, because 
every solution of~\eqref{eq:GRADFLOW}
has the form~${x(t)=g(t)^{-1}x}$, 
where~${g:\R\to\rG^c}$ satisfies 
the differential equation
\begin{equation}\label{eq:KNFLOW}
g^{-1}\dot g = \bi\mu(g^{-1}x).
\end{equation}
Equation~\eqref{eq:KNFLOW} is the negative gradient flow 
of the {\bf Kempf--Ness function}~${\Phi_x:\rG^c/\rG\to\R}$.
The homogeneous space~$\rG^c/\rG$ is simply connected 
and complete with nonpositive sectional curvature, 
and the Kempf--Ness function is Morse--Bott and 
is convex along geodesics (Theorem~\ref{thm:KNF});
its critical manifold may be empty.
Moreover, the {\bf Kempf--Ness Theorem}  
characterizes the stability conditions in terms 
of the properties of the Kempf--Ness function
(Theorem~\ref{thm:KN}); for example a point~${x\in X}$ 
is $\mu$-semistable, i.e.\ the closure of 
its~$\rG^c$-orbit intersects the zero set of the 
moment map, if and only if the Kempf--Ness 
function~$\Phi_x$ is bounded below.

The moment map squared is in general 
far from being Morse--Bott and may
have very complicated critical points. 
However, the aforementioned theorems 
(Kirwan--Ness Inequality, Ness Uniqueness, 
Moment Limit Theorem) exhibit a structure 
of the gradient flow that resembles
the stratification by stable manifolds 
associated to a Morse--Bott function.
More precisely, an element~${x\in X}$
is a critical point of the moment map squared
if it satisfies the equation~${L_x\mu(x)=0}$.
Critical points come in~$\rG$-orbits
and the theorems of Ness and Kirwan--Ness
show that the stable manifold 
of such a critical orbit~$\rG(x)$ 
is a union of~$\rG^c$-orbits, i.e.\ 
\begin{equation}\label{eq:Wsx}
W^s(\rG(x)) = \left\{y\in X\,\bigg|\,x\in\overline{\rG^c(y)},\,
\abs{\mu(x)} = \inf_{\rG^c(y)}\abs{\mu}\right\}.
\end{equation}
This stratification was used by Kirwan~\cite{KIRWAN} to prove 
that the canonical ring homomorphism~${\kappa:H^*_\rG(X)\to H^*(X\dslash\rG)}$
(called the {\bf Kirwan homomorphism})\index{Kirwan homomorphism}
from the equivariant cohomology of~$X$ 
to the cohomology of the Marsden--Weinstein quotient~$X\dslash\rG$
is surjective. Kirwan's theorem is not included in our treatment.
The main motivation for the present book comes 
from infinite-dimensional analogues of GIT
in various areas of geometry. 

One such infinite-dimensional analogue of geometric invariant theory 
is the Donaldson--Uhlenbeck--Yau correspondence between stable holomorphic
vector bundles and Hermitian Yang--Mills connections over K\"ahler manifolds.
(This is a special case of the Kobayashi--Hitchin correspondence.)
In this theory the space of Hermitian connections on a Hermitian 
vector bundle over a closed K\"ahler manifold with curvature of type $(1,1)$ 
is viewed as an infinite-dimensional symplectic manifold, 
the group of unitary gauge transformatons acts on it 
by Hamiltonian symplectomorphisms, 
the moment map assigns to a connection the component of 
the curvature parallel to the symplectic form,  and the zero set 
of the moment map is the space of Hermitian Yang--Mills connections.
For bundles over Riemann surfaces the analogue 
of the Hilbert--Mumford numerical criterion is the 
correspondence between stable bundles and flat connections,
established by Narasimhan--Seshadri~\cite{NS}.
It can be viewed as an extension to higher rank bundles 
of the identification of the Jacobian with a torus.
The moment map picture in this setting was exhibited 
by Atiyah--Bott~\cite{AB} and they proved 
the analogue of the moment-weight inequality. 
Another proof of the Narasimhan--Seshadri theorem 
was given by Donaldson~\cite{DON-NS}.  In dimension four the Hermitian
Yang--Mills conections are anti-self-dual instantons over K\"ahler surfaces.
In this setting the correspondence between stable bundles and 
ASD instantons was established by Donaldson~\cite{DON-asd}
and used to prove nontriviality of the Donaldson invariants for
K\"ahler surfaces. Donaldson's theorem was extended to 
higher-dimensional K\"ahler manifolds by Uhlenbeck--Yau~\cite{UY}.

Another infinite-dimensional analogue of GIT is Donaldson's program 
for the study of constant scalar curvature K\"ahler (cscK) metrics.  
It was noted by Donaldson~\cite{DON-K1} and Fujiki~\cite{FUJIKI} 
that the scalar curvature can be interpreted as a moment map 
for the action of the group of Hamiltonian symplectomorphisms 
of a symplectic manifold $(V,\om_0)$ on the space  $\cJ_0$ 
of $\om_0$-compatible integrable complex structures.  
It was also noted by Donaldson~\cite{DON-K1} that 
in this setting the Futaki invariants~\cite{FUTAKI} are analogues of 
the Mumford numerical invariants, the space of K\"ahler potentials 
is the analogue of~$\rG^c/\rG$, the Mabuchi functional~\cite{MABUCHI1,MABUCHI2}
is the analogue of the the log-norm function in Kempf--Ness~\cite{KN},
and that Tian's notion of K-stability~\cite{TIAN} can be understood as 
an analogue of stability in GIT.  This led to Donaldson's conjecture
relating K-stability to the existence of cscK metrics~\cite{DON-K2,DON-K3} 
and refining earlier conjectures by Yau~\cite{YAU2} and Tian~\cite{TIAN}.  
The Yau--Tian--Donaldson conjecture is the analogue of the 
Hilbert--Mumford criterion for $\mu$-polystability
(Theorem~\ref{thm:HMps}), with the $\mu$-weights 
replaced by Donaldson's generalized Futaki invariants. 
The earlier conjecture of Yau applies to Fano manifolds,
relates K-stability to the existence of K\"ahler--Einstein 
metrics, and has been confirmed in 2013 by 
Chen--Donaldson--Sun~\cite{DS,CDS0,CDS1,CDS2,CDS3}.
The Yau--Tian--Donaldson conjecture, in the case where the 
first Chern class is not a multiple of the K\"ahler class, 
is still open.  The moment-weight inequality in this setting was 
proved by Donaldson~\cite{DON-K4} and Chen~\cite{CIII,CIV}.  

In this situation the duality between the positive curvature 
manifold $\rG$ and the negative curvature manifold $\rG^c/\rG$ 
is particularly interesting.  The analogue of $\rG$ is the 
group $\cG_0$ of Hamiltonian symplectomorphisms of $(V,\om_0)$ 
with the $L^2$ inner product on the Lie algebra of Hamiltonian functions, 
and the analogue of $\rG^c/\rG$ is the space $\cH_0$
of K\"ahler potentials on $(V,J_0,\om_0)$ with an $L^2$ Riemannian metric.
On the one hand the distance function associated to the $L^2$ metric 
on $\cG_0$ is trivial by a result of Eliashberg--Polterovich~\cite{EP}.
On the other hand Chen~\cite{CI} proved that 
in $\cH_0$ any two points are joined 
by a unique geodesic of class $C^{1,1}$
(a solution of the Monge--Amp\`ere equation),
that~$\cH_0$ is a genuine metric space,
and that cscK metrics with negative first Chern class 
are unique in their K\"ahler class.
The latter result was extended by Donaldson~\cite{DON-Kunique} 
to all cscK metrics with discrete automorphism groups.
Calabi--Chen~\cite{CCII} proved that~$\cH_0$ is negatively 
curved in the sense of Alexandrov and, 
assuming all geodesics are smooth, that extremal 
metrics are unique up to holomorphic diffeomorphism.
(Extremal metrics are the analogues of critical points 
of the moment map squared and their uniqueness up
to holomorphic diffeomorphism is the analogue of the 
Ness Uniqueness Theorem~\ref{thm:NESS1}.)
Chen--Tian~\cite{CT} removed the hypothesis on the geodesics.
They also proved that cscK metrics minimize the Mabuchi functional.  
This was independently proved by Donaldson~\cite{DON-K5} in the projective case.
As noted by Chen~\cite{CIV} and Chen--Sun~\cite{CSV}, several straight forward 
statements in GIT have analogues in the cscK  setting that are open questions. 
At the time one of these was convexity of the 
Mabuchi functional along $C^{1,1}$ geodesics.
This has now been settled by Berman--Berndtsson~\cite{BB}
and Chen--Li--Paun~\cite{CLP}.
In~\cite{STPAUL} S.T.~Paul introduced
another notion of stability and in~\cite{STPAUL1} he proved 
that it is equivalent to a properness condition on certain 
finite-dimensional approximations of the Mabuchi functional.
For Fano manifolds with discrete automorphism groups 
it has been shown that P-stability is equivalent to the 
aforementioned K-stability condition of Yau--Tian--Donaldson,
by combining the work of Chen--Donaldson--Sun and S.T.~Paul
with a partial $C^0$ estimate of G\'abor Sz\'ekelyhidi~\cite{S1}
(see Chen--Sun--Wang~\cite{CSW}).

The emphasis in the present book is on self 
contained proofs in the finite-dimensional setting. 
For other expositions of geometric invariant theory and 
the moment map see the papers by Thomas~\cite{RPTHOMAS},
which includes a detailed discussion of the 
cscK analogue, and by Woodward~\cite{W},
which contains many finite-dimensional examples.

Here is a brief description of the content of the book.  
Two preliminary chapters introduce the basic setup 
and the moment map (Chapter~\ref{ch:MOMENT}) 
and examine the negative gradient flow 
of the moment map squared (Chapter~\ref{ch:MU}).
Chapter~\ref{ch:KNF} introduces the Kempf--Ness function 
and establishes its basic properties.
The Mumford numerical invariants are defined
in Chapter~\ref{ch:WEIGHTS} and shown to be 
invariant under the~$\rG^c$-action and under 
the Mumford equivalence relation on the set 
of toral generators. Chapter~\ref{ch:MW1} 
establishes the moment-weight inequality
and derives several consequences 
such as the Kirwan--Ness Inequality, 
the Moment Limit Theorem,
and the Ness Uniqueness Theorem.
The $\mu$-stability notions are introduced in
Chapter~\ref{ch:STABLE} which also proves 
the Kempf--Ness Theorem.  Chapter~\ref{ch:STAG}
deals with the classical algebraic geometric setting
of linear actions on projective varieties by 
reductive groups and shows how it fits into 
the symplectic setup. 
Chapter~\ref{ch:RAT} deals with the converse question 
and examines under which rationality conditions
on the symplectic form and the moment map 
the general symplectic setting of the present book
reduces to the classical algebro geometric setting.
Chapter~\ref{ch:DOMINANT} is devoted to the Kempf 
Existence Theorem and shows that in the 
$\mu$-unstable case the moment-weight inequality 
is actually an equality.  It also shows that the 
supremum in~\eqref{eq:MW} is always attained.
That the supremum over~${\cg\setminus\{0\}}$
in~\eqref{eq:MW} agrees with the supremum 
over~$\Lambda$ requires continuous dependence 
of the weight~$w_\mu(x,\xi)$ on~$\xi$ for torus
actions and this is the subject of Chapter~\ref{ch:TORUS}.
The Hilber--Mumford criterion is proved in Chapter~\ref{ch:HM}.
Chapter~\ref{ch:CRIT} explains a criterion by 
Gabor Sz\'ekelyhidi for points in~$X$ whose $\rG^c$-orbits 
contain critical points of the moment map squared.  
The criterion takes the form of polystability 
with respect to the action of a suitable subgroup.  
Several examples are discussed in Chapter~\ref{ch:EX}.

Five appendices deal with relevant background material. 
Appendix~\ref{app:NONPOS} establishes some properties 
of complete simply connected Riemannian manifolds 
with non\-positive sectional curvature and  
contains a proof of Cartan's Fixed Point Theorem.
Appendix~\ref{app:Gc} establishes the existence of a 
complexification~$\rG^c$ of a compact Lie group~$\rG$
and shows how it is characterized by a universality pro\-perty.
Appendix~\ref{app:GcG} shows that the homogeneous
space~$\rG^c/\rG$ is a complete simply connected Riemannian
manifold with nonpositive sectional curvature.
Appendix~\ref{app:Lambda} introduces parabolic subgroups
and the Mumford equivalence relation on the space 
of toral generators and Appendix~\ref{app:BOREL}
is devoted to the proof that each element of~$\rG^c$
factorizes as a product of an element of a given 
parabolic subgroup and an element of~$\rG$. 


\chapter{The moment map}\label{ch:MOMENT}

Throughout $(X,\om,J)$ denotes a closed K\"ahler manifold,
i.e.\ $X$ is a compact manifold without boundary, 
$\om$  is a symplectic form on $X$, 
and $J$ is an integrable complex structure on $X$,
compatible with $\om$ so
$
\inner{\cdot}{\cdot} := \om(\cdot,J\cdot)
$
is a Riemannian metric. Denote by $\nabla$ the 
corresponding Levi-Civita connection. 

{\bf The moment map}. \index{moment map}
Let $\rG\subset\rU(n)$ be a compact Lie group acting on~$X$ 
by K\"ahler isometries so that the action of $\rG$ 
preserves all three structures $\inner{\cdot}{\cdot}$, 
$\om$, $J$. Denote the group action by 
$$
\rG\times X\to X:(u,x)\mapsto ux
$$
and the infinitesimal action of the Lie 
algebra~${\cg:=\Lie(\rG)\subset\cu(n)}$ by 
$$
\cg\to\Vect(X):\xi\mapsto v_\xi.
$$
We assume that $\cg$ is equipped with an invariant 
inner product and that the group action is Hamiltonian. 
Let $\mu:X\to\cg$ be a moment map for the action, 
i.e.\ $v_\xi$ is a Hamiltonian vector field with Hamiltonian function 
$$
H_\xi:=\inner{\mu}{\xi}
$$
so $\iota(v_\xi)\om=dH_\xi$ or, equivalently, 
\begin{equation}\label{eq:mu1}
\inner{d\mu(x)\xhat}{\xi} = \om(v_\xi(x),\xhat)
\end{equation}
for $x\in X$, $\xhat\in T_xX$, and $\xi\in\cg$.
We assume that the moment map is equivariant, i.e.\ 
\begin{equation}\label{eq:mu2}
\mu(ux) = u\mu(x)u^{-1},\qquad
\inner{\mu(x)}{[\xi,\eta]} = \om(v_\xi(x),v_\eta(x))
\end{equation}
for $x\in X$, $u\in\rG$, and $\xi,\eta\in\cg$.
The two equations in~\eqref{eq:mu2} are equivalent 
whenever $\rG$ is connected. 

\bigbreak

{\bf The complexified group}.\index{Lie group!complexification}
The map\index{complexified group}
$$
\rU(n)\times\cu(n)\to\GL(n,\C):(u,\eta)\mapsto\exp(\i\eta)u
$$ 
is a diffeomorphism, by polar decomposition, 
and the image of $\rG\times\cg$ under 
this diffeomorphism is denoted by 
\begin{equation}\label{eq:Gc}
\rG^c := \left\{hu\,|\,u\in\rG,\,h=\exp(\i\eta),\,\eta\in\cg\right\}.
\end{equation}
This is a Lie subgroup of $\GL(n,\C)$,\footnote{
For us a Lie subgroup is closed.}
called the {\bf complexification of $\rG$}
(see Theorem~\ref{thm:Gc3}). 
It contains $\rG$ as a maximal compact subgroup,
the quotient $\rG^c/\rG$ is connected,
and the Lie algebra of $\rG^c$ is the complexification
$$
\cg^c := \Lie(\rG^c) = \cg\oplus \i\cg
$$
of the Lie algebra of $\rG$ (see Theorem~\ref{thm:Gc1}). 
We will consistently use the notations 
$$
\zeta=\xi+\i\eta,\qquad
\Re(\zeta):=\xi,\qquad \Im(\zeta):=\eta.
$$
for the elements $\zeta\in\cg^c$ and their real 
and imaginary parts $\xi,\eta\in\cg$.
The reader is cautioned that the eigenvalues
of $\xi,\eta$ are imaginary and the eigenvalues 
of $\i\xi,\i\eta$ are real. A complex Lie group 
is called {\bf reductive} if it is the 
complexification\index{Lie group!reductive}
of a compact Lie group.\index{reductive Lie group}

{\bf The action of the complexified group}.
Let $\rG^c\subset\GL(n,\C)$ be the complexification of $\rG$
and let $\cg^c=\cg+\i\cg$ be its Lie algebra.
Then every Lie group homomorphism from $\rG$ 
to a complex Lie group extends uniquely to 
a homomorphism from $\rG^c$ to that 
complex Lie group (see Theorem~\ref{thm:Gc1}). 
Taking the target group to be the group of 
holomorphic automorphisms of $X$ one obtains 
a holomorphic group action of $\rG^c$ on $X$. 
Denote the group action by
$$
\rG^c\times X\to X:(g,x)\mapsto gx
$$ 
and the infinitesimal action of the Lie algebra by
$$
\cg^c\to\Vect(X):\zeta\mapsto v_\zeta:=v_\xi+Jv_\eta.
$$
Here $v_\xi$ is the Hamiltonian vector field of the function $H_\xi$ 
and~${Jv_\eta=\nabla H_\eta}$ is the gradient vector field of the 
function~$H_\eta$.  Since~${v_\zeta}$ is a holomorphic vector field 
and~${\nabla J=0}$, we have 
\begin{equation}\label{eq:Jvxi}
[v_\zeta,Jw] = J[v_\zeta,w],\qquad
\Nabla{Jw}v_\zeta = J\Nabla{w}v_\zeta
\end{equation}
for every $\zeta\in\cg^c$ and every 
vector field\index{Lie bracket!sign convention}
$w\in\Vect(X)$.\footnote
{
We use the sign convention 
$[v,w]:=\Nabla{w}v-\Nabla{v}w$ 
for the Lie bracket of two vector fields
so that the infinitesimal action is a homomorphism
$\cg^c\to\Vect(X)$ of Lie algebras 
(and not an anti-homomorphism).
}

\bigbreak

{\bf Alternative notation for the infinitesimal action}.
For each $x\in X$ we use the alternative notation
$L_x:\cg\to T_xX$ and $L_x^c:\cg^c\to T_xX$ for the 
infinitesimal action of the Lie algebras $\cg$ and $\cg^c$. 
Thus
\begin{equation}\label{eq:Lx}
L_x\xi := v_\xi(x),\qquad
L_x^c\zeta := v_\zeta(x) = L_x\xi+JL_x\eta,
\end{equation}
for $\zeta=\xi+\i\eta\in\cg^c$.  
Then equations~\eqref{eq:mu1} and~\eqref{eq:mu2}
take the form
\begin{equation}\label{eq:mu3}
{L_x}^* = d\mu(x)J,\qquad d\mu(x)^*=JL_x,\qquad
d\mu(x)L_x\xi = -[\mu(x),\xi].
\end{equation}
for $x\in X$ and $\xi\in\cg$. 

\begin{lem}\label{le:ZERO}
Let $x_0,x_1\in\mu^{-1}(0)$.  If $x_1\in\rG^c(x_0)$ then $x_1\in\rG(x_0)$.
In fact, if $\eta\in\cg$ and $u\in\rG$ satisfy 
$
\exp(\i\eta)ux_0=x_1
$ 
then $ux_0=x_1$ and $L_{x_1}\eta=0$. 
\end{lem}

\begin{proof}
Choose $g_0\in\rG^c$ such that $x_1=g_0x_0$ and define
$\eta\in\cg$ and $u\in\rG$ by 
$$
g_0=:\exp(\i\eta)u.
$$  
Define the curve $x:[0,1]\to X$ by 
$x(t) := \exp(\i t\eta)ux_0$.
Then $x(0)=ux_0$, $x(1)=x_1$, and $\dot x = JL_x\eta$. 
Hence, by equation~\eqref{eq:mu1},
\begin{equation}\label{eq:ZERO}
\begin{split}
\frac{d}{dt}\inner{\mu(x)}{\eta} 
&= 
\inner{d\mu(x)\dot x}{\eta} \\
&= 
\om(L_x\eta,\dot x)  \\
&= 
\om(L_x\eta,JL_x\eta) \\
&= 
\abs{L_x\eta}^2  \\
&\ge 
0.
\end{split}
\end{equation}
Since $\inner{\mu(ux_0)}{\eta}=\inner{\mu(x_1)}{\eta}=0$
it follows that $L_{x(t)}\eta=0$ for all $t$.  
Hence $x(t)$ is constant and hence $x_1=ux_0$. 
This proves Lemma~\ref{le:ZERO}.
\end{proof}

Lemma~\ref{le:ZERO} asserts that two points in the zero set of the moment map
are equivalent under $\rG^c$ if and only if they are equivalent under $\rG$.
In notation commonly used in symplectic geometry it says that the 
Marsden--Weinstein quotient\index{Marsden--Weinstein quotient}
is homeomorphic to the algebraic geometry quotient, i.e.\ 
$$
X\dslash\rG \simeq X^\ps/\rG^c, 
$$
where $X\dslash\rG := \mu^{-1}(0)/\rG$ and
$$
X^\ps:=\left\{x\in X\,\big|\,\rG^c(x)\cap\mu^{-1}(0)\ne\emptyset\right\}.
$$
These spaces can be highly singular. 
Lemma~\ref{le:ONE} below asserts that an element ${x\in\mu^{-1}(0)}$ 
is a regular point for the moment map if and only if its isotropy 
subgroup 
$$
\rG_x:=\left\{u\in\rG\,|\,ux=x\right\}
$$ 
is discrete. 
Lemma~\ref{le:isotropy} below asserts that the complex isotropy 
subgroup
$$
\rG^c_x:=\left\{g\in\rG^c\,|\,gx=x\right\}
$$ 
is the complexification of~$\rG_x$ whenever $\mu(x)=0$.  
This means that there is an isomorphism of orbifolds
$$
X^\s\dslash\rG\cong X^\s/\rG^c,
$$
where $X^\s\subset X^\ps$ is the open subset of 
all points $x\in X^\ps$ with discrete isotropy subgroup $\rG^c_x$.
In general, $X^\ps$ is not an open subset of $X$. 
In geometric invariant theory one studies 
the quotient $X^\ss\dslash\rG^c$, where $X^\ss$ is the 
open subset of all points $x\in X$ such that 
$\overline{\rG^c(x)}\cap\mu^{-1}(0)\ne\emptyset$,
and two such points~${x,x'\in X^\ss}$ are equivalent 
iff $\mu^{-1}(0)\cap\overline{\rG^c(x)}\cap\overline{\rG^c(x')}\ne\emptyset$.

\begin{lem}\label{le:ONE}
Let $x\in X$ such that $\mu(x)=0$. 
The following are equivalent. 

\smallskip\noindent{\bf (i)}
$d\mu(x):T_xX\to\cg$ is onto.

\smallskip\noindent{\bf (ii)}
$L_x:\cg\to T_xX$ is injective.

\smallskip\noindent{\bf (iii)}
$L^c_x:\cg^c\to T_xX$ is injective.
\end{lem}

\begin{proof}
The equation $d\mu(x)J={L_x}^*$ in~\eqref{eq:mu3}
shows that~(i) is equivalent to~(ii).  That~(iii) implies~(ii)
is obvious.  Now assume~(ii) and
choose an element 
$
\zeta=\xi+\i\eta\in\ker\,L^c_x.
$  
Then 
\begin{equation*}
\begin{split}
0
&= 
d\mu(x)\bigl(L_x\xi+JL_x\eta\bigr)  \\
&=
-[\mu(x),\xi] + {L_x}^*L_x\eta \\
&=
{L_x}^*L_x\eta
\end{split}
\end{equation*}
by~\eqref{eq:mu3}.
Hence $L_x\eta=0$, hence 
$$
L_x\xi=L_x^c\zeta-JL_x\eta=0,
$$
and hence $\xi=\eta=0$ by~(ii). Thus~(ii) implies~(iii) 
and this proves Lemma~\ref{le:ONE}.
\end{proof}

The hypothesis $\mu(x)=0$ in Lemma~\ref{le:ONE} cannot be removed.
For example a pair of nonantipodal points on the $2$-sphere has 
a trivial isotropy subgroup in $\SO(3)$ but a nontrivial isotropy
subgroup in the complexified group $\PSL(2,\C)$.
(When the points are antipodal the moment map is zero.)

\begin{lem}\label{le:isotropy}
{\bf (i)}
Let $x\in X$ such that $\mu(x)=0$.  
Then $\rG^c_x$ is the complexification of $\rG_x$, 
i.e.\ 
$$
\rG^c_x= \left\{\exp(\i\eta)u\,|\,u\in\rG_x,\,\eta\in\ker\,L_x\right\}.
$$

\smallskip\noindent{\bf (ii)}
Let $\xi\in\cg$.  Then the isotropy subgroup of $\xi$ in~$\rG^c$ 
is the complexification of the isotropy subgroup in~$\rG$, 
i.e.\ if $\eta\in\cg$, $u\in\rG$, and $g:=\exp(\i\eta)u$, then
$$
g\xi g^{-1}=\xi 
\qquad\iff\qquad
u\xi u^{-1}=\xi
\quad\mbox{and}\quad 
[\xi,\eta]=0.
$$
\end{lem}

\begin{proof}
We prove~(i). 
If $u\in\rG$ and $\eta\in\cg$ are such that $g:=\exp(\i\eta)u\in\rG^c_x$,
then Lemma~\ref{le:ZERO} asserts that $ux=x$
and $L_x\eta=0$. This proves~(i).

We prove~(ii).  Assume $g\xi g^{-1}=\xi$, 
abbreviate 
$$
h:=\exp(-\i\eta), 
$$
and let 
$$
\C^n=V_1\oplus\cdots\oplus V_k
$$
be the eigenspace decomposition of $\xi$.  Since 
$$
h\xi h^{-1}=u\xi u^{-1},
$$
it follows that $hV_i\perp hV_j$ for $i\ne j$.  
Since $h=h^*>0$ this implies 
$$
h^2V_i\subset V_i, 
$$
hence $hV_i\subset V_i$, and hence 
$$
uV_i\subset V_i,\qquad \eta V_i\subset V_i
$$ 
for all $i$.  Thus $u$ and $\eta$ commute with $\xi$. 
This proves Lemma~\ref{le:isotropy}.
\end{proof}


\chapter{The moment map squared}\label{ch:MU}

Throughout this chapter we denote by $f:X\to\R$ 
the function defined by
\begin{equation}\label{eq:f}
f(x) := \tfrac12\abs{\mu(x)}^2
\qquad\mbox{for }x\in X.
\end{equation}

\begin{lem}\label{le:MUCRIT}
The gradient of $f$ is given by 
$$
\nabla f(x)=JL_x\mu(x)
$$
for $x\in X$. Hence $x\in X$ is a critical point 
of $f$ if and only if $L_x\mu(x)=0$.
\end{lem}

\begin{proof}
By equation~\eqref{eq:mu1} we have 
$\inner{d\mu(x)\xhat}{\xi} = \om(L_x\xi,\xhat)$.
Take $\xi=\mu(x)$. Then 
$$
df(x)\xhat
= \inner{d\mu(x)\xhat}{\mu(x)}
= \om(L_x\mu(x),\xhat)
= \inner{JL_x\mu(x)}{\xhat}
$$
for $\xhat\in T_xX$.  This proves
Lemma~\ref{le:MUCRIT}.
\end{proof}

By Lemma~\ref{le:MUCRIT} the negative gradient flow line of $f$
through $x_0\in X$ is the solution of the differential equation.
\begin{equation}\label{eq:KN1}
\dot x = -JL_x\mu(x),\qquad x(0)=x_0.
\end{equation}
The complexified group orbits are invariant under the gradient flow. 

\begin{lem}\label{le:GRADFLOW}
Let $x_0\in X$, let $x:\R\to X$ be the unique solution 
of the differential equation~\eqref{eq:KN1}, 
and let~${g:\R\to\rG^c}$ be the unique solution 
of the differential equation
\begin{equation}\label{eq:KN2}
g(t)^{-1}\dot g(t) = \i\mu(x(t)),\qquad g(0)=\one.
\end{equation}
Then 
$$
x(t)=g(t)^{-1}x_0
$$ 
for all $t\in\R$.
\end{lem}

\begin{proof}
Define $y:\R\to X$ by 
$$
y(t) :=  g(t)^{-1}x_0.
$$
Since~${\tfrac{d}{dt}g^{-1}=-g^{-1}\dot gg^{-1}}$
and~${g^{-1}\dot g = \i\mu(x)}$, it follows that
$$
\dot y = - g^{-1}\dot gg^{-1}x_0 
= - L^c_{g^{-1}x_0}(g^{-1}\dot g)
= - JL_y\mu(x)
$$
and~${y(0)=x_0}$.  Hence~${x(t)=y(t)=g(t)^{-1}x_0}$ for every~${t\in\R}$.
This proves Lemma~\ref{le:GRADFLOW}.
\end{proof}

\begin{thm}[{\bf Convergence Theorem}]\label{thm:XINFTY}
\ 

\noindent
There\index{Convergence Theorem}
exist positive constants
${C,c,\eps}$ and~${\tfrac{1}{2}<\alpha<1}$ 
with the following significance. Let~${x_0\in X}$ and 
let~${x:\R\to X}$ be the unique solution of~\eqref{eq:KN1}.
Then the limit 
$$
x_\infty := \lim_{t\to\infty}x(t)
$$
exists and satisfies ${L_{x_\infty}\mu(x_\infty)=0}$.
Moreover, there exists a constant $T\in\R$ 
such that, for all $t>T$,
\begin{equation}\label{eq:Loja1}
\begin{split}
d(x(t),x_\infty) 
&\le 
\int_t^\infty\abs{\dot x(s)}\,ds \\
&\le 
\frac{C}{1-\alpha}\bigl(f(x(t))-f(x_\infty)\bigr)^{1-\alpha} \\
&\le 
\frac{c}{(t-T)^\eps}.
\end{split}
\end{equation}
\end{thm}

The proof of Theorem~\ref{thm:XINFTY} is based on the 
Lojasiewicz gradient inequality~\cite{LOJA,HARAUX}, 
which holds for general analytic gradient flows.  
That it also applies to the moment map squared
was noted by Duistermaat (see Lerman~\cite{LERMAN} 
and also Chen--Sun~\cite[Corollary~3.2]{CS}).  
The result was carried over to certain infinite-dimensional settings 
by Simon~\cite{SIMON} and Morgan--Mrowka--Rubermann~\cite{MMR}.
We include a proof for completeness of the exposition. 

\begin{proof}[Proof of Theorem~\ref{thm:XINFTY}]
Using the Marle--Guillemin--Sternberg local normal form 
(see~\cite[Theorem~2.1]{LERMAN}), one can show that 
the moment map is locally real analytic in suitable coordinates.  
This implies that 
$
f=\tfrac12\abs{\mu}^2
$ 
satisfies the {\bf Lojasiewicz gradient inequality}:
There\index{Lojasiewicz gradient inequality}
exist constants $\delta>0$, $C>0$, and $1/2<\alpha<1$ 
such that, for every critical value $a$ of $f$ and every $x\in X$, 
\begin{equation}\label{eq:Loja2}
\abs{f(x)-a}<\delta\qquad\implies\qquad
\abs{f(x)-a}^\alpha\le C\abs{\nabla f(x)}.
\end{equation}
For a proof see Bierstone--Milman~\cite[Theorem~6.4]{BM}.

\bigbreak

Let~${x:\R\to X}$ be a nonconstant negative gradient 
flow line of~$f$. Then 
$$
a:=\lim_{t\to\infty}f(x(t))
$$ 
is a critical value of $f$. Choose a constant $T\in\R$ 
such that 
$$
a<f(x(t))<a+\delta
$$ 
for $t\ge T$.   Then, for $t\ge T$, 
$$
-\frac{d}{dt}\left(f(x)-a\right)^{1-\alpha}
= 
(1-\alpha)(f(x)-a)^{-\alpha}\abs{\nabla f(x)}^2 \\
\ge
\frac{1-\alpha}{C}\abs{\dot x}.
$$
Integrate this inequality over the interval $[t,\infty)$ to obtain
\begin{equation}\label{eq:Loja3}
\int_t^\infty\abs{\dot x(s)}\,ds
\le \frac{C}{1-\alpha}\bigl(f(x(t))-a\bigr)^{1-\alpha}
\qquad\mbox{for }t\ge T.
\end{equation}
This proves the existence of the limit
$$
x_\infty:=\lim_{t\to\infty}x(t).
$$
This limit is a critical point of $f$ and hence satisfies~${L_{x_\infty}\mu(x_\infty)=0}$
by Lemma~\ref{le:MUCRIT}.

Now abbreviate
$$
\rho(t) := (f(x(t))-a)^{-(2\alpha-1)}. 
$$
Then 
$$
\dot\rho(t)
= (2\alpha-1)\bigl(f(x(t))-a\bigr)^{-2\alpha}\abs{\nabla f(x(t))}^2
\ge \frac{2\alpha-1}{C^2}
\qquad\mbox{for }t\ge T.
$$
This implies 
$$
\rho(t)\ge \frac{2\alpha-1}{C^2}(t-T)
\qquad\mbox{for }t\ge T
$$
and hence
$$
\bigl(f(x(t))-a\bigr)^{1-\alpha} 
= 
\rho(t)^{-\frac{1-\alpha}{2\alpha-1}} \\
\le
\left(\frac{2\alpha-1}{C^2}(t-T)\right)^{-\frac{1-\alpha}{2\alpha-1}}
$$
for $t>T$. Thus
$$
\frac{C}{1-\alpha}\bigl(f(x(t))-a\bigr)^{1-\alpha}
\le \frac{c}{(t-T)^\eps},\;\;
\eps:=\frac{1-\alpha}{2\alpha-1},\;\;
c:= \frac{C}{1-\alpha}
\left(\frac{C^2}{2\alpha-1}\right)^\eps,
$$
for $t>T$. By~\eqref{eq:Loja3}, this 
proves Theorem~\ref{thm:XINFTY}.
\end{proof}

\medskip
We close this section with three lemmas about the Hessian 
of the moment map squared, culminating in
Lijing Wang's inequality in Lemma~\ref{le:WANG}. 
They do not play any role elsewhere in this paper.  

\begin{lem}\label{le:vxi}
Let $x\in X$, $\xhat\in T_xX$, and $\xi,\eta\in\cg$.
If $L_x\xi=0$ then
\begin{equation}\label{eq:dmuv}
d\mu(x)\Nabla{\xhat}v_\xi(x) = -[d\mu(x)\xhat,\xi]
\end{equation}
and
\begin{equation}\label{eq:L*L}
{L_x}^*L_x[\xi,\eta] = [\xi,{L_x}^*L_x\eta].
\end{equation}
\end{lem}

\begin{proof}
Differentiate the function 
$$
\inner{[\mu,\xi]}{\eta} 
= \inner{\mu}{[\xi,\eta]} 
= \om(v_\xi,v_\eta) 
= \inner{Jv_\xi}{v_\eta}
$$
at $x$ in the direction $\xhat$ and use the equations
$v_\xi(x)=0$ and $\nabla J=0$ to obtain
\begin{eqnarray*}
\inner{[d\mu(x)\xhat,\xi]}{\eta}
&=& 
\inner{J\Nabla{\xhat}v_\xi(x)}{v_\eta(x)} \\
&=& 
-\om(v_\eta(x),\Nabla{\xhat}v_\xi(x)) \\
&=& 
-dH_\eta(x)\Nabla{\xhat}v_\xi(x) \\
&=&
-\inner{d\mu(x)\Nabla{\xhat}v_\xi(x)}{\eta}.
\end{eqnarray*}
This proves~\eqref{eq:dmuv}.

It follows from equations~\eqref{eq:Jvxi}, 
\eqref{eq:mu3}, and~\eqref{eq:dmuv} that
\begin{eqnarray*}
[\xi,{L_x}^*L_x\eta]
&=&
-[d\mu(x)JL_x\eta,\xi] \\
&=&
d\mu(x)\Nabla{Jv_\eta}v_\xi(x) \\
&=&
d\mu(x)[v_\xi,Jv_\eta](x) \\
&=& 
d\mu(x)J[v_\xi,v_\eta](x) \\
&=&
d\mu(x)Jv_{[\xi,\eta]}(x) \\
&=&
{L_x}^*L_x[\xi,\eta].
\end{eqnarray*}
This proves~\eqref{eq:L*L} and Lemma~\ref{le:vxi}.
\end{proof}

\begin{lem}\label{le:Hessian}
The covariant Hessian of 
$
f=\tfrac12\abs{\mu}^2
$
at a point~${x\in X}$ is the quadratic 
form~${d^2f_x:T_xX\to\R}$ given by  
\begin{equation}\label{eq:Hess1}
d^2f_x(\xhat) = \abs{d\mu(x)\xhat}^2
+ \inner{\Nabla{J\xhat}v_\xi(x)}{\xhat},\qquad
\xi:=\mu(x),
\end{equation}
for $\xhat\in T_xX$.  If $x$ is a critical point of $f$ then
\begin{equation}\label{eq:Hess2}
d^2f_x(JL_x\eta) = \abs{{L_x}^*L_x\eta}^2-\abs{[\mu(x),\eta]}^2
\end{equation}
for every $\eta\in\cg$.
\end{lem}

\begin{proof}
Since $\nabla f(x)=JL_x\mu(x)=Jv_{\mu(x)}(x)$ and $\nabla J=0$,
we have
$$
\Nabla{\xhat}\nabla f(x) 
= JL_xd\mu(x)\xhat + J\Nabla{\xhat}v_\xi(x),\qquad
\xi:=\mu(x).
$$
Hence the covariant Hessian of $f$ at $x$ is given by 
$$
d^2f_x(\xhat) 
:= \inner{\Nabla{\xhat}\nabla f(x)}{\xhat} 
= \abs{d\mu(x)\xhat}^2
+ \inner{\Nabla{J\xhat}v_\xi(x)}{\xhat}.
$$
Here the last step uses the identity $JL_x=d\mu(x)^*$ in~\eqref{eq:mu3}.  
Take $\xhat:=JL_x\eta$ and assume $L_x\mu(x)=0$ to obtain 
\begin{eqnarray*}
d^2f_x(JL_x\eta) 
&=& 
\abs{d\mu(x)JL_x\eta}^2 - \inner{\Nabla{L_x\eta}v_\xi(x)}{JL_x\eta} \\
&=&
\abs{{L_x}^*L_x\eta}^2 - \inner{d\mu(x)\Nabla{L_x\eta}v_\xi(x)}{\eta} \\
&=&
\abs{{L_x}^*L_x\eta}^2 + \inner{[d\mu(x)L_x\eta,\xi]}{\eta}  \\
&=&
\abs{{L_x}^*L_x\eta}^2 + \inner{d\mu(x)L_x\eta}{[\xi,\eta]}  \\
&=& 
\abs{{L_x}^*L_x\eta}^2 - \abs{[\mu(x),\eta]}^2.
\end{eqnarray*}
Here the third step follows from equation~\eqref{eq:dmuv} in Lemma~\ref{le:vxi}
with $\xhat=L_x\eta$ and the last step follows from equation~\eqref{eq:mu3}.
This proves Lemma~\ref{le:Hessian}. 
\end{proof}

\begin{lem}[{\bf Lijing Wang's Inequality}]\label{le:WANG}
Let $x\in X$ be a\index{Lijing Wang's inequality}  
critical point of the moment map squared. 
Then, for every $\eta\in\cg$,
\begin{equation}\label{eq:WANG}
\abs{[\mu(x),\eta]}\le \abs{{L_x}^*L_x\eta}.
\end{equation} 
\end{lem}

\begin{proof}
The proof is taken from Lijing Wang's paper~\cite[Theorem~3.8]{WANG}.
Define the linear maps $A,B:\cg\to\cg$ by 
$$
A\xi := {L_x}^*L_x\xi,\qquad B\xi := [\mu(x),\xi].
$$
Thus $A$ is self-adjoint and $B$ is skew-adjoint. 
Moreover, $A$ and $B$ commute by Lemma~\ref{le:vxi}.
Identify $\cg^c$ with $\cg\oplus\cg$ and define 
$P^\pm:\cg^c\to\cg^c$ by 
$$
P^+ := \left(\begin{array}{rr}
A & B \\ -B & A \end{array}\right),\qquad
P^- := \left(\begin{array}{rr}
A & -B \\ B & A \end{array}\right).
$$
Then the operators $P^+$ and $P^-$ are self-adjoint, 
they commute, and
$$
\inner{\zeta}{P^\pm\zeta}=\abs{L_x\xi\pm JL_x\eta}^2.
$$
Hence the operator
$$
Q:=P^+P^-
$$
is self-adjoint and nonnegative. Hence
$$
0\le \inner{\xi}{Q\xi} = \abs{{L_x}^*L_x\xi}^2-\abs{[\mu(x),\xi]}^2
$$
for $\xi\in\cg\subset\cg^c$.  This proves Lemma~\ref{le:WANG}.
\end{proof} 


\chapter{The Kempf--Ness function}\label{ch:KNF}

The present chapter introduces the Kempf--Ness 
function
$$
\Phi_x:\rG^c/\rG\to\R
$$ 
on the Hadamard space~$\rG^c/\rG$, associated to an element~${x\in X}$.
In particular, it shows that every gradient flow line of the 
moment map squared in~$\rG^c(x)$ gives rise to
a gradient flow line of the Kempf--Ness function.
It also shows that the Kempf--Ness function is 
convex along geodesics and that it is a Morse--Bott function.

Equip $\rG^c$ with the unique left invariant 
Riemannian metric which agrees with the inner product 
$$
\inner{\xi_1+\i\eta_1}{\xi_2+\i\eta_2}_{\cg^c}
=\inner{\xi_1}{\xi_2}_\cg+\inner{\eta_1}{\eta_2}_\cg
$$
on the tangent space $\cg^c$ to $\rG^c$ at the identity. 
This metric is invariant under the right $\rG$-action. Let 
\begin{equation}\label{eq:piM}
\pi:\rG^c\to M,\qquad   M:=\rG^c/\rG,
\end{equation}
be the projection onto the right cosets of $\rG$.
This is  a principal $\rG$-bundle over a contractible manifold.
The (orthogonal) splitting 
$$
\cg^c=\cg\oplus\i\cg
$$
extends to a left invariant principal connection on $\pi$.
The projection from the horizontal bundle 
(i.e.\ the summand corresponding to $\i\cg$) 
defines a $\rG^c$-invariant Riemannian metric 
of nonpositive curvature on $M$ (Appendix~\ref{app:GcG}).
The geodesic on $M$ through $p=\pi(g)$ in the direction
$v=d\pi(g)g\i\xi\in T_pM$ has the form $\gamma(t)=\pi(g\exp(\i t\xi))$.  
Thus $M$ is complete and, by Hadamard's theorem, 
diffeomorphic to~$\i\cg$ (see also equation~\eqref{eq:Gc}). 

\begin{thm}[{\bf Conjugacy Theorem}]\label{thm:CONJUGACY}
Fix an element $x\in X$.\index{Conjugacy Theorem}

\smallskip\noindent{\bf (i)}
There exists a unique smooth function $\Phi_x:\rG^c\to\R$ such that
\begin{equation}\label{eq:Phix}
d\Phi_x(g)\ghat:=-\inner{\mu(g^{-1}x)}{\Im(g^{-1}\ghat)},\qquad
\Phi_x(u)=0,
\end{equation}
for all $g\in\rG^c$, all $\ghat\in T_g\rG^c$, and all $u\in\rG$.

\smallskip\noindent{\bf (ii)}
Define a map $\psi_x:\rG^c\to\rG^c(x)\subset X$ by~${\psi_x(g):=g^{-1}x}$.
Then $\psi_x$ intertwines the gradient vector field
$\nabla\Phi_x\in\Vect(\rG^c)$ and the gradient vector field
${\nabla f\in\Vect(X)}$, i.e.\ for all $g\in\rG^c$
\begin{equation}\label{eq:CONJUGACY}
d\psi_x(g) \nabla\Phi_x(g) = \nabla f(\psi_x(g)).
\end{equation}
\end{thm}

Assertion~(ii) of Theorem~\ref{thm:CONJUGACY} is a reformulation 
of Lemma~\ref{le:GRADFLOW} and shows again that $\nabla f$ 
is tangent to the $\rG^c$-orbits. Moreover, when the isotropy 
subgroup of $x$ is discrete, equations~\eqref{eq:Phix} 
and~\eqref{eq:CONJUGACY} are equivalent. 
So in this case the function $\Phi_x$ is uniquely 
determined by~\eqref{eq:CONJUGACY} and the 
normalization condition $\Phi_x(u)=0$ for all $u\in\rG$.  
(If $\rG$ is connected it suffices to impose the condition $\Phi_x(\one)=0$.)
In the opposite case, when $x$ is a fixed point of the group action,
equation~\eqref{eq:CONJUGACY} carries no information about $\Phi_x$. 

\begin{defn}\label{def:KN}
The function $\Phi_x:\rG^c\to\R$ in Theorem~\ref{thm:CONJUGACY} 
is called the {\bf lifted Kempf--Ness function} based at $x$.
It is $\rG$-invariant\index{Kempf--Ness function}
and hence descends to a function~${\Phi_x:M\to\R}$
denoted by the same symbol and called 
the {\bf Kempf--Ness function}.
\end{defn}

\begin{proof}[Proof of Theorem~\ref{thm:CONJUGACY}]
Define a  vector field $v_x\in\Vect(\rG^c)$ 
and a $1$-form $\alpha_x$ on $\rG^c$ by 
\begin{equation}\label{eq:alphax}
v_x(g) := -g\i\mu(g^{-1}x), \qquad
\alpha_x(g)\ghat:=-\inner{\mu(g^{-1}x)}{\Im(g^{-1}\ghat)}
\end{equation}
for $g\in\rG^c$ and $\ghat\in T_g\rG^c$. 
The vector field $v_x$ is horizontal since 
$\mu(x)\in\cg$ and is right $\rG$-equivariant, i.e.\ 
$v_x(gu)=v_x(g)u$ for $g\in\rG^c$ and $u\in\rG$.
We must prove the following.

\medskip\noindent {\bf Step 1.} 
The map  $\psi_x$ intertwines the vector fields 
$v_x$ and $\nabla f$, i.e.\ for~${g\in\rG^c}$,
\begin{equation}\label{eq:CON1}
d\psi_x(g)v_x(g)=\nabla f(\psi_x(g)).
\end{equation}

\medskip\noindent{\bf Step 2.} 
If $\Phi_x:\rG^c\to\R$ satisfies $d\Phi_x=\alpha_x$
then its gradient $\nabla \Phi_x$ is $v_x$, i.e.\ 
\begin{equation}\label{eq:CON2}
\alpha_x(g)\ghat=\inner{v_x(g)}{\ghat}_g
\end{equation}
for $g\in\rG^c$ and $\ghat\in T_g\rG^c$.
The inner product on the right is the left-invariant 
Riemannian metric on $\rG^c$.

\medskip\noindent {\bf Step 3.} 
There exists a unique function $\Phi_x:\rG^c\to\R$ 
such that
\begin{equation}\label{eq:CON3}
d\Phi_x=\alpha_x,\qquad \Phi_x|_\rG = 0.
\end{equation}
\bigbreak

\medskip
We prove Step~1.
If $\hat{g}\in T_g\rG^c$ then 
$$
d\psi_x(g)\ghat=-g^{-1}\ghat g^{-1}x = - L^c_{g^{-1}x}(g^{-1}\ghat).
$$
So taking $\ghat=v_x(g)=-g\i\mu(g^{-1}x)$ we get
$$
d\psi_x(g)v_x(g) 
= L^c_{g^{-1}x}(\i\mu(g^{-1}x))
= JL_{g^{-1}x}\mu(g^{-1}x)
= \nabla f(g^{-1}x).
$$
Here the last equality follows from Lemma~\ref{le:MUCRIT}.
This proves~\eqref{eq:CON1} and Step~1.

We prove Step~2. Let $g\in\rG^c$ and $\ghat\in T_g\rG^c$. Then
$$
\inner{v_x(g)}{\ghat}_g
=-\inner{\i\mu(g^{-1}x)}{g^{-1}\ghat}_{\cg^c}
=-\inner{\mu(g^{-1}x)}{\Im(g^{-1}\ghat)}_{\cg}
= \alpha_x(g)\ghat
$$
Here the first and last equations follow 
from the definitions of $v_x$ and $\alpha_x$ 
in~\eqref{eq:alphax}.  
This proves~\eqref{eq:CON2} and Step~2.

We prove Step~3.
By definition $\alpha_x$ is basic, i.e.\ it vanishes 
on the tangent vectors $\ghat=g\xi$ (for $\xi\in\cg$) 
of the group orbit $\pi(g)=g\rG$ and it is invariant under 
the right action of $\rG$ on $\rG^c$.  Hence it descends 
to a $1$-form on $M$. Since $M$ is connected and 
simply connected, it suffices to prove that $\alpha_x$ 
is closed, and for this it suffices to prove that 
the $1$-form  $g^*\alpha_x$ is closed
for every smooth map $g:\R^2\to\rG^c$.

Let $s$ and~$t$ be the standard coordinates on $\R^2$ and 
let $g:\R^2\to\rG^c$ be a smooth map. Define the functions 
$z:\R^2\to X$ and $\zeta_s,\zeta_t:\R^2\to\cg^c$ by 
$$
z:=g^{-1}x,\qquad 
\zeta_s:=g^{-1}\p_sg,\qquad
\zeta_t:=g^{-1}\p_tg.
$$
They satisfy~${\p_sz=-L_z^c\zeta_s}$,
${\p_tz = - L_z^c\zeta_t}$, 
and~${\p_t\zeta_s-\p_s\zeta_t=[\zeta_s,\zeta_t]}$.
Now denote~${\xi_s:=\Re(\zeta_s)}$, ${\eta_s:=\Im(\zeta_s)}$,
${\xi_t:=\Re(\zeta_t)}$, ${\eta_t:=\Im(\zeta_t)}$. Then
\begin{equation}\label{eq:KNF4}
\begin{split}
\p_sz &= -L_z\xi_s-JL_z\eta_s,\qquad
\p_t\xi_s-\p_s\xi_t = [\xi_s,\xi_t] - [\eta_s,\eta_t],\\
\p_tz  &= -L_z\xi_t-JL_z\eta_t,\qquad
\p_t\eta_s-\p_s\eta_t = [\xi_s,\eta_t] + [\eta_s,\xi_t].
\end{split}
\end{equation}
Thus the pullback of $\alpha_x$ under $g$ 
is the $1$-form 
$$
g^*\alpha_x = -\inner{\mu(z)}{\eta_s}\,ds 
- \inner{\mu(z)}{\eta_t}\,dt.
$$
The $1$-form  $g^*\alpha_x$ is closed if and only 
if $\p_t\inner{\mu(z)}{\eta_s}=\p_s\inner{\mu(z)}{\eta_t}$.
Indeed, 
\begin{equation*}
\begin{split}
&
\p_t\inner{\mu(z)}{\eta_s}
- \p_s\inner{\mu(z)}{\eta_t} \\
&=
\inner{\mu(z)}{\p_t\eta_s-\p_s\eta_t} 
+ \inner{d\mu(z)\p_tz}{\eta_s}
- \inner{d\mu(z)\p_sz}{\eta_t} \\
&=
\inner{\mu(z)}{[\eta_s,\xi_t]} 
- \inner{d\mu(z)(L_z\xi_t+JL_z\eta_t)}{\eta_s} \\
&\quad\;
+\,\inner{\mu(z)}{[\xi_s,\eta_t]} 
+ \inner{d\mu(z)(L_u\xi_s+JL_u\eta_s)}{\eta_t} 
= 0.
\end{split}
\end{equation*}
Here the second step follows from~\eqref{eq:KNF4}
and the last step follows from~\eqref{eq:mu3}.
Thus $\alpha_x$ is closed, as claimed,
and this proves Step~3
and Theorem~\ref{thm:CONJUGACY}.
\end{proof}

\bigbreak

\begin{thm}[{\bf Properties of the Kempf--Ness Function}]\label{thm:KNF}
\ 

\smallskip\noindent{\bf (i)}
The Kempf--Ness function
$$
\Phi_x:M\to\R
$$
is Morse--Bott and is convex along geo\-desics.

\smallskip\noindent{\bf (ii)}
The critical set of $\Phi_x$ is a (possibly empty) 
closed connected submanifold of $M$.  It is given by
\begin{equation}\label{eq:KNcrit}
\Crit(\Phi_x) = \left\{\pi(g)\in M\,|\,\mu(g^{-1}x)=0\right\}.
\end{equation}

\smallskip\noindent{\bf (iii)}
If the critical manifold of $\Phi_x$ is nonempty, then it consists 
of the absolute minima of $\Phi_x$ and every negative 
gradient flow line of $\Phi_x$ converges exponentially 
to a critical point. 

\smallskip\noindent{\bf (iv)}
Even if the critical manifold of $\Phi_x$ is empty,
every negative gradient flow line $\gamma:\R\to M$ 
of $\Phi_x$ satisfies
\begin{equation}\label{eq:KNliminf}
\lim_{t\to\infty}\Phi_x(\gamma(t)) = \inf_M\Phi_x.
\end{equation}
(The infimum may be minus infinity.)

\smallskip\noindent{\bf (v)}
The covariant Hessian of $\Phi_x$ at a point $\pi(g)\in M$ 
is the quadratic form 
$$
T_{\pi(g)}M\to\R:d\pi(g)\ghat\mapsto\abs{L_{g^{-1}x}\Im(g^{-1}\ghat)}^2.
$$

\smallskip\noindent{\bf (vi)}
Let $g:\R\to\rG^c$ be a smooth curve.
Then 
$$
\gamma:=\pi\circ g:\R\to M
$$
is a negative gradient flow line of $\Phi_x$ if 
and only if $g$ satisfies the differential equation
\begin{equation}\label{eq:KNF}
\Im(g^{-1}\dot g) = \mu(g^{-1}x).
\end{equation}

\smallskip\noindent{\bf (vii)}
The Kempf--Ness functions satisfy
$$
\Phi_{h^{-1}x}(h^{-1}g) = \Phi_x(g)-\Phi_x(h)
$$
for $x\in X$ and $g,h\in\rG^c$.

\smallskip\noindent{\bf (viii)}
Assume the critical manifold of $\Phi_x$ is nonempty.
Let $g_i\in\rG^c$ be a sequence such that 
$$
\sup_i\Phi_x(\pi(g_i))<\infty.   
$$
Then there exists a sequence~$h_i$ in the identity component~$\rG^c_{x,0}$  
of~$\rG^c_x$ such that~$h_ig_i$ has a convergent subsequence.
\end{thm}

\begin{proof}
It follows directly from the definition that
\begin{equation}\label{eq:KNgrad}
\nabla\Phi_x(\pi(g))=-d\pi(g)g\i\mu(g^{-1}x)
\end{equation}
for every $g\in\rG^c$.  Hence the negative gradient flow lines 
of $\Phi_x$ lift to solutions of equation~\eqref{eq:KNF}
and this proves part~(vi).  It also follows from~\eqref{eq:KNgrad}
that the critical set of $\Phi_x$ is given by~\eqref{eq:KNcrit}.

We compute the covariant Hessian of $\Phi_x$. 
Choose a curve $g:\R\to\rG^c$
and consider the composition 
$\gamma:=\pi\circ g:\R\to M$.  
We compute the covariant derivative 
of the vector field $\nabla\Phi_x$ 
along this curve, using  the formula for the
Levi-Civita connection on $M$ in Appendix~\ref{app:GcG}.
It is given by 
$$
\Nabla{t}\nabla\Phi_x(\pi(g)) = d\pi(g)g\zeta,\quad
\zeta := \frac{d}{dt}\bigl(-\i\mu(g^{-1}x)\bigr)
+ [\Re(g^{-1}\dot g),-\i\mu(g^{-1}x)].
$$
Thus $\zeta=\zeta(t)=\i\eta(t)$, where $\eta(t)$ is given by
\begin{equation*}
\begin{split}
\eta
&=
-\frac{d}{dt}\mu(g^{-1}x) - [\Re(g^{-1}\dot g),\mu(g^{-1}x)] \\
&=
-d\mu(g^{-1}x)\frac{d}{dt}g^{-1}x - [\Re(g^{-1}\dot g),\mu(g^{-1}x)] \\
&=
d\mu(g^{-1}x)g^{-1}\dot gg^{-1}x - [\Re(g^{-1}\dot g),\mu(g^{-1}x)] \\
&=
d\mu(g^{-1}x)L_{g^{-1}x}^cg^{-1}\dot g - [\Re(g^{-1}\dot g),\mu(g^{-1}x)] \\
&=
d\mu(g^{-1}x)L_{g^{-1}x}\Re(g^{-1}\dot g) 
+ d\mu(g^{-1}x)JL_{g^{-1}x}\Im(g^{-1}\dot g) \\
&\quad
-\, [\Re(g^{-1}\dot g),\mu(g^{-1}x)] \\
&=
{L_{g^{-1}x}}^*L_{g^{-1}x}\Im(g^{-1}\dot g).
\end{split}
\end{equation*}
Here the last equation follows from~\eqref{eq:mu3}.
Take the inner product of the vector fields~${\Nabla{t}\nabla\Phi_x(\pi(g))}$ 
and~${d\pi(g)\dot g}$ along the curve $\pi(g)$ in $M$ to obtain
$$
d^2_{\pi(g)}\Phi_x(d\pi(g)\dot g) 
= \Abs{L_{g^{-1}x}\Im(g^{-1}\dot g)}^2.
$$
This proves the formula for the Hessian of $\Phi_x$ in part~(v).

We prove that $\Crit(\Phi_x)$ is a submanifold
of $M$ with tangent spaces 
\begin{equation}\label{eq:Tcrit}
T_{\pi(g)}\Crit(\Phi_x) 
= \left\{d\pi(g)g\i\eta\,|\,\eta\in\ker\,L_{g^{-1}x}\right\}.
\end{equation}
To see this, choose an element~${g\in\rG^c}$ 
such that
$$
\mu(g^{-1}x)=0.
$$
By Hadamard's theorem the 
map~${\cg\to M:\eta\mapsto \pi(g\exp(\i\eta))}$ 
is a diffeomorphism. Moreover, for every~${\eta\in\cg}$, 
the following are equivalent.

\smallskip\noindent{\bf (a)}
$\pi(g\exp(\i\eta))\in\Crit(\Phi_x)$.

\smallskip\noindent{\bf (b)}
$\mu(\exp(-\i\eta)g^{-1}x)=0$.

\smallskip\noindent{\bf (c)}
$\exp(-\i\eta)g^{-1}x\in\rG(g^{-1}x)$.

\smallskip\noindent{\bf (d)}
$L_{g^{-1}x}\eta=0$.

\medskip\noindent
The equivalence of~(a) and~(b) follows from the 
formula~\eqref{eq:KNgrad} for the gradient 
of the Kempf--Ness function.  The equivalence of~(b)
and~(c) follows from Lemma~\ref{le:ZERO}.
The equivalence of~(c) and~(d) follows from the fact
that the isotropy subgroup~${\rG^c_{g^{-1}x}}$ is the 
complexification of~${\rG_{g^{-1}x}}$ by Lemma~\ref{le:isotropy}.
Thus we have proved that the set~${\Crit(\Phi_x)}$ is the image 
of the linear subspace~${\ker L_{g^{-1}x}\subset\cg}$
under the diffeomorphism~${\cg\to M:\eta\mapsto\pi(g\exp(\i\eta))}$.
Hence~${\Crit(\Phi_x)}$ is a closed connected submanifold 
of~$M$ with the tangent space~\eqref{eq:Tcrit} at~$\pi(g)$.  
This proves part~(ii) and, by part~(v), that~$\Phi_x$ 
is a Morse--Bott function.

Part~(v) also shows that the covariant Hessian of $\Phi_x$ is 
everywhere nonnegative.  Hence $\Phi_x$ is convex along geodesics.
Here is an alternative argument.  Let $g_0\in\rG^c$ and $\xi\in\cg$ 
and define $g:\R\to\rG^c$ and $y:\R\to X$ by 
$$
g(t) := g_0\exp(-\i t\xi),\qquad
y(t) := g(t)^{-1}x = \exp(\i t\xi)g_0^{-1}x.
$$
Then~${g^{-1}\dot g=-\i\xi}$ and~${\dot y=JL_y\xi}$.
Moreover, the curve~${\gamma:=\pi\circ g:\R\to M}$ 
is a geodesic and
\begin{equation}\label{eq:d1Phi}
\frac{d}{dt}(\Phi_x\circ\gamma)
= -\inner{\mu(g^{-1}x)}{\Im(g^{-1}\dot g)}
= \inner{\mu(y)}{\xi}.
\end{equation}
Hence, as in equation~\eqref{eq:ZERO},
\begin{equation}\label{eq:d2Phi}
\frac{d^2}{dt^2}(\Phi_x\circ\gamma)
= \frac{d}{dt}\inner{\mu(y)}{\xi}
= \inner{d\mu(y)JL_y\xi}{\xi}
= \abs{L_y\xi}^2\ge 0.
\end{equation}
This shows again that the Kempf--Ness 
function is convex along geodesics.
Thus we have proved assertions~(i), (ii), (v), (vi).

We prove parts~(iii) and~(iv).
Let~${\gamma_0,\gamma_1:\R\to M}$ 
be negative gradient flow lines of~$\Phi_x$. 
Then there exist solutions~${g_0,g_1:\R\to\rG^c}$ 
of the differential 
equation~${g_i^{-1}\dot g_i = \i\mu(g_i^{-1}x)}$
such that~${\gamma_0=\pi\circ g_0}$
and~${\gamma_1=\pi\circ g_1}$. 
Define~${\eta:\R\to\cg}$ and~${u:\R\to\rG}$ 
by
$$
g_1(t) =: g_0(t)\exp(\i\eta(t))u(t).
$$
Then the curve $\beta_t(s):=\pi(g_0(t)\exp(\i s\eta(t)))$ 
for $0\le s\le1$ is the unique geodesic 
connecting~$\gamma_0(t)$ to~$\gamma_1(t)$.  Hence
$$
\rho(t) := d_M(\gamma_0(t),\gamma_1(t)) = \abs{\eta(t)}.
$$
This function is nonincreasing by Lemma~\ref{le:gradflow}.

\bigbreak

To prove part~(iii) assume that $\mu(g_0(0)^{-1}x)=0$.
Then $\gamma_0$ is constant and it follows
that $\gamma_1([0,\infty))$ is contained in
a compact subset of $M$. Since~$\Phi_x$ is a Morse--Bott function,
this implies that $\gamma_1$ converges exponentially to a critical 
point of $\Phi_x$.  This proves part~(iii).

To prove part~(iv) we argue by contradiction and assume
that 
$$
a:=\lim_{t\to\infty}\Phi_x(\gamma_0(t)) > \inf_M\Phi_x.
$$
Then $a>-\infty$ and we can choose $\gamma_1$ such that 
\begin{equation}\label{eq:a}
\Phi_x(\gamma_1(0))<a.
\end{equation}
Since the function $\rho=\abs{\eta}:\R\to\R$ is nonincreasing,
there exists a constant $C>0$ such that 
$\abs{\eta(t)}\le C$ for all $t\ge 0$.
This implies
\begin{equation*}
\begin{split}
\left.\frac{d}{ds}\right|_{s=0}\Phi_x(\beta_t(s))
&=
d\Phi_x(\gamma_0(t))\dot\beta_t(0) \\
&= 
- \inner{\mu(g_0(t)^{-1}x)}{\eta(t)} \\
&\ge 
- \abs{\mu(g_0(t)^{-1}x)} \abs{\eta(t)} \\
&\ge 
-C\abs{\mu(g_0(t)^{-1}x)}.
\end{split}
\end{equation*}
Since the function $\Phi_x\circ\beta_t:[0,1]\to\R$ 
is convex it follows that 
\begin{equation*}
\begin{split}
\Phi_x(\gamma_1(t)) 
&=
\Phi_x(\beta_t(1))  \\
&\ge
\Phi_x(\beta_t(0)) - C\abs{\mu(g_0(t)^{-1}x)} \\
&=x
\Phi_x(\gamma_0(t)) - C\abs{\mu(g_0(t)^{-1}x)}.
\end{split}
\end{equation*}
Since $\Phi_x\circ\gamma_0$ is bounded below and 
$$
\frac{d}{dt}(\Phi_x\circ\gamma_0)(t)=-\abs{\mu(g_0(t)^{-1}x}^2,
$$
there exists a sequence $t_i\to\infty$ such that 
$\lim_{i\to\infty}\abs{\mu(g_0(t_i)^{-1}x)}^2 = 0$.
It follows that 
$$
\lim_{i\to\infty}\Phi_x(\gamma_1(t_i)) 
\ge \lim_{i\to\infty}\Phi_x(\gamma(t_i)) = a.
$$
This contradicts the assumption~\eqref{eq:a}.
Thus we have proved part~(iv).

We prove part~(vii). The $1$-forms $\alpha_x$ and 
$\alpha_{h^{-1}x}$ on $\rG^c$ in~\eqref{eq:alphax}
satisfy
$$
\alpha_x(g)\ghat
= -\inner{\mu(g^{-1}x)}{\Im(g^{-1}\ghat)}
= \alpha_{h^{-1}x}(h^{-1}g)h^{-1}\ghat
$$
for all $g\in\rG^c$ and $\ghat\in T_g\rG^c$.
Thus the pullback of the $1$-form
${\alpha_{h^{-1}x}\in\Om^1(\rG^c)}$ under the 
diffeomorphism~${\rG^c\to\rG^c:g\mapsto h^{-1}g}$
agrees with~$\alpha_x$.  Hence the pullback
of the function~${\Phi_{h^{-1}x}:\rG^c\to\R}$ under the 
same diffeomorphism differs from~$\Phi_x$
by a constant on each connected component of~$\rG^c$.  
Hence assertion~(vii) follows from the normalization 
condition~${\Phi_{h^{-1}x}(u)=\Phi_{h^{-1}x}(\one)=0}$.

\bigbreak

We prove part~(viii) in five steps. 

\medskip\noindent{\bf Step~1.}
{\it For all $g\in\rG^c$ the set
$N_g := \{\pi(hg)\,|\,h\in\rG^c_{x,0}\}\subset M$
is closed.}

\medskip\noindent
Choose a sequence $h_i\in\rG^c_{x,0}$ 
and an element $\tg\in\rG^c$ such that the 
sequence~${\pi(h_ig)}$ converges to~${\pi(\tg)}$ in~$M$. 
Then there exists a sequence~${u_i\in\rG}$
such that~$h_igu_i$ converges to~$\tg$ in~$\rG^c$.
Pass to a subsequence so that the limit~${u:=\lim_{i\to\infty}u_i}$ 
exists in~$\rG$. Then~$h_i$ converges to
$$
h:=\tg u^{-1}g^{-1} \in \rG^c_{x,0}
$$
and hence~${\pi(\tg)=\pi(hg)\in N_g}$. 
This proves Step~1.

\medskip\noindent{\bf Step~2.}
{\it If $\mu(x)=0$ then $\Phi_x$ is constant on $N_g$
for all $g\in\rG^c$.}

\medskip\noindent
It follows from~\eqref{eq:CON1} that
$$
\frac{d}{dt}\Phi_x(\exp(t\zeta))=-\inner{\mu(x)}{\Im(\zeta)}=0
$$
for $\zeta\in\ker\,L_x^c$. Hence $\Phi_x(h)=0$ for all 
$h\in\rG^c_{x,0}$ and hence, by part~(vii),
$$
\Phi_x(h^{-1}g) = \Phi_{h^{-1}x}(h^{-1}g)
= \Phi_x(g) - \Phi_x(h) = \Phi_x(g)
$$
for all $g\in\rG^c$ and all $h\in\rG^c_{x,0}$.
This proves Step~2.

\medskip\noindent{\bf Step~3.}
{\it If $\mu(x)=0$ then there is a constant $\delta>0$
such that, for every $\eta\in\cg$,}
\begin{equation}\label{eq:KNFdelta}
\eta\perp\ker L_x,\;\abs{\eta}\ge 1
\qquad\implies\qquad\Phi_x(\exp(\i\eta))\ge\delta\abs{\eta}.
\end{equation}

\medskip\noindent
If $\eta\in\cg$ satisfies 
$L_x\eta\ne0$, then $\Phi_x(\exp(\i\eta))>0$. 
This follows from the fact that the 
function~${\phi_\eta(t):=\Phi_x(\exp(\i t\eta))}$ 
is convex and satisfies 
$$
\phi_\eta(0)=0,\qquad
\dot\phi_\eta(0)=0,\qquad
\ddot\phi_\eta(0) = \abs{L_x\eta}^2>0.
$$
Now define 
$$
\delta := \min\left\{\Phi_x(\exp(\i\eta))\,|\,
\eta\in\cg,\,\eta\perp\ker\,L_x,\,\abs{\eta}=1\right\}.
$$
Then~${\delta>0}$ and~${\phi_\eta(1)\ge\delta}$ for 
every element~${\eta\in(\ker\,L_x)^\perp}$ of norm one.
Hence the inequality~\eqref{eq:KNFdelta} follows from the 
convexity of the functions $\phi_\eta$. 
This proves Step~3.

\bigbreak

\medskip\noindent{\bf Step~4.}
{\it If $\mu(x)=0$ then part~(viii) holds.}

\medskip\noindent
Choose a sequence $g_i\in\rG^c$ such 
that~${c:=\sup_i\Phi_x(g_i) < \infty}$.
Since $N_{g_i}$ is a closed subset of $M$ by Step~1, 
there exists a sequence $h_i\in\rG^c_{x,0}$ such that
\begin{equation}\label{eq:hN}
r_i := d_M(\pi(h_ig_i),\pi(\one)) 
= \inf_{h\in\rG^c_{x,0}}d_M(\pi(hg_i),\pi(\one)).
\end{equation}
Choose $\eta_i\in\cg$ and $u_i\in\rG$ such that
$$
h_ig_i=\exp(\i\eta_i)u_i.
$$
If $\xi\in\ker\,L_x$ then $\exp(-\i\xi)\in\rG^c_{x,0}$,
and hence $\pi(\exp(-\i\xi)\exp(\i\eta_i))\in N_{g_i}$.
Thus it follows from~\eqref{eq:hN} that, 
for $\xi\in\ker\,L_x$,
\begin{eqnarray*}
r_i
&=&
d_M(\pi(\exp(\i\eta_i)),\pi(\one)) \\
&\le&
d_M(\pi(\exp(-\i\xi)\exp(\i\eta_i)),\pi(\one)) \\
&=&
d_M(\pi(\exp(\i\eta_i)),\pi(\exp(\i\xi))).
\end{eqnarray*}
In other words, for every $\xi\in\ker\,L_x$,
the geodesic 
$$
\gamma_\xi(t):=\pi(\exp(\i t\xi))
$$ 
in~$M$ has minimal distance to the point~${\pi(\exp(\i\eta_i))=\pi(h_ig_i)}$
at $t=0$, and this implies $\inner{\eta_i}{\xi}=0$.
Thus~${\eta_i\perp\ker L_x}$ and, if~${\abs{\eta_i}\ge1}$, 
it follows from Step~2 and Step~3 that
$$
c\ge\Phi_x(g_i) = \Phi_x(h_ig_i)
= \Phi_x(\exp(\i\eta_i))\ge \delta\abs{\eta_i}.
$$
Hence the sequence $\eta_i$ is bounded,
and so the sequence~${h_ig_i=\exp(\i\eta_i)u_i}$ 
has a convergent subsequence.
This proves Step~4.

\medskip\noindent{\bf Step~5.}
{\it We prove part~(viii).}

\medskip\noindent
Let $x\in X$ such that the critical manifold 
of $\Phi_x$ is nonempty. Then there is a $g\in\rG^c$ 
such that $\mu(g^{-1}x)=0$.
Choose a sequence $g_i\in\rG^c$ such that the sequence
$\Phi_x(g_i)$ is bounded. Then by part~(vii) so is the sequence 
$$
\Phi_{g^{-1}x}(g^{-1}g_i)=\Phi_x(g_i)-\Phi_x(g).
$$
Hence, by Step~4, there exists a sequence 
$$
h_i\in\rG^c_{g^{-1}x,0}
$$
such that $h_ig^{-1}g_i$ has a convergent subsequence.
Thus 
$$
\widetilde{h}_i:=gh_ig^{-1}\in\rG^c_{x,0}
$$
for all $i$ and the sequence~${\widetilde{h}_ig_i=gh_ig^{-1}g_i}$ 
has a convergent subsequence. 
This proves Step~5 and Theorem~\ref{thm:KNF}.
\end{proof}


\chapter{$\mu$-Weights}\label{ch:WEIGHTS}

The purpose of the present chapter is to introduce Mumford's
numerical invariants~${w_\mu(x,\zeta)}$ associated to an 
element~${x\in X}$ and a toral generator~${\zeta\in\sT^c}$.
The main result is Mumford's Theorem~\ref{thm:MUMFORD1} which 
establishes the invariance of the $\mu$-weights under the diagonal
action of the complexified Lie group~$\rG^c$ on~${X\times\sT^c}$ 
and under Mumford's equivalence relation on~$\sT^c$. 
Toral generators and Mumford's equivalence relation 
are explained in Appendix~\ref{app:Lambda}.

Introduce\index{weight@$\mu$-weight|(} 
the notations 
$$
\sT^c:=\left\{g\xi g^{-1}\,|\,g\in\rG^c,\xi\in\cg\setminus\{0\}\right\}
$$
and
$$
\Lambda := \left\{\xi\in\cg\setminus\{0\}\,|\,\exp(\xi)=1\right\},\quad
\Lambda^c := \left\{\zeta\in\cg^c\setminus\{0\}\,|\,\exp(\zeta)=1\right\}.
$$

\begin{defn}\label{def:WEIGHT}
The {\bf $\mu$-weight of a pair $(x,\zeta)\in X\times\sT^c$}
is the real number
\begin{equation}\label{eq:WEIGHT}
w_\mu(x,\zeta) := \lim_{t\to\infty}\inner{\mu(\exp(\i t\zeta)x)}{\Re(\zeta)}.
\end{equation}
For $\zeta=\xi\in\cg\setminus\{0\}$ the existence of the limit follows 
from the fact that the function $t\mapsto\inner{\mu(\exp(\i t\xi)x)}{\xi}$
is nondecreasing by~\eqref{eq:ZERO}.  For general elements $\zeta\in\sT^c$
the existence of the limit follows from Lemma~\ref{le:WEIGHT2} below.  
\end{defn}

For $\zeta\in\Lambda^c$ the 
geometric significance of the $\mu$-weight in terms 
of a lift of the~$\rG^c$ action to a line bundle over $X$, 
under a suitable rationality hypothesis, 
is explained in Theorem~\ref{thm:LIFTc} below.
For $\xi\in\cg\setminus\{0\}$ the next lemma shows that
the $\mu$-weight~${w_\mu(x,\xi)}$ is the asymptotic slope 
of the Kempf--Ness function~$\Phi_x$ along the geodesic 
ray~${t\mapsto[\exp(-\i t\xi)]}$ as $t$ tends to~$\infty$.  

\begin{lem}\label{le:SLOPE}
Fix a point~${x\in X}$ and an element~${\xi\in\cg\setminus\{0\}}$.  
Then the function~${t\mapsto t^{-1}\Phi_x(\exp(-\i t\xi))}$ 
is nondecreasing and
\begin{equation}\label{eq:SLOPE}
w_\mu(x,\xi) = \lim_{t\to\infty}\frac{\Phi_x(\exp(-\i t\xi))}{t}.
\end{equation}
\end{lem}

\begin{proof}
By equation~\eqref{eq:d1Phi} we have 
$$
\Phi_x(\exp(-\i t\xi))=\int_0^t\inner{\mu(\exp(\i s\xi)x)}{\xi}\,ds
\qquad\mbox{for all }t>0.
$$
Hence it follows from the definition of the weight in~\eqref{eq:WEIGHT} that
$$
w_\mu(x,\xi) 
= \lim_{t\to\infty}\frac{1}{t}\int_0^t\inner{\mu(\exp(\i s\xi)x)}{\xi}\,ds
= \lim_{t\to\infty}\frac{\Phi_x(\exp(-\i t\xi))}{t}.
$$
That the function $[0,\infty)\to\R:t\mapsto t^{-1}\Phi_x(\exp(-\i t\xi))$
is nondecreasing follows from the fact that the function 
$t\to\Phi_x(\exp(-\i t\xi))$ is convex and vanishes at $t=0$.
This proves Lemma~\ref{le:SLOPE}.
\end{proof}

\begin{thm}[{\bf Mumford}]\label{thm:MUMFORD1}
{\bf (i)} 
The function 
$
w_\mu:X\times\sT^c\to\R
$ 
is $\rG^c$-invariant,
i.e.\ for all $x\in X$, $\zeta\in\sT^c$, and $g\in\rG^c$
$$
w_\mu(gx,g\zeta g^{-1}) = w_\mu(x,\zeta).
$$

\smallskip\noindent{\bf (ii)}
For every $x\in X$ the function $\sT^c\to\R:\zeta\mapsto w_\mu(x,\zeta)$
is constant on the equivalence classes in Theorem~\ref{thm:MUMFORD}, 
i.e.\ for all $\zeta\in\sT^c$ and $p\in\rP(\zeta)$
$$
w_\mu(x,p\zeta p^{-1}) = w_\mu(x,\zeta).
$$
\end{thm}

\begin{proof}
See page~\pageref{proof:MUMFORD1}.
\end{proof}

The proof of Theorem~\ref{thm:MUMFORD1} is based on 
Lemma~\ref{le:WEIGHT2} and Lemma~\ref{le:WEIGHT4} below.
The next lemma establishes the existence of the limit in~\eqref{eq:WEIGHT}.

\begin{lem}\label{le:WEIGHT2}
Let $x_0\in X$ and $\zeta\in\sT^c$.  
Then the limits
\begin{equation}\label{eq:xpm}
x^\pm:=\lim_{t\to\pm\infty}\exp(\i t\zeta)x_0
\end{equation}
exist, the convergence is exponential in $t$, 
and $L_{x^\pm}^c\zeta = 0$.
\end{lem}

\begin{proof} 
Assume first that $\zeta=\xi\in\cg\setminus\{0\}$
and define the function $x:\R\to X$ by~${x(t):=\exp(\i t\xi)x_0}$. 
Then $\dot x=JL_x\xi=\nabla H_\xi(x)$ and so
$x$ is a gradient flow line of the Hamiltonian function
$H_\xi=\inner{\mu}{\xi}$. Since $H_\xi$ is a Morse--Bott 
function the limits~${x^\pm:=\lim_{t\to\pm\infty}\exp(\i t\xi)x_0}$
exist, the convergence is exponential in $t$, 
and the limit points satisfy $L_{x^\pm}\xi=0$.
This proves the lemma for $\zeta=\xi\in\cg\setminus\{0\}$.

Now let $\zeta\in\sT^c$.  By Lemma~\ref{le:Lambda},
there is a $g\in\rG^c$ such that ${g\zeta g^{-1}\in\cg}$.
By what we have just proved the 
limits~${\tx^\pm:=\lim_{t\to\pm\infty}\exp(\i tg\zeta g^{-1})gx}$
exist, the convergence is exponential in~$t$, 
and~${L_{\tx^\pm}^c(g\zeta g^{-1})=0}$.  
Hence the function
$
\exp(\bi t\zeta)x=g^{-1}\exp(\bi tg\zeta g^{-1})gx
$
converges to~${x^\pm:=g^{-1}\tx^\pm}$
as~$t$ tends to~$\pm\infty$ 
and~${L^c_{x^\pm}\zeta = g^{-1}L^c_{\tx^\pm}(g\zeta g^{-1})=0}$.
This proves Lemma~\ref{le:WEIGHT2}.
\end{proof}

\begin{remk}\label{rmk:WEIGHT1}\rm
Let $x_0\in X$ and $\xi\in\cg\setminus\{0\}$.
Define $x(t):=\exp(\i t\xi)x_0$ as in the proof of Lemma~\ref{le:WEIGHT2} 
and let $x^\pm:=\lim_{t\to\pm\infty}\exp(\i t\xi)x_0$.  
Then, as in~\eqref{eq:ZERO}, we have
$
\frac{d}{dt}\inner{\mu(x)}{\xi}
= \abs{L_x\xi}^2 
= \abs{\dot x}^2.
$
Integrate this equation to obtain the energy identity
$$
E(x):=\int_{-\infty}^\infty \abs{\dot x(t)}^2\,dt
= \inner{\mu(x^+)}{\xi} - \inner{\mu(x^-)}{\xi}
= w_\mu(x_0,\xi)+w_\mu(x_0,-\xi).
$$
In particular, $w_\mu(x_0,\xi)+w_\mu(x_0,-\xi)\ge 0$.
\end{remk}

\begin{remk}\label{rmk:WEIGHT2}\rm
Here is another proof of Lemma~\ref{le:WEIGHT2} 
for $\zeta\in\Lambda^c$. Let $x_0\in X$ and
define the map~${z:\bbS:=\R/\Z\times\R\to X}$ by
$$
z(s,t) := \exp((s+\i t)\zeta)x_0.
$$
We prove that $z$ is a finite energy holomorphic curve and 
\begin{equation}\label{eq:Eweight}
E(z) 
:= \int_\bbS\abs{\p_sz}^2
= w_\mu(x_0,\zeta)+w_\mu(x_0,-\zeta).
\end{equation}
To see this, let $\xi:=\Re(\zeta)$ and $\eta:=\Im(\zeta)$.
Then 
\begin{equation}\label{eq:dsdtz}
\p_sz = L_z\xi+JL_z\eta,\qquad 
\p_tz = JL_z\xi-L_z\eta,
\end{equation}
so $\p_sz+J\p_tz=0$ and $z$ is holomorphic. 
It follows also from~\eqref{eq:dsdtz} that
\begin{equation*}
\begin{split}
\p_t\inner{\mu(z)}{\xi} 
&= \inner{d\mu(z)(JL_z\xi-L_z\eta)}{\xi} 
= \abs{L_z\xi}^2-\inner{\mu(z)}{[\xi,\eta]},\\
\p_s\inner{\mu(z)}{\eta} 
&= \inner{d\mu(z)(L_z\xi+JL_z\eta)}{\eta} 
= \abs{L_z\eta}^2-\inner{\mu(z)}{[\xi,\eta]}.
\end{split}
\end{equation*}
Since $\inner{\mu(z)}{[\xi,\eta]}=\om(L_z\xi,L_z\eta)
=-\inner{L_z\xi}{JL_z\eta}$, this implies
$$
\p_t\inner{\mu(z)}{\xi} + \p_s\inner{\mu(z)}{\eta} 
= \abs{L_z\xi+JL_z\eta}^2 = \abs{\p_sz}^2.
$$
Integrate this identity over $0\le s\le1$ to obtain
$$
\frac{d}{dt} \int_0^1\inner{\mu(z)}{\xi}\,ds
= \int_0^1\abs{L_z\xi+JL_z\eta}^2\,ds
= \int_0^1\abs{\p_sz}^2\,ds.
$$
Integrate this identity over $-\infty<t<\infty$ to obtain
\begin{equation}\label{eq:Eu}
E(z)
= \int_\bbS\abs{\p_sz}^2 \\
= \lim_{t\to\infty}\int_0^1\inner{\mu(z)}{\xi}\,ds
- \lim_{t\to-\infty}\int_0^1\inner{\mu(z)}{\xi}\,ds.
\end{equation}
The limits on the right exist because the function 
$t\mapsto \int_0^1\inner{\mu(z(s,t))}{\xi}\,ds$ 
is nondecreasing and bounded. Hence $z$ has finite energy.  
Hence it follows from the removable singularity 
theorem for holomorphic curves in~\cite[Theorem~4.1.2]{MS2}
that the limits 
$
x^\pm := \lim_{t\to\pm\infty}\exp((s+\i t)\zeta)x_0
$ 
exist and the convergence is uniform in $s$ and exponential in~$t$.  
Since the limit is independent of~$s$ it follows
again that~${\exp(s\zeta)x^\pm=x^\pm}$ for all~${s\in\R}$
and so~${L^c_{x^\pm}\zeta=0}$.  Moreover,
$
\lim_{t\to\pm\infty}\int_0^1\inner{\mu(z)}{\xi}\,ds
= \inner{\mu(x^\pm)}{\xi} = \pm w_\mu(x_0,\pm\zeta)
$
and hence equation~\eqref{eq:Eweight} follows from~\eqref{eq:Eu}.
\end{remk}

\begin{lem}\label{le:WEIGHT3}
{\bf (i)}
For all $\zeta=\xi+\i\eta\in\cg^c$ and all $g\in\rG^c$
\begin{equation}\label{eq:normcg}
\begin{split}
\abs{\Re(g\zeta g^{-1})}^2-\abs{\Im(g\zeta g^{-1})}^2
&= \abs{\xi}^2-\abs{\eta}^2,\\
\inner{\Re(g\zeta g^{-1})}{\Im(g\zeta g^{-1})}
&= \inner{\xi}{\eta}.
\end{split}
\end{equation}

\smallskip\noindent{\bf (ii)}
If $\zeta=\xi+\i\eta\in\sT^c$ then 
$\abs{\xi}>\abs{\eta}$ and $\inner{\xi}{\eta}=0$.
\end{lem}

\begin{proof}
For $g=u\in\rG$ equation~\eqref{eq:normcg} follows 
from the invariance of the inner product on $\cg$. 
Hence it suffices to assume~${g=\exp(\i\etahat)}$
for some element~${\etahat\in\cg}$.
Let $\zeta_0\in\cg^c$ and define the functions 
$g:\R\to\rG^c$ and $\zeta:\R\to\cg^c$ by 
$$
g(t) := \exp(\i t\etahat),\qquad
\zeta(t) := g(t)\zeta_0g(t)^{-1}.
$$
Then $\dot gg^{-1} = \i\etahat$ and $\dot\zeta = [\i\etahat,\zeta]$.
Define 
$
\xi(t):=\Re(\zeta(t))
$
and
$ 
\eta(t):=\Im(\zeta(t))
$
so that 
$$
\dot\xi = -[\etahat,\eta],\qquad \dot\eta = [\etahat,\xi].
$$
Then 
$$
\frac{d}{dt}\frac{\abs{\xi}^2}{2} 
= \langle\dot\xi,\xi\rangle 
= - \inner{[\etahat,\eta]}{\xi} 
= \inner{[\etahat,\xi]}{\eta} 
= \inner{\dot\eta}{\eta}
= \frac{d}{dt}\frac{\abs{\eta}^2}{2} 
$$
and 
$$
\frac{d}{dt}\inner{\xi}{\eta} 
= \langle\dot\xi,\eta\rangle + \langle\xi,\dot\eta\rangle
= -\inner{[\etahat,\eta]}{\eta} + \inner{\xi}{[\etahat,\xi]}
= 0.
$$
Hence the functions $t\mapsto\abs{\xi(t)}^2-\abs{\eta(t)}^2$ and 
$t\mapsto\inner{\xi(t)}{\eta(t)}$ are constant.
This proves part~(i). 

\bigbreak

Now let $\zeta=\xi+\i\eta\in\sT^c$.
By Lemma~\ref{le:Lambda} there exists an element 
$g\in\rG^c$ such that 
$
g\zeta g^{-1}\in\cg\setminus\{0\}.
$
Hence $\Im(g\zeta g^{-1})=0$ and it follows from part~(i) that 
$
\abs{\xi}^2-\abs{\eta}^2=\abs{g\zeta g^{-1}}^2>0
$
and
$
\inner{\xi}{\eta}=0.
$
This proves part~(ii) and Lemma~\ref{le:WEIGHT3}.
\end{proof}

\begin{lem}\label{le:WEIGHT4}
{\bf (i)} If $x\in X$ and $\zeta\in\cg^c$ satisfy
$L_x^c\zeta=0$ then, for all $g\in\rG^c$,
\begin{equation}\label{eq:mucg}
\begin{split}
\inner{\mu(gx)}{\Re(g\zeta g^{-1})}
= \inner{\mu(x)}{\Re(\zeta)},\\
\inner{\mu(gx)}{\Im(g\zeta g^{-1})}
= \inner{\mu(x)}{\Im(\zeta)}.
\end{split}
\end{equation}

\smallskip\noindent{\bf (ii)}
If $x\in X$ and $\zeta\in\sT^c$ satisfy
$L_x^c\zeta=0$ then $\inner{\mu(x)}{\Im(\zeta)}=0$. 
\end{lem}

\begin{proof}
For $g=u\in\rG$ equation~\eqref{eq:mucg} follows from the 
invariance of the inner product on $\cg$ and the 
$\rG$-equivariance of the moment map. Hence it suffices
to assume that~${g=\exp(\i\etahat)}$ for sme~${\etahat\in\cg}$.
Fix two elements $x_0\in X$ and $\zeta_0\in\cg^c$, 
define the functions $g:\R\to\rG^c$, $x:\R\to X$, 
and $\zeta:\R\to\cg^c$ by 
$$
g(t) := \exp(\i t\etahat),\qquad
x(t) := g(t)x_0,\qquad
\zeta(t) := g(t)\zeta_0g(t)^{-1},
$$
and denote $\xi(t):=\Re(\zeta(t))$ and $\eta(t):=\Im(\zeta(t))$.
Then $L_x\xi+JL_x\eta=0$ and, as in the proof 
of Lemma~\ref{le:WEIGHT3},
$$
\dot\xi = -[\etahat,\eta],\qquad 
\dot\eta = [\etahat,\xi],\qquad
\dot x = JL_x\etahat.
$$
Hence, by equations~\eqref{eq:mu1} and~\eqref{eq:mu2},
\begin{equation*}
\begin{split}
\frac{d}{dt}\inner{\mu(x)}{\xi}
&=
\inner{d\mu(x)\dot x}{\xi} + \langle \mu(x),\dot\xi\rangle \\
&= 
\om(L_x\xi,\dot x) + \langle \mu(x),\dot\xi\rangle  \\
&= 
-\om(JL_x\eta,\dot x)  + \langle \mu(x),\dot\xi\rangle  \\
&= 
-\om(JL_x\eta,JL_x\etahat)  - \inner{\mu(x)}{[\etahat,\eta]}  \\
&=
\om(L_x\etahat,L_x\eta) - \inner{\mu(x)}{[\etahat,\eta]} \\
&=
0
\end{split}
\end{equation*}
and 
\begin{equation*}
\begin{split}
\frac{d}{dt}\inner{\mu(x)}{\eta}
&= 
\inner{d\mu(x)\dot x}{\eta} + \inner{\mu(x)}{\dot\eta} \\
&= 
\om(L_x\eta,\dot x) + \inner{\mu(x)}{\dot\eta}  \\
&= 
\om(JL_x\xi,\dot x)  + \inner{\mu(x)}{\dot\eta}  \\
&= 
\om(JL_x\xi,JL_x\etahat) + \inner{\mu(x)}{[\etahat,\xi]}    \\
&=
\om(L_x\xi,L_x\etahat) - \inner{\mu(x)}{[\xi,\etahat]}   \\
&=
0.
\end{split}
\end{equation*}
Hence the functions $t\mapsto\inner{\mu(x(t))}{\xi(t)}$ and 
$t\mapsto\inner{\mu(x(t))}{\eta(t)}$ are constant. 
This proves~(i).  To prove~(ii), let $x\in X$
and $\zeta\in\sT^c$ such that ${L_x^c\zeta=0}$.  
Then it follows from Lemma~\ref{le:Lambda} that 
there exists an element~${g\in\rG^c}$ 
such that~${g\zeta g^{-1}\in\cg}$.  By~(i) this implies
${\inner{\mu(x)}{\Im(\zeta)}=\inner{\mu(gx)}{\Im(g\zeta g^{-1})}=0}$.
This proves Lemma~\ref{le:WEIGHT4}.
\end{proof}

\begin{proof}[Proof of Theorem~\ref{thm:MUMFORD1}]
\label{proof:MUMFORD1}
Let $x\in X$, $\zeta\in\sT^c$, and $g\in\rG^c$.
By Lemma~\ref{le:WEIGHT2}, the limit
$$
x^+ := \lim_{t\to\infty}\exp(\i t\zeta)x
$$
exists and satisfies
$$
L^c_{x+}\zeta=0.
$$
Moreover, 
$$
gx^+ = \lim_{t\to\infty}\exp(\i tg\zeta g^{-1})gx.
$$ 
Hence, by Definition~\ref{def:WEIGHT} and Lemma~\ref{le:WEIGHT4}, 
$$
w_\mu(x,\zeta) = \inner{\mu(x^+)}{\Re(\zeta)} 
= \inner{\mu(gx^+)}{\Re(g\zeta g^{-1})}
= w_\mu(gx,g\zeta g^{-1}).
$$
This proves part~(i) of Theorem~\ref{thm:MUMFORD1}.   

Now assume, in addition, that $g\in\rP(\zeta)$. 
Then the limit
$$
g^+ := \lim_{t\to\infty}\exp(\i t\zeta)g\exp(-\i t\zeta)
$$
exists in~$\rG^c$.  It satisfies 
$
\exp(\i s\zeta)g^+\exp(-\i s\zeta)=g^+
$
for $s\in\R$.  Differentiate this equation to obtain
$\zeta g^+-g^+\zeta=0$ and hence
$$
\zeta=g^+\zeta(g^+)^{-1}.
$$
Moreover,
\begin{equation*}
\begin{split}
g^+x^+ 
&= 
\lim_{t\to\infty}\exp(\i t\zeta)g\exp(-\i t\zeta)\cdot 
\lim_{t\to\infty}\exp(\i t\zeta)x \\
&=
\lim_{t\to\infty}\exp(\i t\zeta)gx.
\end{split}
\end{equation*}
Hence it follows from Definition~\ref{def:WEIGHT}  
and part~(i) of Theorem~\ref{thm:MUMFORD1}
(already proved) that
\begin{eqnarray*}
w_\mu(x,g^{-1}\zeta g) 
&=& 
w_\mu(gx,\zeta) \\
&=& 
\inner{\mu(g^+x^+)}{\Re(\zeta)} \\
&=& 
\inner{\mu(g^+x^+)}{\Re(g^+\zeta(g^+)^{-1})} \\
&=&  
\inner{\mu(x^+)}{\Re(\zeta)} \\
&=& 
w_\mu(x,\zeta).
\end{eqnarray*}
Here the fourth equality follows from Lemma~\ref{le:WEIGHT4}
and the fact that $L^c_{x^+}\zeta=0$.
Since~${P(\zeta)}$ is a group (take $p=g^{-1}$) 
this proves\index{weight@$\mu$-weight|)}
part~(ii) of Theorem~\ref{thm:MUMFORD1}.
\end{proof}


\chapter{The moment-weight inequality}\label{ch:MW1}

This chapter is devoted to the proof of the moment-weight inequality
which relates the Mumford numerical invariants~$w_\mu(x,\xi)$
to the norm of the moment map on the complexified group
orbit of~$x$.  We begin by proving the moment-weight inequality 
in a special case and hope that the proof might be of some interest 
in its own right.  The general moment-weight inequality is proved 
in Theorem~\ref{thm:MW2} below.

\begin{thm}[{\bf Restricted Moment-Weight Inequality}]\label{thm:MW1}
Let~${x\in X}$\index{moment-weight inequality!restricted}
and fix a toral generator~${\zeta=\xi+\i\eta\in\sT^c}$ 
such that
\begin{equation}\label{eq:Lx0}
L_x\xi+JL_x\eta=0. 
\end{equation}
Then $\inner{\mu(x)}{\eta}=0$, $\abs{\xi}>\abs{\eta}$, 
and
\begin{equation}\label{eq:MW1}
\frac{\inner{\mu(x)}{\xi}^2}{\abs{\xi}^2-\abs{\eta}^2} 
\le \abs{\mu(gx)}^2\qquad
\mbox{for all }g\in\rG^c.
\end{equation}
\end{thm}

\begin{proof}
The equation~${\inner{\mu(x)}{\eta}=0}$ was proved in Lemma~\ref{le:WEIGHT4} 
and the inequality~${\abs{\xi}>\abs{\eta}}$ was proved in Lemma~\ref{le:WEIGHT3}. 
This shows that the quotient on the left in~\eqref{eq:MW1} is well defined. 
The estimate~\eqref{eq:MW1} holds obviously when~${\mu(x)=0}$. 
So assume~${\mu(x)\ne0}$. Under this assumption we 
prove the inequality~\eqref{eq:MW1} in three steps.

\medskip\noindent{\bf Step~1.}
{\it The inequality~\eqref{eq:MW1} holds when 
$L_x\xi=L_x\eta=0$ and $g=1$.}

\medskip\noindent
By assumption and equation~\eqref{eq:mu3},
$$
[\mu(x),\xi] = -d\mu(x)L_x\xi=0,\qquad
[\mu(x),\eta] = -d\mu(x)L_x\eta=0.
$$
Thus $\zeta=\xi+\i\eta$ commutes with $\mu(x)$.
Hence $\zeta-\lambda\mu(x)\in\sT^c$ for all $\lambda\in\R$
(such that $\zeta\ne\lambda\mu(x)$).
By part~(ii) of Lemma~\ref{le:WEIGHT3}, this implies
$$
\abs{\xi-\lambda\mu(x)}^2\ge\abs{\eta}^2
$$
for all $\lambda\in\R$. 
The term on the left is minimized at
$
\lambda=\abs{\mu(x)}^{-2}\inner{\mu(x)}{\xi}.
$
Hence
$
\abs{\xi}^2 - \abs{\mu(x)}^{-2}\inner{\mu(x)}{\xi}^2
\ge \abs{\eta}^2
$
and this is equivalent to~\eqref{eq:MW1}.

\medskip\noindent{\bf Step~2.}
{\it The inequality~\eqref{eq:MW1} holds 
when $L_x^c\zeta=0$ and $g=1$.}

\medskip\noindent
Choose $x_0\in X$ and 
$\zeta_0=\xi_0+\i\eta_0\in\sT^c$
such that $L_{x_0}\xi_0+JL_{x_0}\eta_0=0$. 
Let~${x:\R\to X}$ be the solution of~\eqref{eq:KN1}
and let $g:\R\to\rG^c$ be the solution of~\eqref{eq:KN2}
so that 
$
g(t)^{-1}\dot g(t)=\i\mu(x(t))
$
and
$
x(t)=g(t)^{-1}x_0
$ 
for every $t\in\R$. Define the function $\zeta:\R\to\sT^c$
by
$$
\zeta(t) := \xi(t)+\i\eta(t) := g(t)^{-1}\zeta_0g(t).
$$
Then
$$
L_x\xi+JL_x\eta=0,\qquad
\dot\zeta=[\zeta,\i\mu(x)],\qquad
\dot\xi = -[\eta,\mu(x)],\qquad
\dot\eta= [\xi,\mu(x)].
$$
By equation~\eqref{eq:mu3}, this implies
\begin{equation*}
\begin{split}
0 
&= 
d\mu(x)(L_x\xi+JL_x\eta)
= 
-[\mu(x),\xi] + {L_x}^*L_x\eta,\\
0 
&= 
d\mu(x)(-L_x\eta+JL_x\xi)
= 
[\mu(x),\eta] + {L_x}^*L_x\xi.
\end{split}
\end{equation*}
Thus~${[\mu(x),\xi]={L_x}^*L_x\eta}$
and~${[\mu(x),\eta]=-{L_x}^*L_x\xi}$, and hence
Thus
\begin{equation*}
\begin{split}
\frac{d}{dt}\frac{\abs{\xi}^2}{2}
&= \langle \dot\xi,\xi \rangle
= - \inner{[\eta,\mu(x)]}{\xi}
= - \inner{\eta}{[\mu(x),\xi]}
= - \abs{L_x\eta}^2,\\
\frac{d}{dt}\frac{\abs{\eta}^2}{2}
&= \langle \dot\eta,\eta \rangle
= \inner{[\xi,\mu(x)]}{\eta}
= \inner{\xi}{[\mu(x),\eta]}
= - \abs{L_x\xi}^2.
\end{split}
\end{equation*}
This shows that the integral 
$\int_0^\infty(\abs{L_x\xi}^2+\abs{L_x\eta}^2)\,dt$
is finite.  Hence there exists a sequence $t_i\to\infty$ 
such that $L_{x(t_i)}\xi(t_i)$ and $L_{x(t_i)}\eta(t_i)$ 
converge to zero.  Passing to a subsequence, 
if necessary, we may assume that the limits 
$$
\xi_\infty:=\lim_{i\to\infty}\xi(t_i),\qquad
\eta_\infty:=\lim_{i\to\infty}\eta(t_i),\qquad
x_\infty:=\lim_{i\to\infty}x(t_i)
$$
exist. These limits satisfy 
$
L_{x_\infty}\xi_\infty=L_{x_\infty}\eta_\infty=0.
$
Since each set of toral generators with given eigenvalues 
and multiplicities is a closed subset of $\cg^c$, we also 
have $\xi_\infty+\i\eta_\infty\in\sT^c$.  
Hence it follows from Lemma~\ref{le:WEIGHT3}, 
Lemma~\ref{le:WEIGHT4}, and Step~1 that
$$
\frac{\inner{\mu(x_0)}{\xi_0}^2}{\abs{\xi_0}^2-\abs{\eta_0}^2} 
= \frac{\inner{\mu(x_\infty)}{\xi_\infty}^2}
{\abs{\xi_\infty}^2-\abs{\eta_\infty}^2} 
\le \abs{\mu(x_\infty)}^2
\le \abs{\mu(x_0)}^2.
$$
This proves Step~2. 

\bigbreak

\medskip\noindent{\bf Step~3.}
{\it The inequality~\eqref{eq:MW1} holds when $L_x^c\zeta=0$.}

\medskip\noindent
Let $x\in X$, $\zeta=\xi+\i\eta\in\sT^c$, and $g\in\rG^c$ 
be as in the hypotheses of the theorem. Then, 
by Lemma~\ref{le:WEIGHT3}, Lemma~\ref{le:WEIGHT4}, 
and Step~2,
\begin{equation*}
\begin{split}
\frac{\inner{\mu(x)}{\xi}^2}{\abs{\xi}^2-\abs{\eta}^2} 
= 
\frac{\inner{\mu(gx)}{\Re(g\zeta g^{-1}}^2}
{\abs{\Re(g\zeta g^{-1})}^2-\abs{\Im(g\zeta g^{-1})}^2} 
\le 
\abs{\mu(gx)}^2.
\end{split}
\end{equation*}
This proves Step~3 and Theorem~\ref{thm:MW1}. 
\end{proof}

\begin{cor}[{\bf Kirwan--Ness Inequality}]\label{cor:KIRNESS1}
Let $x\in X$ be a critical point of the moment map squared.  
Then\index{Kirwan--Ness inequality}
\begin{equation}\label{eq:KIRNESS}
\abs{\mu(x)} \le \abs{\mu(gx)}
\end{equation}
for all~${g\in\rG^c}$.
\end{cor}

\begin{proof}
Assume $\mu(x)\ne 0$. Then $\xi:=\mu(x)\in\sT^c$ and $L_x\xi=0$.
Hence it follows from Theorem~\ref{thm:MW1} that
\begin{equation*}
\begin{split}
\abs{\mu(x)}^2 
= 
\frac{\inner{\mu(x)}{\xi}^2}{\abs{\xi}^2} 
\le
\abs{\mu(gx)}^2
\end{split}
\end{equation*}
for all $g\in\rG^c$. This proves Corollary~\ref{cor:KIRNESS1}.
\end{proof}

The inequality~\eqref{eq:KIRNESS} is implicitly 
contained in the work of Kirwan~\cite{KIRWAN}. 
For linear actions on projective space 
it is proved in Ness~\cite[Theorem~6.2]{NESS2}.
The Kirwan--Ness inequality implies that the Hessian of 
$$
f=\tfrac{1}{2}\abs{\mu}^2
$$ 
is nonnegative on the subspace
$$
\im\,JL_x\subset T_xX
$$ 
for every critical point of $f$. This is equivalent to 
Lijing Wang's inequality in Lemma~\ref{le:WANG}.

\begin{thm}[{\bf First Ness Uniqueness Theorem}]\label{thm:NESS1}
Let~${x_0,x_1\in X}$\index{Ness Uniqueness Theorem} 
be critical points of the moment map squared.
Then
$$
x_1\in\rG^c(x_0)\qquad\implies\qquad x_1\in\rG(x_0).
$$ 
\end{thm}

For linear actions on projective space this is 
Theorem~7.1 in Ness~\cite{NESS2}.  Her proof
carries over to the general case and is reproduced below
as proof~1.  Proof~2 is due to Calabi--Chen~\cite[Corollary~4.1]{CCII}
who used this argument to establish a uniqueness result 
for extremal K\"ahler metrics.  The first proof uses the 
Kirwan--Ness inequality in Corollary~\ref{cor:KIRNESS1}, 
the second proof does not.

\begin{proof}[Proof~1]
The proof has three steps.  By Lemma~\ref{le:ZERO}  
we may assume that $\mu(x_0)$ and~$\mu(x_1)$ are nonzero. 

\medskip\noindent{\bf Step~1.}
{\it We may assume without loss of generality 
that there exists an element $q\in\rG^c$ such that}
$$
x_1=qx_0,\qquad
q\mu(x_0)q^{-1}=\mu(x_0).
$$

\medskip\noindent
Choose an element $g\in\rG^c$ such that 
$$
x_1=gx_0.
$$
By Theorem~\ref{thm:BOREL} there exists an
element~${p\in\rP(-\mu(x_0))}$ such that
$$
u:=pg^{-1}\in\rG
$$
and hence $g=u^{-1}p$. Assume without loss of generality 
that $u=\one$. (Replace $x_1$ by $ux_1$ if necessary.)
Thus~${x_1=px_0}$ and the limit
$$
q:=\lim_{t\to\infty}\exp(-\i t\mu(x_0))p\exp(\i t\mu(x_0))
$$
exists in $\rG^c$ and commutes with $\mu(x_0)$.  
Define the function $x:\R\to X$ by 
$$
x(t) := \exp(-\i t\mu(x_0))x_1.
$$
Since $L_{x_0}\mu(x_0)=0$, we have 
$$
\lim_{t\to\infty}x(t) 
= \lim_{t\to\infty}\exp(-\i t\mu(x_0))p\exp(\i t\mu(x_0))x_0
= qx_0
$$
and hence 
\begin{eqnarray*}
w_\mu(x_1,-\mu(x_0)) 
&=& 
\lim_{t\to\infty}\inner{\mu(x(t))}{-\mu(x_0)} \\
&=& 
- \inner{\mu(qx_0)}{\mu(x_0)} \\
&=& 
- \inner{\mu(qx_0)}{\Re(q\mu(x_0)q^{-1})} \\
&=&
-\abs{\mu(x_0)}^2.
\end{eqnarray*}
The last equation uses Lemma~\ref{le:WEIGHT4}.
Since $t\mapsto\inner{\mu(x(t))}{-\mu(x_0)}$ 
is nondecreasing and $x(0)=x_1$, this implies
$-\inner{\mu(x_1)}{\mu(x_0)}\le-\abs{\mu(x_0)}^2$
and thus
$$
\inner{\mu(x_0)}{\mu(x_1)}
\ge\abs{\mu(x_0)}^2=\abs{\mu(x_1)}^2.
$$
(Here the last equation follows from Corollary~\ref{cor:KIRNESS1}.)
Hence $\mu(x_0)=\mu(x_1)$.  Since $L_{x_1}\mu(x_1)=0$
it follows that the function $x(t)$ is constant
and hence 
$$
x_1=x(0)=\lim_{t\to\infty}x(t)=qx_0.
$$
This proves Step~1.

\bigbreak

\medskip\noindent{\bf Step~2.}
{\it We may assume without loss of generality 
that there exists an element $\eta\in\cg$ such that
$x_1=\exp(\i\eta)x_0$ and $[\eta,\mu(x_0)]=0$.}

\medskip\noindent
Let $q\in\rG^c$ be as in Step~1 and choose elements
$u\in\rG$ and $\eta\in\cg$ such that 
$$
q=\exp(\i\eta)u.
$$
Since 
$
q\mu(x_0)q^{-1}=\mu(x_0)
$ 
it follows from part~(ii) of Lemma~\ref{le:isotropy} that 
$$
u\mu(x_0)u^{-1}=\mu(x_0),\qquad [\eta,\mu(x_0)]=0.
$$ 
Replacing $x_1$ by $u^{-1}x_1$ and $\eta$ by $u^{-1}\eta u$, 
if necessary, we may assume without loss of generality 
that $u=\one$.  This proves Step~2.

\medskip\noindent{\bf Step~3.}
{\it We prove Theorem~\ref{thm:NESS1}.}

\medskip\noindent
Let $\eta\in\cg$ be as in Step~2.
Define the function $x:\R\to X$ by 
$$
x(t) := \exp(\i t\eta)x_0.
$$
Then 
$
\dot x=JL_x\eta
$ 
and, as in equation~\eqref{eq:ZERO},
\begin{equation}\label{eq:ness1}
\frac{d}{dt}\inner{\mu(x)}{\eta} 
= \inner{d\mu(x)\dot x}{\eta}
= \abs{L_x\eta}^2 \ge 0.
\end{equation}
Define the vector field $\xhat(t)\in T_{x(t)}X$ along $x$ by 
$$
\xhat(t):=v_{\mu(x_0)}(x(t))=L_{x(t)}\mu(x_0).
$$
Since $[\eta,\mu(x_0)]=0$ by Step~2,
the Lie bracket $[v_\eta,v_{\mu(x_0)}]$ vanishes,
and hence it follows from~\eqref{eq:Jvxi} that
$$
\Nabla{t}\xhat = \Nabla{\dot x}v_{\mu(x_0)}(x)
= J\Nabla{v_\eta(x)}v_{\mu(x_0)}(x)
= J\Nabla{\xhat}v_\eta(x) ,\qquad
\xhat(0)=0.
$$
Hence 
$
L_{x(t)}\mu(x_0)=0
$ 
for all $t$ and hence 
$$
\frac{d}{dt}\inner{\mu(x)}{\mu(x_0)}
= \inner{d\mu(x)JL_{x}\eta}{\mu(x_0)}
= \inner{L_x\eta}{L_x\mu(x_0)}
= 0.
$$ 
Since $x(0)=x_0$ and $x(1)=x_1$ this implies
\begin{equation}\label{eq:ness2}
\inner{\mu(x_1)}{\mu(x_0)} = \abs{\mu(x_0)}^2 = \abs{\mu(x_1)}^2.
\end{equation}
Here the last equation follows from Corollary~\ref{cor:KIRNESS1}.
It follows from~\eqref{eq:ness2} that 
$$
\mu(x_1)=\mu(x_0)
$$
and hence
$
\inner{\mu(x_1)}{\eta}=\inner{\mu(x_0)}{\eta}.
$
By~\eqref{eq:ness1} this implies
$$
L_{x(t)}\eta=0
$$
for all~$t$, hence~${L_{x_0}\eta=0}$, and hence~${x_1=x_0}$.
This proves Theorem~\ref{thm:NESS1}.
\end{proof}

\begin{proof}[Proof~2]
Choose $g_0\in\rG^c$ such that
\begin{equation}\label{eq:unique1}
x_1=g_0^{-1}x_0.
\end{equation}
Since $x_0$ and $x_1$ are critical points of the moment map 
squared, they satisfy
\begin{equation}\label{eq:unique2}
L_{x_0}\mu(x_0)=0,\qquad L_{x_1}\mu(x_1)=0.
\end{equation}
Define the curves $g,\tg:\R\to\rG^c$ by
$$
g(t) := \exp(\i t\mu(x_0)),\qquad 
\tg(t) := g_0\exp(\i t\mu(x_1)).
$$
Thus the curves $\gamma:=\pi\circ g:\R\to M$ and 
$\tga:=\pi\circ\tg:\R\to M$ are both geodesics.
Moreover, $g(t)^{-1}x_0=x_0$ and 
$\tg(t)^{-1}x_0=x_1$ for all $t$ by~\eqref{eq:unique1}
and~\eqref{eq:unique2}.  Thus $g$ and $\tg$ satisfy 
the differential equation $g^{-1}\dot g=\i\mu(g^{-1}x_0)$.
Hence it follows from Theorem~\ref{thm:KNF}
that $\gamma$ and $\tga$ are negative gradient flow lines 
of the Kempf--Ness function~$\Phi_{x_0}$.
Now define $\eta(t)\in\cg$ and $u(t)\in\rG$ 
by~${g(t)\exp(\i\eta(t))u(t) := \tg(t)}$. Then 
\begin{equation}\label{eq:unique3}
x_1 = \tg(t)^{-1}g(t)x_0 
= u(t)^{-1}\exp(-\i\eta(t))x_0 
\end{equation}
and, by Theorem~\ref{thm:GcG},
$$
\rho(t) := d_M\left(\gamma(t),\tga(t)\right) =  \abs{\eta(t)}
$$
for all $t$. By Lemma~\ref{le:gradflow} this function is nonincreasing. 
If $\rho\equiv0$ then $x_1\in\rG(x_0)$ by assumption.  
Hence assume $\rho\not\equiv0$ and, for each $t$, 
denote by $\gamma(\cdot,t):[0,1]\to M$ the unique geodesic 
from $\gamma(0,t)=\gamma(t)$ to $\gamma(1,t)=\tga(t)$.  
Then 
$$
\gamma(s,t)=\pi(g(t)\exp(\i s\eta(t)))
\qquad\mbox{for }0\le s\le 1.
$$
Hence equation~\eqref{eq:rhodot} in Lemma~\ref{le:gradflow} 
asserts that
\begin{equation}\label{eq:unique4}
\begin{split}
\dot\rho(t) 
&= - \frac{1}{\rho(t)}\int_0^1\frac{\p^2}{\p s^2}\Phi_{x_0}(g(t)\exp(\i s\eta(t))\,ds \\
&= - \frac{1}{\rho(t)}\int_0^1\abs{L_{\exp(-\i s\eta(t))x_0}\eta(t)}^2\,ds.
\end{split}
\end{equation}
Choose a sequence $t_i\to\infty$ such that the limits
$$
\lim_{i\to\infty}\dot\rho(t_i) = 0,\qquad
\lim_{i\to\infty}\eta(t_i) =: \eta_\infty,\qquad
\lim_{i\to\infty}u(t_i) =: u_\infty
$$
exist.  Then $L_{x_0}\eta_\infty=0$ by~\eqref{eq:unique4}.  
Hence it follows from~\eqref{eq:unique3} that 
$$
x_1 = \lim_{i\to\infty} u(t_i)^{-1}\exp(-\i\eta(t_i))x_0
= u_\infty^{-1}\exp(-\i\eta_\infty)x_0
= u_\infty^{-1}x_0.
$$
Thus $x_1\in\rG(x_0)$ and this completes 
the second proof of Theorem~\ref{thm:NESS1}.
\end{proof}

\begin{thm}[{\bf Moment Limit Theorem}]\label{thm:MLT}
Let $x_0\in X$\index{Moment Limit Theorem} 
and let $x:\R\to X$ be the solution of~\eqref{eq:KN1}. 
Define $x_\infty := \lim_{t\to\infty}x(t)$.
Then 
$$
\abs{\mu(x_\infty)} = \inf_{g\in\rG^c}\abs{\mu(gx_0)}.
$$
Moreover, the $\rG$-orbit of $x_\infty$ depends 
only on the $\rG^c$-orbit of $x_0$.
\end{thm}

\begin{proof}
The limit $x_\infty$ exists by Theorem~\ref{thm:XINFTY}. 
Moreover, by Lemma~\ref{le:GRADFLOW}, 
the solution~${x:\R\to X}$ of equation~\eqref{eq:KN1} 
is given by 
$$
x(t)=g(t)^{-1}x_0, 
$$
where~${g:\R\to\rG^c}$ 
is the solution of~\eqref{eq:KN2}. Fix an element~${g_0\in\rG^c}$ 
and let~${\tx:\R\to X}$ and~${\tg:\R\to\rG^c}$ be the solutions 
of the differential equations
$$
\dot\tx = -JL_\tx\mu(\tx),\qquad \tx(0)=g_0^{-1}x_0,
$$
and
$$
\tg^{-1}\dot\tg = \i\mu(\tx),\qquad \tg(0)=g_0.
$$
Define $\eta:\R\to\cg$ and $u:\R\to\rG$ by 
$$
\tg(t) =: g(t)\exp(\i\eta(t))u(t).
$$
Then, by Lemma~\ref{le:GRADFLOW}, 
$$
\tx(t) = \tg(t)^{-1}x_0 = u(t)^{-1}\exp(-\i\eta(t))x(t)
$$ 
for all $t\in\R$.
Denote by $d_M:M\times M\to[0,\infty)$ the distance
function of the Riemannian metric on the homogeneous space
$M$.  By Theorem~\ref{thm:KNF} the curves~${\gamma:=\pi\circ g}$
and~${\tga:=\pi\circ\tg}$ are gradient flow lines of 
the Kempf--Ness function $\Phi_{x_0}:M\to\R$.
Theorem~\ref{thm:KNF} also asserts that
$\Phi_{x_0}$ is convex along geodesics.
Since $M$ is simply connected with 
nonpositive sectional curvature, the function
${\R\to\R:t\mapsto d_M(\gamma(t),\tga(t)) = \abs{\eta(t)}}$
is nonincreasing (see Lemma~\ref{le:gradflow}).  
Hence there exists a sequence $t_\nu\to\infty$ 
such that the limits
$$
\eta_\infty := \lim_{\nu\to\infty}\eta(t_\nu),\qquad
u_\infty := \lim_{\nu\to\infty}u(t_\nu)
$$ 
exist. Hence
$$
\tx_\infty 
:=
\lim_{t\to\infty}\tx(t) 
= \lim_{t\to\infty}u(t)^{-1}\exp(-\i\eta(t))x(t)
= u_\infty^{-1}\exp(-\i\eta_\infty)x_\infty.
$$
This shows that $x_\infty$ and $\tx_\infty$ are critical points of 
the moment map squared belonging to the same $\rG^c$-orbit.
Hence they belong to the same $\rG$-orbit 
by Theorem~\ref{thm:NESS1}, and hence
${
\abs{\mu(x_\infty)} 
= \abs{\mu(\tx_\infty)} 
\le \abs{\mu(g_0^{-1}x_0)}}.
$
This proves Theorem~\ref{thm:MLT}.
\end{proof}

\begin{thm}[{\bf Second Ness Uniqueness Theorem}]\label{thm:NESS2}
Let $x_0\in X$ and\index{Ness Uniqueness Theorem} 
$$
m:=\inf_{g\in\rG^c}\abs{\mu(gx_0)}.  
$$
Then 
$$
x,y\in\overline{\rG^c(x_0)},\quad
\abs{\mu(x)}=\abs{\mu(y)}=m
\qquad\implies\qquad y\in\rG(x).
$$ 
\end{thm}

For linear actions on projective space (and for $x,y\in\rG^c(x_0)$) 
this is Theorem~6.2~(ii) in Ness~\cite{NESS2}.  
For $m=0$ it is Theorem~4.1 in Chen--Sun~\cite{CS}. 
Theorem~\ref{thm:NESS2} shows that every element 
$x\in\overline{\rG^c(x_0)}$ with 
$$
\abs{\mu(x)}=m
$$ 
is a critical point of the moment map squared and is 
the limit point of a negative gradient flow line 
of the moment map squared in $\rG^c(x_0)$.
In our proof of Theorem~\ref{thm:NESS2} 
we follow the argument of Chen--Sun.  

\begin{proof}
Let $x_0\in X$, let $x:\R\to X$ be the unique solution 
of~\eqref{eq:KN1}, and define 
$$
x_\infty:=\lim_{t\to\infty}x(t).
$$  
Then 
$$
x_\infty\in\overline{\rG^c(x_0)},\qquad
\abs{\mu(x_\infty)}=m
$$ 
by Theorem~\ref{thm:MLT}.  
Now let $x\in\overline{\rG^c(x_0)}$ such that $\abs{\mu(x)}=m$.
We must prove that $x\in\rG(x_\infty)$.
To see this, choose a sequence $g_i\in\rG^c$ such that
$$
x=\lim_{i\to\infty}g_i^{-1}x_0
$$
and define $y_i:\R\to X$ and $x_i\in X$ by
$$
\dot y_i = -JL_{y_i}\mu(y_i),\qquad 
y_i(0) = g_i^{-1}x_0,\qquad
x_i:=\lim_{t\to\infty}y_i(t).
$$
Then it follows from the estimate~\eqref{eq:Loja1}
in Theorem~\ref{thm:XINFTY} that there exists 
a constant $c>0$ such that, for $i$ sufficiently large,
$$
d(x_i,g_i^{-1}x_0) 
\le \int_0^\infty\abs{\dot y_i(t)}\,dt 
\le c\left(\abs{\mu(g_i^{-1}x_0)}^2-m^2\right)^{1-\alpha}.
$$
Since 
$$
m=\abs{\mu(x)}=\lim_{i\to\infty}\abs{\mu(g_i^{-1}x_0)},
$$ 
this implies~${x = \lim_{i\to\infty}x_i}$.
Moreover, we have
$$
x_i\in\rG(x_\infty)
$$ 
for all~$i$ by Theorem~\ref{thm:MLT}.  Hence~${x\in\rG(x_\infty)}$ 
because the group orbit~$\rG(x_\infty)$ is compact. 
This proves Theorem~\ref{thm:NESS2}.
\end{proof}

In general, the moment map squared is far from a Morse--Bott 
function and may have very complicated critical points. 
However, it follows from the Kirwan--Ness Inequality,
the Moment Limit Theorem, and the Ness Uniqueness Theorem
that the stable manifolds of this gradient flow
exhibit a structure that resembles a stratification 
by stable manifolds of a Morse--Bott function.
More precisely, let~$x\in X$ be a critical point of 
the moment map squared and define 
the stable manifold of the critical set~$\rG(x)$ by 
\begin{equation}\label{eq:Wsx1}
W^s(\rG(x)) := \left\{y_0\in X\,\Bigg|\,
\begin{array}{l}
\mbox{the unique solution }y:\R\to X\mbox{of }\\
\dot y=-JL_y\mu(y),\,y(0)=y_0\mbox{ satisfies}\\
\lim_{t\to\infty}y(t)=ux\mbox{ for some }u\in\rG
\end{array}\right\}.
\end{equation}
By Theorem~\eqref{thm:XINFTY} $X$ is the union of these stable 
manifolds, and each stable manifold is a union of~$\rG^c$-orbits 
by Theorems~\ref{thm:MLT} and~\ref{thm:NESS2}.

\begin{cor}\label{cor:KIRNESS2}
For~${x\in\Crit(f)}$, let~${W^s(\rG(x))\subset X}$  
be the stable manifold in~\eqref{eq:Wsx1}.
Then the following holds.

\smallskip\noindent{\bf (i)}  
$X=\bigcup_{x\in\Crit(f)}W^s(\rG(x))$.

\smallskip\noindent{\bf (ii)}
Let~${x\in\Crit(f)}$ and~$y_0\in X$.  
Then~$y_0\in W^s(\rG(x))$ if and only if
\begin{equation}\label{eq:WS1}
x\in\overline{\rG^c(y_0)},\qquad
\abs{\mu(x)}=\inf_{g\in\rG^c}\abs{\mu(gy_0)}
\end{equation}

\smallskip\noindent{\bf (iii)}
Let~${x\in\Crit(f)}$. Then~$W^s(\rG(x))$ is a union of $\rG^c$-orbits.
\end{cor}

\begin{proof}
Part~(i) follows directly from the Convergence Theorem~\ref{thm:XINFTY}.
To prove part~(ii), let~${y_0\in X}$, let~${y:\R\to X}$ 
be the unique solution of the initial value 
problem~${\dot y=-JL_y\mu(y)}$ with~${y(0)=y_0}$, 
and define~${y_\infty:=\lim_{t\to\infty}y(t)}$ (Theorem~\ref{thm:XINFTY}). 
Then, by Lemma~\ref{le:GRADFLOW} and Theorem~\ref{thm:MLT}, we have
\begin{equation}\label{eq:WS2}
y_\infty\in\overline{\rG^c(y_0)},\qquad
\abs{\mu(y_\infty)}=\inf_{g\in\rG^c}\abs{\mu(gy_0)}.
\end{equation}
Moreover, it follows from the definition of~$W^s(\rG(x))$
that~$y_0\in W^s(\rG(x))$ 
if and only if~${y_\infty\in\rG(x)}$.
Thus~$y_0\in W^s(\rG(x))$ implies~\eqref{eq:WS1}.
Conversely, if~$y_0$ satisfies~\eqref{eq:WS1}
then it follows from~\eqref{eq:WS2} and the 
Second Ness Uniqueness Theorem~\ref{thm:NESS2} 
that $y_\infty\in\rG(x)$ and hence $y_0\in W^s(\rG(x))$.  
This proves~(ii).  It follows from~(ii) and the 
Moment Limit Theorem~\ref{thm:MLT} that~$W^s(\rG(x))$ 
is a (possibly infinite) union of~$\rG^c$ orbits. 
This proves~(iii) and Corollary~\ref{cor:KIRNESS2}.
\end{proof}

The\index{Kirwan homomorphism} 
stratification in Corollary~\ref{cor:KIRNESS2}
was used by Kirwan~\cite{KIRWAN} to prove that the canonical 
ring homomorphism (the {\bf Kirwan homomorphism})
$$
\kappa:H^*_\rG(X)\to H^*(X\dslash\rG)
$$
from the equivariant cohomology of~$X$ to the cohomology 
of the Marsden--Weinstein quotient~${X\dslash\rG}$ is surjective. 

\bigbreak

\begin{thm}[{\bf General Moment-Weight Inequality}]\label{thm:MW2}
For every element~${x\in X}$,\index{moment-weight inequality} 
every~${\xi\in\cg\setminus\{0\}}$, and every~${g\in\rG^c}$,
\begin{equation}\label{eq:MW2}
\frac{-w_\mu(x,\xi)}{\abs{\xi}} 
\le \abs{\mu(gx)}.
\end{equation}
\end{thm}

We give two proofs of Theorem~\ref{thm:MW2}.
The first proof is due to Mumford~\cite{MUMFORD} 
and Ness~\cite[Lemma~3.1]{NESS2} 
and is based on Theorem~\ref{thm:MUMFORD}.
The second proof is due to Chen~\cite{CIII,CIV}.  
His~methods were developed in the infinite-dimensional 
setting of K-stability for K\"ahler--Einstein metrics 
(Tian~\cite{TIAN}, Donaldson~\cite{DON-K2,DON-K3,DON-K4}).
Chen's infinite-dimensional argument carries over
to the finite-dimensional setting. 

\begin{proof}[Proof~1]
We first prove~\eqref{eq:MW2} for $g=1$.
Choose~${x_0\in X}$ and~${\xi\in\cg\setminus\{0\}}$. 
Define~${x:\R\to X}$ by~${x(t) := \exp(\i t\xi)x_0}$ 
as in Lemma~\ref{le:WEIGHT2}. Then the 
function~${t\mapsto \inner{\mu(x(t))}{\xi}}$ 
is nondecreasing.  Hence 
\begin{equation*}
\begin{split}
\inner{\mu(x_0)}{\xi}
= 
\inner{\mu(x(0))}{\xi} 
\le
\lim_{t\to\infty}\inner{\mu(x(t))}{\xi} 
= 
w_\mu(x_0,\xi).
\end{split}
\end{equation*}
Hence, by the Cauchy--Schwarz inequality,
$$
-w_\mu(x_0,\xi) \le -\inner{\mu(x_0)}{\xi}
\le \abs{\mu(x_0)}\abs{\xi}
$$
and this implies~\eqref{eq:MW2} with $g=1$. 

Now let~${x\in X}$,~${\xi\in\cg\setminus\{0\}}$, 
and~${g\in\rG^c}$.  Then
$$
\zeta:=g\xi\ g^{-1}\in\sT^c.
$$ 
Hence it follows from Theorem~\ref{thm:MUMFORD}
that there exists an element~${p\in\rP(\zeta)}$
such that~${p\zeta p^{-1}\in\cg\setminus\{0\}}$. 
Thus it follows from part~(ii) of 
Theorem~\ref{thm:MUMFORD1} that 
${w_\mu(gx,p\zeta p^{-1}) = w_\mu(gx,\zeta)}$.
Now apply the first step of the proof to
the pair $(gx,p\zeta p^{-1})\in X\times(\cg\setminus\{0\})$. 
Then
\begin{equation*}
\begin{split}
\abs{\mu(gx)}
&\ge 
\frac{-w_\mu(gx,p\zeta p^{-1})}{\abs{p\zeta p^{-1}}} \\
&= 
\frac{-w_\mu(gx,\zeta)}{\sqrt{\abs{\Re(\zeta)}^2-\abs{\Im(\zeta)}^2}} \\
&= 
\frac{-w_\mu(x,\xi)}{\abs{\xi}}.
\end{split}
\end{equation*}
Here the last step follows from part~(i) of Theorem~\ref{thm:MUMFORD1},
Lemma~\ref{le:WEIGHT3}, and the fact that $\zeta=g\xi g^{-1}$.
This completes the first proof of Theorem~\ref{thm:MW2}.
\end{proof}

\begin{proof}[Proof~2]
Let $x\in X$, $\xi\in\cg\setminus\{0\}$,
and $g\in\rG^c$.\index{moment-weight inequality!Chen's proof}
For $t\ge0$ choose $\eta(t)\in\cg$ 
and $u(t)\in\rG$ such that
\begin{equation}\label{eq:xieta1}
\exp(\i\eta(t))u(t) = \exp(\i t\xi)g^{-1}.
\end{equation}
We prove that
\begin{equation}\label{eq:xieta2}
\lim_{t\to\infty}\frac{\eta(t)}{\abs{\eta(t)}} 
= \frac{\xi}{\abs{\xi}}.
\end{equation}
To see this, note that 
$$
\exp(-\i\eta(t))\exp(\i t\xi)g^{-1}\in\rG.
$$
Hence it follows from Lemma~\ref{le:estimate}
that there exists a constant $c>0$ such 
that~${\abs{t\xi-\eta(t)}\le c}$ for all~${t\ge 0}$.
Hence
\begin{equation*}
\begin{split}
\Abs{\frac{\xi}{\abs{\xi}}-\frac{\eta(t)}{\abs{\eta(t)}}}
&\le
\frac{\abs{t\xi-\eta(t)}}{t\abs{\xi}}
+ \abs{\eta(t)}\Abs{\frac{1}{t\abs{\xi}}-\frac{1}{\abs{\eta(t)}}} \\
&=
\frac{\abs{t\xi-\eta(t)}}{t\abs{\xi}}
+ \frac{\abs{\abs{t\xi}-\abs{\eta(t)}}}{t\abs{\xi}} \\
&\le
\frac{2c}{t\abs{\xi}}.
\end{split}
\end{equation*}
This proves~\eqref{eq:xieta2}.
Now abbreviate $u:=u(t)$ and $\eta:=\eta(t)$.
Then, by~\eqref{eq:xieta1},
$$
\exp(\i u^{-1}\eta u)g = u^{-1}\exp(\i t\xi).
$$
Moreover, by~\eqref{eq:ZERO}, the function
$
s\mapsto\inner{\mu(\exp(\i su^{-1}\eta u)gx)}{u^{-1}\eta u}
$
is nondecreasing. Hence
\begin{equation*}
\begin{split}
-\abs{\mu(gx)}
&\le
\abs{\eta}^{-1}
\inner{\mu(gx)}{u^{-1}\eta u} \\
&\le
\abs{\eta}^{-1}
\inner{\mu(\exp(\i u^{-1}\eta u)gx)}{u^{-1}\eta u} \\
&=
\abs{\eta}^{-1}
\inner{\mu(u^{-1}\exp(\i t\xi)x)}{u^{-1}\eta u} \\
&=
\abs{\eta}^{-1}
\inner{\mu(\exp(\i t\xi)x)}{\eta} \\
&=
\abs{\xi}^{-1}\inner{\mu(\exp(\i t\xi)x)}{\xi} +
\inner{\mu(\exp(\i t\xi)x)}
{\abs{\eta}^{-1}\eta-\abs{\xi}^{-1}\xi}.
\end{split}
\end{equation*}
Take the limit $t\to\infty$ and use~\eqref{eq:xieta2}
to obtain
$$
-\abs{\mu(gx)}\le 
\lim_{t\to\infty}
\frac{\inner{\mu(\exp(\i t\xi)x)}{\xi}}{\abs{\xi}}
= \frac{w_\mu(x,\xi)}{\abs{\xi}}.
$$
This completes the second proof 
of Theorem~\ref{thm:MW2}.
\end{proof}

\begin{figure}[htp] 
\centerline{\psfig{figure=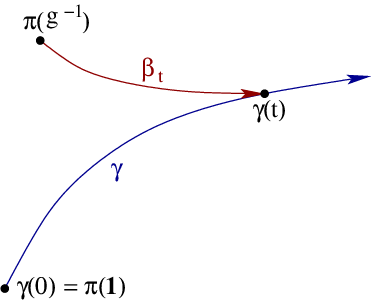,height=1.7in}} 
\caption{Chen's proof of the moment-weight inequality.}
\label{fig:chen}      
\end{figure}

Here is the geometric picture behind Chen's proof
in~\cite{CIII} of the moment-weight inequality (proof~2 above).
The homogeneous space~${M=\rG^c/\rG}$ is a complete, connected,
simply connected Riemannian manifold with nonpositive
sectional curvature (Appendix~\ref{app:GcG}).  
The curve
$$
\gamma(t):=\pi(\exp(-\i t\xi))
$$ 
is a geodesic through the point~${p_0:=\pi(\one)}$.  
The formula~\eqref{eq:xieta1} asserts that, 
for each~$t$, the curve
$$
\beta_t(s):=\pi(g^{-1}\exp(-\i su(t)^{-1}\eta(t)u(t)))
$$
is the unique geodesic connecting~${p_1:=\pi(g^{-1})=\beta_t(0)}$ 
to~${\gamma(t)=\beta_t(1)}$  (see Figure~\ref{fig:chen}).  
Equation~\eqref{eq:xieta2} asserts that the angle between 
the geodesics~$\beta_t$ and~$\gamma$ at the point~$\gamma(t)$ 
where they meet tends to zero as~$t$ tends to infinity.   
(This is closely related to the sphere at infinity discussed in Donaldson's
paper~\cite{DON-lie}.)  Now the main point of Chen's proof 
is a comparison between the slope of the Kempf--Ness
function~${\Phi_x:M\to\R}$ along the 
geodesics~$\gamma$ and~$\beta_t$ 
at the point~$\gamma(t)$. 


\chapter{Stability in symplectic geometry}\label{ch:STABLE}

This chapter introduces the $\mu$-stability conditions
for elements of~$X$ and characterizes them in terms of 
the properties of the Kempf--Ness function~$\Phi_x$.
This is the content of the generalized Kempf--Ness 
Theorem~\ref{thm:KN}.   

\begin{defn}\label{def:mustable}
An element $x\in X$ is called

\smallskip{\bf $\mu$-unstable}\index{unstable}
iff $\overline{\rG^c(x)}\cap\mu^{-1}(0)=\emptyset$,

\smallskip{\bf $\mu$-semistable}\index{semistable}
 iff $\overline{\rG^c(x)}\cap\mu^{-1}(0)\ne\emptyset$,

\smallskip{\bf $\mu$-polystable}\index{polystable}
 iff $\rG^c(x)\cap\mu^{-1}(0)\ne\emptyset$,

\smallskip{\bf $\mu$-stable}\index{stable}
 iff $\rG^c(x)\cap\mu^{-1}(0)\ne\emptyset$ and 
 $\rG^c_x:=\left\{g\in\rG^c\,|\,gx=x\right\}$ is discrete.

\medskip\noindent
Denote the sets of $\mu$-semistable, $\mu$-polystable, 
and $\mu$-stable points by 
\begin{equation}\label{eq:Xpss}
\begin{split}
X^\ss &:=  \left\{x\in X\,|\,x\mbox{ is $\mu$-semistable}\right\},\\
X^\ps &:=  \left\{x\in X\,|\,x\mbox{ is $\mu$-polystable}\right\},\\
X^\s &:=  \left\{x\in X\,|\,x\mbox{ is $\mu$-stable}\right\}.
\end{split}
\end{equation}
\end{defn}

With this terminology Theorem~\ref{thm:NESS2} asserts for $m=0$
that the closure of each $\mu$-semistable $\rG^c$-orbit 
in $X$ contains a unique $\mu$-polystable $\rG^c$-orbit.

\begin{thm}[{\bf $\mu$-Stability Theorem}]\label{thm:STABILITY}
Let $x_0\in X$ and let $x:\R\to X$ be the solution of~\eqref{eq:KN1}.
Define\index{Stability Theorem}
$
x_\infty:=\lim_{t\to\infty}x(t). 
$
Then the following~holds.

\smallskip\noindent{\bf (i)}
$x_0\in X^\ss$ if and only if $\mu(x_\infty)=0$.

\smallskip\noindent{\bf (ii)}
$x_0\in X^\ps$ if and only if $\mu(x_\infty)=0$ and $x_\infty\in\rG^c(x_0)$.

\smallskip\noindent{\bf (iii)}
$x_0\in X^\s$ if and only if the isotropy subgroup $\rG_{x_\infty}$ is discrete.

\smallskip\noindent
Moreover, $X^\ss$ and $X^\s$ are open subsets of $X$.
\end{thm}

\begin{proof}
By definition $x_0\in X^\ss$
if and only if~${\inf_{g\in\rG^c}\abs{\mu(gx_0)}=0}$.  
Hence assertion~(i) follows from Theorem~\ref{thm:MLT}.

We prove that $X^\ss$ is open.
By the Lojasiewicz gradient inequality~\eqref{eq:Loja2} for 
${f=\tfrac12\abs{\mu}^2}$ there exists 
a constant $\delta>0$ such that every~${x\in X}$ satisfies
$$
f(x)<\delta,\quad \nabla f(x)=0\qquad\implies\qquad f(x)=0.
$$
Hence it follows from~(i) that
$
U := \{x\in X\,|\,\tfrac12\abs{\mu(x)}^2<\delta\}\subset X^\ss.
$
Now let~${x_0\in X^\ss}$.  Then there is a~${g_0\in\rG^c}$ 
such that~${g_0x_0\in U}$. Hence there exists a neighborhood~${V\subset X}$ 
of~$x_0$ such that~${g_0V\subset U}$, and thus~${V\subset X^\ss}$.
This shows that~$X^\ss$ is an open subset of $X$.  

We prove~(ii).  If $\mu(x_\infty)=0$ and $x_\infty\in\rG^c(x_0)$ then 
$\rG^c(x_0)\cap\mu^{-1}(0)\ne\emptyset$ and hence $x_0\in X^\ps$.
Conversely, suppose that $x_0\in X^\ps$.  Then there exists an 
element $g_0\in\rG^c$ such that $\mu(g_0x_0)=0$.  
In particular, $g_0x_0$ is a critical point of $f=\tfrac12\abs{\mu}^2$
and thus the negative gradient flow line of $f$ with the initial value 
$g_0x_0$ is constant.  Hence it follows from Theorem~\ref{thm:MLT} 
that $g_0x_0\in\rG(x_\infty)$.  Hence $x_\infty\in\rG^c(x_0)$
and $\mu(x_\infty)=0$.  This proves~(ii).

We prove that $X^\s$ is open.
By Lemma~\ref{le:ONE} the set
$$
Z^\s:=\{x\in X\,|\,\mu(x)=0,\,\ker\,L_x=0\}
$$
is a smooth submanifold of $X$ with tangent spaces 
$T_xZ^\s = \ker\,d\mu(x)$. It follows also from 
Lemma~\ref{le:ONE} and the definition of $\mu$-stability 
that~${X^\s = \rG^cZ^\s}$. Next we prove that there exists 
an open set $U^\s\subset X$ such that 
$$
Z^\s\subset U^\s\subset X^\s.
$$
Define the map $\psi:Z^\s\times\cg\to X$ by 
$$
\psi(x,\eta):=\exp(\i\eta)x
$$
for $x\in X$ and $\eta\in\cg$.   
Its derivative at an element~${(x,0)\in Z^\s\times\{0\}}$ 
is the linear map~${d\psi(x,0):\ker\,d\mu(x)\times\cg\to T_xX}$ 
given by
$$
d\psi(x,0)(\xhat,\etahat) = \xhat +JL_x\etahat
$$
for $\xhat\in T_xX$ and $\etahat\in\cg$. 
We prove that this map is bijective for every~${x\in Z^\s}$.
If~${\xhat\in\ker d\mu(x)}$, ${\etahat\in\cg}$
satisfy~${\xhat +JL_x\etahat=0}$, 
then~${0=d\mu(x)JL_x\etahat = {L_x}^*L_x\etahat}$ 
by~\eqref{eq:mu3}, hence~${\etahat=0}$, and so~${\xhat=0}$.  
Moreover, the map~$d\mu(x)$ is surjective and 
so~${\dim\,\ker d\mu(x)=\dim\,X-\dim\,\cg}$. 
This shows that $d\psi(x,0)$ is bijective for each $x\in Z^\s$.  
Hence $\psi$ restricts to a diffeomorphism from an open 
neighborhood of $Z^\s\times\{0\}$ in $Z^\s\times\cg$ 
to an open neighborhood $U^\s\subset X$ of~$Z^\s$.  
Since~${Z^\s\subset U^\s\subset X^\s}$, 
it follows that~${X^\s=\bigcup_{g\in\rG^c}gU^\s}$
is an open subset of $X$.

We prove~(iii). Let $x_0\in X$, let $x:\R\to X$ be the unique
solution of~\eqref{eq:KN1}, and define~${x_\infty:=\lim_{t\to\infty}x(t)}$. 
Assume first that $x_0\in X^\s$.  Then it follows from~(ii) 
and the definition of $\mu$-stability that
$\mu(x_\infty)=0$, $x_\infty\in\rG^c(x_0)$, and $\ker L^c_{x_0}=0$.
Hence $\ker L^c_{x_\infty}=0$, thus $\ker\,L_{x_\infty}=0$
by Lemma~\ref{le:ONE}, and so~$\rG_{x_\infty}$ is discrete. 
Conversely, assume $\rG_{x_\infty}$ is discrete and 
so $\ker\,L_{x_\infty}=0$.  Since $L_{x_\infty}\mu(x_\infty)=0$ 
by Theorem~\ref{thm:XINFTY}, this implies~${\mu(x_\infty)=0}$.  
Thus~$x_\infty$ is $\mu$-stable.  Since~$X^\s$ is open, 
$x(t)$ is $\mu$-stable for $t$ sufficiently large and, 
since~${x(t)\in\rG^c(x_0)}$ for all $t$ by Lemma~\ref{le:GRADFLOW}, 
$x_0$ is $\mu$-stable.  This proves~(iii) and Theorem~\ref{thm:STABILITY}.
\end{proof}

The generalized Kempf--Ness theorem characterizes the stability 
condition of~$x$ in terms of the properties 
of the {\it Kempf--Ness function}~$\Phi_x$.  
We will derive the original Kempf--Ness theorem~\cite{KN} 
as a corollary in Theorem~\ref{thm:VECTOR}.

\begin{thm}[{\bf Generalized Kempf--Ness Theorem}]
\label{thm:KN}
\ 

\noindent
Let $x\in X$ and let $\Phi_x$ be the Kempf--Ness function of $x$.  
Then\index{Kempf--Ness Theorem!generalized}

\smallskip\noindent{\bf (i)}
$x$ is $\mu$-unstable if and only if $\Phi_x$ is unbounded below,

\smallskip\noindent{\bf (ii)}
$x$ is $\mu$-semistable if and only if $\Phi_x$ is bounded below,

\smallskip\noindent{\bf (iii)}
$x$ is $\mu$-polystable if and only if $\Phi_x$ has a critical point,

\smallskip\noindent{\bf (iv)}
$x$ is $\mu$-stable if and only if $\Phi_x$ is bounded below and proper.
\end{thm}

\begin{proof}
Fix a point $x_0\in X$ and denote by $\Phi:M\to\R$ 
the Kempf--Ness function of $x_0$.  Throughout the proof
${x:\R\to X}$ denotes the solution of~\eqref{eq:KN1}, 
${g:\R\to\rG^c}$ denotes the solution of~\eqref{eq:KN2} 
so that~${x(t)=g(t)^{-1}x_0}$ for all $t$, and 
${\gamma:\R\to M}$ denotes the composition of $g$
with the projection from~$\rG^c$ to~$M$. 
The limit ${x_\infty := \lim_{t\to\infty}x(t)}$
exists by Theorem~\ref{thm:XINFTY}, 
and Theorem~\ref{thm:MLT} asserts that
$\abs{\mu(x_\infty)}=\inf_{g\in\rG^c}\abs{\mu(gx_0)}$.

We prove necessity in~(i).  Assume $x_0$
is $\mu$-unstable. Then $\mu(x_\infty)\ne 0$.
By definition of the Kempf--Ness function
$$
\tfrac{d}{dt}(\Phi\circ\gamma)
= - \inner{\mu(g^{-1}x_0)}{\Im(g^{-1}\dot g)}
= - \abs{\mu(x)}^2
\le - \abs{\mu(x_\infty)}^2.
$$
Thus $\Phi(\gamma(t))\le -t\abs{\mu(x_\infty)}^2$ 
for all $t\ge 0$, and so $\Phi$ is unbounded below.

We prove necessity in~(ii).
Assume $x_0$ is $\mu$-semistable. Then ${\mu(x_\infty)=0}$.
By the Lojasiewicz gradient inequality for $f=\tfrac{1}{2}\abs{\mu}^2$
in~\eqref{eq:Loja2} there exist positive constants~$t_0$ 
and~$c$, and a constant~${1/2<\alpha<1}$, such that
$$
\abs{\mu(x)}^2 
= 2\abs{f(x)}
\le 2\abs{f(x)}^\alpha
\le c\abs{\nabla f(x)}
= c\abs{JL_{x}\mu(x)} 
= c\abs{\dot x}
$$
for $t\ge t_0$.  By Theorem~\ref{thm:XINFTY} 
the function $t\mapsto\abs{\dot x(t)}$ is integrable
over the positive real axis and hence, so is the function
$t\mapsto\abs{\mu(x(t))}^2=-\tfrac{d}{dt}(\Phi\circ\gamma)(t)$.
Hence the limit $a:=\lim_{t\to\infty}\Phi(\gamma(t))$ exists in~$\R$, 
hence~${\inf_M\Phi=a>-\infty}$ by Theorem~\ref{thm:KNF},
and so~$\Phi$ is bounded below.

\bigbreak

Thus we have proved that the conditions 
on the Kempf--Ness function are necessary in~(i) and~(ii).
Since necessity in~(i) is equivalent to sufficiency in~(ii)
and vice versa, this proves~(i) and~(ii).
Assertion~(iii) follows from the fact that $\pi(g)$ is a critical point
of $\Phi$ if and only if $\mu(g^{-1}x_0)=0$.

We prove~(iv). 
Assume first that $\Phi$ is bounded below and proper.
Then $x_0$ is $\mu$-semistable by part~(ii).  
Moreover, by definition of the Kempf--Ness 
function $\Phi$ and the negative gradient flow 
line~$\gamma$,  we have
$$
\Phi(\gamma(0))=0,\qquad
-\infty<a:=\inf_M\Phi\le 0,
$$
and the function $\Phi\circ\gamma:\R\to\R$ is 
nonincreasing. Since $\Phi$ is proper the 
set~${\Phi^{-1}([a,0])}$ is compact 
and contains~${\gamma(t)}$ for every~${t\ge 0}$.  
Hence there exists a sequence $t_i\to\infty$
such that $\gamma(t_i)$ converges. 
Hence $\Phi$ has a critical point. 
By~(iii) this implies that~$x_0$ is $\mu$-polystable.
Hence there exists an element~${g_0\in\rG^c}$ 
such that~${\mu(g_0x_0)=0}$. 
Assume by contradiction that~$x_0$ is not $\mu$-stable.
Then~${\ker\,L_{x_0}^c\ne0}$, hence~${\ker\,L_{g_0x_0}^c\ne0}$,
and hence $\ker\,L_{g_0x_0}\ne0$ by Lemma~\ref{le:ONE}.  
Choose~${\eta\in\cg\setminus\{0\}}$ such that~${L_{g_0x_0}\eta=0}$. 
Then $\exp(\i s\eta)g_0x_0=g_0x_0$
and hence $\mu(\exp(\i s\eta)g_0x_0)=0$.
Hence the curve~${\beta(s):=\pi(g_0^{-1}\exp(-\i s\eta))}$
consists of critical points of $\Phi$ and so,
by part~(iii) of Theorem~\ref{thm:KNF},
$
\Phi(\beta(s))=a
$
for all~$s$. Thus~${\Phi^{-1}(a)}$ is not compact 
in contradiction to the assumption
that~$\Phi$ is proper. Thus~$x_0$ is $\mu$-stable.

Conversely, assume $x_0$ is $\mu$-stable.
Then $\Phi$ is bounded below by part~(ii)
and we must prove that $\Phi$ is proper.
By part~(vii) of Theorem~\ref{thm:KNF}
it suffices to assume that~${\mu(x_0)=0}$. 
We prove first that 
$
p_0:=\pi(\one)\in M
$ 
is the unique point at which $\Phi$ attains 
its minimum $\min_M\Phi=0$. To see this, 
let~${p=\pi(g)\in M\setminus\{p_0\}}$.
Choose $\eta\in\cg$ and $u\in\rG$ such that 
$g=\exp(\i\eta)u$, and define
$
y(t):=\exp(\i t\eta)x_0
$
and
$ 
\phi(t):=\Phi(\pi(\exp(-\i t\eta)).
$
Then $\eta\ne 0$, the function $\phi:\R\to\R$ is convex, and 
$$
\dot\phi(t) = \inner{\mu(y(t))}{\eta},\qquad
\ddot\phi(t) = \abs{L_{y(t)}\eta}^2.
$$
Since $y(0)=x_0$ we obtain $\phi(0)=0$, $\dot\phi(0)=0$,
and $\ddot\phi(0)=\abs{L_{x_0}\eta}^2>0$.
This implies~${\phi(t)>0}$ for every $t\in\R\setminus 0$
and, in particular, $\Phi(p)=\phi(-1)>0$. 

For every~${r>0}$ denote by~${B_r\subset M}$ 
the ball of radius~$r$ centered at~$p_0$. 
Define the number~${\delta := \inf_{\p B_1}\Phi > 0}$
Then it follows from convexity that
$$
d(p_0,p)\ge 1\qquad\implies\qquad
\Phi(p) \ge \delta d_M(p_0,p)
$$
for all~${p\in M}$.   Hence, for every $c>0$, 
$$
c\ge\delta\qquad\implies\qquad
\Phi^{-1}([0,c]) \subset B_{c/\delta}.
$$
Hence, for every $c>0$, the set 
$\Phi^{-1}([0,c])$ is closed and bounded,
and hence compact because $M$ is complete.
This proves~(iv) and Theorem~\ref{thm:KN}.
\end{proof}

\bigbreak

A central result in geometric invariant theory is the
Hilbert--Mumford criterion, which characterizes the  
$\mu$-stability conditions in terms of the $\mu$-weights
introduced in Chapter~\ref{ch:WEIGHTS}.  The necessity of 
these conditions follows directly from the general moment-weight 
inequality~\eqref{eq:MW2} in Theorem~\ref{thm:MW2}.

\begin{thm}[{\bf Stability and Weights}]
\label{thm:HMnec}
Let~${x_0\in X}$.  Then the following holds.

\smallskip\noindent{\bf (i)}
If $x_0$ is $\mu$-semistable and $\zeta\in\sT^c$ 
then $w_\mu(x_0,\zeta)\ge0$.

\smallskip\noindent{\bf (ii)}
If $x_0$ is $\mu$-polystable and $\zeta\in\sT^c$ 
then $w_\mu(x_0,\zeta)\ge0$ and
$$
w_\mu(x_0,\zeta)=0
\qquad\iff\qquad
\lim_{t\to\infty}\exp(\i t\zeta)x_0\in\rG^c(x_0).
$$

\smallskip\noindent{\bf (iii)}
If $x_0$ is $\mu$-stable and $\zeta\in\sT^c$ 
then $w_\mu(x_0,\zeta)>0$.
\end{thm}

\begin{proof}
See page~\pageref{proof:HMnec}.
\end{proof}

The Hilbert--Mumford criterion asserts that the necessary 
conditions for $\mu$-semistability, $\mu$-polystability, 
and $\mu$-stability in Theorem~\ref{thm:HMnec} are in fact
necessary and sufficient.  This is proved in
Chapter~\ref{ch:HM} below. 

\begin{lem}\label{le:POLYWEIGHT1}
Let $x_0\in X$ such that $\mu(x_0)=0$ 
and let~${\xi\in\cg\setminus\{0\}}$.  
Then the following are equivalent.

\smallskip\noindent{\bf (i)} 
$w_\mu(x_0,\xi)=0$.

\smallskip\noindent{\bf (ii)} 
$L_{x_0}\xi=0$.

\smallskip\noindent{\bf (iii)} 
$\lim_{t\to\infty}\exp(\i t\xi)x_0\in\rG^c(x_0)$.
\end{lem}

\begin{proof}
Assume~(i) and define $x(t):=\exp(\i t\xi)x_0$. 
Then, by~(i) and~\eqref{eq:ZERO}, 
$$
\inner{\mu(x(0))}{\xi}=0,\qquad
\lim_{t\to\infty}\inner{\mu(x(t))}{\xi}=0,\qquad
\frac{d}{dt}\inner{\mu(x(t))}{\xi}=\abs{L_{x(t)}\xi}^2
$$
for all~$t$. Hence $L_{x(t)}\xi=0$ for all $t$
and this shows that~(i) implies~(ii). 
That~(ii) implies~(iii) follows from the fact 
that $\exp(\i t\xi)x_0=x_0$ for~${t\in\R}$ 
and~${\xi\in\ker L_{x_0}}$.  Now assume~(iii) 
and define 
$
x^+:=\lim_{t\to\infty}\exp(\i t\xi)x_0.
$
By~(iii) there exists an element $g\in\rG^c$ 
such that $x^+=gx_0$. Since $L_{x^+}\xi=0$ by 
Lemma~\ref{le:WEIGHT2}, we have
\begin{equation*}
\begin{split}
w_\mu(x_0,\xi)
&= 
\inner{\mu(x^+)}{\xi} 
= 
w_\mu(x^+,\xi) \\ 
&= 
w_\mu(gx_0,\xi) 
= 
w_\mu(x_0,g^{-1}\xi g) \\
&= 
\inner{\mu(x_0)}{\Re(g^{-1}\xi g)} 
= 0.
\end{split}
\end{equation*}
Here we have used Theorem~\ref{thm:MUMFORD1} 
and the fact that $\exp(\i t g^{-1}\xi g)x_0=x_0$ for all~$t$.  
Thus~(iii) implies~(i) and this proves Lemma~\ref{le:POLYWEIGHT1}.
\end{proof}

\begin{lem}\label{le:POLYWEIGHT2}
Let $x_0\in X$ such that $\mu(x_0)=0$ 
and let~$\zeta\in\sT^c$.  Then
$$
w_\mu(x_0,\zeta)=0
\qquad\iff\qquad
x^+:=\lim_{t\to\infty}\exp(\i t\zeta)x_0\in\rG^c(x_0).
$$
\end{lem}

\begin{proof}
By Theorem~\ref{thm:MUMFORD} there exists a~${p\in\rP(\zeta)}$ 
such that~${\xi:=p\zeta p^{-1}\in\cg}$.  Then the limit 
$
p^+ := \lim_{t\to\infty}\exp(\i t\zeta)p\exp(-\i t\zeta)
$
exists in $\rG^c$.  It satisfies
\begin{equation*}
\begin{split}
x^+
&=
\lim_{t\to\infty}\exp(\i t\zeta)x_0 \\
&=
\lim_{t\to\infty}\bigl(\exp(\i t\zeta)p\exp(-\i t\zeta)\bigr)
\lim_{t\to\infty}\exp(\i t\zeta)p^{-1}x_0 \\
&= 
p^+p^{-1}\lim_{t\to\infty}p\exp(\i t\zeta)p^{-1}x_0 \\
&=
p^+p^{-1}\lim_{t\to\infty}\exp(\i t\xi)x_0. 
\end{split}
\end{equation*}
This shows that $x^+\in\rG^c(x_0)$ if and only if
$\lim_{t\to\infty}\exp(\i t\xi)x_0\in\rG^c(x_0)$.
Moreover, we have~${w_\mu(x_0,\xi)=w_\mu(x_0,\zeta)}$ 
by part~(ii) of Theorem~\ref{thm:MUMFORD1}.
Hence it follows from Lemma~\ref{le:POLYWEIGHT1} 
that $w_\mu(x_0,\zeta)=0$ if and only if $x^+\in\rG^c(x_0)$,
and this proves Lemma~\ref{le:POLYWEIGHT2}.
\end{proof}

\begin{proof}[Proof of Theorem~\ref{thm:HMnec}]\label{proof:HMnec}
We prove part~(i) by an indirect argument.  Suppose that there exists
a toral generator $\zeta\in\sT^c$ such that $w_\mu(x_0,\zeta)<0$. 
Then, by Theorem~\ref{thm:MUMFORD} and part~(ii) 
of Theorem~\ref{thm:MUMFORD1}, there exists a~${\xi\in\cg\setminus\{0\}}$ 
such that~${w_\mu(x_0,\xi)<0}$.  Hence~${\inf_{g\in\rG^c}\Abs{\mu(gx_0)}>0}$ 
by Theorem~\ref{thm:MW2} and so $x_0$ is $\mu$-unstable.
This proves part~(i).  

We prove part~(ii).   Suppose $x_0$ is $\mu$-polystable 
and fix an element~${\zeta\in\sT^c}$.  
Then~${w_\mu(x_0,\zeta)\ge0}$ by part~(i).
Now choose~${g\in\rG^c}$ such that~${\mu(gx_0)=0}$ 
and define~${x^+:=\lim_{t\to\infty}\exp(\i t\zeta)x_0}$.
Then
$$
\lim_{t\to\infty}\exp(\i tg\zeta g^{-1})gx_0=gx^+.
$$
Hence~${w_\mu(gx_0,g\zeta g^{-1})=0}$ if and only
if~${gx^+\in\rG^c(gx_0)}$, by Lemma~\ref{le:POLYWEIGHT2}.
Moreover, we have~${w_\mu(gx_0,g\zeta g^{-1})=w_\mu(x_0,\zeta)}$
by part~(i) of Theorem~\ref{thm:MUMFORD1}, and this shows 
that~${w_\mu(x_0,\zeta)=0}$ if and only if~${x^+\in\rG^c(x_0)}$.
This proves part~(ii). 

We prove part~(iii).  Assume that $x_0$ is $\mu$-stable.
Then $w_\mu(x_0,\zeta)\ge0$ for all $\zeta\in\sT^c$ by part~(i).
Moreover, there exists an element $g\in\rG^c$ such that 
\begin{equation}\label{eq:HMstab}
\mu(gx_0)=0,\qquad \ker L_{gx_0}=\{0\}.
\end{equation}
Suppose, by contradiction, that there is 
a~${\zeta\in\sT^c}$ such that~${w_\mu(x_0,\zeta)=0}$.
Then~${w_\mu(gx_0,g\zeta g^{-1})=0}$
by part~(i) of Theorem~\ref{thm:MUMFORD1}. 
Hence it follows from Theorem~\ref{thm:MUMFORD} 
and part~(ii) of Theorem~\ref{thm:MUMFORD1}, that there exists 
an element~${\xi\in\cg\setminus\{0\}}$ such that~${w_\mu(gx_0,\xi)=0}$.  
Hence $L_{gx_0}\xi=0$ by Step~1 above, in contradiction 
to~\eqref{eq:HMstab}. This proves part~(iii) and Theorem~\ref{thm:HMnec}.
\end{proof}

\begin{remk}\label{rmk:POLYWEIGHT}\rm
{\bf (i)} In Lemma~\ref{le:POLYWEIGHT1}
the hypothesis~${\mu(x_0)=0}$ cannot be replaced
by the hypothesis that~$x_0$ is $\mu$-polystable. 
For example, consider the diagonal action of~${\rG:=\SO(3)}$ 
on~${X:=S^2\times S^2}$, where~${S^2\subset\R^3}$ is the 
unit sphere with its standard symplectic form. Let
$$
x=(-\xi,\eta) \in S^2\times S^2
$$
be a pair of distinct non-antipodal points on $S^2$, 
so that 
$$
\eta\ne\pm\xi.
$$  
Then~$x$ is $\mu$-polystable,~${\mu(x)\ne 0}$, 
and the isotropy subgroup of~$x$ 
in~$\SO(3)$ is trivial. Identify the Lie 
algebra~${\cg:=\cso(3)}$ with~$\R^3$ 
so that~${\xi\in S^2\subset\cg}$. Then
\begin{equation}\label{eq:polynec1}
\lim_{t\to\infty}\exp(\i t\xi)x= (-\xi,\xi)\in \rG^c(x),\qquad
w_\mu(x,\xi)=0,\qquad L_x\xi\ne 0.
\end{equation}

\smallskip\noindent{\bf (ii)} 
In Lemma~\ref{le:POLYWEIGHT1} the element~${\xi\in\cg\setminus\{0\}}$ 
cannot be replaced by~${\zeta\in\sT^c}$.  
For example, let~${(x,\xi)\in X\times(\cg\setminus\{0\})}$ 
be as in part~(i) so that~$x$ is $\mu$-polystable 
and~\eqref{eq:polynec1} holds. Choose an element~${g\in\rG^c}$ 
such that~${\mu(gx)=0}$ and define 
$$
x_0 := gx,\qquad \zeta:=g\xi g^{-1}\in\sT^c.
$$  
Then, by part~(i) of Theorem~\ref{thm:MUMFORD1}, we have
\begin{equation}\label{eq:polynec2}
\mu(x_0)=0,\qquad w_\mu(x_0,\zeta)=0,\qquad L_{x_0}^c\zeta\ne 0.
\end{equation}

\smallskip\noindent{\bf (iii)} 
Let~${x_0\in X}$ be $\mu$-polystable
and fix an element $\zeta\in\sT^c$.  Then
\begin{equation}\label{eq:polynec3}
L_{x_0}^c\zeta=0\qquad\implies\qquad w_\mu(x_0,\zeta)=0.
\end{equation}
To see this, assume that $L_{x_0}^c\zeta=0$ and 
choose $g\in\rG^c$ such that $\mu(gx_0)=0$.  
Then~${L_{gx_0}^c(g\zeta g^{-1})=0}$ and hence it follows 
from part~(i) of Theorem~\ref{thm:MUMFORD1} 
that
$$
w_\mu(x_0,\zeta)
= w_\mu(gx_0,g\zeta g^{-1})
= \inner{\mu(gx_0)}{\Re(g\zeta g^{-1})}
= 0.
$$
This proves~\eqref{eq:polynec3}.
By part~(ii) the converse does not hold, 
i.e.~${w_\mu(x_0,\zeta)=0}$ does not imply $L_{x_0}^c\zeta=0$, 
even in the case $\mu(x_0)=0$.

\smallskip\noindent{\bf (iv)} 
The hypothesis that $x_0$ is $\mu$-polystable  
cannot be removed in~\eqref{eq:polynec3}. 
For example, consider the standard action of~${\rG=\SO(3)}$ 
on~${X=S^2}$, fix an element~${x_0\in S^2}$, 
and choose~${\xi:=x_0\in S^2\subset\cg}$.  
Then $x_0$ is $\mu$-unstable and~${L_{x_0}\xi=0}$,
however, the Mumford weights 
are~${w_\mu(x_0,\pm\xi)=\pm1}$.
\end{remk}


\chapter{Stability in algebraic geometry}\label{ch:STAG}

This chapter examines linear group actions on projective space,
which is the classical setting of geometric invariant 
theory~\cite{KEMPF,KN,MUMFORD,NESS1,NESS2}.

Let $V$ be a finite-dimensional complex vector space,
let $\rG\subset\rU(n)$ be a compact Lie group equipped with 
a representation $\rG\to\GL(V)$, and let~${\rG^c\subset\GL(n,\C)}$ 
be its complexification and
$$
\rG^c\times V\to V:(g,v)\mapsto gv.
$$
be the complexified representation.

\begin{defn}\label{def:stable}
A nonzero vector $v\in V$ is called

\smallskip{\bf unstable}\index{unstable}
iff $0\in\overline{\rG^c(v)}$,

\smallskip{\bf semistable}\index{semistable}
iff $0\notin\overline{\rG^c(v)}$,

\smallskip{\bf polystable}\index{polystable}
iff $\rG^c(v)=\overline{\rG^c(v)}$,

\smallskip{\bf stable}\index{stable}
iff $\rG^c(v)=\overline{\rG^c(v)}$ and the
isotropy subgroup $\rG^c_v$ is discrete.
\end{defn}

The linear action of $\rG^c$ on $V$ induces an
action on projective space $\bbP(V)$. Fix a 
$\rG$-invariant Hermitian structure on $V$
and denote by~${\inner{\cdot}{\cdot}}$ the
associated real inner product on~$V$.
Choose a scaling factor\footnote
{
Think of $\hbar$ as Planck's constant.  
If the symplectic form $\omega$ is replaced by
a positive real multiple $c\,\omega$ and $\hbar$ is replaced 
by $c\,\hbar$ all the formulas which follow remain correct.
Our choice of $\hbar$ is consistent with the physics notation 
in that $\om,\mu,\hbar$ have the units of action and
the Hamiltonian $H_\xi=\inner{\mu}{\xi}$ 
has the units of energy.
}
${\hbar>0}$ and restrict the symplectic structure 
to the sphere of radius~${r:=\sqrt{2\hbar}}$. 
It is $S^1$-invariant and descends 
to $\bbP(V)$. (This is $2\hbar$ times the Fubini--Study form. 
Its integral over over the positive generator 
of $\pi_2(\bbP(V))$ is $2\pi\hbar$.)

\begin{lem}\label{le:Vmu}
A moment map for the action of $\rG$ on $\bbP(V)$ 
is given by
\begin{equation}\label{eq:Vmu}
\inner{\mu(x)}{\xi} 
= \hbar \frac{\inner{v}{\i\xi v}}{\abs{v}^2},\qquad
x:=[v]\in\bbP(V).
\end{equation}
\end{lem}


\begin{proof}
Define the map~${\mu:\bbP(V)\to\cg^*}$ by~\eqref{eq:Vmu}.
Fix an element~${x\in\bbP(V)}$ and a tangent vectors
${\xhat\in T_x\bbP(V)}$, and choose~${v,\vhat\in V}$
such that
$$
x := [v],\qquad \xhat := [\vhat],\qquad
\abs{v}=\sqrt{2\hbar},\qquad \inner{v}{\vhat}=0.
$$
The infinitesimal action~${L_x:\cg\to T_x\bbP(V)}$ is given 
by~${L_x\xi=[\xi v]}$ for~${\xi\in\cg}$ and the symplectic form
at~$x$ is
$$
\om_x(L_x\xi,\xhat)=\inner{\i\xi v}{\vhat}.
$$
Now choose a smooth path~${\R\to V:t\mapsto v(t)}$
such that
$$
v(0)=v,\qquad \dot v(0)=\vhat,\qquad
\abs{v(t)}=\sqrt{2\hbar}
$$
for all~$t$.   Differentiate the curve~${t\mapsto\inner{\mu(x(t))}{\xi}}$ 
at~${t=0}$ to obtain
\begin{equation*}
\begin{split}
\inner{d\mu(x)\xhat}{\xi}
&= 
\left.\frac{d}{dt}\right|_{t=0}\inner{\mu(x(t))}{\xi} \\
&=
\left.\frac{d}{dt}\right|_{t=0}\tfrac{1}{2} \inner{\i\xi v(t)}{v(t)}  \\
&=
\inner{\i\xi v}{\vhat}  \\
&=
\om_x(L_x\xi,\xhat).  \\
\end{split}
\end{equation*}
This proves Lemma~\ref{le:Vmu}.
\end{proof}

\begin{lem}\label{le:VKN}
Let $\mu:\bbP(V)\to\cg^*$ be the moment map in Lemma~\ref{le:Vmu}.
Then the Kempf--Ness function $\Phi_x:\rG^c\to\R$ associated to
$x=[v]$ is given by
\begin{equation}\label{eq:VKN}
\Phi_x(g) = \hbar \Bigl(\log\abs{g^{-1}v}-\log\abs{v}\Bigr).
\end{equation}
\end{lem}

\begin{proof}
Define $\Phi_x:\rG^c\to\R$ by~\eqref{eq:VKN}.
Then~${\Phi_x(u)=0}$ for all $u\in\rG$ and
\begin{equation*}
\begin{split}
d\Phi_x(g)\ghat
&=
\hbar \frac{\inner{g^{-1}v}{-g^{-1}\ghat g^{-1}v}}
{\abs{g^{-1}v}^2} \\
&=
- \hbar\frac{\inner{g^{-1}v}{\i\Im(g^{-1}\ghat)g^{-1}v}}
{\abs{g^{-1}v}^2} \\
&=
- \inner{\mu(g^{-1}x)}{\Im(g^{-1}\ghat)}
\end{split}
\end{equation*}
for $g\in\rG^c$ and $\ghat\in T_g\rG^c$.
Thus $\Phi_x$ is the Kempf--Ness function. 
\end{proof}

\begin{lem}\label{le:Vweight}
Let $\mu:\bbP(V)\to\cg^*$ be the moment map in Lemma~\ref{le:Vmu}
and fix an element~${x=[v]\in\bbP(V)}$.  Then the following holds.

\smallskip\noindent{\bf (i)}
Let~${\zeta\in\sT^c}$, denote by 
$$
\lambda_1<\cdots<\lambda_k
$$ 
the eigenvalues 
of~$\i\zeta$ (understood  as a linear operator on~$V$), 
and denote by~${V_i\subset V}$ the corresponding 
eigenspaces.  Write
$$
v=\sum_{i=1}^kv_i,\qquad
v_i\in V_i. 
$$
Then
\begin{equation}\label{eq:Vweight}
w_\mu(x,\zeta) = \hbar\max_{v_i\ne 0}\lambda_i.
\end{equation}

\smallskip\noindent{\bf (ii)}
Assume $v$ is semistable, i.e.\ $0\notin\overline{\rG^c(v)}$.
If $\zeta\in\cg^c$ and $v$ is an eigenvector of $\zeta$
then $\zeta v=0$. 
\end{lem}

\begin{proof}
Write~${\zeta=\xi+\i\eta\in\sT^c}$ with~${\xi,\eta\in\cg}$.  Then
\begin{equation*}
\begin{split}
w_\mu(x,\zeta)
&:=
\lim_{t\to\infty}\inner{\mu(\exp(\i t\zeta)x}{\xi} \\
&\phantom{:}=
\hbar\lim_{t\to\infty} 
\frac{\inner{\exp(\i t\zeta)v}{\i\xi\exp(\i t\zeta)v}}
{\abs{\exp(\i t\zeta)v}^2} \\
&\phantom{:}=
\hbar\lim_{t\to\infty} 
\frac{\inner{\exp(\i t\zeta)v}{\i\zeta\exp(\i t\zeta)v}}
{\abs{\exp(\i t\zeta)v}^2} \\
&\phantom{:}=
\hbar\lim_{t\to\infty} 
\frac{\inner{\sum_{i=1}^ke^{\lambda_it}v_i}
{\sum_{j=1}^k\lambda_je^{\lambda_jt}v_j}}
{\abs{\exp(\i t\zeta)v}^2} \\
&\phantom{:}=
\hbar\lim_{t\to\infty} 
\frac{\sum_{i,j=1}^k\lambda_je^{(\lambda_i+\lambda_j)t}\inner{v_i}{v_j}}
{\sum_{i,j=1}^ke^{(\lambda_i+\lambda_j)t}\inner{v_i}{v_j}} \\
&\phantom{:}=
\hbar\max_{v_i\ne 0}\lambda_i.
\end{split}
\end{equation*}
This proves~(i).  Now assume $0\notin\overline{\rG^c(v)}$
and let $\zeta\in\cg^c$ and $\lambda\in\C$ 
such that
$$
\zeta v=\lambda v.
$$
Then the curves
$$
e^{\lambda t}v=\exp(t\zeta)v,\qquad
e^{\i\lambda t}v=\exp(\i t\zeta)v
$$
in~$\rG^c(v)$ cannot 
contain~$0$ in their closure, and hence~${\lambda=0}$.
This proves~(ii) and Lemma~\ref{le:Vweight}.
\end{proof}

The next result is the original Kempf--Ness theorem in~\cite{KN}.

\begin{thm}[{\bf Kempf--Ness}]\label{thm:VECTOR}
Let $\mu:\bbP(V)\to\cg^*$ be the moment map in Lemma~\ref{le:Vmu}
and fix an element $x=[v]\in\bbP(V)$.\index{Kempf--Ness Theorem}  
Then the following holds.

\smallskip\noindent{\bf (i)}
$v$ is unstable if and only if
$x$ is $\mu$-unstable.

\smallskip\noindent{\bf (ii)}
$v$ is semistable if and only if
$x$ is $\mu$-semistable.

\smallskip\noindent{\bf (iii)}
$v$ is polystable if and only if
$x$ is $\mu$-polystable.

\smallskip\noindent{\bf (iv)}
$v$ is stable if and only if
$x$ is $\mu$-stable.
\end{thm}

\begin{proof}
The Kempf--Ness function 
$\Phi_x(\pi(g))=\hbar(\log\abs{g^{-1}v}-\log\abs{v})$
is unbounded below if and only if $0\in\overline{\rG^c(v)}$.
By Theorem~\ref{thm:KN}, $\Phi_x$ is unbounded below 
if and only if $x$ is $\mu$-unstable. This proves~(i) and~(ii).

To prove~(iii) assume first that $v$ is polystable.
Thus $\rG^c(v)$ is a closed subset of $V$
and hence~${c:=\inf_{g\in\rG^c}\abs{g^{-1}v}>0}$. 
Choose a sequence~${g_i\in\rG^c}$ such that
the sequence~${\abs{g_i^{-1}v}}$ converges to~$c$.  
Then the sequence~${g_i^{-1}v\in V}$ is bounded and hence
has a convergent subsequence. Since~$\rG^c(v)$ is closed, 
the limit of this subsequence has the form~$g^{-1}v$ 
for some~${g\in\rG^c}$. Thus~${\abs{g^{-1}v}=c}$ 
and so~${\pi(g)\in\rG^c/\rG}$ is a critical point of~$\Phi_x$. 
Hence~$x$ is $\mu$-polystable by part~(iii) of Theorem~\ref{thm:KN}.  

Conversely, assume $x$ is $\mu$-polystable.
Choose a sequence $g_i\in\rG^c$ such that
the limit~${w := \lim_{i\to\infty} g_i^{-1}v}$
exists. Since $\Phi_x$ has a critical point, there is
a sequence $h_i\in\rG^c_{x,0}$ such that $h_ig_i$ 
has a convergent subsequence
(see part~(viii) of Theorem~\ref{thm:KNF}). 
Pass to a subsequence so that 
the limit
$$
g:= \lim_{i\to\infty} h_ig_i
$$ 
exists.  
Since~${h_i^{-1}v=v}$ for all~$i$,
by part~(ii) of Lemma~\ref{le:Vweight},
we have
$$
w = \lim_{i\to\infty} g_i^{-1}v 
= \lim_{i\to\infty} g_i^{-1}h_i^{-1}v 
= g^{-1}v\in\rG^c(v).
$$
Thus $\rG^c(v)$ is closed and so $v$ is polystable.
This proves~(iii). 

Since the identity componenets of the isotropy subgroups 
of~$x$ and~$v$ agree in the semistable case, 
by part~(ii) of Lemma~\ref{le:Vweight}, part~(iv) follows from~(iii). 
This proves Theorem~\ref{thm:VECTOR}.
\end{proof}

\begin{remk}\label{rmk:projective}\rm
In Definition~\ref{def:stable} the space $\bbP(V)$ 
can be replaced by any $\rG^c$-invariant closed 
complex submanifold $X\subset\bbP(V)$.
Theorem~\ref{thm:VECTOR} then continues to hold 
for the complex and symplectic structures and moment 
map on~$X$ obtained by restriction.
\end{remk}


\chapter{Rationality}\label{ch:RAT}

The definition of stability extends to actions of $\rG^c$ 
by holomorphic automorphisms on a holomorphic line 
bundle~$E$ over a closed complex manifold~$X$.  
In Definition~\ref{def:stable} replace the complement 
of the origin in~$V$ by the complement of the zero section 
in~$E$ and replace the relation~${x=[v]}$ 
by~${v\in E_x\setminus\{0\}}$. 
In the relevant applications the dual 
bundle~${E^{-1}=\Hom(E,\C)}$ is ample.  

\begin{defn}\label{def:linearization}
A {\bf linearization}\index{linearization of a group action} 
of a holomorphic action of $\rG^c$
on a complex manifold $X$ consists 
of a holomorphic line bundle 
$$
E\to X
$$ 
and a lift of the action of~$\rG^c$ on~$X$ to an action on~$E$ 
by holomorphic line bundle automorpisms.
\end{defn}

A linearization of the $\rG^c$-action is required for
the definition of stability in the intrinsic algebraic geometric 
setting. From the symplectic viewpoint the choice of the line bundle
$E$ is closely related to the symplectic form~$\om$ and 
the choice of the lift of the $\rG^c$-action on~$X$ to
a  $\rG^c$-action on~$E$ is closely related to the choice 
of the moment map~$\mu$.  To obtain a linearization as in 
Definition~\ref{def:linearization} we will need to impose some
additional conditions on the triple $(X,\om,\mu)$. 
Throughout we denote by 
$$
{\bbD:=\left\{z\in\C\,|\,\abs{z}\le 1\right\}}
$$
the closed unit disc in the complex plane
and by 
$$
{\rZ(\cg):=\left\{\xi\in\cg\,|\,
[\xi,\eta]=0\;\forall\,\eta\in\cg\right\}}
$$
the center of $\cg$. 

\bigbreak

Let $x_0\in X$ and let $u:\R/\Z\to\rG$ be a smooth loop.
Then the loop 
$$
u^{-1}x_0:\R/\Z\to X
$$
is contractible.
(It is homotopic to a loop of the form ${t\mapsto\exp(-t\xi)x_0}$
with $\xi\in\Lambda$.  This is a periodic orbit of a time 
independent periodic Hamiltonian system and the corresponding
Hamiltonian function has a critical point. Connect $x_0$ to that 
critical point by a curve to obtain a homotopy to a 
constant loop.)\label{contractible}
If~${\ox:\bbD\to X}$ is a smooth map satisfying 
\begin{equation}\label{eq:bc}
\ox(e^{2\pi\bi t})=u(t)^{-1}x_0,
\end{equation}
define the {\bf equivariant symplectic action} of the 
triple~${(x_0,u,\ox)}$ by\index{equivariant symplectic action} 
\begin{equation}\label{eq:muaction}
\cA_\mu(x_0,u,\ox) := - \int_\bbD \ox^*\om 
+ \int_0^1\inner{\mu(u^{-1}x_0)}{u^{-1}\dot u}\,dt. 
\end{equation}
The equivariant symplectic action depends only 
on the homotopy class of the loop~$u:\R/\Z\to\rG$
(see part~(i) of Theorem~\ref{thm:RAT} below).

\begin{defn}\label{def:RAT}
Fix a constant $\hbar>0$. The triple $(X,\om,\mu)$
is called {\bf rational (with factor $\hbar$)} if the 
following\index{rational!triple $(X,\om,\mu)$} 
holds.\index{moment map!rational}

\smallskip\noindent{\bf (A)} 
The cohomology class of $\om$ lifts
to $H^2(X;2\pi\hbar\Z)$. 

\smallskip\noindent{\bf (B)} 
$\cA_\mu(x_0,u,\ox)\in 2\pi\hbar\Z$
for  all triples $(x_0,u,\ox)$ satisfying~\eqref{eq:bc}.
\end{defn}

\begin{thm}\label{thm:LIFT}
Let $(X,\om,J)$ be a closed K\"ahler manifold
equip\-ped with a Hamiltonian group action by 
a compact Lie group $\rG$ which is generated by 
an equivariant moment $\mu:X\to\cg$. 
Then the following holds.

\smallskip\noindent{\bf (i)}
Condition~(A) holds if and only if 
there exists a Hermitian line bundle $E\to X$ 
and a Hermitian connection $\nabla$ on $E$
such that the curvature of $\nabla$ is
\begin{equation}\label{eq:Fom}
F^\nabla = \frac{\i}{\hbar}\om.
\end{equation}

\smallskip\noindent{\bf (ii)}
Assume condition~(A) and let $(E,\nabla)$ be as in~(i).
Assume further that~$\rG$ is connected.
Then condition~(B) holds if and only if 
there is a lift of the $\rG$-action on~$X$ 
to a $\rG$-action on~$E$ such that, 
for all $x_0\in X$, $v_0\in E_{x_0}$, and 
$u:\R\to\rG$, the section~${v:=u^{-1}v_0\in E_x}$ 
along the path $x:=u^{-1}x_0$ satisfies 
\begin{equation}\label{eq:LIFT}
\Nabla{t}v = \frac{\i}{\hbar}\inner{\mu(x)}{u^{-1}\dot u}v.
\end{equation}
\end{thm}

\begin{proof}
See page~\pageref{proof:LIFT}.
\end{proof}

\begin{remk}\label{rmk:E}\rm
The complex line bundle $E\to X$ in 
Theorem~\ref{thm:LIFT} has the real first Chern class
$$
c_1(E)=\left[-\frac{\om}{2\pi\hbar}\right] \in H^2(X;\R).
$$ 
Since the curvature of~$\nabla$ is a $(1,1)$-form, 
the Cauchy--Riemann operator
$$
\bar\p^\nabla:\Om^{0,0}(X,E)\to\Om^{0,1}(X,E)
$$
defines a holomorphic structure on~$E$.  Conversely,
for every holomorphic structure on~$E$ there exists
a complex gauge transformation~${g:X\to\C^*}$ 
and a Hermitian connection~$\nabla$ on~$E$ 
such that~${\bar\p^\nabla=g\circ\bar\p\circ g^{-1}}$
and~${F^\nabla = \frac{\i}{\hbar}\om}$. 
In other words, the map~${\nabla\mapsto\bar\p^\nabla}$ 
descends to an isomorphism from the moduli space of
unitary gauge equivalence classes of Hermitian connections 
satisfying~\eqref{eq:Fom}, which is isomorphic to
the torus $H^1(X;\i\R)/H^1(X;2\pi\i\Z)$, to the 
Jacobian of holomorphic structures on~$E$.  
This isomorphism is analogous to the correspondence 
between symplectic and complex quotients in~GIT.
\end{remk}

\begin{cor}\label{cor:RAT}\rm
Assume $\bbP(V)$ is equipped with a $\rG^c$-action,
Fubini--Study form $\om$, and moment map
$\mu$ as in Chapter~\ref{ch:STAG} and let $X$ 
be a $\rG$-invariant complex submanifold of $\bbP(V)$.
Then $(X,\om,\mu)$ is rational.
\end{cor}

\begin{proof}
Condition~(A) holds 
as was noted in the text immediately after
Definition~\ref{def:stable}.  That condition~(B) holds
follows by inserting the formula~\eqref{eq:Vmu} 
for the moment map in equation~\eqref{eq:muaction} 
and choosing $u:\R/\Z\to\rG\subset\rU(V)$ to be a 
loop of the form~${u(t)=\exp(t\xi)}$ with~${\xi\in\Lambda}$.
\end{proof}

\begin{proof}[Proof of Theorem~\ref{thm:LIFT}]\label{proof:LIFT}
Part~(i) follows directly from Chern--Weil theory and the 
classification of complex line bundles by their first Chern class.
The proof of part~(ii) has two steps.

\medskip\noindent{\bf Step~1.}
{\it Let $x_0\in X$, $u:\R/\Z\to\rG$, $\ox:\bbD\to X$ 
with~${\ox(e^{2\pi\i t}) = u(t)^{-1}x_0}$.
Define~${x(t):=u(t)^{-1}x_0}$ and let~${v(t)\in E_{x(t)}}$ 
be a horizontal section. Then}
$$
v(1) = \exp\left(-\frac{\i}{\hbar}\int_\bbD \ox^*\om\right)v(0).
$$

\medskip\noindent
Assume without loss of generality that
$v_0:=v(0)\ne0$, that $\ox(s)=x_0$ for ${0\le s\le1}$,
and that $\ox$ is constant near the origin. Define 
$
z(s,t):=\ox(se^{2\pi\i t})
$ 
for $0\le s,t\le1$ and let $v(s,t)\in E_{z(s,t)}$ 
be the solution of 
$$
\Nabla{t}v=0,\qquad 
v(s,0)=v_0.
$$ 
Since~${v\ne 0}$, there is a unique function~${\lambda:[0,1]^2\to\R}$ 
such that~${\Nabla{s}v=\lambda v}$. It satisfies~${\lambda(s,0)=0}$ and
$$
(\p_t\lambda)v 
= 
\Nabla{t}(\lambda v) 
= 
\Nabla{t}\Nabla{s}v  
= 
\Nabla{s}\Nabla{t}v - F^\nabla(\p_sz,\p_tz)v 
= 
- \frac{\i}{\hbar}\om(\p_sz,\p_tz)v.
$$
Hence
$$
\lambda(s,1) = -\frac{\i}{\hbar}
\int_0^1\om(\p_sz,\p_tz)\,dt
$$
and hence 
$$
\Nabla{s}v(s,1) = \left(-\frac{\i}{\hbar}
\int_0^1\om(\p_sz,\p_tz)\,dt\right) v(s,1).
$$
Thus $v(1,1)=e^{\i\alpha}v(0,1)=e^{\i\alpha}v_0$, where
$$
\i\alpha = -\frac{\i}{\hbar}
\int_0^1\int_0^1\om(\p_sz,\p_tz)\,dt\,ds
= -\frac{\i}{\hbar}\int_\bbD\ox^*\om.
$$
This proves Step~1.

\medskip\noindent{\bf Step~2.}
{\it Let $x_0\in X$ and $u:\R/\Z\to\rG$, 
and define the loop $x:\R/\Z\to X$ by~${x(t):=u(t)^{-1}x_0}$.
If $v(t)\in E_{x(t)}$ is a section of $E$ along $x$ 
satisfying~\eqref{eq:LIFT} then}
$$
v(1) = \exp\left(\frac{\i}{\hbar}\cA_\mu(x_0,u,\ox)\right)v(0).
$$ 

\medskip\noindent
Let $v_0(t)\in E_{x(t)}$ be the horizontal lift
satisfying~${v_0(0)=v(0)=:v_0}$. 
Then we have~${\Nabla{t}v_0\equiv0}$ 
and there exists a function $\lambda:[0,1]\to\R$
such that $v=\lambda v_0$. It satisfies
$$
\dot\lambda v_0
= \Nabla{t}(\lambda v_0)
= \Nabla{t}v
= \frac{\i}{\hbar}\inner{\mu(x)}{u^{-1}\dot u}v
= \frac{\i}{\hbar}\inner{\mu(x)}{u^{-1}\dot u}\lambda v_0.
$$
Hence 
$$
\lambda^{-1}\dot\lambda
=\frac{\i}{\hbar}\inner{\mu(x)}{u^{-1}\dot u},\qquad
\lambda(0)=1,
$$ 
and hence
$$
\lambda(1) = \exp\left(\frac{\i}{\hbar}
\int_0^1\inner{\mu(u^{-1}x_0)}{u^{-1}\dot u}\,dt\right).
$$
By Step~1 and~(B) this implies 
$
v(1)=\lambda(1)v_0(1)
= \exp(\frac{\i}{\hbar}\cA_\mu(x_0,u,\ox))v(0). 
$
This proves Step~2. 

Since $\rG$ is connected it follows from
Step~2 that the $\rG$-action on $X$ lifts to a 
$\rG$-action on $E$ via~\eqref{eq:LIFT}
if and only if $\cA_\mu(x_0,u,\ox)\in 2\pi\hbar\Z$
for  all triples~${(x_0,u,\ox)}$.  
This proves Theorem~\ref{thm:LIFT}.
\end{proof}

The following theorem shows that the equivariant 
symplectic action~\eqref{eq:muaction} is a homotopy 
invariant and gives conditions on~$\om$ 
under which the moment map~$\mu$ can be chosen 
so that the triple~${(X,\om,\mu)}$ is rational.

\begin{thm}\label{thm:RAT}
Assume $\rG$ is connected.
Then the following holds.

\smallskip\noindent{\bf (i)}
The action integral~\eqref{eq:muaction} 
is invariant under homotopy. 

\smallskip\noindent{\bf (ii)}
There exists an integer~${N\in\N}$ such that
each torsion class~${\alpha\in H_1(\rG;\Z)}$
satisfies~${N\alpha=0}$. 

\smallskip\noindent{\bf (iii)}
Let $N\in\N$ be as in~(ii). 
If $\inner{\om}{\pi_2(X)}\subset 2\pi\hbar N\Z$ then 
there exists a central element~${\tau\in\rZ(\cg)}$ 
such that the moment map~${\mu+\tau}$ satisfies~(B).

\smallskip\noindent{\bf (iv)}
Assume~$(X,\om,\mu)$ is rational with factor $\hbar$ 
and let $\tau\in\rZ(\cg)$ so $\mu+\tau$ is an equivariant 
moment map.  Then $(X,\om,\mu+\tau)$ is rational with factor $\hbar$
if and only if~${\inner{\tau}{\xi}\in 2\pi\hbar\Z}$ for every $\xi\in\Lambda$.
\end{thm}

\begin{proof}
Choose functions $x_0:\R\to X$, $u:\R\times\R/\Z\to\rG$,
$\ox:\R\times\bbD\to X$ such that 
$\ox(s,e^{2\pi\i t})=u(s,t)^{-1}x_0(s)=:x(s,t)$
and define~${\xi_s:=u^{-1}\p_su}$ and~${\xi_t:= u^{-1}\p_tu}$.
Then~${\p_s\xi_t-\p_t\xi_s+[\xi_s,\xi_t]=0}$
and~${\p_tx=-L_x\xi_t}$. Differentiate the function 
$$
a(s):=\cA_\mu(x(s),u(s,\cdot),\ox(s,\cdot))
= -\int_\bbD\ox(s,\cdot)^*\om
+ \int_0^1\inner{\mu(x)}{\xi_t}\,dt.
$$
Then
\begin{equation*}
\begin{split}
\dot a
&=
\int_0^1\Bigl(
- \om(\p_sx,\p_tx)
+ \inner{d\mu(x)\p_sx}{\xi_t}
+ \inner{\mu(x)}{\p_s\xi_t} 
\Bigr)\,dt \\
&=
\int_0^1\Bigl(
\om(\p_tx,\p_sx)
+ \om(L_x\xi_t,\p_sx)
+ \inner{\mu(x)}{\p_s\xi_t}
\Bigr)\,dt \\
&=
\int_0^1\inner{\mu(x)}{\p_s\xi_t}\,dt \\
&=
\int_0^1\inner{\mu(x)}{\p_t\xi_s-[\xi_s,\xi_t]}\,dt \\
&=
\int_0^1\Bigl(
\inner{\mu(x)}{\p_t\xi_s}
-\inner{d\mu(x)L_x\xi_t}{\xi_s}
\Bigr)\,dt \\
&=
\int_0^1\Bigl(
\inner{\mu(x)}{\p_t\xi_s}
+\inner{d\mu(x)\p_tx}{\xi_s}
\Bigr)\,dt \\
&=
\int_0^1\p_t\inner{\mu(x)}{\xi_s}\,dt 
=
0.
\end{split}
\end{equation*}
This proves~(i).

Assertion~(ii) follows from the fact that the fundamental
group of $\rG$ is abelian and finitely generated. 
We prove~(iii). Let $H^1(\rG)$ denote the space 
of harmonic $1$-forms on $\rG$ with respect to the Riemannian
metric induced by the invariant inner product on $\cg$.  
Then there is a vector space isomorphism 
$
\rZ(\cg)\to H^1(\rG):\tau\mapsto\alpha_\tau
$
which assigns to each element~${\tau\in\rZ(\cg)}$
the harmonic $1$-form~${\alpha_\tau\in\Om^1(\rG)}$, 
given by~${\alpha_\tau(u,\uhat) := \inner{\tau}{u^{-1}\uhat}}$.
That the map $\tau\mapsto\alpha_\tau$ is surjective 
follows from the fact that the Ricci tensor of $\rG$ 
is nonnegative, so every harmonic $1$-form 
$\alpha$ is parallel.  

Now let $\xi\in\Lambda$ and define the loop 
$u_\xi:\R/\Z\to\rG$ by $u_\xi(t):=\exp(t\xi)$.
An iterate of $u_\xi$ is contractible if and only
if 
$$
\int_{\R/\Z} u_\xi^*\alpha_\tau=\inner{\tau}{\xi}=0
\qquad\mbox{for all }\tau\in\rZ(\cg), 
$$
or equivalently~${\xi\in\rZ(\cg)^\perp=[\cg,\cg]}$.  
Thus the torsion classes in $\pi_1(\rG)$ 
correspond to~${\Lambda\cap\rZ(\cg)^\perp}$ and 
the free part corresponds to~${\Lambda\cap\rZ(\cg)}$.

Choose an integral basis 
${\xi_1,\dots,\xi_k\in\rZ(\cg)\cap\Lambda}$ 
of $\rZ(\cg)$ and a point ${x_0\in X}$. 
Then each loop~${\R/\Z\to X:t\mapsto \exp(-t\xi_j)x_0}$
is contractible (see page~\pageref{contractible}).
Hence there exist functions~${\ox_j:\bbD\to X}$ such 
that~${\ox_j(e^{2\pi\i t})=\exp(-t\xi_j)x_0}$ 
for~${t\in\R}$ and~${j=1,\dots,k}$.  
Define~${\tau\in\rZ(\cg)}$ by
\begin{equation}
\inner{\tau}{\xi_j}
:= \int_\bbD \ox_j^*\om - \inner{\mu(x_0)}{\xi_j}
= -\cA_\mu(x_0,u_{\xi_j},\ox_j),\qquad
j=1,\dots,k.
\end{equation}
We claim that $\mu+\tau$ satisfies~(B). 
By the definition of $\tau$, 
$$
\cA_{\mu+\tau}(x_0,u_{\xi_j},\ox_j)=0,\qquad
j=1,\dots,k.
$$
Since $\inner{\om}{\pi_2(X)}\subset 2\pi\hbar N\Z$ 
this implies that 
${\cA_{\mu+\tau}(x_0,u_{\xi_j},\ox)\in 2\pi\hbar N\Z}$
for every $j$ and every smooth map $\ox:\bbD\to X$ 
such that $\ox(e^{2\pi\i t})=\exp(-t\xi_j)x_0$.
Now let $\xi\in\Lambda\cap\rZ(\cg)^\perp$.
Then the loop $\R/\Z\to X:t\mapsto\exp(-t\xi)x_0$
is contractible and hence there 
is a smooth map $\ox:\bbD\to X$ such that 
$$
{\ox(e^{2\pi\i t})=u_\xi(t)^{-1}x_0=\exp(-t\xi)x_0}
$$ 
for all $t$.  Define~${\ox_N(z):=\ox(z^N)}$. 
Since the action integral is invariant under homotopy, 
the loop~$u_{N\xi}$ is contractible, 
and~${\inner{\om}{\pi_2(X)}\subset 2\pi\hbar N\Z}$,
it follows that~${\cA_\mu(x_0,u_{N\xi},\ox_N)\in 2\pi\hbar N\Z}$. 
Thus
$$
\cA_{\mu+\tau}(x_0,u_\xi,\ox)
= \cA_\mu(x_0,u_\xi,\ox)
= \frac{1}{N}\cA_\mu(x_0,u_{N\xi},\ox_N)\in 2\pi\hbar\Z.
$$
Since the action integral is invariant under homotopy and
additive under catenation it follows that 
$\cA_{\mu+\tau}(x_0,u,\ox)\in2\pi\hbar\Z$ for every triple $(x_0,u,\ox)$.
This proves~(iii). Assertion~(iv) follows directly from the 
definitions and this proves Theorem~\ref{thm:RAT}.
\end{proof}

\begin{thm}\label{thm:LIFTc}
Assume~$(X,\om,\mu)$ is rational with factor $\hbar$
and $\rG$ is connected.  Let~${(E,\nabla)}$ 
be as in part~(i) of Theorem~\ref{thm:LIFT}.
Then the following holds.

\smallskip\noindent{\bf (i)} 
The $\rG$-action on $E$ in part~(ii) of Theorem~\ref{thm:LIFT}
extends uniquely to a $\rG^c$-action on~$E$ by holomorphic 
vector bundle automorphisms.  Moreover,
\begin{equation}\label{eq:LIFTc}
\Nabla{t}v = \frac{\i}{\hbar}\inner{\mu(x)}{\xi}v 
- \frac{1}{\hbar}\inner{\mu(x)}{\eta}v,\qquad
\zeta:=\xi+\i\eta :=g^{-1}\dot g,
\end{equation}
for all $x_0\in X$, $v_0\in E_{x_0}$, $g:\R\to\rG^c$
with $x:=g^{-1}x_0$ and ${v:=g^{-1}v_0\in E_x}$.

\smallskip\noindent{\bf (ii)}
Let $x\in X$ and $v\in E_x\setminus\{0\}$.
The Kempf--Ness function of $x$ is 
$$
\Phi_x(g)=\hbar(\log\abs{g^{-1}v}-\log\abs{v}).
$$

\smallskip\noindent{\bf (iii)}  
Let $x\in X$ and $\zeta\in\Lambda^c$
and define~${x^+:=\lim_{t\to\infty}\exp(\i t\zeta)x}$.
Then 
\begin{equation}\label{eq:wZ}
w_\mu(x,\zeta) = \inner{\mu(x^+)}{\Re(\zeta)}\in 2\pi\hbar\Z
\end{equation}
and the action of $\zeta$ on $E_{x^+}$ is given by
\begin{equation}\label{eq:zetax+}
\exp(-t\zeta)v^+ = e^{\frac{\bi t}{\hbar}w_\mu(x,\zeta)}v^+.
\end{equation}
\end{thm}

\begin{proof}
The Hermitian connection determines a $\rG$-invariant 
integrable complex structure on $E$ via the splitting of the 
tangent bundle into horizontal and vertical subbundles.
Hence the $\rG$-action on $E$ extends to a $\rG^c$-action 
and this proves~(i).  
Part~(ii) follows by computing the derivative of the function 
$$
\Phi_x(g):=\hbar(\log\abs{g^{-1}v} - \log\abs{v}).
$$
Equation~\eqref{eq:zetax+} in~(iii) follows from~\eqref{eq:LIFTc},
and~\eqref{eq:wZ} follows from~\eqref{eq:zetax+}.
This proves Theorem~\ref{thm:LIFTc}.
\end{proof}

\begin{thm}[{\bf Kempf--Ness}]\label{thm:RATSTABLE}
Assume that~$(X,\om,\mu)$ is rational with factor $\hbar$
and that $\rG$ is connected.\index{Kempf--Ness Theorem}  
Let $(E,\nabla)$ be as in part~(i) of Theorem~\ref{thm:LIFT} and let 
the $\rG^c$-action on $E$ be as in Theorem~\ref{thm:LIFTc}.
Let $Z\subset E$ denote the zero section and
let $x\in X$ and 
$$
v\in E_x\setminus\{0\}. 
$$
Then

\smallskip\noindent{\bf (i)}
$x$ is $\mu$-unstable if and only if
$\overline{\rG^c(v)}\cap Z\ne\emptyset$,

\smallskip\noindent{\bf (ii)}
$x$ is $\mu$-semistable if and only if
$\overline{\rG^c(v)}\cap Z=\emptyset$,

\smallskip\noindent{\bf (iii)}
$x$ is $\mu$-polystable if and only if
$\rG^c(v)$ is a closed subset of $E$,

\smallskip\noindent{\bf (iv)}
$x$ is $\mu$-stable if and only if
$\rG^c(v)$ is closed and $\rG^c_v$ is discrete.
\end{thm}

\begin{proof}
The proof is verbatim the same 
as that of Theorem~\ref{thm:VECTOR}.
\end{proof}

\bigbreak

\begin{defn}\label{def:ratinner}
Fix a constant~${\hbar>0}$.\index{rational!inner product}
An invariant inner product on the Lie algebra~$\cg$ 
is called {\bf rational (with factor $\hbar$)} if
\begin{equation}\label{eq:ratinner}
\xi,\eta\in\Lambda,\quad[\xi,\eta]=0
\qquad\implies\qquad
\inner{\xi}{\eta}\in2\pi\hbar\Z.
\end{equation}
If~$\rG$ is a Lie subgroup of~$\rU(n)$ then an example of such an 
inner product is given by 
$$
\inner{\xi}{\eta}:=-2\pi\hbar\frac{\trace(\xi\eta)}{4\pi^2}
$$
for $\xi,\eta\in\cg\subset\cu(n)$.
\end{defn}

\begin{thm}\label{thm:mucrit}
Assume that the triple~$(X,\om,\mu)$ and the inner product on~$\cg$
are rational with factor $\hbar$.   Then, for every critical 
point $x\in X$ of the square of the moment map, there exists 
an integer $k$ such that $k\mu(x)\in\Lambda$, i.e.\ 
$$
L_x\mu(x)=0\qquad\implies\qquad \mu(x)\in\Q\Lambda.
$$
\end{thm}

\begin{proof}
First observe that 
\begin{equation}\label{eq:muxi}
x\in X,\quad\xi\in\Lambda,\quad L_x\xi=0
\qquad\implies\qquad \inner{\mu(x)}{\xi}\in 2\pi\hbar\Z.
\end{equation}
To see this fix an element $x\in X$ and an element 
$\xi\in\Lambda$ such that $L_x\xi=0$.  
Then~${\exp(\i t\xi)x=x}$ for all $t$ and hence 
$\inner{\mu(x)}{\xi}=w_\mu(x,\xi)\in2\pi\hbar\Z$ 
by part~(iii) of Theorem~\ref{thm:LIFTc}. 
This proves~\eqref{eq:muxi}.

Now let $x\in X$ be a critical point of the square 
of the moment map, so 
$$
L_x\mu(x)=0.
$$
Let 
$$
\rT:=\overline{\{\exp(t\mu(x))\,|\,t\in\R\}}\subset\rG
$$ 
be the torus generated by $\mu(x)$.  Then the Lie algebra 
$\ct:=\Lie(\rT)$ is contained in the kernel of $L_x$.  
Choose an integral basis $\xi_1,\dots,\xi_k$ of $\ct\cap\Lambda$. 
Since $\mu(x)\in\ct$ there exist real numbers 
$\lambda_1,\dots,\lambda_k$ such that
$$
\mu(x) = \sum_{i=1}^k \lambda_i\xi_i.
$$
Since $L_x\xi_j=0$, it then follows from~\eqref{eq:muxi} that
$$
\sum_{i=1}^k\lambda_i\frac{\inner{\xi_i}{\xi_j}}{2\pi\hbar} 
= \frac{\inner{\mu(x)}{\xi_j}}{2\pi\hbar} \in \Z
$$
for all $j$.  Since the integer matrix 
$((2\pi\hbar)^{-1}\inner{\xi_i}{\xi_j})_{i,j=1}^k$
is invertible, it follows that~${\lambda_i\in\Q}$ for all $i$.
Hence $\mu(x)\in\Q\Lambda$ and this proves
Theorem~\ref{thm:mucrit}.
\end{proof}


\chapter{The dominant $\mu$-weight}\label{ch:DOMINANT}


The\index{weight@$\mu$-weight!dominant|(}
moment-weight inequality in Theorem~\ref{thm:MW2}
can be stated in the form
\begin{equation}\label{eq:MW3}
\sup_{0\ne\xi\in\cg}\frac{-w_\mu(x,\xi)}{\abs{\xi}} 
\le \inf_{g\in\rG^c}\abs{\mu(gx)}.
\end{equation}
The main results of the present chapter assert that the supremum 
on the left in~\eqref{eq:MW3} is always attained (Theorem~\ref{thm:SUPMAX}),
that it is attained at a unique element~$\xi_0$ up to scaling 
whenever~$x$ is $\mu$-unstable (Theorem~\ref{thm:KEMPF1}),
that the inequality is strict in the $\mu$-stable case (Theorem~\ref{thm:MODSTAB})
and that equality holds in~\eqref{eq:MW3} 
in the $\mu$-unstable case (Theorem~\ref{thm:KEMPF}).
That equality also holds when $x$ is $\mu$-semistable, 
but not $\mu$-stable, will follow from the Hilbert--Mumford 
criterion (Corollary~\ref{cor:M}).
In~\cite{KEMPF} Kempf proved that the supremum in~\eqref{eq:MW3}  
is attained at a unique element~${\xi_0\in\Lambda}$ up to scaling, 
whenever $x$ is $\mu$-unstable and the triple $(X,\om,\mu)$ and 
the inner product on~$\cg$ are rational (Corollary~\ref{cor:KEMPF}). 
In general, the ray in~${\cg\setminus\{0\}}$ along which the
supremum is attained need not intersect $\Lambda$. 

\begin{thm}\label{thm:SUPMAX}
Let~${x_0\in X}$.  Then there exists an element~${\xi_0\in\cg}$ such that
\begin{equation}\label{eq:SUPMAX}
\Abs{\xi_0}=1,\qquad 
-w_\mu(x_0,\xi_0) = \sup_{0\ne\xi\in\cg}\frac{-w_\mu(x_0,\xi)}{\abs{\xi}}.
\end{equation}
\end{thm}

\begin{proof}
See page~\pageref{proof:SUPMAX}.
\end{proof}

\begin{thm}[{\bf Generalized Kempf Uniqueness Theorem}]\label{thm:KEMPF1}
Assume that~${x_0\in X}$ is $\mu$-unstable.  Then $\xi_0\in\cg$ 
is uniquely determined by~\eqref{eq:SUPMAX}.
\end{thm}

\begin{proof}
See page~\pageref{proof:KEMPF1}.
\end{proof}

\begin{thm}\label{thm:MODSTAB}
Assume $x_0$ is $\mu$-stable. Then 
\begin{equation}\label{eq:MODSTAB2}
0 < \lambda(x_0)
:= \inf_{0\ne\xi\in\cg}\frac{w_\mu(x_0,\xi)}{\abs{\xi}}
\le \inf_{x\in\overline{\rG^c(x_0)}\setminus\rG^c(x_0)}\Abs{\mu(x)},
\end{equation}
the set~${\overline{\rG^c(x_0)}\setminus\rG^c(x_0)}$ is compact,
every~${x\in\overline{\rG^c(x_0)}\setminus\rG^c(x_0)}$ 
is $\mu$-un\-stable, and~${\inf_{g\in\rG^c}\Abs{\mu(gx^+)}\ge\lambda(x_0)}$
for~${\zeta\in\sT^c}$ and~${x^+:=\lim_{t\to\infty}\exp(\i t\zeta)x_0}$.
\end{thm}

\begin{proof}
See page~\pageref{proof:MODSTAB}.
\end{proof}

It follows from Theorem~\ref{thm:MUMFORD1} 
and Lemma~\ref{le:WEIGHT3} that
\begin{equation}\label{eq:MODSTAB3}
\lambda(x_0)
:= 
\inf_{0\ne\xi\in\cg}\frac{w_\mu(x_0,\xi)}{\abs{\xi}}
= 
\inf_{\zeta\in\sT^c}
\frac{w_\mu(x,\zeta)}{\sqrt{\abs{\Re(\zeta)}^2-\abs{\Im(\zeta)}^2}}
\end{equation}
for all $x_0\in X$. 
Thus the function $\lambda:X\to\R$ is $\rG^c$-invariant.
By Theorem~\ref{thm:MODSTAB} the number $\lambda(x_0)$ 
is positive whenever $x_0\in X$ is~{$\mu$-stable}.  
In the work of Sz\'ekelyhidi~\cite{S} the number $\lambda(x_0)$ 
is called the {\bf modulus of stability}\index{modulus of stability} 
in the $\mu$-stable case.  

Another proof of Kempf Uniqueness for torus actions is 
contained in Lemma~\ref{le:TORUS} in the next section.  
In the $\mu$-unstable case an alternative proof of 
Theorem~\ref{thm:SUPMAX} is given in the following theorem,
which also establishes equality in~\eqref{eq:MW3} and is 
the main result of this chapter.

\begin{thm}[{\bf Generalized Kempf Existence Theorem}]\label{thm:KEMPF}
Assume that~${x_0\in X}$\index{Kempf Existence Theorem}
is $\mu$-unstable so that
\begin{equation}\label{eq:m}
m := \inf_{g\in\rG^c}\abs{\mu(gx_0)} > 0.
\end{equation}
Let $x:\R\to X$ be the unique solution of~\eqref{eq:KN1} 
and define~${x_\infty:=\lim_{t\to\infty}x(t)}$.
Let $g:\R\to\rG^c$ be the solution of~\eqref{eq:KN2} 
so that ${x(t)=g(t)^{-1}x_0}$ for all~$t$.
Define the functions $\R\to\cg:t\mapsto\xi(t)$ 
and $\R\to\rG:t\mapsto u(t)$ by
\begin{equation}\label{eq:gxih}
g(t) =: \exp(-\i\xi(t))u(t).
\end{equation}
Then the limit
\begin{equation}\label{eq:xinfty1}
\xi_\infty := \lim_{t\to\infty}\frac{\xi(t)}{t}
\end{equation}
exists and satisfies
\begin{equation}\label{eq:xinfty2}
w_\mu(x_0,\xi_\infty) = -m^2,\qquad 
\abs{\xi_\infty} = m.
\end{equation}
Moreover, there exists an element $u_\infty\in\rG$ such that
\begin{equation}\label{eq:xinfty3}
\xi_\infty = -u_\infty^{-1}\mu(x_\infty)u_\infty.
\end{equation}
\end{thm}

\begin{proof}
See page~\pageref{proof:KEMPF}.
\end{proof}

\bigbreak

The proof of Theorem~\ref{thm:KEMPF} is essentially due 
to Chen--Sun~\cite[Theorems~4.4 and~4.5]{CS}. 
They use the same argument to establish the existence 
of negative weights for linear actions on projective space 
under the assumption~\eqref{eq:m}.  Their~$\chi$ is our~$\xi_\infty$
and their~$\gamma$ is the negative gradient flow line of the 
Kempf--Ness function going by the same name 
in the proof below.

\begin{figure}[htp] 
\centering 
\includegraphics{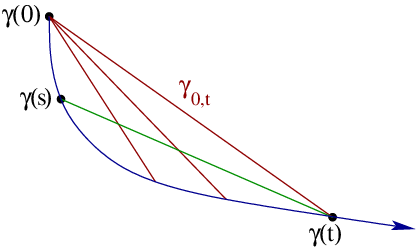} 
\caption{{A negative gradient flow line and geodesics.}}
\label{fig:geograd}      
\end{figure}

\begin{proof}[Proof of Theorem~\ref{thm:KEMPF}]\label{proof:KEMPF}
Let~$\gamma$\index{Kempf Existence Theorem!Chen--Sun's proof|(} 
denote the image of $g$ in $M=\rG^c/\rG$, i.e.\ 
$$
\gamma(t) := \pi(g(t)) = \pi(\exp(-\i\xi(t)))
$$
for $t\in\R$. Let $\nabla$ be the Levi-Civita connection on $M$. 
For $0\le s<t$ let 
$$
\gamma_{s,t}:[s,t]\to M
$$  
be the geodesic connecting the points 
$$
\gamma_{s,t}(s)=\gamma(s),\qquad 
\gamma_{s,t}(t)=\gamma(t).  
$$
It is given by
\begin{equation}\label{eq:gast}
\gamma_{s,t}(r):=\pi\left(
g(s)\exp\left(-\i\frac{r-s}{t-s}\xi(s,t)\right)
\right)
\qquad\mbox{for }s\le r\le t,
\end{equation}
where $\xi(s,t)\in\cg$ and $u(s,t)\in\rG$ are chosen such that
\begin{equation}\label{eq:xist}
g(s)\exp\bigl(-\i\xi(s,t)\bigr)u(s,t) = g(t).
\end{equation}
(See Figure~\ref{fig:geograd}.)  For $0\le s<t$ define 
the function 
$
\rho_{s,t}:[s,t]\to[0,\infty)
$ 
by 
\begin{equation}\label{eq:rhost}
\rho_{s,t}(r):=d_M(\gamma_{s,t}(r),\gamma(r))
\qquad\mbox{for }s\le r\le t.
\end{equation}
With this notation in place, we prove the assertions in nine steps.

\bigbreak

\medskip\noindent{\bf Step~1.}
{\it For every $t>0$ we have}
$$
\abs{\nabla\dot\gamma(t)}
= \abs{{L_{x(t)}}^*L_{x(t)}\mu(x(t))}.
$$

\medskip\noindent
By~\eqref{eq:KN2}, we have 
$
g(t)^{-1}\dot g(t)=\i\mu(x(t))\in\i\cg
$ 
for all $t$. Hence Theorem~\ref{thm:GcG} asserts that 
$
\nabla\dot\gamma
= d\pi(g)g\i d\mu(x)\dot x
= -d\pi(g)g\i{L_x}^*L_x\mu(x).
$
This proves Step~1.

\medskip\noindent{\bf Step~2.}
{\it There exist positive constants $c$ and $\eps$ 
such that, for all $t>0$,}
$$
\int_t^\infty \abs{{L_{x(r)}}^*L_{x(r)}\mu(x(r))}\,dr
\le \frac{c}{t^\eps}.
$$

\medskip\noindent
This follows directly from Theorem~\ref{thm:XINFTY}.

\medskip\noindent{\bf Step~3.}
{\it Let~$c$ and~$\eps$ be as in Step~2 and fix three real numbers $r_0,s,t$. Then}
\begin{equation}\label{eq:rst}
0\le s<r_0<t,\quad \rho_{s,t}(r_0)\ne0
\qquad\implies\qquad
\dot\rho_{s,t}(r_0)\le\frac{c}{r_0^\eps}.
\end{equation}

\medskip\noindent
Assume $s<r_0<t$ and $\rho_{s,t}(r_0)\ne0$. 
Let~$r_1$ be the smallest real number bigger 
than~$r_0$ such that $\rho_{s,t}(r_1)=0$.  
Thus 
$$
r_0<r_1\le t,\qquad \rho_{s,t}(r_1)=0, 
$$
and
$$
\rho_{s,t}(r)\ne 0\qquad\mbox{for }r_0\le r<r_1.  
$$
Hence by Lemma~\ref{le:distance}
$$
- \ddot\rho_{s,t}(r)
\le \abs{\nabla\dot\gamma(r)}
= \abs{{L_{x(r)}}^*L_{x(r)}\mu(x(r))}
$$
for $r_0\le r<r_1$. Here the equality follows from Step~1. 
Integrate this inequality over the interval~${r_0\le r<r_1}$ 
and use Lemma~\ref{le:distance} to obtain
\begin{equation*}
\begin{split}
\dot\rho_{s,t}(r_0) 
&=
\frac{d\rho_{s,t}}{dr^-}(r_1) 
- \int_{r_0}^{r_1}\ddot\rho_{s,t}(r)\,dr \\
&= 
- \abs{\dot\gamma(r_1)-\dot\gamma_{s,t}(r_1)}
- \int_{r_0}^{r_1}\ddot\rho_{s,t}(r)\,dr \\
&\le
\int_{r_0}^{r_1}\abs{{L_{x(r)}}^*L_{x(r)}\mu(x(r))}\,dr \\
&\le
\frac{c}{r_0^\eps}.
\end{split}
\end{equation*} 
Here the last inequality follows from Step~2.
This proves Step~3.

\bigbreak

\medskip\noindent{\bf Step~4.}
{\it Let~${\eps>0}$ be as in Step~2. Then there exists a constant $C>0$ such that
$$
\rho_{s,t}(r)\le C\left(r^{1-\eps}-s^{1-\eps}\right)
$$
for all real numbers $r,s,t$ such that $0\le s\le r\le t$.}

\medskip\noindent
Fix a real number $r_1$ such that $s<r_1<t$.  Suppose without loss
of generality that~${\rho_{s,t}(r_1)\ne 0}$ and choose~$r_0$ such 
that~${s\le r_0<r_1}$, ${\rho_{s,t}(r_0)=0}$, and ${\rho_{s,t}(r)\ne0}$ 
for~${r_0<r\le r_1}$.  Now integrate the inequality 
$$
\dot\rho_{s,t}(r)\le cr^{-\eps}
$$ 
in Step~3 over the interval~${r_0<r\le r_1}$ to obtain
$$
\rho_{s,t}(r_1) = \int_{r_0}^{r_1}\dot\rho_{s,t}(r)\,dr
\le \int_s^{r_1}\frac{cdr}{r^\eps}
= \frac{c}{1-\eps}\left({r_1}^{1-\eps}-s^{1-\eps}\right).
$$
This proves Step~4 with $C:=c/(1-\eps)$.

\medskip\noindent{\bf Step~5.} 
{\it Let $s\ge 0$. Then
$$
\Abs{\frac{\xi(s,t')}{t'-s}-\frac{\xi(s,t)}{t-s}} 
\le \frac{C}{t^\eps}
\qquad\mbox{for all }t'\ge t\ge s+1
$$
Thus the limit 
$$ 
\xi_\infty(s) := \lim_{t\to\infty} \frac{\xi(s,t)}{t-s}
$$
exists in $\cg$.}

\medskip\noindent
The geodesics $\gamma_{s,t}$ and $\gamma_{s,t'}$ intersect at
$\gamma(s)=\gamma_{s,t}(s)=\gamma_{s,t'}(s)$ (see Figure~\ref{fig:geograd}).
Hence it follows from equation~\eqref{eq:gast}, Lemma~\ref{le:expand},
and Step~4 that
\begin{equation*}
\begin{split}
\Abs{\frac{\xi(s,t')}{t'-s}-\frac{\xi(s,t)}{t-s}}
&=
\Abs{\dot\gamma_{s,t}(s)-\dot\gamma_{s,t'}(s)} \\
&\le 
\frac{d_M\left(\gamma_{s,t}(t),\gamma_{s,t'}(t)\right)}{t-s} \\
&= 
\frac{\rho_{s,t'}(t)}{t-s} \\
&\le 
C\frac{t^{1-\eps}-s^{1-\eps}}{t-s} \\
&\le
\frac{C}{t^\eps}
\end{split}
\end{equation*}
for~${t'\ge t\ge s+1}$.  This proves Step~5. 

\bigbreak

\medskip\noindent{\bf Step~6.} 
{\it $w_\mu(x(s),\xi_\infty(s))=-m^2$ for all $s\ge0$.}

\medskip\noindent
By Theorem~\ref{thm:MLT}, we have $\abs{\mu(x_\infty)}=m$.
Fix a real number $s\ge0$ and define the geodesic 
$\gamma_{s,\infty}:[s,\infty)\to M$ by 
$$
\gamma_{s,\infty}(r):=\pi(g(s)\exp(-\i (r-s)\xi_\infty(s)))
= \lim_{t\to\infty}\gamma_{s,t}(r)
$$
for $r\ge s$.  By Step~4, we have
$$
d_M\left(\gamma(r),\gamma_{s,t}(r)\right) 
= \rho_{s,t}(r)
\le C\bigl(r^{1-\eps}-s^{1-\eps}\bigr)
$$
for~${s\le r\le t}$. Take the limit $t\to\infty$ to obtain
\begin{equation}\label{eq:dist1}
d_M\left(\gamma(r),\gamma_{s,\infty}(r)\right) 
\le C\bigl(r^{1-\eps}-s^{1-\eps}\bigr)
\end{equation}
for $r\ge s\ge 0$.  Now the Kempf--Ness function 
is globally Lipschitz continuous with Lipschitz 
constant~${L := \sup_{g\in\rG^c}\abs{\mu(gx_0)}}$.
Hence it follows from~\eqref{eq:dist1} that
\begin{equation}\label{eq:dist2}
\Abs{\Phi_{x_0}(\gamma(t)) - \Phi_{x_0}(\gamma_{s,\infty}(t))}
\le LC\bigl(t^{1-\eps}-s^{1-\eps}\bigr)
\qquad\mbox{for }t\ge s\ge0.
\end{equation}
Integrate the equation 
$$
\frac{d}{dr}\Phi_{x_0}(\gamma(r))=-\abs{\mu(x(r))}^2
$$
to obtain
\begin{equation}\label{eq:dist3}
\Phi_{x_0}(\gamma(t)) 
= \Phi_{x_0}(\gamma(s)) 
- \int_s^t\abs{\mu(x(r))}^2\,dr.
\end{equation}
By Lemma~\ref{le:SLOPE} and~\eqref{eq:dist2}
and~\eqref{eq:dist3} we have
\begin{equation*}
\begin{split}
w_\mu(x(s),\xi_\infty(s))
&=
\lim_{t\to\infty} 
\frac{\Phi_{x_0}(\gamma_{s,\infty}(s+t))}{t} \\
&=
\lim_{t\to\infty} 
\frac{\Phi_{x_0}(\gamma_{s,\infty}(t))-\Phi_{x_0}(\gamma(s))}{t-s} \\
&=
\lim_{t\to\infty} 
\frac{\Phi_{x_0}(\gamma(t))-\Phi_{x_0}(\gamma(s))}{t-s} \\
&=
-\lim_{t\to\infty} 
\frac{1}{t-s}\int_s^t\abs{\mu(x(r))}^2\,dr \\
&=
- \abs{\mu(x_\infty)}^2 \\
&=
- m^2.
\end{split}
\end{equation*}
This proves Step~6.

\medskip\noindent{\bf Step~7.} 
{\it $\abs{\xi_\infty(s)}=m$ for all $s\ge0$.}

\medskip\noindent
By definition of $\xi(s,t)$ in~\eqref{eq:xist} we have 
\begin{equation*}
\begin{split}
\frac{\abs{\xi(s,t)}}{t-s}
&=
\frac{d_M(\gamma(s),\gamma(t))}{t-s} \\
&\le
\frac{1}{t-s}\int_s^t\abs{\dot\gamma(r)}\,dr \\
&=
\frac{1}{t-s}\int_s^t\abs{\mu(x(r))}\,dr.
\end{split}
\end{equation*}
Take the limit $t\to\infty$.  Then
\begin{equation*}
\begin{split}
\abs{\xi_\infty(s)} 
\le 
\lim_{t\to\infty}
\frac{1}{t-s}\int_s^t\abs{\mu(x(r))}\,dr 
= 
\abs{\mu(x_\infty)} 
= 
m.
\end{split}
\end{equation*}
Moreover, it follows from the moment-weight 
inequality in Theorem~\ref{thm:MW2} that 
\begin{equation*}
\begin{split}
m^2 
= 
-w(x(s),\xi_\infty(s))  
\le 
\abs{\xi_\infty(s)}\inf_{g\in\rG^c}\abs{\mu(gx_0)} 
= 
m\abs{\xi_\infty(s)}
\end{split}
\end{equation*}
and hence $\abs{\xi_\infty(s)}\ge m$.
This proves Step~7.

\medskip\noindent{\bf Step~8.} 
{\it For every $s\ge0$ there exists an element $u_\infty(s)\in\rG$
such that} 
$$
\xi_\infty(0)=u_\infty(s)^{-1}\xi_\infty(s)u_\infty(s).
$$

\medskip\noindent
The geodesics $\gamma_{s,t}$ and $\gamma_{0,t}$ intersect 
at the point~${\gamma(t)=\gamma_{s,t}(t)=\gamma_{0,t}(t)}$ 
(see Figure~\ref{fig:geograd}).
Hence it follows from~\eqref{eq:gast}, \eqref{eq:xist} 
and Lemma~\ref{le:expand} that 
\begin{equation*}
\begin{split}
\Abs{u(s,t)^{-1}\frac{\xi(s,t)}{t-s}u(s,t)
- u(0,t)^{-1}\frac{\xi(0,t)}{t}u(0,t)}
&=
\abs{\dot\gamma_{s,t}(t)-\dot\gamma_{0,t}(t)} \\
&\le
\frac{d_M(\gamma_{s,t}(s),\gamma_{0,t}(s))}{t-s}  \\
&= 
\frac{\rho_{0,t}(s)}{t-s}  \\
&\le
\frac{Cs^{1-\eps}}{t-s}
\end{split}
\end{equation*}
for $t\ge s+1\ge0$.
Here the last inequality follows from Step~4.
Now choose a sequence $t_i\to\infty$ such that
the limit~${u_\infty(s) := \lim_{i\to\infty}u(s,t_i)u(0,t_i)^{-1}}$
exists.  Then, by Step~5, 
we have~${\xi_\infty(0) = \lim_{i\to\infty}t_i^{-1}\xi(0,t_i)}$
and hence
\begin{equation*}
\begin{split}
\xi_\infty(0)
&= 
\lim_{i\to\infty}
u(0,t_i)u(s,t_i)^{-1}\frac{\xi(s,t_i)}{t_i-s}u(s,t_i)u(0,t_i)^{-1} \\
&= 
u_\infty(s)^{-1}\xi_\infty(s)u_\infty(s).
\end{split}
\end{equation*}
This proves Step~8.

\bigbreak

\medskip\noindent{\bf Step~9.} 
{\it $\lim_{s\to\infty}\xi_\infty(s)=-\mu(x_\infty)$.}

\medskip\noindent
By Step~5 the inequality
$$
\Abs{\frac{\xi(s,t')}{t'-s}-\frac{\xi(s,t)}{t-s}} 
\le \frac{C}{t^\eps} \le\frac{C}{s^\eps}
$$
holds for $t'\ge t\ge s+1$. Take the limit $t'\to\infty$ to obtain
$$
\Abs{\xi_\infty(s)-\xi(s,s+1)} \le \frac{C}{s^\eps}
\qquad\mbox{for all }s\ge 0.
$$
Now suppose~${\xi(s,s+1)+\mu(x(s))\ne 0}$.
Then~${\dot\gamma(s)\ne\dot\gamma_{s,s+1}(s)}$, 
so there is a real number~$t$ such that~${s<t\le s+1}$,
${\rho_{s,s+1}(t)=0}$, and~${\rho_{s,s+1}(r)\ne 0}$ for~${s<r<t}$.
By Step~3 this implies~${\dot\rho_{s,s+1}(r)\le cr^{-\eps}}$ 
for~${s<r<t}$ and hence, by Lemma~\ref{le:distance}, 
\begin{equation*}
\begin{split}
\abs{\mu(x(s))+\xi(s,s+1)}
= 
\abs{\dot\gamma(s)-\dot\gamma_{s,s+1}(s)}  
= 
\lim_{r\searrow s}\dot\rho_{s,s+1}(r)  
\le 
\frac{c}{s^\eps}.
\end{split}
\end{equation*}
Hence~${\Abs{\xi_\infty(s)+\mu(x(s))}\le\frac{c+C}{s^\eps}}$
for all~${s\ge 0}$ and this proves Step~9.

The existence of the limit in~\eqref{eq:xinfty1} follows from Step~5,
that it satisfies~\eqref{eq:xinfty2} follows from Steps~6 and~7, 
and that it satisfies~\eqref{eq:xinfty3} follows from Steps~8 and~9.
This\index{Kempf Existence Theorem!Chen--Sun's proof|)} 
proves Theorem~\ref{thm:KEMPF}.
\end{proof}

In preparation for the proofs of  the remaining theorems in this chapter
we establish a convergence result for the Kempf--Ness function.

\begin{lem}\label{le:KNSPHERE}
Let~${x_0\in X}$, define
\begin{equation}\label{eq:lambda}
\lambda(x_0):=\inf_{0\ne\xi\in\cg}\frac{w_\mu(x_0,\xi)}{\abs{\xi}},
\end{equation}
denote by~${\Phi_{x_0}:\rG^c\to\R}$ the 
lifted Kempf--Ness function, and define
$$
S_t := \left\{\exp(-\i t\xi)\,\big|\,\xi\in\cg,\,\Abs{\xi}=1\right\}
$$
for~${t>0}$.  Then~${\lambda(x_0)>-\infty}$ and the following holds.

\smallskip\noindent{\bf (i)}
$\frac{1}{t}\inf_{S_t}\Phi_{x_0}\le\lambda(x_0)$ for all~${t>0}$.

\smallskip\noindent{\bf (ii)}
There exists a sequence $t_i>0$ of positive real numbers
and a convergent sequence $\xi_i\in\cg$ such that
$$
\Abs{\xi_i}=1,\qquad 
\Phi_{x_0}(\exp(-\i t_i\xi_i)) = \inf_{S_{t_i}}\Phi_{x_0},\qquad
\lim_{i\to\infty}t_i=\infty.
$$
The limit~${\xi_0:=\lim_{i\to\infty}\xi_i}$
of any such sequence satisfies~\eqref{eq:SUPMAX}.

\smallskip\noindent{\bf (iii)}
$\lim_{t\to\infty}t^{-1}\inf_{S_t}\Phi_{x_0}=\lambda(x_0)$.
\end{lem}

\begin{proof}
We prove part~(i).  The number~$\lambda(x_0)$ 
in~\eqref{eq:lambda} is finite by the moment-weight 
inequality~\eqref{eq:MW2} in Theorem~\ref{thm:MW2}. 
Fix two real numbers~${t>0}$ and~${\eps>0}$, 
and choose~${\xi\in\cg}$ such that
$$
\abs{\xi}=1,\qquad 
w_\mu(x_0,\xi)<\lambda(x_0)+\eps. 
$$
Since the function~${s\mapsto\inner{\mu(\exp(\i s\xi)x_0)}{\xi}}$ is 
nondecreasing and converges to the weight~${w_\mu(x_0,\xi)}$
as~$s$ tends to infinity, we have 
\begin{equation*}
\begin{split}
\frac{1}{t}\Phi_{x_0}(\exp(-\i t\xi))
= 
\frac{1}{t}\int_0^t\inner{\mu(\exp(\i s\xi)x_0)}{\xi}\,ds 
\le 
w_\mu(x_0,\xi) 
< 
\lambda(x_0)+\eps.
\end{split}
\end{equation*}
Here the first equality follows from~\eqref{eq:d1Phi}
and the fact that~${\Phi_{x_0}(\one)=0}$.
Since~${\eps>0}$ was chosen arbitrary, 
this proves~(i).

We prove part~(ii).
Choose any sequence $t_i>0$ that tends to infinity.
Since $\Phi_{x_0}$ is continuous and~$S_{t_i}$ is compact, 
the axiom of countable choice asserts that there exists 
a sequence~${\xi_i\in\cg}$ such that 
$$
\abs{\xi_i}=1,\qquad
\Phi_{x_0}(\exp(-\i t_i\xi_i)) = \inf_{S_{t_i}}\Phi_{x_0}
$$
for all~$i$.  Passing to a subsequence, we may assume 
that the sequence~$(\xi_i)_{i\in\N}$ converges.  
Denote the limit by~${\xi_0:=\lim_{i\to\infty}\xi_i}$. 
Then 
$$
\Abs{\xi_0}=\lim_{i\to\infty}\Abs{\xi_i}=1.
$$  
Moreover, by Lemma~\ref{le:SLOPE}, 
the function~${t\mapsto t^{-1}\Phi_{x_0}(\exp(-\i t\xi_i))}$
is nondecreasing and hence by part~(i) we have
$$
\frac{1}{t}\Phi_{x_0}(\exp(-\i t\xi_i))
\le 
\frac{1}{t_i}\inf_{S_{t_i}}\Phi_{x_0}
\le 
\lambda(x_0)
\qquad\mbox{for }0<t\le t_i. 
$$
Take the limit $i\to\infty$ to obtain 
$$
\frac{1}{t}\Phi_{x_0}(\exp(-\i t\xi_0))\le \lambda(x_0)
$$
for all $t>0$. By Lemma~\ref{le:SLOPE}, this implies 
$$
w_\mu(x_0,\xi_0)
= 
\lim_{t\to\infty}\frac{1}{t}\Phi_{x_0}(\exp(-\i t\xi_0))
\le
\lambda(x_0)
$$
and hence
$$
w_\mu(x_0,\xi_0)=\lambda(x_0)
$$
by definition of~$\lambda(x_0)$ in~\eqref{eq:lambda}.
This proves~(ii).

\bigbreak

We prove part~(iii).
Choose $t_i$, $\xi_i$, and~${\xi_0=\lim_{i\to\infty}\xi_i}$ 
as in~(ii). Then
$$
\lambda(x_0)
= w_\mu(x_0,\xi_0)
= \lim_{t\to\infty}\frac{1}{t}\Phi_{x_0}(\exp(-\i t\xi_0))
$$ 
by Lemma~\ref{le:SLOPE}. Fix a constant $\eps>0$ and choose
$t_0>0$ such that 
$$
\frac{1}{t_0}\Phi_{x_0}(\exp(-\i t_0\xi_0)) > \lambda(x_0)-\eps.
$$
Then there exists a constant $i_0>0$ such that,
for all $i\in\N$,
$$
i\ge i_0\qquad\implies\qquad
t_i\ge t_0\quad\mbox{and}\quad
\frac{1}{t_0}\Phi_{x_0}(\exp(-\i t_0\xi_i)) > \lambda(x_0)-\eps.
$$
Hence, for every $i\in\N$ with $i\ge i_0$, we have
\begin{equation*}
\begin{split}
\lambda(x_0)-\eps
&<
\frac{1}{t_0}\Phi_{x_0}(\exp(-\i t_0\xi_i)) \\
&\le
\frac{1}{t_i}\Phi_{x_0}(\exp(-\i t_i\xi_i)) \\
&=
\frac{1}{t_i}\inf_{S_{t_i}}\Phi_{x_0} \\
&\le 
\lambda(x_0).
\end{split}
\end{equation*}
Here the second inequality holds 
because~${t\mapsto t^{-1}\Phi_{x_0}(\exp(-\i t\xi_i))}$ 
is nondecreasing by Lemma~\ref{le:SLOPE},
the equality follows from the choice of the sequence~$\xi_i$
in part~(ii), and the last inequality follows from part~(i). 
Thus
\begin{equation}\label{eq:MODSTAB5}
\lim_{i\to\infty}\frac{1}{t_i}\inf_{S_{t_i}}\Phi_{x_0}=\lambda(x_0).
\end{equation}
Now Lemma~\ref{le:SLOPE} asserts that
$$
\frac{1}{t_i}\Phi_{x_0}(\exp(-\i t_i\xi))\le \frac{1}{t}\Phi_{x_0}(\exp(-\i t\xi))
$$
for all $i\in\N$, all $t\ge t_i$, and all $\xi\in\cg$ with $\Abs{\xi}=1$.  
Take the infimum over all~${\xi\in\cg}$ with $\abs{\xi}=1$ and use 
part~(i) to obtain
$$
\frac{1}{t_i}\inf_{S_{t_i}}\Phi_{x_0}\le \frac{1}{t}\inf_{S_t}\Phi_{x_0}\le\lambda(x_0)
\qquad\mbox{for all }i\in\N\mbox{ and all }t\ge t_i.
$$
By~\eqref{eq:MODSTAB5} this proves~(iii) and Lemma~\ref{le:KNSPHERE}.
\end{proof}

\begin{proof}[Proof of Theorem~\ref{thm:SUPMAX}]
\label{proof:SUPMAX}
This follows from part~(ii) of Lemma~\ref{le:KNSPHERE}.
\end{proof}

\begin{figure}[htp] 
\centerline{\psfig{figure=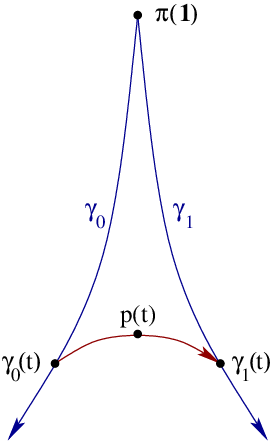,height=2.4in}} 
\caption{Proof of Kempf uniqueness.}
\label{fig:kempf}      
\end{figure}

\begin{proof}[Proof of Theorem~\ref{thm:KEMPF1}]
\label{proof:KEMPF1}
Assume $x_0$ is $\mu$-unstable.  Then
$$
\lambda(x_0) 
= \inf_{0\ne\xi\in\cg}\frac{w_\mu(x_0,\xi)}{\abs{\xi}} 
\le \frac{w_\mu(x_0,\xi_\infty)}{\abs{\xi_\infty}}
< 0 
$$ 
by Theorem~\ref{thm:KEMPF}.  Define~${m:=-\lambda(x_0)}$,
let~${\xi_0,\xi_1\in\cg}$ such that
\begin{equation}\label{eq:xi01}
\Abs{\xi_0}=\Abs{\xi_1}=1,\qquad
w_\mu(x_0,\xi_0) = w_\mu(x_0,\xi_1) = \lambda(x_0) = -m,
\end{equation}
and consider the geodesics
$$
\gamma_0(t):= \pi\bigl(\exp(-\i t\xi_0)\bigr),\qquad
\gamma_1(t):= \pi\bigl(\exp(-\i t\xi_1)\bigr).
$$
(See Figure~\ref{fig:kempf}.)
For $t>0$ choose $\eta(t)\in\cg$ and $u(t)\in\rG$ such that
$$
\exp(-\i t\xi_0)\exp(\i\eta(t)) = \exp(-\i t\xi_1)u(t)
$$
and define 
$$
p(t) := \pi\bigl(\exp(-\i t\xi_0)\exp(\i\tfrac12\eta(t))\bigr)
\in M = \rG^c/\rG.
$$
Then $p(t)$ is the midpoint of the geodesic 
joining~$\gamma_0(t)$ and~$\gamma_1(t)$. 
Hence it follows from Lemma~\ref{le:ALEXANDROV} that
\begin{equation*}
\begin{split}
d(\pi(\one),p(t))^2
&\le
\frac{d(\pi(\one),\gamma_0(t))^2 + d(\pi(\one),\gamma_1(t))^2}{2}
- \frac{d(\gamma_0(t),\gamma_1(t))^2}{4} \\
&=
t^2 - \frac{d(\gamma_0(t),\gamma_1(t))^2}{4} \\
&\le
t^2\left(1 - \frac{\abs{\xi_0-\xi_1}^2}{4}\right).
\end{split}
\end{equation*}
The last inequality holds for $t\ge 1$, by Lemma~\ref{le:expand}.
Thus 
\begin{equation}\label{eq:rt}
\frac{r(t)}{t} \le \sqrt{1 - \frac{\abs{\xi_0-\xi_1}^2}{4}},\qquad
r(t):=d(\pi(\one),p(t)),
\end{equation}
for $t\ge 1$.
Moreover, 
$$
\Phi_{x_0}(\gamma_0(t))\le -tm,\qquad
\Phi_{x_0}(\gamma_1(t))\le -tm
$$
for all~${t\ge 0}$ by equation~\eqref{eq:xi01} and Lemma~\ref{le:SLOPE}.  
Since $\Phi_{x_0}$ is convex along geodesics, this implies 
\begin{equation}\label{eq:Phimt}
\Phi_{x_0}(p(t))\le -tm
\end{equation}
for all~${t\ge 0}$. In particular, the function $r(t)=d(\pi(\one),p(t))$
diverges to infinity as $t$ tends to infinity.  Hence
$$
\lim_{t\to\infty}\frac{1}{r(t)}\inf_{S_{r(t)}}\Phi_{x_0} = -m
$$
by part~(iii) of Lemma~\ref{le:KNSPHERE} and, for $t\ge 1$, we have
$$
\frac{m}{\sqrt{1-\frac{\abs{\xi_0-\xi_1}^2}{4}}}
\le 
\frac{tm}{r(t)} 
\le 
- \frac{\Phi_{x_0}(p(t))}{r(t)} 
\le
- \frac{1}{r(t)}\inf_{S_{r(t)}}\Phi_{x_0}.
$$
Here the first inequality follows from~\eqref{eq:rt}, 
the second from~\eqref{eq:Phimt},
and the last from the fact that $p(t)\in S_{r(t)}$.
Take the limit $t\to\infty$ to obtain 
$$
\frac{m}{\sqrt{1-\frac{\abs{\xi_0-\xi_1}^2}{4}}} \le m
$$
and hence $\xi_0=\xi_1$.  
This proves Theorem~\ref{thm:KEMPF1}.
\end{proof}

\begin{proof}[Proof of Theorem~\ref{thm:MODSTAB}]
\label{proof:MODSTAB}
Assume $x_0$ is $\mu$-stable and define~$\lambda(x_0)$ by~\eqref{eq:lambda}.
Then $\lambda(x_0)>0$ by Theorem~\ref{thm:SUPMAX} 
and Theorem~\ref{thm:HMnec}.  Now let 
$$
x\in\overline{\rG^c(x_0)}\setminus\rG^c(x_0).
$$
Then there exist sequences $\xi_i\in\cg$, $t_i>0$, and $u_i\in\rG$
such that $t_i\to\infty$, $\Abs{\xi_i}=1$ for all $i$, 
and~${x=\lim_{i\to\infty}u_i\exp(\i t_i\xi_i)x_0}$.
Hence 
\begin{equation*}
\begin{split}
\frac{1}{t_i}\inf_{S_{t_i}}\Phi_{x_0}
&\le 
\frac{1}{t_i}\Phi_{x_0}(\exp(-\i t_i\xi_i)) 
= 
\frac{1}{t_i}\int_0^{t_i}\inner{\mu(\exp(\i t\xi_i)x_0)}{\xi_i}\,dt \\
&\le 
\inner{\mu(\exp(\i t_i\xi_i)x_0)}{\xi_i} \\
&\le
\Abs{\mu(u_i\exp(\i t_i\xi_i)x_0)}
\end{split}
\end{equation*}
for all $i\in\N$.  Take the limit $i\to\infty$ 
and use part~(iii) of Lemma~\ref{le:KNSPHERE}
to obtain $\lambda(x_0)\le\Abs{\mu(x)}$.
This proves~\eqref{eq:MODSTAB2}. 

\bigbreak

It follows from~\eqref{eq:MODSTAB2} that every
element~${x\in\overline{\rG^c(x_0)}\setminus\rG^c(x_0)}$
satisfies
$$
\inf_{g\in\rG^c}\abs{\mu(gx)}\ge\lambda(x_0) > 0
$$
and hence is $\mu$-unstable.   Since the set 
$X^\us\subset X$ of $\mu$-unstable points is closed
by Theorem~\ref{thm:STABILITY}, so is the set
$$
\overline{\rG^c(x_0)}\setminus\rG^c(x_0) 
= \overline{\rG^c(x_0)}\cap X^\us,
$$
and hence this set is compact.  
Now let $\zeta\in\sT^c$.  Then the limit point
$$
x^+ := \lim_{t\to\infty}\exp(\i t\zeta)x_0
$$
satisfies $L_{x^+}^c\zeta=0$ by Lemma~\ref{le:WEIGHT2},
and hence cannot belong to $\rG^c(x_0)$. Thus 
$$
x^+\in\overline{\rG^c(x_0)}\setminus\rG^c(x_0)
$$
and so
$
\inf_{g\in\rG^c}\abs{\mu(gx^+)}\ge\lambda(x_0).
$
This proves Theorem~\ref{thm:MODSTAB}.
\end{proof}

\begin{cor}\label{cor:KEMPF}
Let $x_0\in X$ be $\mu$-unstable.
Then there exists a unique element $\xi_0\in\cg$ such that
$\abs{\xi_0}=1$ and 
\begin{equation}\label{eq:w2}
-w_\mu(x_0,\xi_0)
= \inf_{g\in\rG^c}\abs{\mu(gx_0)}
= \sup_{0\ne\xi\in\cg}
\frac{-w_\mu(x_0,\xi)}{\abs{\xi}}.
\end{equation}
If the triple $(X,\om,\mu)$ and the inner product on $\cg$ are rational
with factor $\hbar$ then there exists a positive integer $\ell$
such that $\sqrt{2\pi\hbar\ell}\xi_0\in\Lambda$.
\end{cor}

\begin{proof}
Let $\xi_\infty$ be as in Theorem~\ref{thm:KEMPF}.
Then the element 
$$
\xi_0:=\frac{\xi_\infty}{\abs{\xi_\infty}}\in\cg
$$ 
satisfies the first equation in~\eqref{eq:w2}.  
The second equation in~\eqref{eq:w2} 
follows from the first~($\le$) and the moment-weight 
inequality in Theorem~\ref{thm:MW2}~($\ge$).
Uniqueness follows from Theorem~\ref{thm:KEMPF1}.
To prove the last assertion, let~${k\in\N}$ 
such that $k\mu(x_\infty)\in\Lambda$ (Theorem~\ref{thm:mucrit})
and take
$$
\ell:=\frac{k^2\abs{\mu(x_\infty)}^2}{2\pi\hbar}.
$$
Then $\ell\in\N$ and
$$
\sqrt{2\pi\hbar\ell}\xi_0
= k\abs{\mu(x_\infty)}\xi_0
= k\xi_\infty
= -u_\infty^{-1}k\mu(x_\infty)u_\infty\in\Lambda.
$$
This\index{weight@$\mu$-weight!dominant|)}
proves Corollary~\ref{cor:KEMPF}.
\end{proof}


\chapter{Torus actions}\label{ch:TORUS}

Throughout this chapter we assume that the Lie group~${\rG=\rT}$ 
is a torus with Lie algebra~${\ct:=\Lie(\rT)}$ and complexification~$\rT^c$. 
The first main result is Theorem~\ref{thm:TORUS}, which
asserts that in the case of a torus action the function
\begin{equation}\label{eq:WEIGHTT}
\ct\setminus\{0\}\to\R:\xi\mapsto w_\mu(x,\xi)
\end{equation}
is continuous for every~${x\in X}$.  The second main result 
(Theorem~\ref{thm:CONVEX}) is not used elsewhere in this book.
It asserts that the closure of the image of a complexified group orbit 
under the moment map is convex in the case of a torus action.

\begin{thm}[{\bf Continuity}]\label{thm:TORUS}
Let~${x\in X}$. Then
\begin{equation}\label{eq:wtorus}
w_\mu(x,\xi)
= \sup_{g\in\rT^c}\inner{\mu(gx)}{\xi}
\end{equation}
for every $\xi\in\ct\setminus\{0\}$ and the 
function~\eqref{eq:WEIGHTT} is continuous.
\end{thm}

\begin{proof}
Let~${x_0\in X}$ and~${\xi\in\ct\setminus\{0\}}$. Then
$$
w_\mu(x_0,\xi)
=
\lim_{t\to\infty}\inner{\mu(\exp(\i t\xi)x_0)}{\xi}
\le
\sup_{g\in\rT^c}\inner{\mu(gx_0)}{\xi}.
$$
Now let $g\in\rT^c$. Then the function 
$t\mapsto\inner{\mu(\exp(\i t\xi)gx_0)}{\xi}$
is nondecreasing by~\eqref{eq:ZERO} and converges 
to $w_\mu(gx_0,\xi)$ as $t$ tends to infinity.
This implies
$$
\inner{\mu(gx_0)}{\xi}
\le 
w_\mu(gx_0,\xi)
= 
w_\mu(gx_0,g\xi g^{-1})
= 
w_\mu(x_0,\xi)
$$ 
by Theorem~\ref{thm:MUMFORD1}.
Hence
$$
\sup_{g\in\rT^c}\inner{\mu(gx_0)}{\xi}
\le w_\mu(x_0,\xi)
$$ 
and this proves~\eqref{eq:wtorus} for~${x=x_0}$.

\bigbreak

Now assume, by contradiction, that the function~\eqref{eq:WEIGHTT} 
with~${x=x_0}$ is not continuous. Then there exists a 
sequence~${\xi_i\in\ct\setminus\{0\}}$ converging to 
an element~${\xi\in\ct\setminus\{0\}}$ 
such that the sequence~${w_\mu(x_0,\xi_i)}$ 
does not converge to~${w_\mu(x_0,\xi)}$.  
Since~$w_\mu(x_0,\xi_i)$ is a bounded sequence,
there is a subsequence, still denoted by $\xi_i$, 
such that the limit 
\begin{equation}\label{eq:wcont1}
c_\infty := \lim_{i\to\infty}w_\mu(x_0,\xi_i)
\end{equation}
exists and is not equal to $w_\mu(x_0,\xi)$. 
Choose $x\in\overline{\rT^c(x_0)}$ such that
$$
w_\mu(x_0,\xi) 
= \sup_{g\in\rT^c}\inner{\mu(gx_0)}{\xi} 
= \inner{\mu(x)}{\xi}.
$$
Then, for every $i\in\N$, we have
\begin{equation*}
\begin{split}
w_\mu(x_0,\xi)  + \inner{\mu(x)}{\xi_i-\xi}
&=
\inner{\mu(x)}{\xi_i} \\
&\le
\sup_{g\in\rT^c}\inner{\mu(gx_0)}{\xi_i} \\
&=
w_\mu(x_0,\xi_i).
\end{split}
\end{equation*}
Take the limit $i\to\infty$ to obtain
$w_\mu(x_0,\xi)\le\lim_{i\to\infty}w_\mu(x_0,\xi_i)=c_\infty$
and so
\begin{equation}\label{eq:wcont2}
w_\mu(x_0,\xi) < c_\infty.
\end{equation}
Choose a sequence $x_i\in\overline{\rT^c(x_0)}$ 
such that 
\begin{equation}\label{eq:wcont3}
w_\mu(x_0,\xi_i)= \inner{\mu(x_i)}{\xi_i}
\qquad\mbox{for all }i\in\N.
\end{equation}
Passing to a subsequence we may assume that the 
limit~${x_\infty := \lim_{i\to\infty}x_i}$ exists.  
Then~${x_\infty\in\overline{\rT^c(x_0)}}$ and,
by~\eqref{eq:wcont1}, \eqref{eq:wcont2},
and~\eqref{eq:wcont3}, we have
\begin{equation*}
\begin{split}
w_\mu(x_0,\xi) 
&<
c_\infty \\
&=
\lim_{i\to\infty}w_\mu(x_0,\xi_i) \\
&=
\lim_{i\to\infty}\inner{\mu(x_i)}{\xi_i} \\
&=
\inner{\mu(x_\infty)}{\xi}.
\end{split}
\end{equation*}
This contradicts equation~\eqref{eq:wtorus} 
and thus completes the proof 
of Theorem~\ref{thm:TORUS}.
\end{proof}

The next lemma makes use of the Generalized 
Kempf Existence Theorem~\ref{thm:KEMPF}.
Denote the unit sphere in~$\ct$ 
by~${S(\ct) := \left\{\xi\in\ct\,|\,\abs{\xi}=1\right\}}$.  

\begin{lem}\label{le:TORUS}
Let $x_0\in X$ such that
$$
\inf_{g\in\rT^c}\abs{\mu(gx_0)}>0,
$$
and let~${\eta\in\overline{\mu(\rT^c(x_0))}}$, 
Then the following are equivalent:
\begin{equation}\label{eq:convex1}
\abs{\eta} = \inf_{g\in\rT^c}\abs{\mu(gx_0)},
\end{equation}
\begin{equation}\label{eq:convex2}
\abs{\eta}^2 = \inf_{g\in\rT^c}\inner{\eta}{\mu(gx_0)}.
\end{equation}
Moreover, there is a unique~${\eta\in\overline{\mu(\rT^c(x_0))}}$ 
satisfying these conditions, 
the function~${S(\ct)\to\R:\xi\mapsto w_\mu(x_0,\xi)}$ 
takes on its minimum at~${\xi_0:=-\abs{\eta}^{-1}\eta}$ 
and only at that point, and
$$
\inf_{\xi\in S(\ct)}w_\mu(x_0,\xi)
= - \abs{\eta}
= - \inf_{\rT^c(x_0)}\abs{\mu}.
$$
\end{lem}

\begin{proof}
Existence is obvious for equation~\eqref{eq:convex1} and
uniqueness is obvious for equation~\eqref{eq:convex2}.  
We prove that there exists an element~${\eta\in\overline{\mu(\rT^c(x_0))}}$
that satisfies~\eqref{eq:convex2}. Let $x:\R\to X$ 
be the solution of~\eqref{eq:KN1} and define
$$
x_\infty := \lim_{t\to\infty}x(t),\qquad
\eta := \mu(x_\infty).
$$
Then~${\eta\in\overline{\mu(\rT^c(x_0))}}$ and, 
by Theorem~\ref{thm:MLT}, we have
$$
\abs{\eta} = \abs{\mu(x_\infty)} = \inf_{g\in\rT^c}\abs{\mu(gx_0)}.
$$
Since~$\rT^c$ is abelian, it follows 
from~\eqref{eq:xinfty2} and~\eqref{eq:xinfty3} 
in Theorem~\ref{thm:KEMPF} that
$$
w_\mu(x_0,-\eta) = - \abs{\eta}^2.
$$
Hence, by part~(i) of Theorem~\ref{thm:MUMFORD1}, we have
$$
w_\mu(gx_0,-\eta) 
= w_\mu(x_0,-\eta) 
= -\abs{\eta}^2
\qquad\mbox{for all }g\in\rT^c.
$$
Since the function $t\mapsto\inner{\mu(\exp(-\i t\eta)gx_0)}{-\eta}$ 
is nondecreasing by~\eqref{eq:ZERO}, and converges 
to $w_\mu(gx_0,-\eta)=-\abs{\eta}^2$ as $t$ tends to infinity, 
it follows (by evaluating at~${t=0}$) that 
$\inner{\mu(gx_0)}{-\eta} \le -\abs{\eta}^2$ and 
hence~${\abs{\eta}^2\le \inner{\mu(gx_0)}{\eta}}$
for all~${g\in\rT^c}$.  This proves the existence of 
an element $\eta\in\overline{\mu(\rT^c(x_0))}$ 
that satisfies~\eqref{eq:convex2}.

If $\eta\in\overline{\mu(\rT^c(x_0))}$ satisfies~\eqref{eq:convex2}
then $\abs{\eta}^2\le\abs{\eta}\abs{\mu(gx_0)}$ 
for all $g\in\rT^c$ by the Cauchy--Schwarz inequality.  
Since $\eta\ne0$ this implies $\abs{\eta}\le\abs{\mu(gx_0)}$
for all $g\in\rT^c$. Hence $\eta$ satisfies~\eqref{eq:convex1}.

Conversely, assume that $\eta\in\overline{\mu(\rT^c(x_0))}$ 
satisfies~\eqref{eq:convex1}.  By what we have proved above 
there is an $\eta_0\in\overline{\mu(\rT^c(x_0))}$ 
that satisfies~\eqref{eq:convex2} and hence 
also~\eqref{eq:convex1}. This implies 
${\abs{\eta}^2=\abs{\eta_0}^2\le\inner{\eta_0}{\eta}}$
and hence~${\eta=\eta_0}$.  Thus~$\eta$ 
satisfies~\eqref{eq:convex2} and is uniquely 
determined by either condition.

Now define
$$
m := \inf_{g\in\rT^c}\Abs{\mu(gx_0)}>0
$$ 
and observe that~${\inner{\mu(gx_0)}{\xi}\ge -\abs{\mu(gx_0)}}$
for all~${g\in\rT^c}$ and all~${\xi\in S(\ct)}$.
Take the supremum over all~${g\in\rT^c}$ to obtain
$$
w_\mu(x_0,\xi)
=
\sup_{g\in\rT^c}\inner{\mu(gx_0)}{\xi}
\ge 
-\inf_{g\in\rT^c}\abs{\mu(gx_0)}
= 
-m
$$
for all~${\xi\in S(\ct})$ by~\eqref{eq:wtorus}.
Now take the infimum over all~${\xi\in S(\ct})$  
to obtain
\begin{equation}\label{eq:xim}
\inf_{\xi\in S(\ct)}w_\mu(x_0,\xi) \ge -m.
\end{equation}
To prove equality, recall from~\eqref{eq:convex1} 
and~\eqref{eq:convex2} that  
\begin{equation}\label{eq:etam}
\inf_{g\in\rT^c} \inner{\mu(gx_0)}{\eta} 
= \abs{\eta}^2 
= m^2.
\end{equation}
Thus $\eta\ne0$ and, by equation~\eqref{eq:wtorus}, we have
$$
w_\mu\left(x_0,-\abs{\eta}^{-1}\eta\right) 
= -\inf_{g\in\rT^c}\inner{\mu(gx_0)}{\abs{\eta}^{-1}\eta}
= -\abs{\eta}  
= -m.
$$
By~\eqref{eq:xim} this shows that the 
function~${S(\ct)\to\R:\xi\mapsto w_\mu(x_0,\xi)}$
attains its minimum at the 
point~${\xi_0 := -\abs{\eta}^{-1}\eta\in S(\ct)}$, i.e.\ 
$$
\inf_{\xi\in S(\ct)}w_\mu(x_0,\xi)
= w_\mu\left(x_0,-\abs{\eta}^{-1}\eta\right) 
= -m.
$$
Now let~$\xi\in S(\ct)$ such that $\xi\ne-\abs{\eta}^{-1}\eta$.  Then 
$\inner{\abs{\eta}^{-1}\eta}{\xi}>-1$ and hence
$$
w_\mu(x_0,\xi) 
=
\sup_{g\in\rT^c}\inner{\mu(gx_0)}{\xi}
\ge 
\inner{\eta}{\xi} 
> 
- \abs{\eta} 
= 
- m.
$$
Here the first step follows from~\eqref{eq:wtorus},
the second step follows from the fact 
that~${\eta\in\overline{\mu(\rT^c(x_0))}}$, 
the third step uses the 
inequality~${\inner{\abs{\eta}^{-1}\eta}{\xi}>-1}$,
and the last step follows from~\eqref{eq:etam}.
This proves Lemma~\ref{le:TORUS}.
\end{proof}

The following theorem asserts that the image of every complexified
group orbit under the moment map has a convex closure.
This is a variant of Atiyah--Guillemin--Sternberg convexity
(see~\cite[\S5\,\&\,Appendix]{NESS2}).

\begin{thm}[{\bf Convexity}]\label{thm:CONVEX}
For every $x\in X$ the set $\overline{\mu(\rT^c(x))}$ is convex.
\end{thm}

\begin{proof}
For every~${\tau\in\ct}$ and every~${r>0}$ 
denote the closed ball of radius~$r$ about~$\tau$ 
by~${B_r(\tau)\subset\ct}$.  
Fix an element~${x\in X}$ and define
$$
\Delta:=\overline{\mu(\rT^c(x))}.
$$
For $\tau\in\ct\setminus\Delta$ define 
$$
d(\tau,\Delta):=\inf_{\xi\in\Delta}\abs{\tau-\xi}.
$$
Then Lemma~\ref{le:TORUS}, with $\mu$ 
replaced by $\mu-\tau$, asserts that
\begin{equation}\label{eq:CONVEX}
\left.\begin{array}{l}
\tau\in\ct\setminus\Delta\mbox{ and}\\
\eta\in B_{d(\tau,\Delta)}(\tau)\cap\Delta
\end{array}\right\}
\quad\implies\quad
d(\tau,\Delta)^2 = \inf_{\eta'\in\Delta}\inner{\eta-\tau}{\eta'-\tau}.
\end{equation}
This implies that $\Delta$ is convex.
To see this, suppose by contradiction that 
there exist elements $\tau_0,\tau_1\in\Delta$ 
and a constant $0<\lambda<1$ such 
that
$$
\tau:=(1-\lambda)\tau_0+\lambda\tau_1\notin\Delta.
$$
Then
$$
r:=d(\tau,\Delta)>0. 
$$
Fix an element~${\eta\in B_r(\tau)\cap\Delta}$. 
Then by~\eqref{eq:CONVEX} we have
$$
\inner{\eta-\tau}{\tau_i-\tau}\ge r^2
\qquad\mbox{for }i=0,1.
$$
However, since~${\tau_0-\tau=\lambda(\tau_0-\tau_1)}$ 
and~${\tau_1-\tau=(1-\lambda)(\tau_1-\tau_0)}$, the inner 
products~${\inner{\eta-\tau}{\tau_0-\tau}}$ 
and~${\inner{\eta-\tau}{\tau_1-\tau}}$ have opposite signs.  
This is a contradiction and proves Theorem~\ref{thm:CONVEX}.
\end{proof}


\chapter{The Hilbert--Mumford criterion}\label{ch:HM}

In this chapter we return to the general case where~$\rG$ 
is any compact Lie group acting on a closed K\"ahler 
manifold~$X$ by K\"ahler isometries, and the action is 
generated by an equivariant moment map $\mu:X\to\cg=\Lie(\rG)$. 

\begin{thm}[{\bf The Mumford Numerical Function}]\label{thm:M}
\ 

\smallskip\noindent{\bf (i)}
For\index{Mumford numerical function}
every $x\in X$ we have
\begin{equation}\label{eq:MUMFORD}
\begin{split}
m_\mu(x) 
&:= 
\inf_{0\ne\xi\in\cg} \frac{w_\mu(x,\xi)}{\abs{\xi}}  
=
\inf_{\zeta\in\sT^c}
\frac{w_\mu(x,\zeta)}{\sqrt{\abs{\Re(\zeta)}^2-\abs{\Im(\zeta)}^2}} \\
&\phantom{:}= 
\inf_{\xi\in\Lambda} \frac{w_\mu(x,\xi)}{\abs{\xi}} 
=
\inf_{\zeta\in\Lambda^c}
\frac{w_\mu(x,\zeta)}{\sqrt{\abs{\Re(\zeta)}^2-\abs{\Im(\zeta)}^2}}.
\end{split}
\end{equation}
The function $m_\mu:X\to\R$ defined by equation~\eqref{eq:MUMFORD} 
for $x\in X$ is called the {\bf Mumford numerical function}.

\smallskip\noindent{\bf (ii)}
The Mumford numerical function $m_\mu:X\to\R$ 
is $\rG^c$-invariant.

\smallskip\noindent{\bf (iii)}
Every $x\in X$ satisfies $m_\mu(x)+\inf_{\rG^c(x)}\abs{\mu}\ge0$.

\smallskip\noindent{\bf (iv)}
If $x\in X$ is $\mu$-unstable then
$0>m_\mu(x)=-\inf_{\rG^c(x)}\abs{\mu}$.

\smallskip\noindent{\bf (v)}
If $x\in X$ is $\mu$-stable then
$0<m_\mu(x)\le\inf_{\overline{\rG^c(x)}\setminus\rG^c(x)}\abs{\mu}$.

\smallskip\noindent{\bf (vi)}
For each $x\in X$ there is a $\xi\in\cg$ such that
$\Abs{\xi}=1$ and ${w_\mu(x,\xi)=m_\mu(x)}$.
\end{thm}

\begin{proof}
By Theorem~\ref{thm:MUMFORD} every $\zeta\in\sT^c$ 
is equivalent to an element $\xi\in\cg$, 
we have~${w_\mu(x,\zeta)=w_\mu(x,\xi)}$ by Theorem~\ref{thm:MUMFORD1},
and~${\abs{\Re(\zeta)}^2-\abs{\Im(\zeta)}^2=\abs{\xi}^2}$
by Lemma~\ref{le:WEIGHT3}.  This proves the second 
and last equalities in~\eqref{eq:MUMFORD}.

It remains to prove that the infimum over $\cg\setminus\{0\}$ 
in~\eqref{eq:MUMFORD} agrees with the infimum over~$\Lambda$.
Let $x_0\in X$ and~${\xi_0\in\cg}$ such that
$$
\Abs{\xi_0}=1,\qquad 
w_\mu(x_0,\xi_0)=m_\mu(x_0)
$$
(Theorem~\ref{thm:MODSTAB}). Consider the torus
$$
\rT := \overline{\left\{\exp(t\xi_0)\,|\,t\in\R\right\}} 
\subset \rG,\qquad
\ct := \Lie(\rT).
$$
Denote the unit sphere in~$\ct$ by
$$
S(\ct):=\{\xi\in\ct\,|\,\abs{\xi}=1\}.
$$  
Then the function ${S(\ct)\to\R:\xi\mapsto w_\mu(x_0,\xi)}$
is continuous by Theorem~\ref{thm:TORUS}, 
the set~${\{\abs{\eta}^{-1}\eta\,|\,\eta\in\ct\cap\Lambda\}}$ 
is dense in $S(\ct)$, and~${\inf_{\xi\in S(\ct)}w_\mu(x_0,\xi)=m_\mu(x_0)}$.
Hence, for each $\eps>0$, there exists an 
element~${\eta\in\ct\cap\Lambda}$ such that 
$$
\abs{\eta}^{-1}w_\mu(x_0,\eta)
= w_\mu(x_0,\abs{\eta}^{-1}\eta)
< m_\mu(x_0) + \eps.
$$
This implies
$$
\inf_{\eta\in\Lambda}(\abs{\eta}^{-1}w_\mu(x_0,\eta))
\le m_\mu(x_0) = \inf_{0\ne\xi\in\cg} (\abs{\xi}^{-1}w_\mu(x_0,\xi)).
$$
The converse inequality is obvious and this proves~(i).

Part~(ii) follows from~(i), Theorem~\ref{thm:MUMFORD1}, and Lemma~\ref{le:WEIGHT3},
part~(iii) is equivalent to the moment-weight inequality in Theorem~\ref{thm:MW2},
part~(iv) follows from Corollary~\ref{cor:KEMPF}, 
and parts~(v) and~(vi) follow from Theorem~\ref{thm:MODSTAB}. 
This proves Theorem~\ref{thm:M}.
\end{proof}

The main results of the present chapter are the Hilbert--Mumford 
numerical criteria for $\mu$-semistability, 
$\mu$-polystability, and $\mu$-stability. 
We begin with the $\mu$-semistable case,
where the Hilbert--Mumford criterion is a direct 
consequence of Theorem~\ref{thm:HMnec} (necessity)
and Theorem~\ref{thm:M} (sufficiency).  

\begin{thm}[{\bf Hilbert--Mumford Criterion: Semistable Case}]\label{thm:HMss}
\ 

\noindent
For every~${x_0\in X}$ the following are 
equivalent.\index{Hilbert--Mumford Criterion!for semistability}

\smallskip\noindent{\bf (i)}
$x_0$ is $\mu$-semistable.

\smallskip\noindent{\bf (ii)}
Every $\xi\in\Lambda$ satisfies $w_\mu(x_0,\xi)\ge 0$.

\smallskip\noindent{\bf (iii)}
Every $\xi\in\cg\setminus\{0\}$ satisfies $w_\mu(x_0,\xi)\ge 0$.

\smallskip\noindent{\bf (iv)}
Every $\zeta\in\Lambda^c$ satisfies $w_\mu(x_0,\zeta)\ge 0$.

\smallskip\noindent{\bf (v)}
Every $\zeta\in\sT^c$ satisfies $w_\mu(x_0,\zeta)\ge 0$.
\end{thm}

\begin{proof}
The equivalence of the assertions~(ii), (iii), (iv) and~(v) follows 
from equation~\eqref{eq:MUMFORD} in Theorem~\ref{thm:M},
that each of these conditions implies~(i) follows from part~(iv)
of Theorem~\ref{thm:M}, and that~(i) implies~(v) was proved in 
part~(i) of Theorem~\ref{thm:HMnec}.  

More precisely, Mumford's Theorems~\ref{thm:MUMFORD} 
and~\ref{thm:MUMFORD1} show that~(ii)$\iff$(iv) and (iii)$\iff$(v).  
The equivalence~(ii)$\iff$(iii) follows from the continuity of the weights 
for torus actions in Theorem~\ref{thm:TORUS}, with the 
argument spelled out in the proof of Theorem~\ref{thm:M}.
That~(i)$\implies$(iii) follows from the moment-weight inequality 
in Theorem~\ref{thm:MW2}, which shows that the existence 
of a negative weight implies that $x_0$ is $\mu$-unstable.
That~(iii)$\implies$(i) follows from the Kempf Existence 
Theorem~\ref{thm:KEMPF}, which produces 
an element~${\xi\in\cg\setminus\{0\}}$ with~${w_\mu(x_0,\xi)<0}$
whenever $x_0$ is $\mu$-unstable. This proves Theorem~\ref{thm:HMss}.
\end{proof}

\begin{thm}[{\bf Hilbert--Mumford Criterion: Unstable Case}]\label{thm:HMus}
\ 

\noindent
For every~${x_0\in X}$ the following are equivalent.  

\smallskip\noindent{\bf (i)}
$x_0$ is $\mu$-unstable.

\smallskip\noindent{\bf (ii)}
There exists a $\xi\in\Lambda$ such that $w_\mu(x_0,\xi)<0$.

\smallskip\noindent{\bf (iii)}
There exists a $\xi\in\cg\setminus\{0\}$ such that $w_\mu(x_0,\xi)<0$.

\smallskip\noindent{\bf (iv)}
There exists a $\zeta\in\Lambda^c$ such that $w_\mu(x_0,\zeta)<0$.

\smallskip\noindent{\bf (v)}
There exists a $\zeta\in\sT^c$ such that $w_\mu(x_0,\zeta)<0$.
\end{thm}

\begin{proof}
This follows directly from the definitions and Theorem~\ref{thm:HMss}.
\end{proof}

The next theorem is the {\bf Hilbert--Mumford numerical criterion}
in its classical form. We derive it as a corollary of 
Theorem~\ref{thm:HMus}.\index{Hilbert--Mumford Criterion!classical}

\begin{thm}[{\bf Hilbert--Mumford Criterion: Classical Case}]
\label{thm:HMc}
\ 

\noindent
Let~${\rG\subset\rU(n)}$ be a compact Lie group, 
let~${\rG^c\subset\GL(n,\C)}$ be its complexification,
and let $V$ be a finite-dimensional complex vector space 
equipped with a holomorphic representation of $\rG^c$.
If $v\in V$ is a nonzero vector such that~${0\in\overline{\rG^c(v)}}$, 
then there exists a~${\xi\in\Lambda}$ such 
that~${\lim_{t\to\infty}\exp(\i t\xi)v = 0}$.
\end{thm}

\begin{proof}
By assumption $\rG$ induces a Hamiltonian group action on
the projective space $X=\bbP(V)$ with the moment map 
of Lemma~\ref{le:Vmu}.  Let $v\in V$ be a nonzero vector such that
$0\in\overline{\rG^c(v)}$. Then $v$ is unstable and thus the 
point~${{x:=[v]}\in\bbP(V)}$ is $\mu$-unstable by Theorem~\ref{thm:VECTOR}.
Hence Theorem~\ref{thm:HMus} asserts that there exists 
a~${\xi\in\Lambda}$ such that $w_\mu(x,\xi)<0$. 
By Lemma~\ref{le:Vweight} this means that~$v$ 
is contained in the direct sum of the negative 
eigenspaces of the Hermitian operator on $V$
determined by~$\i\xi$. Hence~${\lim_{t\to\infty}\exp(\i t\xi)v=0}$ 
and this proves Theorem~\ref{thm:HMc}.
\end{proof}

\begin{thm}[{\bf Hilbert--Mumford Criterion: Polystable Case}]\label{thm:HMps}
\ 

\noindent
For every $x_0\in X$ the following are 
equivalent.\index{Hilbert--Mumford Criterion!for polystability}

\smallskip\noindent{\bf (i)}
$x_0$ is $\mu$-polystable.

\smallskip\noindent{\bf (ii)}
Every $\xi\in\Lambda$ satisfies $w_\mu(x_0,\xi)\ge 0$ and
$$
w_\mu(x_0,\xi)=0\qquad\implies\qquad
\lim_{t\to\infty}\exp(\i t\xi)x_0\in\rG^c(x_0).
$$

\smallskip\noindent{\bf (iii)}
Every $\xi\in\cg\setminus\{0\}$ satisfies $w_\mu(x_0,\xi)\ge 0$ and
$$
w_\mu(x_0,\xi)=0\qquad\implies\qquad
\lim_{t\to\infty}\exp(\i t\xi)x_0\in\rG^c(x_0).
$$

\smallskip\noindent{\bf (iv)}
Every $\zeta\in\Lambda^c$ satisfies $w_\mu(x_0,\zeta)\ge 0$ and
$$
w_\mu(x_0,\zeta)=0\qquad\implies\qquad
\lim_{t\to\infty}\exp(\i t\zeta)x_0\in\rG^c(x_0).
$$

\smallskip\noindent{\bf (v)}
Every $\zeta\in\sT^c$ satisfies $w_\mu(x_0,\zeta)\ge 0$ and
$$
w_\mu(x_0,\zeta)=0\qquad\implies\qquad
\lim_{t\to\infty}\exp(\i t\zeta)x_0\in\rG^c(x_0).
$$
\end{thm}


\begin{proof}
That~(i)$\implies$(v) was proved in part~(ii) of Theorem~\ref{thm:HMnec}, 
and the implications~(v)$\implies$(iv)$\implies$(ii) 
and~(v)$\implies$(iii)$\implies$(ii) follow directly from the definitions.  
Thus it remains to prove that~(ii) implies~(i). 

The proof of sufficiency of the $\mu$-weight condition for $\mu$-polystability 
is due to Chen--Sun~\cite[Theorem~4.7]{CS}.  Here is their argument.  
Fix an element~${x_0\in X}$ that is $\mu$-semistable 
but not $\mu$-polystable, let 
$
x:\R\to X
$
be the solution of the differential equation
$$
\dot x=-JL_x\mu(x),\qquad x(0)=x_0,
$$
and define 
$
x_\infty:=\lim_{t\to\infty}x(t).
$
Then, by Theorem~\ref{thm:STABILITY},  
\begin{equation}\label{eq:NP1}
x_\infty:=\lim_{t\to\infty}x(t)\notin\rG^c(x_0),\qquad \mu(x_\infty)=0,\qquad
\ker\,L_{x_\infty}\ne 0.
\end{equation}
We prove in seven steps that there exists 
an element $\xi\in\Lambda$ such that 
\begin{equation}\label{eq:NP2}
w_\mu(x_0,\xi)=0,\qquad
 \lim_{t\to\infty}\exp(\i t\xi)x_0\notin\rG^c(x_0).
\end{equation}

\noindent{\bf Step~1.}
{\it  The complex isotropy subgroup $\rG^c_{x_\infty}$
is the complexification of $\rG_{x_\infty}$.
Moreover, there exists a $\rG_{x_\infty}$-equivariant 
local holomorphic coordinate chart $\psi:U_\infty\to X$
on a $G_{x_\infty}$-invariant open 
neighborhood~${U_\infty\subset T_{x_\infty}X}$
of the origin such that $\psi(0)=x_\infty$ and $d\psi(0)=\id$.}

\medskip\noindent
Since $\mu(x_\infty)=0$, it follows from Lemma~\ref{le:isotropy}
that $\rG^c_{x_\infty}$ is the complexification of $\rG_{x_\infty}$.
Now let 
$
\phi:(T_{x_\infty}X,0)\to(X,x_\infty)
$ 
be any holomorphic coordinate chart,
defined in a neighborhood of the origin in~$T_{x_\infty}X$,
such that~${\phi(0)=x_\infty}$ and~${d\phi(0)=\id}$.  
Let~$\dvol_\infty$ denote the Haar measure on $\rG_{x_\infty}$ 
and define a map~$\chi$ from an open neighborhood of~$x_\infty$
in~$X$ to an open neighborhood of the origin in~${T_{x_\infty}X}$ by
$$
\chi(x) := \frac{1}{\Vol(\rG_{x_\infty})}
\int_{\rG_{x_\infty}} u^{-1}\phi^{-1}(ux)\dvol_\infty(u)
$$
for $x\in X$ sufficiently close to $x_\infty$.
Then $\chi(x_\infty)=0$, $d\chi(x_\infty)=\id$, 
and $\chi$ is holomorphic and $\rG_{x_\infty}$-equivariant.  
Hence, by the inverse function theorem, 
it restricts to a $\rG_{x_\infty}$-equivariant holomorphic diffeomorphism 
from a $\rG_{x_\infty}$-invariant open neighborhood 
of~$x_\infty$ in~$X$ to a $\rG_{x_\infty}$-invariant open 
neighborhood $U_\infty\subset T_{x_\infty}X$ of the origin.   
The inverse $\psi:=\chi^{-1}$ of this restriction 
satisfies the requirements of Step~1.

\medskip\noindent{\bf Step~2.}
{\it   Recall the notation~${v_\xi(x)=L_x\xi}$
for the infinitesimal action of~${\xi\in\cg}$ on~$X$,
and let~${\psi:U_\infty\to X}$ be the holomorphic coordinate chart
in Step~1.  Then there exists a $\delta>0$ such
that ${B_\delta(x_\infty)
:=\left\{\xhat\in T_{x_\infty}X\,\big|\,\abs{\xhat}<\delta\right\}
\subset U_\infty}$ and,
for all~${\xhat\in B_\delta(x_\infty)\cap\im(L_{x_\infty}^c)^\perp}$
and all ${\yhat\in T_{x_\infty}X}$, ${\zeta=\xi+\bi\eta\in\cg^c}$,}
\begin{equation}\label{eq:AL}
\begin{array}{l}
d\psi(\xhat)\yhat = L_{\psi(\xhat)}^c\zeta\\
\mbox{\it{and} }\yhat\perp\im(L_{x_\infty}^c)
\end{array}
\qquad\iff\qquad
\begin{array}{l}
\yhat = \Nabla{\xhat}v_\xi(x_\infty) + J\Nabla{\xhat}v_\eta(x_\infty)\\
\mbox{\it{and} } L_{x_\infty}^c\zeta=0.
\end{array}
\end{equation}

\medskip\noindent
The action of~$\rG_{x_\infty}$ on~$X$ by K\"ahler isometries
gives rise to a unitary action of~$\rG_{x_\infty}$ 
on the tangent space~$T_{x_\infty}X$ at the fixed point~$x_\infty$.
The infinitesimal action takes the form of a Lie algebra 
homomorphism
$$
\cg_{x_\infty}:=\Lie(\rG_{x_\infty})=\ker(L_{x_\infty})\to\cu(T_{x_\infty}X):
\xi\mapsto A_\xi = \nabla v_\xi(x_\infty).
$$
Thus
$
A_\xi\xhat := \left.\tfrac{d}{dt}\right|_{t=0}\exp(t\xi)\xhat 
= \Nabla{\xhat}v_\xi(x_\infty)
$
for~${\xhat\in T_{x_\infty}X}$ and~${\xi\in\ker(L_{x_\infty})}$.
Differentiate the identity $\psi(u\xhat)=u\psi(\xhat)$ 
for~${\xhat\in U_\infty}$ and ${u\in\rG_{x_\infty}}$,
to obtain for all~${\xhat\in U_\infty}$ and all~${\xi\in\cg}$,
\begin{equation}\label{eq:AL1}
L_{x_\infty}\xi=0\qquad\implies\qquad
d\psi(\xhat)A_\xi\xhat = L_{\psi(\xhat)}\xi.
\end{equation}
Since~${uL_{x_\infty}\eta=L_{x_\infty}(u\eta u^{-1})}$
for~${u\in\rG_{x_\infty}}$ and~${\eta\in\cg}$,
we also have for~${\xi,\eta\in\cg}$,
\begin{equation}\label{eq:AL2}
L_{x_\infty}\xi=0\qquad\implies\qquad
A_\xi L_{x_\infty}\eta = L_{x_\infty}[\xi,\eta].
\end{equation}
By~\eqref{eq:AL2} the subspaces~$\im(L_{x_\infty}^c)$ 
and~$\im(L_{x_\infty}^c)^\perp$ are invariant under the 
endomorphism~$A_\xi$ for every~${\xi\in\ker(L_{x_\infty})}$.  
Hence the implication ``$\Longleftarrow$'' in~\eqref{eq:AL} 
for all~${\xhat\in U_\infty\cap\im(L_{x_\infty}^c)^\perp}$
follows from~\eqref{eq:AL1}, because~$d\psi(\xhat)$ 
is complex linear and~$L_{x_\infty}^c\zeta=0$ 
if and only if~${L_{x_\infty}\xi=L_{x_\infty}\eta=0}$.

\bigbreak

For~${\xhat\in U_\infty}$ define the linear operator
$
\sL_{\xhat}:\cg^c\times(\im L_{x_\infty}^c)^\perp\to T_{\psi(\xhat)}X
$
by 
$$
\sL_{\xhat}(\zeta,\yhat) := L^c_{\psi(\xhat)}\zeta-d\psi(\xhat)\yhat.
$$
for~$\zeta\in\cg^c$ and~${\yhat\in\im(L_{x_\infty}^c)^\perp}$.   
The index of this operator (the dimension of the source minus 
the dimension of the target) is~${2k}$, where~${k:=\dim(\rG_{x_\infty})}$.
Moreover, if~${\xhat\in U_\infty\cap \im(L_{x_\infty}^c)^\perp}$,
then it follows from the implication ``$\Longleftarrow$'' in~\eqref{eq:AL}
(already proved)  that the $2k$-dimensional subspace
$$
\sZ_\xhat := \left\{(\xi+\bi\eta,A_\xi\xhat+JA_\eta\xhat)\,\big|\,\xi,\eta\in\ker(L_{x_\infty})\right\}
\subset\cg^c\times\im(L_{x_\infty}^c)^\perp
$$
is contained in the kernel of the operator~$\sL_{\xhat}$.
Since~$\sL_\xhat$ is surjective for~${\xhat=0}$,
there exists a constant $\delta>0$ such that $B_\delta(x_\infty)\subset U_\infty$
and $\sL_{\xhat}$ is surjective for every~${\xhat\in B_\delta(x_\infty)}$.
Hence~${\dim(\ker(\sL_{\xhat}))=2k}$ 
and so~${\ker(\sL_{\xhat}) = \sZ_\xhat}$
for every~${\xhat\in B_\delta(x_\infty)\cap\im(L_{x_\infty}^c)^\perp}$.
This proves Step~2. 

\medskip\noindent{\bf Step~3.}
{\it Let $\delta>0$ be as inStep~2.
Then there exists a $t_0>0$ and smooth curves 
$\xi,\eta:[t_0,\infty)\to(\ker L_{x_\infty})^\perp$,
$\xhat:[t_0,\infty)\to (\im L_{x_\infty}^c)^\perp$ 
such that
$$
x(t) = \exp(\i\eta(t))\exp(\xi(t))\psi(\xhat(t)),\qquad
\abs{\xhat(t)}<\delta,
$$
for every $t\ge t_0$.}

\medskip\noindent
Define the map 
$f:\ker(L_{x_\infty})^\perp\times\ker(L_{x_\infty})^\perp
\times(B_\delta(x_\infty)\cap\im(L_{x_\infty}^c)^\perp)\to X$ 
by~${f(\xi,\eta,\xhat) := \exp(\i\eta)\exp(\xi)\psi(\xhat)}$.
Then the derivative of $f$ at the origin is bijective.  
Hence $f$ restricts to a diffeomorphism
from an open neighborhood of the origin 
in~${\ker(L_{x_\infty})^\perp\times\ker(L_{x_\infty})^\perp
\times\im(L_{x_\infty}^c)^\perp}$
onto an open neighborhood of~$x_\infty$ in~$X$. 
This proves Step~3.

\medskip\noindent{\bf Step~4.}
{\it Let $t_0,\xi,\eta,\xhat$ be as in Step~3 and let
$g:\R\to\rG^c$ be the unique solution of the 
equation $g^{-1}\dot g=\i\mu(x)$ with $g(0)=\one$.
For $t\ge t_0$ define
$$
h(t) := \exp(\i\eta(t))\exp(\xi(t)),\qquad
g_\infty(t) := h(t_0)^{-1}g(t_0)^{-1}g(t)h(t).
$$
Then, for all $t\ge t_0$,}
$$
g_\infty(t_0)=\one,\quad
g_\infty(t)\in\rG_{x_\infty}^c,\quad
\xhat(t)=g_\infty(t)^{-1}\xhat(t_0),\quad
x(t)=h(t)\psi(\xhat(t)).
$$

\medskip\noindent
By Lemma~\ref{le:GRADFLOW} and Step~3, we have
$$
g(t)^{-1}x_0 = x(t) = h(t)\psi(\xhat(t))
$$
for $t\ge t_0$.  Hence,
for all~${t\ge t_0}$,
\begin{equation*}
\begin{split}
\psi(\xhat(t)) 
&= 
h(t)^{-1}g(t)^{-1}x_0 \\
&= 
h(t)^{-1}g(t)^{-1}g(t_0)h(t_0)\psi(\xhat(t_0)) \\
&=
g_\infty(t)^{-1}\psi(\xhat(t_0)).
\end{split}
\end{equation*}
Differentiate this equation to obtain 
$$
d\psi(\xhat)\p_t\xhat = L^c_{\psi(\xhat)}\zeta_\infty,\qquad
\zeta_\infty := \xi_\infty+\i\eta_\infty := -g_\infty^{-1}\dot g_\infty.
$$
Since $\abs{\xhat(t)}<\delta$ 
and $\p_t\xhat(t)\perp\im(L_{x_\infty}^c)$, 
it follows from Step~2 that 
$$
\zeta_\infty(t)\in\ker(L_{x_\infty}^c),\qquad
\p_t\xhat(t) = \Nabla{\xhat(t)}v_{\xi_\infty(t)}(x_\infty) 
+ J\Nabla{\xhat(t)}v_{\eta_\infty(t)}(x_\infty)
$$
for all $t\ge t_0$.  Since $g_\infty(t_0)=\one$, we obtain
$$
g_\infty(t)\in\rG^c_{x_\infty},\qquad
\xhat(t)=g_\infty(t)^{-1}\xhat(t_0)
$$
for all $t\ge t_0$.  This proves Step~4.

\medskip\noindent{\bf Step~5.}
{\it There exists an element $\xi\in\Lambda$ such that 
$L_{x_\infty}\xi=0$ and}
$$
\lim_{t\to\infty}\exp(\i t\xi)h(t_0)^{-1}x(t_0)=x_\infty,\qquad
w_\mu(h(t_0)^{-1}x(t_0),\xi)=0.
$$

\medskip\noindent
By Step~4, $\xhat(t)\in\rG_{x_\infty}^c(\xhat(t_0))$ for every
$t\ge t_0$ and $\lim_{t\to\infty}\xhat(t)=0$. Moreover, 
the compact Lie group $\rG_{x_\infty}$
acts on $\im(L^c_{x_\infty})^\perp$ by unitary 
automorphisms.  Hence, by Theorem~\ref{thm:HMc},
there exists an element $\xi\in\Lambda$ such that
$$
\lim_{t\to\infty}\exp(\bi t\xi)\xhat(t_0)=0,\qquad 
L_{x_\infty}\xi=0.
$$
Since $\i\xi$ acts on the tangent space $T_{x_\infty}X$ 
by a Hermitian endomorphism, it follows that the 
function $t\mapsto\abs{\exp(\bi t\xi)\xhat(t_0)}$ 
is decreasing.  Hence the vector 
$$
\exp(\bi t\xi)\xhat(t_0)\in \im(L^c_{x_\infty})^\perp
$$ 
is contained in the domain $U_\infty$ of the 
holomorphic coordinate chart $\psi$ for all $t\ge 0$.
Hence
$$
\psi(\exp(\bi t\xi)\xhat(t_0))
= \exp(\bi t\xi)\psi(\xhat(t_0)) 
= \exp(\bi t\xi)h(t_0)^{-1}x(t_0)
$$
for $t\ge0$, hence
$$
\lim_{t\to\infty}\exp(\bi t\xi)h(t_0)^{-1}x(t_0) = \psi(0) = x_\infty
$$
and hence 
$$
w_\mu(h(t_0)^{-1}x(t_0),\xi) = \inner{\mu(x_\infty)}{\xi} = 0.
$$
This proves Step~5.

\bigbreak

\medskip\noindent{\bf Step~6.}
{\it There exists an element $\zeta\in\Lambda^c$ 
such that}
$$
w_\mu(x_0,\zeta)=0,\qquad
\lim_{t\to\infty}\exp(\i t\zeta)x_0\notin\rG^c(x_0).
$$

\noindent
Let $g,h,\xi$ be as in Steps~4 and~5 and define
$$
\zeta := g(t_0)h(t_0)\xi h(t_0)^{-1}g(t_0)^{-1}.
$$
Then, by part~(i) of Theorem~\ref{thm:MUMFORD1},
\begin{equation*}
\begin{split}
w_\mu(x_0,\zeta)
&= 
w_\mu(h(t_0)^{-1}g(t_0)^{-1}x_0,
h(t_0)^{-1}g(t_0)^{-1}\zeta g(t_0)h(t_0)) \\
&=
w_\mu(h(t_0)^{-1}x(t_0),\xi) \\
&= 
0
\end{split}
\end{equation*}
and
\begin{equation*}
\begin{split}
\lim_{t\to\infty}\exp(\i t\zeta)x_0
&=
\lim_{t\to\infty}\exp\Bigl(\i tg(t_0)h(t_0)\xi h(t_0)^{-1}g(t_0)^{-1}\Bigr)x_0 \\
&= 
g(t_0)h(t_0)\lim_{t\to\infty}\exp(\i t\xi)h(t_0)^{-1}x(t_0)  \\
&=
g(t_0)h(t_0)x_\infty.
\end{split}
\end{equation*}
Hence 
$
\lim_{t\to\infty}\exp(\i t\zeta)x_0\notin\rG^c(x_0)
$
by~\eqref{eq:NP1} and this proves Step~6.

\medskip\noindent{\bf Step~7.}
{\it There exists an element $\xi\in\Lambda$ 
that satisfies~\eqref{eq:NP2}, i.e.}
$$
w_\mu(x_0,\xi)=0,\qquad
\lim_{t\to\infty}\exp(\i t\xi)x_0\notin\rG^c(x_0).
$$

\medskip\noindent
Let $\zeta\in\Lambda^c$ be as in Step~6. 
By Theorem~\ref{thm:MUMFORD} there exist elements 
$$
p,p^+\in\rP(\zeta)
$$ 
such that 
$$
\xi:=p^{-1}\zeta p\in\Lambda,\qquad 
p^+=\lim_{t\to\infty}\exp(\i t\zeta)p\exp(-\i t\zeta).
$$ 
Hence,  by part~(ii) of Theorem~\ref{thm:MUMFORD1},
$$
w_\mu(x_0,\xi)=w_\mu(x_0,\zeta)=0
$$ 
and
\begin{equation*}
\begin{split}
\lim_{t\to\infty}\exp(\i t\xi)x_0
&=
\lim_{t\to\infty} p^{-1}
\exp(\i t\zeta)p\exp(-\i t\zeta)\exp(\i t\zeta)x_0 \\
&= 
p^{-1}p^+\lim_{t\to\infty}\exp(\i t\zeta)x_0 \\
&\notin
\rG^c(x_0).
\end{split}
\end{equation*}
The last assertion follows from Step~6.
This proves~\eqref{eq:NP2}, Step~7, and Theorem~\ref{thm:HMps}.
\end{proof}

\begin{thm}[{\bf Hilbert--Mumford Criterion: Stable Case}]\label{thm:HMs}
\ 

\noindent
For every $x_0\in X$ the following are 
equivalent.\index{Hilbert--Mumford Criterion!for stability}

\smallskip\noindent{\bf (i)}
$x_0$ is $\mu$-stable.

\smallskip\noindent{\bf (ii)}
Every $\xi\in\Lambda$ satisfies $w_\mu(x_0,\xi)>0$.

\smallskip\noindent{\bf (iii)}
Every $\xi\in\cg\setminus\{0\}$ satisfies $w_\mu(x_0,\xi)>0$.

\smallskip\noindent{\bf (iv)}
Every $\zeta\in\Lambda^c$ satisfies $w_\mu(x_0,\zeta)>0$.

\smallskip\noindent{\bf (v)}
Every $\zeta\in\sT^c$ satisfies $w_\mu(x_0,\zeta)>0$.
\end{thm}

\begin{proof}
That~(i)$\implies$(v) was proved in part~(iii) of Theorem~\ref{thm:HMnec}, 
and the implications~(v)$\implies$(iv)$\implies$(ii) 
and~(v)$\implies$(iii)$\implies$(ii) follow directly from the definitions.  
Thus it remains to prove that~(ii) implies~(i).

Assume $w_\mu(x_0,\xi)>0$ for all $\xi\in\Lambda$.
Then $x_0$ is $\mu$-polystable by Theorem~\ref{thm:HMps}.  
Hence there exists an element~${g\in\rG^c}$ such that~${\mu(gx_0)=0}$.  
Assume, by contradiction, that~${\ker\,L_{gx_0}\ne\{0\}}$.
Then~${\ker L_{gx_0}}$ is a nontrivial Lie subalgebra of~$\cg$
and hence~${\Lambda\cap\ker\,L_{gx_0}\ne\emptyset}$.
Let~${\xi_0\in\Lambda\cap\ker\,L_{gx_0}}$.
Then~${w_\mu(gx_0,\xi_0)=0}$ by Lemma~\ref{le:POLYWEIGHT1}
and so~${w_\mu(x_0,g^{-1}\xi_0g)=0}$ 
by part~(i) of Theorem~\ref{thm:MUMFORD1}.
Since~${g^{-1}\xi_0g\in\Lambda^c}$, it follows from 
Theorem~\ref{thm:MUMFORD} and part~(ii) of 
Theorem~\ref{thm:MUMFORD1} that there exists an 
element~${\xi\in\Lambda}$ such that~${w_\mu(x_0,\xi)=0}$,
in contradiction to our assumption. 
This shows that~${\ker L_{gx_0}=\{0\}}$,
hence we have~${\ker L_{gx_0}^c=\{0\}}$ by Lemma~\ref{le:ONE},
and therefore~${\ker\,L_{x_0}^c=\{0\}}$. 
Thus~$x_0$ is $\mu$-stable and this proves Theorem~\ref{thm:HMs}.
\end{proof}

\begin{cor}\label{cor:M}
For every $x_0\in X$ the following holds.

\smallskip\noindent{\bf (i)}
$x_0$ is $\mu$-unstable if and only if $m_\mu(x_0)<0$.

\smallskip\noindent{\bf (ii)}
$x_0$ is $\mu$-semistable if and only if $m_\mu(x_0)\ge0$.

\smallskip\noindent{\bf (iii)}
$x_0$ is $\mu$-stable if and only if $m_\mu(x_0)>0$.

\smallskip\noindent{\bf (iv)}
$x_0$ is $\mu$-stable if and only if the moment-weight 
inequality~\eqref{eq:MW3} is strict, i.e.\ 
$m_\mu(x_0)+\inf_{\rG^c(x_0)}\abs{\mu}>0$.
\end{cor}

\begin{proof}
Part~(i) follows directly from Theorem~\ref{thm:HMus}
and the definition of the Mumford numerical 
function~$m_\mu$ in~\eqref{eq:MUMFORD}.  
Part~(ii) follows~from~(i).   
Part~(iii) follows from Theorem~\ref{thm:HMs}
and Theorem~\ref{thm:MODSTAB}.
To prove part~(iv), observe that
$
m_\mu(x_0)>0=\inf_{\rG^c(x_0)}\abs{\mu}
$ 
whenever~$x_0$ is $\mu$-stable (by~(iii)), 
that
$
-m_\mu(x_0)=\inf_{\rG^c(x_0)}\abs{\mu}
$ 
whenever~$x_0$ is $\mu$-unstable 
(by Corollary~\ref{cor:KEMPF}),
and that
$
m_\mu(x_0)=0=\inf_{\rG^c(x_0)}\abs{\mu}
$ 
whenever~$x_0$ is $\mu$-semistable, 
but not $\mu$-stable (by parts~(ii) and~(iii)).
This proves Corollary~\ref{cor:M}.
\end{proof}


\chapter{Critical orbits}\label{ch:CRIT}

The results of this chapter are based on the work of
G\'abor Sz\'ekelyhidi~\cite{S}. We assume throughout that~$\rG$ 
is a compact Lie group whose Lie algebra~${\cg=\Lie(\rG)}$ 
is equipped with an invariant inner product,
that it acts on a closed K\"ahler manifold~$X$ by K\"ahler isometries, 
and that the action is generated by an equivariant moment map~${\mu:X\to\cg}$. 
The goal of this section is to examine complexified group orbits 
that contain critical points of the moment map squared.  
In the notation of Chapter~\ref{ch:MOMENT} the problem 
can be rephrased as that of finding a solution of the 
equation 
$
L_x\mu(x)=0
$ 
as $x$ ranges over a $\rG^c$-orbit in $X$.  
This is the finite-dimensional analogue of finding an {\it extremal metric} 
in the cscK setting described in Chapter~\ref{ch:INTRO}.
The main result is the generalized Sz\'ekelyhidi criterion 
in Theorem~\ref{thm:SZEK1} for the existence of a critical point 
in the complexified group orbit, in terms of polystability 
with respect to the action of a suitable quotient group.  
As a warmup we begin with the following criterion in terms 
of the gradient flow of the square of the moment map, 
which is analogous to Theorem~\ref{thm:STABILITY}. 

\begin{thm}[{\bf Critical Orbits}]\label{thm:CRIT}
Let $x_0\in X$, let $x:\R\to X$ be the unique solution of~\eqref{eq:KN1},
and define $x_\infty:=\lim_{t\to\infty}x(t)$.  Then the following are equivalent.

\smallskip\noindent{\bf (i)}
$\rG^c(x_0)$ contains a critical point of the square of the moment map.

\smallskip\noindent{\bf (ii)} 
$x_\infty\in\rG^c(x_0)$. 
\end{thm}

\begin{proof} 
We prove that~(i) implies~(ii).   Assume that
there exists an~${x\in\rG^c(x_0)}$ such that~${L_x\mu(x)=0}$. 
Then~${\Abs{\mu(x)}=\inf_{g\in\rG^c}\mu(gx_0)=\Abs{\mu(x_\infty)}}$
by Corollary~\ref{cor:KIRNESS1} and Theorem~\ref{thm:MLT}.
Hence~${x_\infty\in\rG(x)\subset\rG^c(x_0)}$ 
by Theorem~\ref{thm:NESS2}.  This shows that~(i) implies~(ii). 
The converse implication follows from the fact 
that~${L_{x_\infty}\mu(x_\infty)=0}$ by Theorem~\ref{thm:XINFTY}. 
This proves Theorem~\ref{thm:CRIT}.
\end{proof}

\bigbreak

In his PhD thesis~\cite{S} G\'abor Sz\'ekelyhidi found a criterion for 
the existence of a critical point of the moment map squared in the 
complexified group orbit in terms of polystability with respect to 
the action of a suitable subgroup. To describe his criterion we 
recall the the notations 
$$
\rG_x := \left\{u\in\rG\,|\,ux=x\right\},\qquad
\rG_x^c := \left\{g\in\rG^c\,|\,gx=x\right\}
$$
for the compact and complex stabilizer subgroups.
Recall also that~$\rG_x^c$ is the complexification 
of~$\rG_x$ whenever~${\mu(x)=0}$ (Lemma~\ref{le:isotropy}).
However, in general, the {\it complexified} stabilizer subgroup~$(\rG_x)^c$ 
is a proper subgroup of the {\it complex} stabilizer subgroup~$\rG_x^c$. 
The latter may not even be reductive. 

Throughout a {\it torus} is a compact connected abelian Lie group.
The Sz\'ekelyhidi criterion requires the choice of a maximal torus 
$$
\rT\subset\rG^c_x
$$
and it may not be possible to choose 
this torus such that it is contained in~$\rG$. 
However, if~${\rT\subset\rG^c}$ 
is any torus and~$\ct:=\Lie(\rT)$ is its Lie algebra,
then~${\ct\setminus\{0\}\subset\sT^c}$
and so the $\rG^c$-invariant pairing 
\begin{equation}\label{eq:innerc}
\inner{\zeta_1}{\zeta_2}_c 
:= \inner{\Re(\zeta_1)}{\Re(\zeta_2)}
- \inner{\Im(\zeta_1)}{\Im(\zeta_2)}
\end{equation}
on $\cg^c$ is positive definite on $\ct$ by 
Lemma~\ref{le:WEIGHT3}. This implies that there exists 
a unique linear projection~${\Pi_\rT:\cg^c\to\ct}$ such that 
\begin{equation}\label{eq:PIT}
\inner{\zeta-\Pi_\rT(\zeta)}{\tau}_c=0
\qquad\mbox{for all }\zeta\in\cg^c
\mbox{ and all }\tau\in\ct.
\end{equation}
With this understood, we introduce the following notation.

Let $x\in X$ and let $\rT\subset\rG^c_x$ be a torus
with the Lie algebra
$$
\ct:=\Lie(\rT)\subset\cg^c.
$$ 
Let $\rG_T^c\subset\rG^c$ be the identity component of the 
centralizer of $\rT$ (the subgroup of all elements 
of~$\rG^c$ that commute with each element of~$\rT$),
i.e.\ 
\begin{equation}\label{eq:GcT}
\begin{split}
\cg_\rT^c 
&:= 
\left\{\zeta\in\cg^c\,\big|\,[\zeta,\tau]=0
\mbox{ for all }\tau\in\ct\right\}, \\
\rG_\rT^c 
&:= 
\left\{g(1)\,\Bigg|\,
\begin{array}{l}
g:[0,1]\to\rG\mbox{ is a smooth path} \\
\mbox{such that }g(0)=\one\mbox{ and } \\
\dot g(t)g(t)^{-1}\in\cg_\rT^c \mbox{ for all }t\in[0,1]
\end{array}\right\}.
\end{split}
\end{equation}
By the Closed Subgroup Theorem $\rG_\rT^c$ is a Lie subgroup of $\rG^c$ 
with the Lie algebra~${\Lie(\rG_\rT^c)=\cg_\rT^c}$. 
For~${x\in X}$ and~${\zeta\in\sT^c\cap\cg_\rT^c\setminus\ct}$
the {\bf $(\mu,\rT)$-weight}\index{weights@$(\mu,T)$-weight} 
of the pair $(x,\zeta)$ is defined by 
\begin{equation}\label{eq:wmuT}
w_{\mu,\rT}(x,\zeta) := \lim_{t\to\infty}
\inner{\mu(\exp(\i t\zeta)x)}{\Re(\zeta-\Pi_\rT(\zeta))}.
\end{equation}
With this terminology in place we are in a position 
to formulate the main results of this section.

\bigbreak


\begin{thm}[{\bf Generalized Sz\'ekelyhidi Criterion}]
\label{thm:SZEK1}
\ 

\noindent
Let~${x\in X}$, let~${\rT\subset\rG^c_x}$ be a maximal 
torus with the Lie algebra~$\ct:=\Lie(\rT)$, and let~${g\in\rG^c}$ 
such that~${g\rT g^{-1}\subset\rG}$.
Then the following are\index{Sz\'ekelyhidi Criterion!generalized}
equivalent.

\smallskip\noindent{\bf (i)}
$\rG^c(x)$ contains a critical point of 
the square of the moment map.

\smallskip\noindent{\bf (ii)} 
There exists an element $h\in\rG_\rT^c$ such that $g^{-1}\mu(ghx)g\in\ct$.
\end{thm}

\begin{proof}
Let~${h\in\rG_\rT^c}$ such that~${g^{-1}\mu(ghx)g\in\ct\subset\ker L_x^c}$. 
Then~${h\tau h^{-1}=\tau}$ and hence~${L_{hx}^c\tau=hL_x^c\tau=0}$ 
for all~${\tau\in\ct}$. Thus
$$
L_{ghx}\mu(ghx)=gL_{hx}^c(g^{-1}\mu(ghx)g)=0.
$$
This shows that~(ii) implies~(i). 
The converse is proved on page~\pageref{proof:SZEK1}.
\end{proof}

\begin{thm}[{\bf Sz\'ekelyhidi Moment-Weight Inequality}]
\label{thm:SZEK2}
Let~${x\in X}$,\index{moment-weight inequality!Sz\'ekelyhidi}
let~${\rT\subset\rG^c_x}$ be a torus with the Lie algebra~$\ct:=\Lie(\rT)$, 
and let $g\in\rG^c$ such that~${g\rT g^{-1}\subset\rG}$.
Then~${\inf_{h\in\rG^c}\abs{\mu(hx)}\ge\abs{\Pi_\rT(g^{-1}\mu(gx)g)}_c}$ 
and
\begin{equation}\label{eq:SMW0}
\sup_{\zeta\in\sT^c\cap\cg^c_\rT\setminus\ct}
\frac{-w_{\mu,\rT}(x,\zeta)}{\sqrt{\abs{\zeta}_c^2 - \abs{\Pi_\rT(\zeta)}_c^2}}
\le
\inf_{h\in\rG^c}\sqrt{\abs{\mu(hx)}^2-\abs{\Pi_\rT(g^{-1}\mu(gx)g)}_c^2}.
\end{equation}
Moreover, the supremum on the left is attained, it remains unchanged 
when taken over all~${\zeta\in\sT^c\cap\cg^c_\rT\setminus\ct}$ 
with~${\exp(\zeta)\in\rT}$, ${g\zeta g^{-1}\in\cg}$, ${\Pi_\rT(\zeta)=0}$,
and equality holds in~\eqref{eq:SMW0} if and only if
the left hand side is nonnegative.
\end{thm}

\begin{proof}
See page~\pageref{proof:SZEK2}.
\end{proof}

\begin{thm}[{\bf Hilbert--Mumford Criterion: Critical Orbits}]
\label{thm:SZEK3}
\ 

\noindent
Let~${x\in X}$ and fix a maximal torus~${\rT\subset\rG_x^c}$ 
with the Lie algebra~${\ct:=\Lie(\rT)}$.
Choose an element~${g\in\rG^c}$ such that~${g\rT g^{-1}\subset\rG}$.
Then the following are\index{Hilbert--Mumford Criterion!for critical orbits} 
equivalent.

\smallskip\noindent{\bf (i)}
$\rG^c(x)$ contains a critical point of 
the square of the moment map.

\smallskip\noindent{\bf (ii)}
If ${\zeta\in\sT^c}$ satisfies~${[\zeta,\ct]=0}$,
${\exp(\zeta)\in\rT}$, ${g\zeta g^{-1}\in\cg}$,
and~${\Pi_\rT(\zeta)=0}$, then~${w_\mu(x,\zeta)\ge0}$ and
$$
w_\mu(x,\zeta)=0\qquad\implies\qquad
\lim_{t\to\infty}\exp(\i t\zeta)x\in\rG_\rT^c(x).
$$

\smallskip\noindent{\bf (iii)}
If $\zeta\in\sT^c\setminus\ct$ satisfies $[\zeta,\ct]=0$,
then $w_{\mu,\rT}(x,\zeta)\ge0$ and
$$
w_{\mu,\rT}(x,\zeta)=0\qquad\implies\qquad
\lim_{t\to\infty}\exp(\i t\zeta)x\in\rG_\rT^c(x).
$$
\end{thm}

\begin{proof}
That~(iii) implies~(ii) follows from the fact 
that~${w_{\mu,\rT}(x,\zeta)=w_\mu(x,\zeta)}$
whenever~${\zeta\in\sT^c}$ satisfies~${[\zeta,\ct]=0}$
and~${\Pi_\rT(\zeta)=0}$ (see equation~\eqref{eq:wmuT}).  
That~(ii) implies~(i) and~(i) implies~(iii) 
will be proved on page~\pageref{proof:SZEK3}.
\end{proof}

We will see below that the Sz\'ekelyhidi criterion
in Theorem~\ref{thm:SZEK1} is a restatement 
of the polystability condition on $x$  
with respect to the action of a suitable quotient group 
on a suitable submanifold of $X$.  In the case
$$
g=\one,
$$
where~$\rT$ is a subgroup of~$\rG$, the Lie group 
in question is the quotient 
$$
\rG_\rT/\rT
$$
of the identity component $\rG_\rT\subset\rG$ of the 
centralizer of $\rT$ by the torus $\rT$, and the submanifold
$X_\rT$, on which it acts, consists of all elements of $X$
that contain the torus $\rT$ in their stabilizer subgroup.
Then the moment map descends to an equivariant moment map 
$$
\mu_\rT:X_\rT\to\cg_\rT/\ct
$$
for this action and, in the case $g=\one$, condition~(ii) 
in Theorem~\ref{thm:SZEK1} says that $x\in X_\rT$ 
is $\mu_\rT$-polystable (see Lemma~\ref{le:GT}).
Theorem~\ref{thm:SZEK2} states the corresponding moment-weight
inequality for $\mu_\rT$ and Theorem~\ref{thm:SZEK3} adapts 
the Hilbert--Mumford criterion for polystability to this setting. 
The proofs of all three theorems reduce to the case~${g=\one}$,
either by replacing $(x,\rT)$ with 
$$
(gx,g\rT g^{-1})
$$ 
or, equivalently, by leaving the pair~${(x,\rT)}$ 
unchanged while replacing the triple~${(\rG,\om,\mu)}$ 
with~${(g^*\rG,g^*\om,g^*\mu)}$, where 
$$
g^*\rG := g^{-1}\rG g,\qquad 
(g^*\om)_x(\xhat_1,\xhat_2) := \om_{gx}(g\xhat_1,g\xhat_2),
$$
and
$$
(g^*\mu)(x) := g^{-1}\mu(gx)g.
$$ 
This is the first ingredient in the proof of the Sz\'ekelyhidi
criterion and is explained in Lemma~\ref{le:CONJUGATE} below.
The second ingredient in the proof is the study of the 
action of the quotient group $\rG_\rT/\rT$ on the 
submanifold~$X_\rT$ and its moment map $\mu_\rT:X_\rT\to\cg_\rT/\ct$
in the case $\rT\subset\rG$. 
This is the content of the Lemma~\ref{le:GT}.
The third ingredient, required for the proof of the 
generalized Sz\'ekelyhidi moment-weight inequality
and of the Hilbert--Mumford criterion for critical orbits,
is the study of the toral generators for the quotient 
group (Lemma~\ref{le:GcT}) and the corresponding Mumford 
numerical invariants.  These are shown in Lemma~\ref{le:TWEIGHT}
to be the $(\mu,\rT)$-weights in~\eqref{eq:wmuT}.
After these preparations we are ready to prove the 
main theorems of this section. The proof for $g=\one$ 
will again use the gradient flow of the moment map 
squared and the Kempf--Ness function as the central
technical ingredients.


\subsubsection*{Conjugation and the balancing condition}

The purpose of the present subsection is to explain why it suffices to 
prove Theorems~\ref{thm:SZEK1}, \ref{thm:SZEK2}, and~\ref{thm:SZEK3} 
in the case~${g=\one}$. 

\begin{lem}\label{le:CONJUGATE}
Fix an element $g\in\rG^c$, define
\begin{equation}\label{eq:tG}
\trG:=g^{-1}\rG g,\qquad
\tcg:=g^{-1}\cg g,
\end{equation}
and define the map $\tmu:X\to\tcg$ by 
\begin{equation}\label{eq:tmu}
\tmu(x) := g^{-1}\mu(gx)g\qquad\mbox{for }x\in X.
\end{equation}
Let 
$$
\tom:=g^*\om
$$ 
be the pullback of $\om$ under the diffeomorphism 
induced by $g$. Then the following holds.

\smallskip\noindent{\bf (i)}
The bilinear form~\eqref{eq:innerc} 
defines an invariant inner product on~$\tcg$. 

\smallskip\noindent{\bf (ii)}
The $2$-form $\tom$ is a K\"ahler form on the complex manifold $(X,J)$.

\smallskip\noindent{\bf (iii)}
The Lie group $\trG$ acts on $(X,\tom,J)$ by K\"ahler isometries.

\smallskip\noindent{\bf (iv)}
The map $\tmu:X\to\tcg$ is $\trG$-equivariant.

\smallskip\noindent{\bf (v)}
The map $\tmu:X\to\tcg$ is a moment map for the $\trG$-action on $(X,\tom)$.

\smallskip\noindent{\bf (vi)}
Let $x\in X$ and let $\rT\subset\rG^c_x$ be a torus with Lie algebra $\ct$.
Then
\begin{equation}\label{eq:wtmu}
w_{\tmu,\rT}(x,\zeta)=w_{\mu,\rT}(x,\zeta)
\end{equation}
for every $\zeta\in\sT^c\cap\cg_\rT^c\setminus\ct$.
\end{lem}

\begin{proof}
We prove part~(i).  The bilinear form~\eqref{eq:innerc} on $\cg^c$
is symmetric by definition and, since $\tcg\setminus\{0\}\subset\sT^c$,
its restriction to $\tcg$ is positive definite by Lemma~\ref{le:WEIGHT3}.
This proves~(i).

Part~(ii) follows from the fact that $\om$ is a K\"ahler form on $(X,J)$
and that the diffeomorphism induced by $g$ preserves the complex structure.

We prove part~(iii).  Fix an element $\tu\in\trG$ and define 
$$
u:=g\tu g^{-1}.
$$
Then $u\in\rG$ by~\eqref{eq:tG} and hence 
$$
\tu^*\tom
= \tu^*g^*\om
= (g\tu)^*\om
= (ug)^*\om
= g^*u^*\om
= g^*\om=\tom.
$$
Moreover, the $\trG$-action preserves the complex structure 
because $\trG\subset\rG^c$. This proves~(iii).

\bigbreak

We prove part~(iv).  Let $x\in X$ and $\tu\in\trG$, 
and define~${u:=g\tu g^{-1}\in\rG}$ as above.  
Then, by~\eqref{eq:tmu},
\begin{equation*}
\begin{split}
\tmu(\tu x) 
&= g^{-1}\mu(g\tu x)g \\
&= g^{-1}\mu(ugx)g \\
&= g^{-1}u\mu(gx)u^{-1}g \\
&= \tu g^{-1}\mu(gx)g\tu^{-1} \\
&= \tu\tmu(x)\tu^{-1}
\end{split}
\end{equation*}
and this proves~(iv).

We prove part~(v). Fix elements $x\in X$, $\xhat\in T_xX$, 
and $\txi\in\tcg$, and define 
$$
\xi:=g\txi g^{-1}.  
$$
Then $\xi\in\cg$ by~\eqref{eq:tG} and
\begin{equation*}
\begin{split}
\tom_x(L^c_x\txi,\xhat)
&= \om_{gx}(gL^c_x\txi,g\xhat) \\
&= \om_{gx}(L_{gx}(g\txi g^{-1}),g\xhat) \\
&= \om_{gx}(L_{gx}\xi,g\xhat) \\
&= \inner{d\mu(gx)g\xhat}{\xi} \\
&= \inner{g^{-1}(d\mu(gx)g\xhat)g}{g^{-1}\xi g}_c \\
&= \langle d\tmu(x)\xhat,\txi\rangle_c.
\end{split}
\end{equation*}
Here the penultimate equality follows from Lemma~\ref{le:WEIGHT3}
and the last equality follows from~\eqref{eq:tmu}.  This proves~(v).

We prove part~(vi).  Let $x\in X$, let $\rT\subset\rG^c_x$ be a torus 
with Lie algebra $\ct$, let~${\zeta\in\sT^c\setminus\ct}$ with $[\zeta,\ct]=0$, 
and define $x^+:=\lim_{t\to\infty}\exp(\i t\zeta)x$.
Then
$$
x^+=\lim_{t\to\infty}\exp(\i t(\zeta-\Pi_\rT(\zeta)))x
$$ 
because~${L_x^c\Pi_\rT(\zeta)=0}$ and~$\zeta$ commutes with~$\Pi_\rT(\zeta)$.
Hence $L_{x^+}^c(\zeta-\Pi_\rT(\zeta))=0$ by Lemma~\ref{le:WEIGHT2}, 
and so
\begin{equation*}
\begin{split}
w_{\mu,\rT}(x,\zeta) 
&= 
\inner{\mu(x^+)}{\zeta-\Pi_\rT(\zeta)}_c \\
&= 
\inner{\mu(gx^+)}{g(\zeta-\Pi_\rT(\zeta))g^{-1})}_c \\
&= 
\inner{g^{-1}\mu(gx^+)g}{\zeta-\Pi_\rT(\zeta)}_c \\
&= 
w_{\tmu,\rT}(x,\zeta). 
\end{split}
\end{equation*}
Here the second equality follows from Lemma~\ref{le:WEIGHT4},
the third equality follows from Lemma~\ref{le:WEIGHT3}, 
and the first and last equalities follow from the 
definition of the relative weights in~\eqref{eq:wmuT}.
This proves~(vi) and Lemma~\ref{le:CONJUGATE}.
\end{proof}

\begin{defn}\label{def:balanced}
An element $x\in X$\index{balanced point@$\mu$-balanced point}
is called {\bf $\mu$-balanced}
if there exists a maximal torus~${\rT\subset\rG^c_x}$ 
such that $\rT\subset\rG$. 
\end{defn}

Here we slightly abuse notation because the balancing condition 
depends only on the maximal compact subgroup of~$\rG^c$ 
but not on the choice of the moment map.
Lemma~\ref{le:CONJUGATE} shows that it suffices to prove 
Theorems~\ref{thm:SZEK1}, \ref{thm:SZEK2}, and~\ref{thm:SZEK3}
under the assumption that~$x$ is $\mu$-balanced and~${g=\one}$.

\begin{exple}\label{ex:BALANCED}\rm
\smallskip\noindent{\bf (i)}
If $g\in\rG^c$ belongs to the center of $\rG^c$, then the
maximal compact subgroup~${\trG=g^{-1}\rG g}$ in 
Lemma~\ref{le:CONJUGATE} is equal to $\rG$, however,
the resulting symplectic form $\tom:=g^*\om$ and the moment 
map $\tmu=g^*\mu$ may well be different from $\om$
and $\mu$, respectively.

\smallskip\noindent{\bf (ii)} 
Every fixed point of the $\rG^c$-action is $\mu$-balanced
and so is every point with discrete isotropy in~$\rG^c$.

\smallskip\noindent{\bf (iii)}
If~$\rG$ is abelian then every torus in~$\rG^c$ is necessarily 
contained in~$\rG$ and so every element of~$X$ is $\mu$-balanced.

\smallskip\noindent{\bf (iv)}
Every $\rG^c$-orbit contains a $\mu$-balanced element.
Let~${x\in X}$ and suppose that~${\rT\subset\rG_x^c}$ is a 
maximal torus. By Lemma~\ref{le:KG} there is a $g\in\rG^c$
such that~${g\rT g^{-1}\subset\rG}$.  Then~${g\rT g^{-1}}$ 
is a maximal torus in~${g\rG^c_x g^{-1}=\rG_{gx}^c}$
and hence~$gx$ is $\mu$-balanced.

\smallskip\noindent{\bf (v)}
Consider the diagonal action of~${\rG=\SO(3)}$ on~${(S^2)^n}$.  
Then~${\rG^c}$ is the group of M\"obius transformations, 
the moment map~${\mu:(S^2)^n\to\R^3}$ 
is given by~${\mu(x)=\sum_{i=1}^nx_i}$ 
for~${x=(x_1,\dots,x_n)\in(S^2)^n}$, and its critical 
points are the $n$-tuples~${x\in(S^2)^n}$ 
that satisfy~${x_i=\pm x_j}$ for all~$i$ and~$j$.
Now let~${x\in(S^2)^n}$ such that~${x_2=\cdots=x_n\ne\pm x_1}$.
Then~${\rG_x=\{\one\}}$ and~${\rG_x^c\cong\C^*}$.  
Thus~$\rG^c_x$ contains a unique maximal 
torus~$\rT\not\subset\rG$ and so~$x$ is not $\mu$-balanced. 
\end{exple} 

\begin{cor}[{\bf Generalized Sz\'ekelyhidi Criterion}]\label{cor:SZEK0}
Let~${x\in X}$ be a $\mu$-balanced element
and let~${\rT\subset\rG_x}$ be a maximal torus 
with the Lie algebra~${\ct:=\Lie(\rT)}$.
Then\index{Hilbert--Mumford Criterion!for critical orbits} 
$
\inf_{g\in\rG^c}\abs{\mu(gx)}\ge\abs{\Pi_\rT(\mu(x))}
$ 
and\index{moment-weight inequality!Sz\'ekelyhidi}
\begin{equation}\label{eq:SMWbalanced}
\sup_{\xi\in\cg\setminus\{0\},\,[\xi,\ct]=0,\,\xi\perp\ct}
\frac{-w_\mu(x,\xi)}{\abs{\xi}}
\le
\inf_{g\in\rG^c}\sqrt{\abs{\mu(gx)}^2-\abs{\Pi_\rT(\mu(x))}^2}.
\end{equation}
Moreover,\index{Sz\'ekelyhidi Criterion!for $\mu$-balanced points}
the following are equivalent.

\smallskip\noindent{\bf (i)}
$\rG^c(x)$ contains a critical point of 
the square of the moment map.

\smallskip\noindent{\bf (ii)} 
There exists an element $h\in\rG_\rT^c$ such that $\mu(hx)\in\ct$.

\smallskip\noindent{\bf (iii)} 
If ${\xi\in\cg\setminus\{0\}}$ satisfies~${[\xi,\ct]=0}$
and~${\xi\perp\ct}$ then~${w_\mu(x,\xi)\ge0}$ and
$$
w_\mu(x,\xi)=0\qquad\implies\qquad
\lim_{t\to\infty}\exp(\i t\xi)x\in\rG_\rT^c(x).
$$
\end{cor}

\begin{proof}
Theorem~\ref{thm:SZEK1}, Theorem~\ref{thm:SZEK2}, and Theorem~\ref{thm:SZEK3}.
\end{proof}


\subsubsection*{Centralizer and quotient group}

For a torus $\rT\subset\rG$ with the Lie algebra~${\ct:=\Lie(\rT)\subset\cg}$
denote by $\rG_T\subset\rG$ the identity component 
of the centralizer of $\rT$ (the subgroup of all elements 
of~$\rG$ that commute with each element of~$\rT$),
i.e.\ 
\begin{equation}\label{eq:GT}
\begin{split}
\cg_\rT 
&:= 
\left\{\xi\in\cg\,\big|\,[\xi,\tau]=0
\mbox{ for all }\tau\in\ct\right\}, \\
\rG_\rT 
&:= 
\left\{u(1)\,\Bigg|\,
\begin{array}{l}
u:[0,1]\to\rG\mbox{ is a smooth path} \\
\mbox{such that }u(0)=\one\mbox{ and } \\
\dot u(t)u(t)^{-1}\in\cg_\rT \mbox{ for all }t\in[0,1]
\end{array}\right\}.
\end{split}
\end{equation}
By the Closed Subgroup Theorem,~$\rG_\rT$ is a Lie subgroup of~$\rG$ 
with the Lie algebra~${\Lie(\rG_\rT)=\cg_\rT}$.  
Moreover,~$\rT$ is a subgroup of the center of~$\rG_\rT$.   
Thus~$\rT$ is a normal subgroup of~$\rG_\rT$ and~$\rG_\rT/\rT$ 
is a compact Lie group with the Lie algebra~${\Lie(\rG_\rT/\rT)=\cg_\rT/\ct}$.
For~${\xi\in\cg_\rT}$ denote by
$$
[\xi]_\rT:=\xi+\ct\in\cg_\rT/\ct
$$
the equivalence class of~$\xi$. 
The quotient group acts on the space 
\begin{equation}\label{eq:XT}
X_\rT := \left\{x\in X\,\big|\,\rT\subset\rG_x\right\}
\end{equation}
of all elements of~$X$ that contain the torus~$\rT$ in
their stabilizer subgroup.  The next lemma shows 
that~$X_\rT$ is a complex submanifold of~$X$ and that 
the~$\rG_\rT/\rT$-action on~$X_\rT$ is Hamiltonian. 

\begin{lem}\label{le:GT}
Let $\rT\subset\rG$ be a torus with the Lie algebra $\ct$,
let $\rG_\rT,\cg_\rT$ be as in~\eqref{eq:GT},
and let $\rG_\rT^c,\cg_\rT^c$ be as in~\eqref{eq:GcT}.
Then the following holds.

\smallskip\noindent{\bf (i)}
The formula
\begin{equation}\label{eq:innerT}
\inner{\left[\xi\right]_\rT}{\left[\eta\right]_\rT}_{\cg_\rT/\ct} 
:= \inner{\xi}{\eta-\Pi_\rT(\eta)}
\end{equation}
for $\xi,\eta\in\cg_\rT$ defines an invariant inner product on $\cg_\rT/\ct$.

\smallskip\noindent{\bf (ii)}
If $\zeta\in\sT^c\cap\cg_\rT^c\setminus\ct$ then 
$\zeta+\tau\in\sT^c\cap\cg_\rT^c\setminus\ct$
for all $\tau\in\ct$ and 
$$
\Im(\zeta)\perp\ct.
$$

\smallskip\noindent{\bf (iii)}
The set $X_\rT$ in~\eqref{eq:XT} is a closed $\rG_\rT$-invariant 
complex submanifold of~$X$ and~${\mu(X_\rT)\subset\cg_\rT}$.   

\smallskip\noindent{\bf (iv)}
The quotient group $\rG_\rT/\rT$ acts on $X_\rT$ by K\"ahler isometries,
and the action is generated by the $\rG_\rT/\rT$-equivariant moment map 
\begin{equation}\label{eq:muT}
\mu_\rT:X_\rT\to\cg_\rT/\ct,\qquad \mu_\rT(x) := [\mu(x)]_\rT = \mu(x)+\ct.
\end{equation}
\end{lem}

\begin{proof}
See page~\pageref{proof:GT}.
\end{proof}

\bigbreak

In preparation for the proof of Lemma~\ref{le:GT} 
we establish some basic properties of the projections~$\Pi_\rT$. 

\begin{lem}\label{le:PIT}
Let $\rT\subset\rG^c$ be a torus and let $\Pi_\rT:\cg^c\to\ct$
be the projection defined by~\eqref{eq:PIT}.
Fix two elements $g\in\rG^c$ and $\zeta\in\cg^c$.
Then
\begin{equation}\label{eq:PITg}
\Pi_{g\rT g^{-1}}(g\zeta g^{-1}) = g\Pi_\rT(\zeta)g^{-1}.
\end{equation}
Moreover, if $x\in X$ satisfies $\rT\subset\rG^c_x$ then 
\begin{equation}\label{eq:PITmug}
\inner{\mu(gx)}{\Re(\Pi_{g\rT g^{-1}}(g\zeta g^{-1}))}
= \inner{\mu(x)}{\Re(\Pi_\rT(\zeta))}
\end{equation}
\end{lem}

\begin{proof}
Let $\ct:=\Lie(\rT^c)\subset\cg^c$ 
and define~${\tau:=\Pi_\rT(\zeta)}\in\ct$.  
Then~${\inner{\zeta-\tau}{\ct}_c=0}$ by~\eqref{eq:PIT} 
and hence, by Lemma~\ref{le:WEIGHT4},
$$
\inner{g\zeta g^{-1}-g\tau g^{-1}}{g\ct g^{-1}}_c=0.
$$
Thus~${g\tau g^{-1}=\Pi_{g\rT g^{-1}}(g\zeta g^{-1})}$
and this proves~\eqref{eq:PITg}.

Now let $x\in X$ such that $\rT\subset\rG^c_x$.
Then $L^c_x\tau=0$ and hence it follows from Lemma~\ref{le:WEIGHT4}
that~${\inner{\mu(x)}{\Re(\tau)}=\inner{\mu(gx)}{\Re(g\tau g^{-1})}}$.
Thus equation~\eqref{eq:PITmug} follows from~\eqref{eq:PITg}
and this proves Lemma~\ref{le:PIT}.
\end{proof}

\begin{proof}[Proof of Lemma~\ref{le:GT}]\label{proof:GT}
We prove part~(i).  Let~${\xi,\eta\in\cg_\rT}$ and~${u\in\rG_\rT}$. 
Then we have~${\Pi_\rT(u\eta u^{-1})=\Pi_\rT(\eta)}$
by Lemma~\ref{le:PIT} and so
\begin{equation*}
\begin{split}
\inner{[u\xi u^{-1}]_\rT}{[u\eta u^{-1}]_\rT}_{\cg_\rT/\ct}
&= 
\inner{u\xi u^{-1}}{u\eta u^{-1}-\Pi_\rT(\eta)} \\
&= 
\inner{\xi }{\eta-\Pi_\rT(\eta)} \\
&= 
\inner{[\xi]_\rT}{[\eta]_\rT}_{\cg_\rT/\ct}.
\end{split}
\end{equation*}
This proves~(i).   

We prove part~(ii). Fix an element
$\zeta\in\sT^c\cap\cg_\rT^c\setminus\ct$ and let~${\tau\in\ct}$.
Then~$\zeta$ and~$\tau$ are two commuting toral generators
and $\zeta+\tau\ne0$.  Hence~${\zeta+\tau}$ is again a toral generator,
commutes with $\ct$, and is not an element of $\ct$.  
Thus~${\zeta+\tau\in\sT^c\cap\cg_\rT^c\setminus\ct}$.
This implies 
$
\inner{\Re(\zeta)+\tau}{\Im(\zeta)}=0
$ 
for all $\tau\in\ct$ by Lemma~\ref{le:WEIGHT3}.  
Hence $\Im(\zeta)\perp\ct$ and this 
proves~(ii).

We prove part~(iii).  That $X_\rT$ is a closed subset of $X$ follows 
directly from the definitions, and that it is a submanifold of $X$ 
with the tangent spaces 
$
T_xX_\rT 
= \left\{\xhat\in T_xX\,|\,
a\xhat=\xhat\mbox{ for all }a\in\rT
\right\}
$
for $x\in X_\rT$ is a general fact about smooth group actions.  
That $X_\rT$ is a complex submanifold follows from the fact 
that the $\rT$-action preserves the complex structure.
If~${x\in X_\rT}$ and~${g\in\rG_\rT}$, then~${ax=x}$ for all~${a\in\rT}$, 
hence~${agx=gax=gx}$ for all~${a\in\rT}$ because~$a$ and~$g$ commute,
and hence~${gx\in X_\rT}$.   Thus~$X_\rT$ is $\rG_\rT$-invariant.
Moreover, the group~$\rT$ acts trivially on~$X_\rT$ by definition, 
and hence the action of~$\rG_\rT$ on~$X_\rT$ by K\"ahler 
isometries descends to an action of the quotient group~$\rG_\rT/\rT$. 

\bigbreak

If $x\in X_\rT$, then $\ct\subset\ker L_x$, hence 
$
[\mu(x),\tau]=-d\mu(x)L_x\tau=0
$
for all $\tau\in\ct$ by~\eqref{eq:mu3}, and so $\mu(x)\in\cg_\rT$.
This shows that
$
\mu(X_\rT)\subset\cg_\rT
$ 
and so the map $\mu_\rT:X_\rT\to\cg_\rT/\ct$ in~\eqref{eq:muT} 
is well defined. That it is $\rG_\rT/\rT$-equivariant
follows directly from the $\rG$-equivariance of map $\mu:X\to\cg$.  
This proves~(iii).

We prove part~(iv). Fix elements $x\in X_\rT$, 
$\xhat\in T_xX_\rT$, and $\xi\in\cg_\rT$. Then
\begin{equation*}
\begin{split}
\inner{d\mu_\rT(x)\xhat}{[\xi]_\rT}_{\cg_\rT/\ct}
&=
\inner{d\mu(x)\xhat}{\xi-\Pi_\rT(\xi)} \\
&=
\inner{d\mu(x)\xhat}{\xi} \\
&=
\om(L_x\xi,\xhat).
\end{split}
\end{equation*}
Here the second equality follows from the equation $d\mu(x)^*=JL_x$
in~\eqref{eq:mu3} and the fact that $\ct\subset\ker L_x$. 
This proves~(iv) and Lemma~\ref{le:GT}.
\end{proof}

By definition, an element $x\in X_\rT$ is $\mu_\rT$-polystable if and only if
there exists an element $h\in\rG_\rT^c$ such that $\mu(hx)\in\ct$.
Thus the generalized Sz\'ekelyhidi criterion  in Theorem~\ref{thm:SZEK1}
asserts that, when $\rT\subset\rG^c_x$ is a maximal torus
and~${\rT\subset\rG}$, the complexified group orbit $\rG^c(x)$
contains a critical point of the square of the moment map 
if and only if $x\in X_\rT$ is $\mu_\rT$-polystable.
If~$\rG$ is a torus, this follows directly definitions.


\subsubsection*{The $(\mu,\rT)$-weights}

This subsection establishes the basic properties of the relative
$(\mu,\rT)$-weights, in preparation for the proof of the Sz\'ekelyhidi 
moment-weight inequality and of the Hilbert--Mumford criterion for critical orbits.
The first step is to characterizate of the toral generators in the 
quotient Lie algebra $\cg_\rT^c/\ct^c$. 

\begin{lem}\label{le:GcT}
Let $\rT\subset\rG^c$ be a torus and let $\rG^c_\rT$ and $\cg_\rT^c$
be as in~\eqref{eq:GcT}.  Then $\rG^c_\rT/\rT^c$ is a reductive 
Lie group with the Lie algebra 
$$
\cg_\rT^c/\ct^c = \Lie(\rG^c_\rT/\rT^c).
$$
Moreover, the set of toral generators in $\cg_\rT^c/\ct^c$ is given by
$$
\sT^c_{\rG^c_\rT/\rT^c}
= \left\{[\zeta]_{\rT^c}\,\big|\,\zeta\in\sT^c\cap\cg^c_\rT\setminus\ct\right\}.
$$
Here~${[\zeta]_{\rT^c}:=\zeta+\ct^c}$ denotes 
the equivalence class of $\zeta\in\cg_\rT^c$ in $\cg_\rT^c/\ct^c$. 
\end{lem}

\begin{proof}
Assume first that $\rT\subset\rG$.   Then $\rG^c_\rT/\rT^c$ is a connected 
complex Lie group which contains $\rG_\rT/\rT$ as a maximal compact subgroup,
and whose Lie algebra is the complexification of $\cg_\rT/\ct=\Lie(\rG_\rT/\rT)$. 
Thus~$\rG^c_\rT/\rT^c$ is the complexification of~$\rG_\rT/\rT$ and this proves~(i).

Now let $\zeta\in\sT^c\cap\cg_\rT^c\setminus\ct$.  
Then the set
$$
\rT_\zeta:=\overline{\{\exp(t\zeta)\,|\,t\in\R\}}
$$
is a compact subgroup of~$\rG_\rT^c$, and so projects 
to a compact subgroup of~$\rG_\rT^c/\rT^c$ 
which is generated by~$[\zeta]_{\rT^c}$.  
Since~${\Im(\Pi_\rT(\zeta))\perp\ct}$ by part~(ii) 
of Lemma~\ref{le:GT}, we have~${[\zeta]_{\rT^c}\ne 0}$ 
and so~${[\zeta]_{\rT^c}}$ is a toral generator.

Conversely, choose an element~${\zeta\in\cg_\rT^c}$ 
such that $[\zeta]_{\rT^c}$ is a toral generator.  
Then~${\zeta\notin\ct^c}$ and we may 
assume without loss of generality that~${\Im(\zeta)\perp\ct}$. 
(If necessary, replace~$\zeta$ by~${\zeta-\i\Im(\Pi_\rT(\zeta))}$ 
without changing the equivalence class~$[\zeta]_{\rT^c}$.)
Since~$[\zeta]_{\rT^c}$ is a toral generator, 
there exists an~${h\in\rG_\rT^c}$ such that~${\Im(h\zeta h^{-1})\in\ct}$.  
Since~${\Im(\zeta)\perp\ct}$, we have~${\Im(h \zeta h^{-1})\perp\ct}$
by Lemma~\ref{le:WEIGHT3} and so~${\Im(h\zeta h^{-1})=0}$.  
Thus~${\zeta\in\sT^c\cap\cg_\rT^c\setminus\ct}$ and this proves 
Lemma~\ref{le:GcT} under the assumption that~${\rT\subset\rG}$.  

To prove the result in general, choose $g\in\rG^c$ 
such that
$$
\trT:=g\rT g^{-1}\subset\rG.
$$
Then $\rG_\rT^c/\rT^c$ is isomorphic to $\rG_\trT^c/\trT^c$
and hence is reductive.  Now let~${\zeta\in\cg_\rT^c}$.  
Then the equivalence class~${[\zeta]_{\rT^c}\in\cg_\rT^c/\ct^c}$ 
is a toral generator if an only if $[g\zeta g^{-1}]_{\trT^c}$ 
is a toral generator, and 
$$
\zeta\in\sT^c\cap\cg_\rT^c\setminus\ct
\qquad\iff\qquad
g\zeta g^{-1}\in\sT^c\cap\cg_\trT^c\setminus\tct.
$$
This proves Lemma~\ref{le:GcT}. 
\end{proof}

\begin{lem}\label{le:TWEIGHT}
Let~${x\in X}$ and let $\rT\subset\rG_x^c$ be a torus
with the Lie algebra~$\ct$.\index{weights@$(\mu,T)$-weight}   
Then the following holds.

\smallskip\noindent{\bf (i)}
If~${\zeta\in\sT^c\cap\cg_\rT^c\setminus\ct}$, then
\begin{equation}\label{eq:TWEIGHT1}
w_{\mu,\rT}(x,\zeta) 
= w_\mu(x,\zeta) - \inner{\mu(x)}{\Re(\Pi_\rT(\zeta))}
\end{equation}

\smallskip\noindent{\bf (ii)}
If~${\zeta\in\sT^c\cap\cg_\rT^c\setminus\ct}$ and~${\tau\in\ct}$,
then~${\zeta+\tau\in\sT^c\cap\cg^c_\rT\setminus\ct}$ and
\begin{equation}\label{eq:TWEIGHT2}
w_{\mu,\rT}(x,\zeta+\tau)
= w_{\mu,\rT}(x,\zeta) 
\end{equation}

\smallskip\noindent{\bf (iii)}
If~${\zeta\in\sT^c\cap\cg_\rT^c\setminus\ct}$ and~${g\in\rG^c}$, then
\begin{equation}\label{eq:TWEIGHT3}
w_{\mu,\rT}(x,\zeta)
= w_{\mu,g\rT g^{-1}}(gx,g\zeta g^{-1})
\end{equation}

\smallskip\noindent{\bf (iv)}
If $\rT\subset\rG_x$ 
and~$\rG_\rT,\cg_\rT$ are as in~\eqref{eq:GT} 
and~${\zeta\in\sT^c\cap\cg^c_\rT\setminus\ct}$,
then
\begin{equation}\label{eq:TWEIGHT4}
w_{\mu,\rT}(x,\zeta) = w_{\mu_\rT}(x,[\zeta]_{\rT^c}).
\end{equation}
Here the term on the right hand side is the weight associated to the 
moment map~\eqref{eq:muT} for the action of $\rG_\rT/\rT$
on the manifold $X_\rT$ in Lemma~\ref{le:GT}.
\end{lem}

\bigbreak

\begin{proof}
We prove~(i).
Let~${\zeta\in\sT^c\cap\cg_\rT^c\setminus\ct}$.
Then the group element
$$
g:=\exp(\i t\zeta)\in\rG_\rT^c
$$
satisfies
$$
g\zeta g^{-1}=\zeta,\qquad g\rT g^{-1}=\rT.
$$
Hence it follows from equation~\eqref{eq:PITmug} 
in Lemma~\ref{le:PIT} that
$$
\inner{\mu(\exp(\i t\zeta)x)}{\Re(\Pi_\rT(\zeta))}
= \inner{\mu(x)}{\Re(\Pi_\rT(\zeta))}
$$ 
for all~${t\in\R}$.  Thus, by~\eqref{eq:wmuT}, we have
\begin{equation*}
\begin{split}
w_{\mu,\rT}(x,\zeta)
&=
\lim_{t\to\infty}\inner{\mu(\exp(\i t\zeta)x)}{\Re(\zeta-\Pi_\rT(\zeta))} \\
&=
\lim_{t\to\infty}\inner{\mu(\exp(\i t\zeta)x)}{\Re(\zeta)} 
- \inner{\mu(x)}{\Re(\Pi_\rT(\zeta))}  \\
&=
w_\mu(x,\zeta) - \inner{\mu(x)}{\Re(\Pi_\rT(\zeta))}.
\end{split}
\end{equation*}
This proves~\eqref{eq:TWEIGHT1} and~(i).

We prove~(ii).
Let~${\zeta\in\sT^c\cap\cg_\rT^c\setminus\ct}$ and~${\tau\in\ct}$.
Then 
$
{\zeta+\tau\in\sT^c\cap\cg^c_\rT\setminus\ct}
$
by part~(ii) of Lemma~\ref{le:GT}.  
Moreover, since $[\zeta,\tau]=0$, we have
$$
\exp(\i t(\zeta+\tau))x=\exp(\i t\zeta)\exp(\i t\tau)x=\exp(\i t\zeta)x
$$
for all~${t\in\R}$, hence
$$
x^+ := \lim_{t\to\infty}\exp(\i t\zeta)x
=\lim_{t\to\infty}\exp(\i t(\zeta+\tau))x,
$$ 
and so
\begin{equation*}
\begin{split}
w_{\mu,\rT}(x,\zeta+\tau)
&= 
\inner{\mu(x^+)}{\Re(\zeta+\tau-\Pi_\rT(\zeta+\tau))} \\
&= 
\inner{\mu(x^+)}{\Re(\zeta-\Pi_\rT(\zeta))} \\
&= 
w_{\mu,\rT}(x,\zeta).
\end{split}
\end{equation*}
by~\eqref{eq:wmuT}.  This proves~\eqref{eq:TWEIGHT2} and~(ii).

We prove~(iii).
Let~${\zeta\in\sT^c\cap\cg^c_\rT\setminus\ct}$ and $g\in\rG^c$. 
Then, by Lemma~\ref{le:PIT}, part~(i), and Theorem~\ref{thm:MUMFORD1}, 
we have
\begin{equation*}
\begin{split}
w_{\mu,\rT}(x,\zeta)
&= 
w_\mu(x,\zeta) 
- \inner{\mu(x)}{\Re(\Pi_\rT(\zeta))} \\
&= 
w_\mu(gx,g\zeta g^{-1}) 
- \inner{\mu(gx)}{\Re(\Pi_{g\rT g^{-1}}(g\zeta g^{-1}))} \\
&=
w_{\mu,g\rT g^{-1}}(gx,g\zeta g^{-1}).
\end{split}
\end{equation*}
This proves~\eqref{eq:TWEIGHT3} and~(iii).

We prove~(iv).
Assume $\rT\subset\rG_x$ and
let~${\zeta\in\sT^c\cap\cg^c_\rT\setminus\ct}$.
Then it follows from part~(ii) of Lemma~\ref{le:GT} that 
$
\Pi_\rT(\Im(\zeta))=0.
$
Now define
$$
x^+:=\lim_{t\to\infty}\exp(\i t\zeta)x.
$$
Then, by~\eqref{eq:wmuT}, we have
\begin{equation*}
\begin{split}
w_{\mu,\rT}(x,\zeta)
&= 
\inner{\mu(x^+)}{\Re(\zeta-\Pi_\rT(\zeta))} \\
&= 
\inner{\mu(x^+)}{\Re(\zeta)-\Pi_\rT(\Re(\zeta))} \\
&= 
\inner{\mu_\rT(x^+)}{[\Re(\zeta)]_\rT} \\
&= 
w_{\mu_\rT}(x,[\zeta]_{\rT^c}).
\end{split}
\end{equation*}
Here the second equality follows from the 
fact that $\ct\subset\cg$ and $\Pi_\rT(\Im(\zeta))=0$,
and the third equality follows from the definition of 
the inner product on $\cg_\rT/\ct$ in~\eqref{eq:innerT}
and the definition of the moment map $\mu_\rT:X_\rT\to\cg_\rT/\ct$
in~\eqref{eq:muT}. This proves~\eqref{eq:TWEIGHT4}, 
part~(iv), and Lemma~\ref{le:TWEIGHT}.
\end{proof}


\subsubsection*{Proof of the generalized Sz\'ekelyhidi criterion}

The proof of Theorem~\ref{thm:SZEK1} is based 
on the following three lemmas.

\begin{lem}\label{le:SZEK1}
Let $x_0\in X$, let $\rT\subset\rG_{x_0}$ be a torus
with the Lie algebra~$\ct$, 
and let $\rG_\rT$ and $\cg_\rT$ be as in~\eqref{eq:GT}.  
Then the following holds for every element~${g\in\rG_\rT^c}$.

\smallskip\noindent{\bf (i)}
$\rT\subset\rG_{gx_0}$ and $\ct\subset\ker L_{gx_0}$.

\smallskip\noindent{\bf (ii)}
$\mu(gx_0)\in\cg_\rT$.

\smallskip\noindent{\bf (iii)}
$\mu(gx_0)-\mu(x_0)\perp\ct$. 
\end{lem}

\begin{proof}
Parts~(i) and~(ii) follow from Lemma~\ref{le:GT}
because $\rG_\rT^c(x_0)\subset X_\rT$. 
To prove part~(iii), choose a smooth path $g:[0,1]\to\rG_\rT^c$
with $g(0)=\one$ and define the paths~${x:[0,1]\to X}$ 
and~${\eta,\xi:[0,1]\to\cg_\rT}$ by 
$$
x(t) := g(t)^{-1}x_0,\qquad 
\xi(t)+\i\eta(t) := g(t)^{-1}\dot g(t)
$$
for $0\le t\le1$.  Then 
$
\dot x=-L_x\xi-JL_x\eta
$
and hence 
$$
\frac{d}{dt}\mu(x) = - d\mu(x)L_x\xi - d\mu(x)JL_x\xi
= [\mu(x),\xi] - L_x^*L_x\eta.
$$
Thus the inner product of~${\tfrac{d}{dt}\mu(x)}$ with any element~${\tau\in\ct}$ 
vanishes by part~(i) and so the path~${t\mapsto\inner{\mu(x(t))}{\tau}}$ 
is constant.  This proves part~(iii) and Lemma~\ref{le:SZEK1}.
\end{proof}

\begin{lem}\label{le:SZEK2}
Let $x_0\in X$, let $x:\R\to X$ be the unique solution of~\eqref{eq:KN1},
and define~${x_\infty:=\lim_{t\to\infty}x(t)}$.  Then the following holds
for every $t\in\R$. 

\smallskip\noindent{\bf (i)}
$\rG_{x(t)}=\rG_{x_0}$ and $\rG_{x_0}\subset\rG_{x_\infty}$.

\smallskip\noindent{\bf (ii)}
$\ker L_{x(t)}=\ker L_{x_0}$ and $\ker L_{x_0}\subset\ker L_{x_\infty}$.

\smallskip\noindent{\bf (iii)}
$[\mu(x(t)),\xi]=0=[\mu(x_\infty),\xi]$ for all $\xi\in\ker L_{x_0}$.

\smallskip\noindent{\bf (iv)}
$\mu(x(t))-\mu(x_0)\perp\ker L_{x_0}$ 
and $\mu(x_\infty)-\mu(x_0)\perp\ker L_{x_0}$. 
\end{lem}

\begin{proof}
Let $u\in\rG$ and $s\in\R$. 
Then the unique solution $y:\R\to X$ of the differential equation 
$\dot y=-JL_y\mu(y)$ with $y(s)=ux(s)$ is~${y(t)=ux(t)}$. 
This implies $\rG_{x(s)}\subset\rG_{x(t)}$ for all $s,t\in\R$.  
Interchange~$s$ and~$t$ to obtain
$$
\rG_{x(s)}=\rG_{x(t)}
$$ 
for all~${s,t\in\R}$. This proves~(i) and~(ii). 
To prove~(iii), let~${\xi\in\ker L_{x_0}}$.  
Then
$$
L_{x(t)}\xi=0
$$ 
for all~${t\in\R}$ by part~(i), and hence
$$
[\mu(x(t)),\xi]=-d\mu(x(t))L_{x(t)}\xi=0
$$
for all~$t$ by~\eqref{eq:mu3}.  This proves~(iii).
To prove~(iv), we use~\eqref{eq:mu3} to compute
$$
\frac{d}{dt}\mu(x(t))
= 
d\mu(x(t))\dot x(t) 
= 
- d\mu(x(t))JL_{x(t)}\mu(x(t)) 
= 
- L_{x(t)}^*L_{x(t)}\mu(x(t)).
$$
Now fix an element~${\xi\in\ker L_{x_0}}$.  Then~${L_{x(t)}\xi=0}$ 
for all~$t$ by part~(ii) and this implies 
$$
\frac{d}{dt}\inner{\mu(x(t))}{\xi}
= 
-\inner{L_{x(t)}\xi}{L_{x(t)}\mu(x(t))}
= 
0.
$$ 
This proves part~(iv) and Lemma~\ref{le:SZEK2}.
\end{proof}

\begin{lem}\label{le:SZEK3}
Let~${x_0\in X}$, let $x:\R\to X$ be the unique solution 
of~\eqref{eq:KN1} and define~${x_\infty:=\lim_{t\to\infty}x(t)}$.
Let~${\rT\subset\rG_{x_0}}$ be a torus with the Lie algebra~$\ct$, 
and let~$\rG_\rT$ and~$\cg_\rT$ be as in~\eqref{eq:GT}.  
Then the following holds.

\smallskip\noindent{\bf (i)}
$\mu(x(t))\in\cg_\rT$ for all $t\in\R$ and $\mu(x_\infty)\in\cg_\rT$. 

\smallskip\noindent{\bf (ii)}
$\mu(x(t))-\mu(x_0)\perp\ct$ for all $t\in\R$ and $\mu(x_\infty)-\mu(x_0)\perp\ct$. 

\smallskip\noindent{\bf (iii)}
$x(t)\in\rG_\rT^c(x_0)$ for all $t\in\R$ and $x_\infty\in\overline{\rG_\rT^c(x_0)}$. 

\smallskip\noindent{\bf (iv)}
$\Abs{\mu(gx_0)}\ge\Abs{\mu(x_\infty)}\ge \Abs{\Pi_\rT(\mu(x_0))}$ 
for all $g\in\rG^c$ and 
\begin{equation}\label{eq:GTG}
\begin{split}
\Abs{\mu(x_\infty)-\Pi_\rT(\mu(x_0))}
&= 
\inf_{g\in\rG_\rT^c}\sqrt{\abs{\mu(gx_0)}^2-\abs{\Pi_\rT(\mu(x_0))}^2} \\
&= 
\inf_{g\in\rG^c}\sqrt{\abs{\mu(gx_0)}^2-\abs{\Pi_\rT(\mu(x_0))}^2}.
\end{split}
\end{equation}
\end{lem}

\begin{proof}
Since $\ct\subset\ker L_{x_0}$, parts~(i) and~(ii) follow directly 
from parts~(iii) and~(iv) in Lemma~\ref{le:SZEK2}.
To prove part~(iii) denote by $g:\R\to\rG^c$ the 
unique solution of~\eqref{eq:KN2}, so that
$$
g^{-1}\dot g = \i\mu(x),\qquad g(0)=\one.
$$
Since $\mu(x(t))\in\cg_\rT$ for all~$t$ by~(i),
this implies $g(t)\in\rG_\rT^c$,
hence 
$$
x(t)=g(t)^{-1}x_0\in\rG_\rT^c(x_0)
$$
for all~$t$ by Lemma~\ref{le:GRADFLOW},
and so
$$
x_\infty=\lim_{t\to\infty}x(t)\in\overline{\rG_\rT^c(x_0)}.
$$
This proves~(iii). To prove part~(iv) observe that
$$
\Pi_\rT(\mu(x_\infty))=\Pi_\rT(\mu(x_0))
$$ 
by~(ii) and hence, by Theorem~\ref{thm:MLT},
$$
\Abs{\Pi_\rT(\mu(x_0))}
= \Abs{\Pi_\rT(\mu(x_\infty))}
\le \Abs{\mu(x_\infty)}
= \inf_{g\in\rG^c}\Abs{\mu(gx_0)}.
$$
Since $x_\infty\in\overline{\rG_\rT^c(x_0)}$ by~(iii),
this implies
\begin{equation*}
\begin{split}
\sqrt{\abs{\mu(x_\infty)}^2-\abs{\Pi_\rT(\mu(x_0))}^2} 
&=
\inf_{g\in\rG^c}\sqrt{\abs{\mu(gx_0)}^2-\abs{\Pi_\rT(\mu(x_0))}^2} \\
&\le
\inf_{g\in\rG_\rT^c}\sqrt{\abs{\mu(gx_0)}^2-\abs{\Pi_\rT(\mu(x_0))}^2} \\
&\le 
\sqrt{\abs{\mu(x_\infty)}^2-\abs{\Pi_\rT(\mu(x_0))}^2}. 
\end{split}
\end{equation*}
This proves~\eqref{eq:GTG} and Lemma~\ref{le:SZEK3}.
\end{proof}

\begin{proof}[Proof of Theorem~\ref{thm:SZEK1} (i)$\implies$(ii)]
\label{proof:SZEK1}
Let $x_0\in X$ such that $\rG^c(x_0)$ contains 
a critical point of the square of the moment map
and let~${\rT\subset\rG^c_{x_0}}$ be a maximal torus.
By Lemma~\ref{le:CONJUGATE} it suffices to assume that~$x_0$ 
is $\mu$-balanced and~${\rT\subset\rG}$. 
We must prove that there exists an element~${h\in\rG_\rT^c}$ 
such that~${\mu(hx_0)\in\ct:=\Lie(\rT)}$.
Let $x:\R\to X$ be the unique solution of~\eqref{eq:KN1}
and define~${x_\infty:=\lim_{t\to\infty}x(t)}$.
Then~${L_{x_\infty}\mu(x_\infty)=0}$ by Theorem~\ref{thm:XINFTY}.  
We prove in five steps that 
\begin{equation}\label{eq:SZEK}
x_\infty\in\rG_\rT^c(x_0),\qquad \mu(x_\infty)\in\ct.
\end{equation}
The $\mu$-balanced condition is used in Step~2.

\medskip\noindent{\bf Step~1.}
{\it Let $g:\R\to\rG^c$ be the unique solution of~\eqref{eq:KN2}.
Then~${g(t)\in\rG_\rT^c}$ and we have~${x(t)=g(t)^{-1}x_0}$ for all $t\in\R$.}

\medskip\noindent
Part~(i) of Lemma~\ref{le:SZEK3} asserts that~${\mu(x(t))\in\cg_\rT}$ 
for all~$t$.  Hence it follows from~\eqref{eq:KN2} 
that
$
g(t)^{-1}\dot g(t)=\i\mu(x(t))\in\cg_\rT^c
$ 
and this implies~${g(t)\in\rG_\rT^c}$ for all~$t$.
The formula~${x(t)=g(t)^{-1}x_0}$ follows from 
Lemma~\ref{le:GRADFLOW} and this proves Step~1.

\medskip\noindent{\bf Step~2.}
{\it $\rT$ is a maximal torus in $\rG^c_{x_\infty}$ 
and there exists an element $g_\infty\in\rG^c$ 
such that~${x_\infty=g_\infty^{-1}x_0}$ 
and $\rT=g_\infty^{-1}\rT g_\infty$.}

\medskip\noindent
By assumption, the complexified group orbit $\rG^c(x_0)$ 
contains a critical point of the square of the moment map. 
Hence~${x_\infty\in\rG^c(x_0)}$ by Theorem~\ref{thm:CRIT}.
Choose an element $g_0\in\rG^c$ such that~${g_0x_0=x_\infty}$.
Since~$\rT$ is a maximal torus in~$\rG_{x_0}^c$ it follows 
that~${g_0\rT g_0^{-1}}$ is a maximal torus in
$$
g_0\rG_{x_0}^cg_0^{-1}=\rG_{g_0x_0}^c=\rG_{x_\infty}^c.
$$
Moreover, since $x_0$ is $\mu$-balanced and~$\rT\subset\rG$,
it follows from part~(i) of Lemma~\ref{le:SZEK2} 
that~${\rT\subset\rG_{x_0}\subset\rG_{x_\infty}\subset\rG_{x_\infty}^c}$.  
Since $\rT$ has the same dimension as~${g_0\rT g_0^{-1}}$, 
this shows that $\rT$ is another maximal torus in $\rG^c_{x_\infty}$.  
Now it follows from the Cartan--Iwasawa--Malcev Theorem 
in~\cite[Thm~14.1.3]{HN} that any two maximal tori in any 
connected Lie group are conjugate. Apply this result 
to the identity component of $\rG^c_{x_\infty}$ to obtain 
an element~${g_1\in\rG^c_{x_\infty}}$ such that
$
g_1g_0\rT g_0^{-1}g_1^{-1} = \rT.
$
Then we also have~${g_1g_0x_0 = g_1x_\infty = x_\infty}$ 
and so the element~${g_\infty:=(g_1g_0)^{-1}\in\rG^c}$
satisfies the requirements of Step~2. 

\medskip\noindent{\bf Step~3.}
{\it $\mu(x_\infty)\in\ct$.}

\medskip\noindent
Note that 
$$
L_{x_\infty}\mu(x_\infty)=0,\qquad
[\mu(x_\infty),\ct]=0
$$ 
by part~(iii) of Lemma~\ref{le:SZEK2}.
Hence 
$$
\ct' := \ct + \R\mu(x_\infty) \subset \ker L_{x_\infty}
$$
is an abelian Lie subalgebra and therefore the set
$$
\rT' := \overline{\exp(\ct')} \subset\rG_{x_\infty}
$$
is a torus.   It contains~$\rT$ and so must be equal to~$\rT$ by Step~2.  
Hence 
$$
\mu(x_\infty)\in\ct'\subset\Lie(\rT')=\ct
$$ 
and this proves Step~3.

\bigbreak

\medskip\noindent{\bf Step~4.}
{\it Let $g:\R\to\rG^c$ be as in Step~1 and let $g_\infty\in\rG^c$ 
be as in Step~2.  Define 
$$
\tau_0:=\Pi_\rT(\mu(x_0))
$$
and choose $u:\R\to\rG$ and $\eta:\R\to\cg$ such that
\begin{equation}\label{eq:gainftyeta}
g(t)\exp(\i\eta(t))u(t) = g_\infty\exp(\i t\tau_0)
\qquad\mbox{for all }t\in\R.
\end{equation}
Then the function $\R\to\R:t\mapsto\Abs{\eta(t)}$ is nonincreasing.}

\medskip\noindent
Recall that $\pi:\rG^c\to M:=\rG^c/\rG$ denotes the canonical projection 
and define the curves~${\gamma:\R\to M}$ and~${\gamma_\infty:\R\to M}$ by
$$
\gamma(t) := \pi\bigl(g(t)\bigr),\qquad
\gamma_\infty(t) := \pi\bigl(g_\infty\exp(\i t\tau_0)\bigr)
$$
for $t\in\R$.
We prove that these are both negative gradient flow lines 
of the Kempf--Ness function $\Phi_{x_0}:M\to\R$. 
For $\gamma$ this follows directly from the definition 
and part~(vi) of Theorem~\ref{thm:KNF}.  
For $\gamma_\infty$ we use the fact that~${\mu(x_\infty)\in\ct}$ 
by Step~3, thus~${\mu(x_\infty)=\Pi_\rT(\mu(x_0))=\tau_0}$
by part~(ii) of Lemma~\ref{le:SZEK3}, and~so
$$
\mu(x_\infty) 
= \tau_0
\in \ct
\subset \ker L_{x_0}
\subset \ker L_{x_\infty}
$$
by part~(ii) of Lemma~\ref{le:SZEK2}.
This implies
\begin{equation*}
\begin{split}
\bigl(g_\infty\exp(\i t\tau_0)\bigr)^{-1}
\frac{d}{dt}g_\infty\exp(\i t\tau_0)  
&= 
\i\tau_0 \\
&=
\i\mu(x_\infty) \\
&=
\i\mu(\exp(-\i t\tau_0)x_\infty) \\
&=
\i\mu\bigl((g_\infty\exp(\i t\tau_0))^{-1}x_0\bigr).
\end{split}
\end{equation*}
Thus $\gamma$ and $\gamma_\infty$
are negative gradient flow lines of the Kempf--Ness 
function~$\Phi_{x_0}$ as claimed.

By equation~\eqref{eq:gainftyeta} and 
part~(ii) of Theorem~\ref{thm:GcG},
the curve
$$
[0,1]\to M:s\mapsto\pi\bigl(g(t)\exp(\i s\eta(t))\bigr)
$$ 
is a geodesic joining~$\gamma(t)$ to~$\gamma_\infty(t)$.
Since $M$ has nonpositive sectional curvature by
part~(iv) of Theorem~\ref{thm:GcG}, geodesics are unique,
and hence
$$
\abs{\eta(t)} = d(\gamma(t),\gamma_\infty(t))
\qquad\mbox{for all }t\in\R.
$$
Since the Kempf--Ness function is convex along geodesics 
by part~(i) of Theorem~\ref{thm:KNF}, it follows from 
Lemma~\ref{le:gradflow} that the function~${t\mapsto\Abs{\eta(t)}}$ 
is nonincreasing.  This proves Step~4.

\bigbreak

\medskip\noindent{\bf Step~5.}
{\it $x_\infty\in\rG_\rT^c(x_0)$.}

\medskip\noindent
Let $g:\R\to\rG^c_\rT$ be as in Step~1, 
let $g_\infty$ be as in Step~2, 
and let 
$$
\tau_0:=\Pi_\rT(\mu(x_0))\in\ct
$$
and $u:\R\to\rG$ and $\eta:\R\to\cg$ be as in Step~4.  Then 
\begin{equation}\label{eq:SZEK1}
\tau_1 := g_\infty \tau_0g_\infty^{-1} \in\ct.
\end{equation}
By~\eqref{eq:gainftyeta} and~\eqref{eq:SZEK1}, we have
\begin{equation}\label{eq:SZEK2}
g(t)\exp(\i\eta(t))u(t) = g_\infty\exp(\i t\tau_0) = \exp(\i t\tau_1)g_\infty
\end{equation}
for all $t\in\R$. Since the function $t\mapsto\Abs{\eta(t)}$ 
is nonincreasing by Step~4, there exists a sequence $t_i\to\infty$
such that the limits 
\begin{equation}\label{eq:SZEK3}
\eta_\infty := \lim_{i\to\infty}\eta(t_i),\qquad 
u_\infty := \lim_{i\to\infty}u(t_i)
\end{equation}
exist.  Define 
\begin{equation}\label{eq:SZEK4}
a_\infty := \exp(\i\eta_\infty)u_\infty.
\end{equation}
Since $\tau_1\in\ct\subset\ker L_{x_0}$, it follows 
from~\eqref{eq:SZEK2}, \eqref{eq:SZEK3}, and~\eqref{eq:SZEK4} that
\begin{equation*}
\begin{split}
x_\infty 
&= 
\lim_{i\to\infty}g(t_i)^{-1}x_0 \\
&= 
\lim_{i\to\infty}g(t_i)^{-1}\exp(\i t_i\tau_1)x_0 \\
&= 
\lim_{i\to\infty}\exp(\i\eta(t_i))u(t_i)g_\infty^{-1}x_0 \\
&= 
\exp(\i\eta_\infty)u_\infty x_\infty \\
&= 
a_\infty x_\infty.
\end{split}
\end{equation*}
Thus $a_\infty\in\rG_{x_\infty}^c$ and hence
$$
x_\infty = a_\infty x_\infty = a_\infty g_\infty^{-1}x_0.
$$
Moreover,
$$
\exp(\i\eta(t))u(t)g_\infty^{-1}=g(t)^{-1}\exp(\i t\tau_1)\in\rG_\rT^c
$$
for all $t\in\R$ by~\eqref{eq:SZEK2} and Step~1. 
Hence, by~\eqref{eq:SZEK3} and~\eqref{eq:SZEK4}, we have
$$
a_\infty g_\infty^{-1} 
= \exp(\i\eta_\infty)u_\infty g_\infty^{-1}
= \lim_{i\to\infty}\exp(\i\eta(t_i))u(t_i)g_\infty^{-1}
\in\rG_\rT^c.
$$
This proves Step~5 and Theorem~\ref{thm:SZEK1}.
\end{proof}

\subsubsection*{Proof of the generalized Sz\'ekelyhidi moment-weight inequality}

The next goal is to establish the generalized Sz\'ekelyhidi 
moment-weight inequality in Theorem~\ref{thm:SZEK2}.  
In the rational case it is due to Sz\'ekelyhidi~\cite[Theorem~1.3.6]{S}.  
We derive it as a corollary of the standard moment-weight inequality
in Theorem~\ref{thm:MW2} for the $\rG_\rT/\rT$-action on~$X_\rT$
in Lemma~\ref{le:GT}.   

\begin{proof}[Proof of Theorem~\ref{thm:SZEK2}]
\label{proof:SZEK2}
Let~${x\in X}$, let $\rT\subset\rG_x$ be a torus with the Lie algebra~$\ct$, 
and let~$\rG_\rT$ and~$\cg_\rT$ be as in~\eqref{eq:GT}.  
Then, by~\eqref{eq:innerT}, we have
\begin{equation}\label{eq:innerquot}
\Abs{\left[\xi\right]_\rT}_{\cg_\rT/\ct}^2 
= \Abs{\xi}^2 - \Abs{\Pi_\rT(\xi)}^2 = \Abs{\xi-\Pi_\rT(\xi)}^2
\end{equation}
for all $\xi\in\cg_\rT$. 
Moreover, by part~(ii) of Lemma~\ref{le:GT}, we have~${\Pi_\rT(\Im(\zeta))=0}$.  
For~${\zeta=\xi+\i\eta\in\sT^c\cap\cg_\rT^c\setminus\ct}$ 
this implies
$$
\Abs{\left[\eta\right]_\rT}_{\cg_\rT/\ct} =\Abs{\eta}
$$ 
and hence
\begin{equation}\label{eq:innerquotc}
\begin{split}
\Abs{[\zeta]_{\rT^c}}_c^2 
&= 
\Abs{[\xi]_\rT}_{\cg_\rT/\ct}^2
- \Abs{[\eta]_\rT}_{\cg_\rT/\ct}^2 \\
&= 
\Abs{\xi}^2 - \Abs{\eta}^2 - \Abs{\Pi_\rT(\xi)}^2 \\
&= 
\Abs{\zeta}_c^2 - \Abs{\Pi_\rT(\zeta)}^2 \\
&= 
\Abs{\zeta-\Pi_\rT(\zeta)}_c^2.
\end{split}
\end{equation}
Combine~\eqref{eq:innerquot} and~\eqref{eq:innerquotc}
with Lemma~\ref{le:TWEIGHT} and Theorem~\ref{thm:M} 
to obtain
\begin{equation}\label{eq:TMW0}
\begin{split}
\sup_{\zeta\in\sT^c\cap\cg_\rT^c\setminus\ct}
\frac{-w_{\mu,\rT}(x,\zeta)}{\sqrt{\abs{\zeta}_c^2 - \abs{\Pi_\rT(\zeta)}^2}}
&= 
\sup_{\zeta\in\sT^c\cap\cg_\rT^c\setminus\ct}
\frac{-w_{\mu_\rT}(x,[\zeta]_{\rT^c})}
{\Abs{[\zeta]_{\rT^c}}_c} \\
&= 
\sup_{\xi\in\cg_\rT\setminus\ct\atop\exp(\xi)\in\rT}
\frac{-w_{\mu_\rT}(x,[\xi]_\rT)}
{\Abs{[\xi]_\rT}_{\cg_\rT/\ct}} \\
&= 
\sup_{\xi\in\cg_\rT\setminus\ct\atop\exp(\xi)\in\rT}
\frac{-w_{\mu,\rT}(x,\xi)}{\Abs{\xi-\Pi_\rT(\xi)}} \\
&=
\sup_{\xi\in\cg_\rT\cap\ct^\perp\setminus\{0\}\atop\exp(\xi)\in\rT}
\frac{-w_\mu(x,\xi)}{\Abs{\xi}}.
\end{split}
\end{equation}
Here the first step follows from~\eqref{eq:innerquotc}
and part~(iv) of Lemma~\ref{le:TWEIGHT}, 
the second from Theorem~\ref{thm:M}
for the $\rG_\rT/\rT$-action on $X_\rT$ 
with the moment map~${\mu_\rT:X_\rT\to\cg_\rT/\ct}$ 
in Lemma~\ref{le:GT}, the third from~\eqref{eq:innerquot} 
and part~(iv) of Lemma~\ref{le:TWEIGHT}, and the last 
from part~(ii) of Lemma~\ref{le:TWEIGHT} by replacing $\xi$ 
with~${\xi-\Pi_\rT(\xi)}$. 

\bigbreak

Since~${\Pi_\rT(\mu(hx))=\Pi_\rT(\mu(x))}$ for all $h\in\rG_\rT^c$ 
by Lemma~\ref{le:SZEK1}, we have
\begin{equation}\label{eq:muquot}
\begin{split}
\Abs{\mu_\rT(hx)}_{\cg_\rT/\ct}
&= 
\sqrt{\abs{\mu(hx)}^2 - \abs{\Pi_\rT(\mu(hx))}^2} \\
&= 
\sqrt{\abs{\mu(hx)}^2 - \abs{\Pi_\rT(\mu(x))}^2} \\
&= 
\Abs{\mu(hx)-\Pi_\rT(\mu(x))}
\end{split}
\end{equation}
for all $h\in\rG_\rT^c$ by~\eqref{eq:muT} and~\eqref{eq:innerquot}. 
Hence
\begin{equation}\label{eq:TMW}
\begin{split}
\sup_{\zeta\in\sT^c\cap\cg_\rT^c\setminus\ct}
\frac{-w_{\mu,\rT}(x,\zeta)}{\sqrt{\abs{\zeta}_c^2 - \abs{\Pi_\rT(\zeta)}^2}}
&= 
\sup_{\zeta\in\sT^c\cap\cg_\rT^c\setminus\ct}
\frac{-w_{\mu_\rT}(x,[\zeta]_{\rT^c})}
{\Abs{[\zeta]_{\rT^c}}_c} \\
&\le
\inf_{h\in\rG_\rT^c}\Abs{\mu_\rT(hx)}_{\cg_\rT/\ct} \\
&=
\inf_{h\in\rG_\rT^c}\sqrt{\abs{\mu(hx)}^2-\abs{\Pi_\rT(\mu(x))}^2} \\
&=
\inf_{h\in\rG^c}\sqrt{\abs{\mu(hx)}^2-\abs{\Pi_\rT(\mu(x))}^2}.
\end{split}
\end{equation}
Here the first step follows from~\eqref{eq:TMW0},
the second step follows from Theorem~\ref{thm:MW2},
the third step follows from~\eqref{eq:muquot}, and the last 
step follows from part~(iv) of Lemma~\ref{le:SZEK3}.
That the supremum on the left is attained follows from Theorem~\ref{thm:MODSTAB},
and that equality holds in~\eqref{eq:TMW} if and only if the left hand
side is nonnegative follows from Corollary~\ref{cor:M} 
for the $\rG_\rT/\rT$-action on~$X_\rT$ with the moment
map $\mu_\rT:X_\rT\to\cg_\rT/\ct$ in Lemma~\ref{le:GT}.
This proves the Sz\'ekelyhidi moment-weight 
inequality~\eqref{eq:SMW0} in the case~${g=\one}$.  

To prove the result in general, replace~$\rG$ by~${\trG:=g^{-1}\rG g}$,
with the inner product~\eqref{eq:innerc} on its Lie algebra $\tcg := g^{-1}\cg g$, 
and use the moment map $\tmu:X\to\tcg$ in~\eqref{eq:tmu}
for the $\trG$-action on $(X,g^*\om)$.  Then Lemma~\ref{le:CONJUGATE} 
shows that the  estimate~\eqref{eq:SMW0} is equivalent 
to the same estimate with~$\mu$ replaced by~$\tmu$.
This proves Theorem~\ref{thm:SZEK2}.
\end{proof}

In the rational case the inequality~\eqref{eq:SMW0} in 
Theorem~\ref{thm:SZEK2} is due to Sz\'ekelyhidi.
In~\cite[Thm~1.3.6]{S} he stated the estimate 
in the following form.

\begin{cor}[{\bf Sz\'ekelyhidi}]\label{cor:SMW}
Let~${x\in X}$, let $\rT\subset\rG_x$ be a torus with the Lie algebra $\ct$, 
and let~$\rG_\rT$ and~$\cg_\rT$ be as in~\eqref{eq:GT}.  
Let $\xi\in\cg\setminus\ct$ such that~${[\xi,\ct]=0}$ 
and~${w_{\mu,\rT}(x,\xi) < 0}$.\index{moment-weight inequality!Sz\'ekelyhidi}  
Then
\begin{equation}\label{eq:SMW}
\abs{\Pi_\rT(\mu(x))}^2  
+ \frac{w_{\mu,\rT}(x,\xi)^2}{\abs{\xi}^2-\abs{\Pi_\rT(\xi)}^2}
\le
\inf_{h\in\rG^c}\abs{\mu(hx)}^2.
\end{equation}
\end{cor}

\begin{proof}
This follows by taking $\zeta=\xi$ and $g=\one$ 
in Theorem~\ref{thm:SZEK2} and squaring the 
inequality~\eqref{eq:SMW0}.
\end{proof}

\subsubsection*{Proof of the Hilbert--Mumford criterion for critical orbits}

By Theorem~\ref{thm:SZEK1} we can use the Hilbert--Mumford 
criterion for polystability in Theorem~\ref{thm:HMps} to characterize
the complexified group orbits that contain critical points of 
the square of the moment map. 

\begin{proof}[Proof of Theorem~\ref{thm:SZEK3}]
\label{proof:SZEK3}
Assume first that $x\in X$ is $\mu$-balanced, 
that $\rT\subset\rG_x$ is a maximal torus, and that $g=\one$.
Then $x\in X_\rT$. Consider the Hamiltonian 
$\rG_\rT/\rT$-action on $X_\rT$ in Lemma~\ref{le:GT}
with the moment $\mu_\rT:X_\rT\to\cg_\rT/\ct$.  

\medskip\noindent{\bf Claim~1.}
{\it Condition~(i) in Theorem~\ref{thm:SZEK3}
holds if and only if $x$ is $\mu_\rT$-polystable.}

\medskip\noindent{\bf Claim~2.}
{\it Condition~(ii) in Theorem~\ref{thm:SZEK3}
is equivalent to condition~(ii) in Theorem~\ref{thm:HMps} 
for the quadruple~${(X_\rT,\rG_\rT/\rT,\mu_\rT,x)}$.}

\medskip\noindent{\bf Claim~3.}
{\it Condition~(iii) in Theorem~\ref{thm:SZEK3}
is equivalent to condition~(v) in Theorem~\ref{thm:HMps} 
for the quadruple~${(X_\rT,\rG_\rT/\rT,\mu_\rT,x)}$.}

\medskip\noindent
Theorem~\ref{thm:SZEK1} with~${g=\one}$ asserts that~$\rG^c(x)$ 
contains a critical point of the square of the moment map 
if and only if there exists an~${h\in\rG^c_\rT}$ 
such that~${\mu(hx)\in\ct}$.  By Lemma~\ref{le:GT} 
this holds if and only if $x$ is $\mu_\rT$-poystable.
This proves Claim~1. Claim~3 follows directly from part~(iv) 
of Lemma~\ref{le:TWEIGHT} and the characterization 
of toral generators in Lemma~\ref{le:GcT}.

We prove Claim~2.  Every equivalence class in $\cg_\rT/\ct$
has a unique representative~${\xi\in\cg_\rT\cap\ct^\perp}$.
The equivalence class~${[\xi]_\rT\in\cg_\rT/\ct}$ is nonzero if and 
only if~${\xi\ne0}$, and it satisfies~${\exp([\xi]_\rT)=[\one]_\rT\in\rG_\rT/\rT}$ 
if and only if~${\exp(\xi)\in\rT}$.   Thus our lattice vector 
in~${\cg_\rT/\ct}$ can be represented by a unique 
element $\xi\in\sT^c$ that satisfies the conditions
$$
[\xi,\ct]=0,\qquad 
\exp(\xi)\in\rT,\qquad 
\xi\in\cg,\qquad 
\Pi_\rT(\xi)=0.
$$
Since~${w_{\mu_\rT}(x,[\xi]_\rT) = w_{\mu,\rT}(x,\xi) = w_\mu(x,\xi)}$
for any such $\xi$, by part~(iv) of Lemma~\ref{le:TWEIGHT}, 
this proves Claim~2. 

Under the assumption $\rT\subset\rG_x$, the assertions 
of Theorem~\ref{thm:SZEK3} follow directly from 
Claim~1, Claim~2, Claim~3, and Theorem~\ref{thm:HMps}.

Now let $\rT$ be a maximal torus in $\rG_x^c$ and let $g\in\rG^c$
such that $g\rT g^{-1}\subset\rG$. 
Consider the subgroup 
$
\trG:=g^{-1}\rG g\subset\rG^c
$
and the moment map~${\tmu:X\to\tcg}$ in~\eqref{eq:tmu}. 
Then~$\rG^c(x)$ contains a critical point of~$\frac12\abs{\mu}^2$ 
if and only if it contains a critical point of~$\frac12\abs{\tmu}^2$.
Moreover,~${w_{\tmu,\rT}(x,\zeta) = w_{\mu,\rT}(x,\zeta)}$
for all~${\zeta\in\sT^c\cap\cg_\rT^c\setminus\ct}$ 
by Lemma~\ref{le:CONJUGATE}.  
Thus each of the conditions~(i),~(ii),~(iii) in Theorem~\ref{thm:SZEK3}
for~$(X,\om,\rG,\mu,g,\rT,x)$ with $\rT\subset\rG^c_x$
is equivalent to the corresponding condition 
for~$(X,\tom,\trG,\tmu,\one,\rT,x)$ with $\rT\subset\trG_x$,
for which the equivalence has already been established.
This proves Theorem~\ref{thm:SZEK3}.
\end{proof}

\subsubsection*{The rational case}

If the inner product on $\cg$ is rational with some factor $\hbar>0$
(Definition~\ref{def:ratinner}), then for every torus $T\subset\rG$ 
with the Lie algebra $\ct=\Lie(\rT)$ there exists a unique connected 
Lie subgroup $\rH\subset\rG$ with the Lie algebra
${\ch = \cg_\rT\cap\ct^\perp}$ (Lemma~\ref{le:H}). 
In this situation the action of the quotient group~$\rG_\rT/\rT$ on~$X_\rT$ 
can be replaced by the action of the subgroup~$\rH$ on all of~$X$ 
and the Sz\'ekelyhidi criterion in Theorem~\ref{thm:SZEK1}
can be restated in terms of $\mu_\rH$-polystability
(Corollary~\ref{cor:SZEK}).

\begin{lem}[{\bf Sz\'ekelyhidi~\cite[Lemma~1.3.2]{S}}]\label{le:H}
Assume the inner product on $\cg$ is rational with factor $\hbar$.
Let~${\rT\subset\rG}$ be a torus with the Lie algebra~${\ct:=\Lie(\rT)}$,
let $\rG_\rT$ and $\cg_\rT$ be as in~\eqref{eq:GT},
and define
\begin{equation}\label{eq:H}
\begin{split}
\ch 
&:= \left\{\xi\in\cg\,\big|\,[\xi,\tau]=0
\mbox{ and }\inner{\xi}{\tau}=0
\mbox{ for all }\tau\in\ct\right\}, \\
\rH 
&:= 
\left\{u(1)\,\Bigg|\,
\begin{array}{l}
u:[0,1]\to\rG\mbox{ is a smooth path} \\
\mbox{such that }u(0)=\one\mbox{ and } \\
\dot u(t)u(t)^{-1}\in\ch \mbox{ for all }t\in[0,1]
\end{array}\right\}.
\end{split}
\end{equation}
Then $\rH$ is a Lie subgroup of $\rG$
with the Lie algebra~${\Lie(\rH)=\ch}$.  Moreover,
\begin{equation}\label{eq:GTH}
\rG_\rT := \rT\rH = \left\{av\,\big|\,a\in\rT,\,v\in\rH\right\},\qquad
\cg_\rT = \ct\oplus\ch,
\end{equation}
and the complexifications are related by $\rG_\rT^c=\rT^c\rH^c$
and $\cg_\rT^c=\ct^c\oplus\ch^c$.  
\end{lem}

\begin{proof}
The linear subspace~${\ch\subset\cg}$ in~\eqref{eq:H} 
is a Lie subalgebra and hence $\rH$ is a subgroup of $\rG$ 
(called the {\it integral subgroup of $\ch$}). 
We must prove that it is closed. By a theorem of 
Malcev~\cite{MALCEV} the subgroup~$\rH$ is closed
if and only if~${\overline{\{\exp(t\eta)\,|\,t\in\R\}}\subset\rH}$
for all~${\eta\in\ch}$ (see also Hilgert--Neeb~\cite[Cor~14.5.6]{HN}).

Fix an element $\eta\in\ch$. Then $\exp(t\eta)$ commutes 
with every element of $\ct$ and hence the linear subspace 
$\ct+\R\eta$ is an abelian subalgebra of $\cg$.  
Hence
$$
T_\eta:=\overline{\exp(\ct+\R\eta)}\subset\rG
$$ 
is a torus.  Denote
$$
\ct_\eta := \Lie(\rT_\eta) \subset \cg.
$$
Then the lattice~${\ct\cap\Lambda}$ spans~$\ct$,
and the lattice~${\ct_\eta\cap\Lambda}$ spans~$\ct_\eta$. 
Consider the orthogonal decomposition
$$
\ct_\eta = \ct\oplus\ct',\qquad
\ct_\eta' := \ct_\eta\cap\ct^\perp.
$$
We claim that the intersection~${\ct_\eta'\cap\Lambda}$ 
spans~$\ct_\eta'$.  To see this, let~${m:=\dim(\ct)}$
and~${n:=\dim(\ct_\eta)}$, choose an integral 
basis~${e_1,\dots,e_m}$ of~${\ct\cap\Lambda}$, 
and extend it to an integral basis~${e_1,\dots,e_n}$ 
of~${\ct_\eta\cap\Lambda}$.  

Use Gram--Schmidt to obtain an orthogonal 
basis~${e_1',\dots,e_n'}$ of~$\ct_\eta$ 
defined recursively by $e_1':=e_1$ and
$$
e_k':=e_k-\sum_{i=1}^{k-1}\frac{\inner{e_k}{e_i'}}{\Abs{e_i'}^2}e_i'
\qquad\mbox{for }k=2,\dots,n.
$$
It follows by induction that~$e_k'$ is a rational linear combination 
of~${e_1,\dots,e_k}$ for each~$k$ and thus satisfies 
$$
\Abs{e_k'}^2\in2\pi\hbar\Q,\qquad
\inner{e_j}{e_k'}\in2\pi\hbar\Q
$$ 
for all $j$.   Since 
$$
\mathrm{span}(e_1',\dots,e_k')=\mathrm{span}(e_1,\dots,e_k)
$$ 
for all~$k$, the vectors $e_{m+1}',\dots,e_n'$ form 
a rational basis of~${\ct_\eta'\cap\Q\Lambda}$, 
so $\ct_\eta'\cap\Lambda$ spans~$\ct_\eta'$ as claimed. 
Thus
$$
\rT_\eta':=\exp(\ct_\eta')
$$ 
is a closed subgroup of~$\rG$
such that
$$
\exp(\R\eta)\subset\rT_\eta'\subset\rH
$$ 
and so~${\overline{\exp(\R\eta)}\subset\rT_\eta'\subset\rH}$.
This shows that
$$
\overline{\exp(\R\eta)}\subset\rH
\qquad\mbox{for all }\eta\in\ch.
$$
Hence it follows from Malcev's theorem~\cite{MALCEV}  
that~$\rH$ is a closed subgroup of~$\rG$ 
and so is a Lie subgroup of~$\rG$. 

To prove~\eqref{eq:GTH}, observe that
$$
\cg_\rT=\ct\oplus\ch,\qquad
\rT\rH\subset\rG_\rT.
$$ 
To prove the converse inclusion, 
choose a smooth path~${u:[0,1]\to\rG}$ 
such that~${u(0)=\one}$ and~${\xi(t):=u(t)^{-1}\dot u(t)\in\cg_\rT}$
for~${0\le t\le1}$. Then~${\xi(t)\in\ct\oplus\ch}$ for all~$t$ 
and we write
$$
\xi(t)=\tau(t)+\eta(t),\qquad
\tau(t)\in\ct,\qquad
\eta(t)\in\ch.
$$  
Define the curves~${a:[0,1]\to\rT}$ and~${v:[0,1]\to\rH}$ by
$$
a^{-1}\dot a=\tau,\qquad 
v^{-1}\dot v=\eta,\qquad
a(0)=v(0)=\one.
$$
Then~${u(t)=a(t)v(t)}$ for all~$t$, because~$\rT$ and~$\rH$
commute.  Tus~${\rG_\rT\subset\rT\rH}$ and so~${\rG_\rT=\rT\rH}$.
The same argument shows that~${\rG_T^c=\rT^c\rH^c}$
and this proves Lemma~\ref{le:H}.
\end{proof}

The following corollary is the Sz\'ekelyhidi criterion
in its original form, as stated in~\cite[Theorem~1.3.4]{S}.

\begin{cor}[{\bf Sz\'ekelyhidi Criterion}]\label{cor:SZEK}
Assume\index{Sz\'ekelyhidi Criterion} 
that the inner product on~$\cg$ 
is rational\index{Sz\'ekelyhidi Criterion!for $\mu$-balanced points}
with factor $\hbar$ and that~${x_0\in X}$ is $\mu$-balanced.
Let~${\rT\subset\rG_{x_0}}$ be a maximal torus with the Lie algebra~$\ct$, 
let $\rH$ and $\ch$ be as in Lemma~\ref{le:H}, and let
$\Pi_\rH:\cg\to\ch$ be the orthogonal projection. 
Then~${\mu_\rH:=\Pi_\rH\circ\mu:X\to\ch}$
is a moment map for the action of $\rH$ on $X$
and the following are equivalent.

\smallskip\noindent{\bf (i)}
$\rG^c(x_0)$ contains a critical point of 
the square of the moment map.

\smallskip\noindent{\bf (ii)}
$x_0$ is $\mu_\rH$-polystable.

\smallskip\noindent{\bf (iii)}
Every $\zeta\in\sT^c\cap\ch^c$ 
satisfies $w_\mu(x_0,\zeta)\ge0$ and
\begin{equation}\label{eq:HMCRIT}
w_\mu(x_0,\zeta)=0\qquad\implies\qquad
\lim_{t\to\infty}\exp(\i t\zeta)x_0\in\rH^c(x_0).
\end{equation}

\smallskip\noindent{\bf (iv)}
Every $\zeta\in\Lambda^c\cap\ch^c$ 
satisfies $w_\mu(x_0,\zeta)\ge0$ and~\eqref{eq:HMCRIT}.

\smallskip\noindent{\bf (v)}
Every $\zeta=\xi\in\ch\setminus\{0\}$ 
satisfies $w_\mu(x_0,\xi)\ge0$ and~\eqref{eq:HMCRIT}.

\smallskip\noindent{\bf (vi)} 
Every $\zeta=\xi\in\Lambda\cap\ch$ 
satisfies $w_\mu(x_0,\xi)\ge0$ and~\eqref{eq:HMCRIT}.
\end{cor}

\begin{proof}
Assertion~(ii) is equivalent to condition~(ii) in Theorem~\ref{thm:SZEK1} 
because
$$
\rG^c_\rT=\rH^c\rT^c
$$ 
by Lemma~\ref{le:H}.  Hence the equivalence of~(i) and~(ii) 
follows from Theorem~\ref{thm:SZEK1}.
That~(ii) is also equivalent to the remaining assertions
was proved in Theorem~\ref{thm:HMps}.  
This proves Corollary~\ref{cor:SZEK}. 
\end{proof}


\chapter{Examples}\label{ch:EX}


\begin{exple}[{\bf Circle actions}]\label{ex:S1actions}\rm
Let $(X,\om)$ be a closed symplectic manifold and let~${H:X\to\R}$
be a smooth function whose Hamiltonian flow is $2\pi$-periodic. 
Denote by~${v_H\in\Vect(X)}$ the Hamiltonian vector field
of~$H$ and by~${\{\varphi^s_H\{_{s\in\R}}$ the flow of~$v_H$.  Thus
$$
\iota(v_H)\om = dH,\qquad
\frac{d}{ds}\varphi^s_H = v_H\circ\varphi_H^s,\qquad
\varphi_H^0=\id,\qquad
\varphi_H^{s+2\pi}=\varphi_H^s
$$
for all $s\in\R$. Then the map
$
S^1\times X\to X:(e^{\i s},x)\mapsto\varphi_H^s(x)
$ 
is a Hamiltonian circle action on~$X$ 
and~${\mu:=\i H:X\to\i\R=\Lie(S^1)}$ is an 
equivariant moment map for this action. 
Choose an $S^1$-invariant and $\om$-compatible
almost complex structure~$J$ on~$X$.   
Then~${\inner{\cdot}{\cdot}:=\om(\cdot,J\cdot)}$
is a Riemannian metric on~$X$, the gradient of~$H$
with respect to this Riemannian metric is given 
by~${\nabla H = Jv_H}$, and the equation $\cL_{v_H}J=0$ implies
$$
[v_H,\nabla H]=0.
$$ 
Hence the flows induced by the vector fields $v_H$ and $\nabla H$
commute and so the circle action extends to a $\C^*$-action 
$\C^*\times X\to X:(g,x)\mapsto g\cdot x$ via
\begin{equation}\label{eq:Cstaraction}
e^{\i(s+\i t)}\cdot x := \varphi_H^s\circ\varphi_{\nabla H}^t(x)
\end{equation}
for $s,t\in\R$ and $x\in X$, where 
$\R\to\Diff(X):t\mapsto\varphi^t_{\nabla H}$
denotes the (positive) gradient flow of $H$. 
Thus, for each $x_0\in X$, the curve $x:\R\to X$, 
defined by 
$$
x(t) := e^{-t}\cdot x_0 = \varphi_{\nabla H}^t(x_0)
$$
for $t\in\R$, is the positive gradient flow line of $H$
through $x_0$, i.e.\ 
\begin{equation}\label{eq:gradflowH}
\dot x(t) = \nabla H(x(t)),\qquad x(0)=x_0.
\end{equation}
Denote the limit points of this gradient flow line by 
\begin{equation}\label{eq:limitH}
x^+:=\lim_{t\to\infty}x(t),\qquad
x^-:=\lim_{t\to-\infty}x(t).
\end{equation}
In spite of the fact that $J$ need not be integrable,
and so $(X,\om,J)$ need not be a K\"ahler manifold,
all the definitions and results in this book 
carry over to the present situation. For example, 
the weights of $x_0$ are 
\begin{equation}\label{eq:weightH}
w_\mu(x_0,\i) = H(x^+),\qquad w_\mu(x_0,-\i) = - H(x^-).
\end{equation}
Thus the weights are all nonnegative if and only if 
$$
H(x^-)\le 0\le H(x^+), 
$$
and this holds if and only 
if the closure 
$$
\overline{\C^*\cdot x_0} = \C^*\cdot x_0\cup\{x^-,x^+\}
$$ 
of the complexified group orbit intersects 
the zero set $H^{-1}(0)$ of the moment map,
i.e.\ $x_0$ is $\mu$-semistable. 
If $x_0$ is $\mu$-polystable, but not $\mu$-stable,  
then the criterion of Theorem~\ref{thm:HMps} translates 
into the condition $x^-=x_0=x^+$, in which case 
the gradient flow line~\eqref{eq:gradflowH} 
is constant and so
$$
H(x^-)=0=H(x^+).
$$
The point $x_0$ is $\mu$-stable if and only if its 
complexified group orbit $\C^*\cdot x_0$ intersects $H^{-1}(0)$ 
and the gradient flow line~\eqref{eq:gradflowH} is nonconstant.
This holds if and only if 
$$
H(x^-) < 0 < H(x^+),
$$
confirming the criterion of Theorem~\ref{thm:HMs}.
In this situation the Kempf--Ness function 
$\Phi_{x_0}:(0,\infty)\cong\C^*/S^1\to\R$ is given by 
\begin{equation}\label{eq:KNcircle}
\Phi_{x_0}(r) = \int_1^r
\frac{H(\varphi^{\log(\rho)}_{\nabla H}(x_0))}{\rho}\,d\rho
\qquad\mbox{for }r>0.
\end{equation}
The metric on the space $M=\C^*/S^1\cong(0,\infty)$ is given by $r^{-1}dr$ 
and the geodesics have the form $\gamma(t)=r_0e^{ct}$.  
Convexity of $\Phi_{x_0}$ along geodesics translates into 
the equation 
$
\tfrac{d^2}{dt^2}\Phi_{x_0}(e^t) 
= \Abs{\nabla H(\varphi_{\nabla H}^t(x_0))}^2 
\ge 0.
$
However, 
$$
r^2\Phi_{x_0}''(r)=\abs{\nabla H(\varphi_{\nabla H}^{\log(r)}(x_0))}^2
- H(\varphi_{\nabla H}^{\log(r)}(x_0))
$$ 
and thus the function $\Phi_{x_0}:(0,\infty)\to\R$ need not be convex. 
The Generalized Sz\'ekelyhidi Criterion for critical orbits
in Theorem~\ref{thm:SZEK1} asserts, in the present case, 
that the complexified group orbit $\C^*\cdot x_0$ 
contains a critical point of the function $\tfrac{1}{2}H^2:X\to\R$
if and only if $x_0$ is either $\mu$-stable (the case $\rT=\{1\}$)
or is a fixed point of the circle action (the case $\rT=S^1$). 
\end{exple}


\begin{exple}[{\bf A circle action on $S^2$}]
\label{ex:S1S2}\rm
Consider the unit sphere 
\begin{equation}\label{eq:S2}
S^2:=\left\{x=(x_1,x_2,x_3)\in\R^3\,\Big|\,
\abs{x}^2=\sum_{i=1}^3x_i^2=1\right\}
\end{equation}
equipped with the standard symplectic and complex structures 
given by 
$$
\sigma_x(\xhat,\yhat):=\inner{x\times\xhat}{\yhat},\qquad
J(x)\xhat :=x\times\xhat
$$ 
for $x\in S^2$ and $\xhat,\yhat\in T_xS^2=x^\perp$.  
The standard circle action is given by rotation 
in the $(x_1,x_2)$-plane and is generated by the Hamiltonian 
function~${H_c(x):=x_3+c}$ for~${x\in S^2}$.  
The complexified group action is given by 
\begin{equation}\label{eq:S2action}
g\cdot x = \frac{1}{(1+x_3)+\abs{g}^2(1-x_3)}
\left(\begin{array}{c}
2(\Re(g)x_1-\Im(g)x_2) \\
2(\Re(g)x_2+\Im(g)x_1) \\
(1+x_3) - \abs{g}^2(1-x_3)
\end{array}\right)
\end{equation}
for $g\in\C^*$ and $x\in S^2$.  This corresponds to the standard 
$\C^*$-action on the Riemann sphere $\overline{\C}=\C\cup\{\infty\}$ 
under the stereographic projection 
\begin{equation}\label{eq:stereo}
S^2\to\overline{\C}:x\mapsto\frac{x_1+\i x_2}{1+x_3},
\end{equation}
which sends the north pole~${N:=(0,0,1)}$ to~${z=0}$ 
and sends the south pole~${S:=(0,0,-1)}$ to~${z=\infty}$.  
The Kempf--Ness function~${\Phi_x:(0,\infty)\to\R}$ 
of an element $x\in S^2$ is given by the explicit formula
\begin{equation}\label{eq:S2KN}
\Phi_x(r) = (c-1)\log(r) 
+ \log\left(\frac{\abs{z}^2+r^2}{\abs{z}^2+1}\right),\qquad
z := \frac{x_1+\i x_2}{1+x_3}.
\end{equation}
The slopes are 
\begin{equation*}
\begin{split}
w_\mu(x,\i) 
= \lim_{t\to\infty}\frac{\Phi(e^t)}{t}
&= \left\{\begin{array}{ll}
c+1,&\mbox{if }x\ne S,\\
c-1,&\mbox{if }x= S,
\end{array}\right. \\
w_\mu(x,-\i) 
= \lim_{t\to\infty}\frac{\Phi(e^{-t})}{t}
&= \left\{\begin{array}{ll}
-(c-1),&\mbox{if }x\ne N,\\
-(c+1),&\mbox{if }x=N.
\end{array}\right. 
\end{split}
\end{equation*}
This is consistent with~\eqref{eq:weightH}, 
Lemma~\ref{le:SLOPE}, and the Hilbert--Mumford criterion. 
Namely, when~${-1<c<1}$ then the elements 
of~${S^2\setminus\{S,N\}}$ 
are stable while~$S$ and~$N$ are unstable, 
and when $\abs{c}>1$ all elements of $S^2$ are unstable.
If $c=1$ then $S$ is polystable, 
the elements of $S^2\setminus\{S,N\}$
are all semistable, and $N$ is unstable. 
If $c=-1$ then $N$ is polystable, 
the elements of $S^2\setminus\{S,N\}$
are all semistable, and $S$ is unstable. 
\end{exple}


\begin{exple}[{\bf The $\SO(3)$-action on $S^2$}]
\label{ex:SO3S2}\rm
The group $\rG:=\SO(3)$ acts on the unit sphere 
$S^2\subset\R^3$ and the action preserves the 
standard K\"ahler structure in Example~\ref{ex:S1S2}. 
Throughout we identify the Lie algebra $\cso(3)$
with~$\R^3$ via the isomorphism
\begin{equation}\label{eq:R3so3}
\R^3\to\cso(3):
\xi=(\xi_1,\xi_2,\xi_3)
\mapsto A_\xi 
:= \left(\begin{array}{rrr}
0 & -\xi_3 & \xi_2 \\
\xi_3 & 0 & -\xi_1 \\
-\xi_2 & \xi_1 & 0
\end{array}\right).
\end{equation}
This isomorphism identifies the cross product with the Lie bracket,
satisfies the equation $A_{u\xi}=uA_\xi u^{-1}$ for $u\in\SO(3)$ and $\xi\in\R^3$,
and the standard inner product on $\R^3$ is invariant and coincides with 
half the trace, i.e.~${\inner{\xi}{\eta}=-\frac12\trace(A_\xi A_\eta)}$ 
for all~${\xi,\eta\in\R^3}$. The $\SO(3)$-action on $S^2$ is Hamiltonian
and, under the above identification of the Lie algebra with $\R^3$,
the moment map $\mu:S^2\to\R^3$ is the canonical inclusion,
i.e.\ $\mu(x)=x$. 

Throughout we identify $\SO(3)$ with $\SU(2)/\{\pm\one\}$ 
via the isomorphism induced by the Lie group 
homomorphism $\SU(2)\to\SO(3)$ given by 
\begin{equation}\label{eq:SU2SO3}
\left(\begin{array}{rr}
a & b \\
- \ob & \oa
\end{array}\right)
\mapsto
\left(\begin{array}{ccc}
\Re(a^2-b^2) & -\Im(a^2-b^2) & \Re(2ab) \\
\Im(a^2-b^2) & \Re(a^2+b^2) & \Im(2ab) \\
-\Re(2\overline{a}b) & -\Im(2\overline{a}b) & \abs{a}^2-\abs{b}^2
\end{array}\right).
\end{equation}
The homomorphism~\eqref{eq:SU2SO3} has been chosen such that,
under the stereographic projection $S^2\to\overline\C$ in~\eqref{eq:stereo}, 
the $\SO(3)$-action on $S^2$ corresponds to the standard
$\SU(2)$-action on the Riemann sphere $\overline{\C}$ via
\begin{equation}\label{eq:SU2Mob}
\left(\begin{array}{rr}
a & b \\
-\ob & \oa
\end{array}\right) \cdot z 
= \frac{az+b}{-\ob z+\oa}
\end{equation}
(see~\cite[Proposition~4.2.2]{G}).  Thus 
$$
\rG^c=\PSL(2,\C)
$$ 
is the group of M\"obius transformations, it acts 
on the Riemann sphere by 
\begin{equation}\label{eq:MobC}
\left(\begin{array}{rr}
a & b \\
c & d
\end{array}\right) \cdot z 
= \frac{az+b}{cz+d},
\end{equation}
and the corresponding action on $S^2$ via the stereographic 
projection~\eqref{eq:stereo} is given by 
\begin{equation}\label{eq:MobS2}
\left(\begin{array}{rr}
a & b \\
c & d
\end{array}\right) \cdot x
= \frac{1}{\abs{az+b}^2+\abs{cz+d}^2}
\left(\begin{array}{c}
2\Re\bigl((az+b)\overline{(cz+d)}\bigr) \\
2\Im\bigl((az+b)\overline{(cz+d)}\bigr) \\
\abs{az+b}^2-\abs{cz+d}^2
\end{array}\right)
\end{equation}
for $x\in S^2$ and $z:=(1+x_3)^{-1}(x_1+\i x_2)$.

\bigbreak

To obtain an explicit formula for the Kempf--Ness 
function, it is convenient to identify the symmetric 
space~${\rG^c/\rG}$ with the hyperbolic $3$-space 
\begin{equation}\label{eq:H3}
\bbH^3 = \left\{y=(y_0,y_1,y_2,y_3)\in\R^4\,\big|\,
y_0>0,\,Q(y,y)=-1\right\},
\end{equation}
where~${Q:\R^4\times\R^4\to\R}$ is 
the quadratic form
$$
Q(y,y'):=-y_0y_0'+y_1y_1'+y_2y_2'+y_3y_3'
$$
for~${y,y'\in\R^4}$ and the Riemannian metric on~$\bbH^3$
is given by the restriction of~$Q$ to the tangent spaces.
The isometry~${\cF:\rG^c/\rG\to\bbH^3}$ sends 
the equivalence class of~${g\in\SL(2,\C)}$ to the 
vector~${(y_0,y_1,y_2,y_3)\in\bbH^3}$ given by
\begin{equation}\label{eq:yg}
\begin{array}{ll}
y_0 &:= (a+d)^2/2 - 1\\
y_1 &:= (a+d)\Re(b),\\
y_2 &:= (a+d)\Im(b),\\
y_3 &:= (d^2-a^2)/2,
\end{array}
\qquad
\left(\begin{array}{cc}
a & b \\ \ob & d 
\end{array}\right)
:= (gg^*)^{1/2}.
\end{equation}
With this identification the lifted 
Kempf--Ness function~${\Phi_x:\PSL(2,\C)\to\R}$
associated to an element~${x\in S^2}$ is given by 
\begin{equation}\label{eq:KNH3}
\Phi_x(g) = \log\left(y_0- x_1y_1 - x_2y_2 - x_3y_3\right)
\end{equation}
for~${g\in\PSL(2,\C)}$, where~${y\in\bbH^3}$
is given by~\eqref{eq:yg} (see~\cite[Thms~4.2.5 and~4.2.6]{G}).
In this example the $\mu$-weight of a 
pair~${(x,\xi)\in S^2\times\R^3}$ is 
\begin{equation}\label{eq:weightS2}
w_\mu(x,\xi) = \left\{\begin{array}{rl}
\abs{\xi},&\mbox{if }\xi\ne-\abs{\xi}x,\\
-\abs{\xi},&\mbox{if }\xi=-\abs{\xi}x.
\end{array}\right.
\end{equation}
The vector~${\xi\in\R^3\setminus\{0\}}$ 
corresponds to the skew-Hermitian matrix
\begin{equation}\label{eq:uxihat}
\uhat_\xi
:= \frac12 \left(\begin{array}{cc}
\i\xi_3 & \xi_2-\i\xi_1 \\
-\xi_2-\i\xi_1 & -\i\xi_3
\end{array}\right)
\in \csu(2)
\end{equation}
via the deriviative of of the 
homomorphism~\eqref{eq:SU2SO3}
and so the geodesic ray~${\gamma_\xi}$ 
in~$\bbH^3$ is given 
by~${\gamma_\xi(t) = \cF(\exp(\i t\uhat_\xi))
= (\cosh(t\abs{\xi}),\sinh(t\abs{\xi})\frac{\xi}{\abs{\xi}})}$
for~${t\ge 0}$ (see~\cite[Lemma~4.2.8]{G}).
The Kempf--Ness function along this ray is
$$
\Phi_x(\gamma_\xi)(t) 
= \log\left(
\cosh(t\abs{\xi})
- \sinh(t\abs{\xi})\frac{\inner{\xi}{x}}{\abs{\xi}}
\right)
$$
and the asymptotic slope of this function is the 
weight~${w_\mu(x,\xi)}$ in~\eqref{eq:weightS2}.
The explicit computation shows that the weight
function~${\xi\mapsto w_\mu(x,\xi)}$ is discontinuous.
Geometrically, the asymptotic slope of the the 
Kempf--Ness function~\eqref{eq:KNH3} is $-1$
along precisely one geodesic ray of unit speed,
and is~$+1$ along all others. This corresponds 
to the fact that all elements~${x\in S^2}$ 
are~$\mu$-unstable because~${\abs{\mu}\equiv 1}$.
\end{exple}


\begin{exple}[{\bf The $\SO(3)$-action on $(S^2)^n$}]
\label{ex:SO3S2n}\rm
Fix a positive integer~$n$ and real numbers~${\lambda_i>0}$
for~${i=1,\dots,n}$. Let~${\sigma\in\Om^2(S^2)}$ 
be the standard symplectic form in Example~\ref{ex:S1S2}
and consider the symplectic form
$$
\om_\lambda := \lambda_1\sigma\oplus
\lambda_2\sigma\oplus\cdots\oplus\lambda_n\sigma
\in\Om^2((S^2)^n)
$$
on~${X:=(S^2)^n}$.  This is a K\"ahler form 
for the standard complex structure~$J$.  
The diagonal action of~${\rG=\SO(3)}$ on~$(S^2)^n$ 
preserves the K\"ahler structure~${(\om_\lambda,J)}$ 
and is generated by the moment 
map~${\mu_\lambda:(S^2)^n\to\R^3\cong\cso(3)}$,
which assigns to each~$n$-tuple~${x=(x_1,\dots,x_n)\in(S^2)^n}$ 
its center of mass
$$
\mu_\lambda(x_1,\dots,x_n) = \lambda_1x_1+\lambda_2x_2+\cdots+\lambda_nx_n
$$
(Example~\ref{ex:SO3S2}). The gradient flow of the 
momen map squared is given by 
\begin{equation}\label{eq:gradflowmusquared}
\dot x_i = \Inner{\sum_{j=1}^n\lambda_jx_j}{x_i}x_i 
- \sum_{j=1}^n\lambda_jx_j,\qquad i=1,\dots,n.
\end{equation}
An element~${x=(x_1,\dots,x_n)\in(S^2)^n}$ is a critical point 
of the mument map squared (i.e.\ satisfies~${L_x\mu_\lambda(x)=0}$)
if and only if either~${\mu_\lambda(x)=0}$ 
or the vectors~${x_1,\dots,x_n}$ are parallel. 
By the computations in Example~\ref{ex:SO3S2} the Kempf--Ness function 
of an element~${x=(x_1,\dots,x_n)\in(S^2)^n}$ is given by 
\begin{equation}\label{eq:KNH3n}
\Phi_x(g) = \sum_{i=1}^n\lambda_i
\log\left(y_0 - x_{i1}y_1-x_{i2}y_2-x_{i3}y_3\right)
\end{equation}
for~${g\in\PSL(2,\C)}$, where~${y\in\bbH^3}$ is given 
by~\eqref{eq:yg}.  The $\mu_\lambda$-weight of
a pair~$(x,\xi)$ with~${x=(x_1,\dots,x_n)\in(S^2)^n}$
and~${\xi\in\R^3\setminus\{0\}}$ is given by
\begin{equation}\label{eq:weightS2n}
\frac{w_{\mu_\lambda}(x,\xi)}{\abs{\xi}}
= \sum_{i=1}^n\lambda_i\frac{w_\mu(x_i,\xi)}{\abs{\xi}}
= \sum_{x_i\ne-\xi/\abs{\xi}}\lambda_i
- \sum_{x_i=-\xi/\abs{\xi}}\lambda_i.
\end{equation}
Fix an element~${x=(x_1,\dots,x_n)\in(S^2)^n}$.
Then the weighted sum of the Dirac measures 
at the points~${x_i}$ defines a 
function~${\nu_x:S^2\to[0,\nu]}$ given by
\begin{equation}\label{eq:nux}
\nu := \sum_{i=1}^n\lambda_i,\qquad
\nu_x(p):=\sum_{x_i=p}\lambda_i
\qquad\mbox{for }p\in S^2.
\end{equation}
Define
\begin{equation}\label{eq:kpm}
\kappa^+(x) := \max_{i=1,\dots,n}\nu_x(x_i),\qquad
\kappa^-(x) := \min_{i=1,\dots,n}\nu_x(x_i).
\end{equation}
With these formulas at hand the stability properties 
of an element~${x\in(S^2)^n}$ are 
determined by the Hilbert--Mumford criterion as follows. 

\bigbreak

When~${\xi\in\R^3}$ has norm~${\abs{\xi}=1}$, it follows 
from~\eqref{eq:weightS2n} and~\eqref{eq:nux} 
that
\begin{equation}\label{eq:weightS2nkappa}
w_{\mu_\lambda}(x,\xi)=\nu-2\nu_x(-\xi).
\end{equation}
By Theorem~\ref{thm:HMus} an element~${x\in(S^2)^n}$
is $\mu_\lambda$-unstable if and only if it admits
a negative weight, and according to~\eqref{eq:weightS2nkappa}
this means that more than half the mass is concentrated
at a single element of~$S^2$.  Thus, by~\eqref{eq:kpm}, 
\begin{equation}\label{eq:S2n-unstable}
x\mbox{ is $\mu_\lambda$-unstable}
\qquad\iff\qquad
\kappa^+(x) > \frac{\nu}{2}.
\end{equation}
Among the $\mu_\lambda$-unstable points are those where
the mass is concentrated at precisely two elements 
of~$S^2$ with greater mass at one of these points.
Those are the elements whose~$\rG^c$-orbits contain 
higher critical points of the moment map squared.
It follows directly from~\eqref{eq:S2n-unstable} that
\begin{equation}\label{eq:S2n-semistable}
x\mbox{ is $\mu_\lambda$-semistable}
\qquad\iff\qquad
\kappa^+(x) \le \frac{\nu}{2}.
\end{equation}
If~${\kappa^+(x)=\nu/2}$ then precisely half the 
mass is located at one element of~$S^2$.
The other half of the mass is also located at a single point 
if and only if the element is $\mu_\lambda$-polystable 
but not $\mu_\lambda$-stable, because any pair of points 
is equivalent to an antipodal pair via a 
M\"obius transformation. Thus
\begin{equation}\label{eq:S2n-polystable}
x\mbox{ is $\mu_\lambda$-polystable and not $\mu_\lambda$-stable}
\qquad\iff\qquad
\kappa^\pm(x) = \frac{\nu}{2}.
\end{equation}
Using the fact that for~$\xi,x_i\in S^2$ we 
have~${\lim_{t\to\infty}\exp(\i t\xi)x_i=\xi}$ 
whenever~${x_i\ne-\xi}$ 
and~${\lim_{t\to\infty}\exp(\i t\xi)x_i=-\xi}$ 
whenever~${x_i=-\xi}$, one can verify that 
this is consistent with the Hilbert--Mumford 
criterion in Theorem~\ref{thm:HMps}.
Finally, it follows from Theorem~\ref{thm:HMs} 
and~\eqref{eq:weightS2nkappa} that
\begin{equation}\label{eq:S2n-stable}
x\mbox{ is $\mu_\lambda$-stable}
\qquad\iff\qquad
\kappa^+(x) < \frac{\nu}{2}.
\end{equation}
In this case at least three points on~$S^2$ have 
positive mass and so the isotropy subgroup is trivial. 
The Hilbert--Mumford criterion asserts in this case
that there exists a M\"obius transformation~${g\in\PSL(2,\C)}$
such that the weighted center of mass of the 
points~${gx_1,\dots,gx_n}$ is zero.  

This example does not satisfy the 
rationality conditions of Chapter~\ref{ch:RAT} 
whenever the~$\lambda_i$ are rationally independent.
The three notions of~$\mu_\lambda$-semi\-stability, 
$\mu_\lambda$-polystability, and~$\mu_\lambda$-stability are 
equivalent whenever~$n$ is odd and~${\lambda_i=1}$,
and then the quotient space
$
{\cM_n := X^\s/\rG^c = \mu^{-1}(0)/\rG = X\dslash\rG}
$
is a smooth manifold,\index{Mumford quotient} 
called the {\bf Mumford quotient}.
It can be viewed as compactification of the configuration 
space~${\cC_n := ((S^2)^n\setminus\Delta_n)/\PSL(2,\C)}$,
where~${\Delta_n\subset(S^2)^n}$ iss the fat diagonal. 
A finer compactification is the Deligne--Mumford 
space~$\ocM_{0,n}$, and the projection~${\pi:\ocM_{0,n}\to\cM_n}$ 
sends a stable curve of genus zero with~$n$ marked points
to its {\it Mumford component}.
\end{exple}

\bigbreak

\begin{exple}[{\bf Normal matrices}]\label{ex:LA}\rm
This\index{normal matrices|(} 
example is taken from a lecture 
by Peter Kronheimer~\cite{KRONHEIMER}.
Consider the vector space~${V:=\C^{n\times n}}$ with the 
Hermitian inner product
$$
\inner{A}{B}:=\Re(\trace(A^*B))
$$ 
for~${A,B\in\C^{n\times n}}$ and the projective space~${X:=\bbP(V)}$
with the Fubini--Study form
\begin{equation}\label{eq:LAom}
\om_A(\Ahat_1,\Ahat_2) := \frac{\inner{\i\Ahat_1}{\Ahat_2}}{\abs{A}^2}
- \frac{\inner{\i\Ahat_1}{A}\inner{\Ahat_2}{A}}{\abs{A}^4} 
+ \frac{\inner{\Ahat_1}{A}\inner{\i\Ahat_2}{A}}{\abs{A}^4}  
\end{equation}
for~${A\in V\setminus\{0\}}$ and~${\Ahat_i\in V}$.  
The special unitary group~${\rG:=\SU(n)}$ 
acts on~$V$ by
$$
u\cdot A := uAu^{-1}
$$
for~${u\in\SU(n)}$ and~${A\in V}$. 
The induced action on the projective space~$\bbP(V)$ 
is Hamiltonian and generated by the moment 
map~${\mu:\bbP(V)\to\csu(n)}$ whose lift to~${V\setminus\{0\}}$ 
is given by 
\begin{equation}\label{eq:LAmu}
\mu(A) = -\tfrac{\i}{2}\frac{[A,A^*]}{\abs{A}^2}
\end{equation}
for~${A\in V\setminus\{0\}}$.  Here the Lie algebra~${\cg:=\csu(n)}$
is equipped with the standard inner product
$$
\inner{\xi}{\eta}:=\Re(\trace(\xi^*\eta))
$$
for~${\xi,\eta\in\csu(n)}$ and the infinitesimal covariant action of~$\cg$
on~$\bbP(V)$ is given by
\begin{equation}\label{eq:LAinfact}
L_A\xi:=[\xi,A] - \frac{\trace([\xi,A]A^*)}{\Abs{A}^2}A
\end{equation}
for~${\xi\in\csu(n)}$ and~${A\in V\setminus\{0\}}$.
The right hand side in~\eqref{eq:LAinfact} is the projection 
of the commutator~${[\xi,A]}$ onto the complex orthogonal 
complement of~$A$. 
Equation~\eqref{eq:LAinfact} implies~\eqref{eq:LAmu}.  
The complexified group action is given by
$$
g\cdot A=gAg^{-1}
$$
for~${g\in\SL(n,\C)}$ and~$A\in\C^{n\times n}$. 
A matrix~${A\in\C^{n\times n}\setminus\{0\}}$ belongs to the zero
set of the moment map if and only if it is normal.
It is $\mu$-unstable if and only if it has a trivial spectrum~${\sigma(A)=\{0\}}$,
is $\mu$-semistable if and only if~${\sigma(A)\ne\{0\}}$, 
and is~$\mu$-polystable if and only if it is diagonalizable
(and hence has at least one nonzero eigenvalue). 
In the case~${n\ge 2}$ there are no $\mu$-stable points 
because the isotropy subgroup is always nontrivial.

The stability conditions in this example are characterized 
by the Mumford $\mu$-weights as in Lemma~\ref{le:Vweight}. 
If~${\xi\in\csu(n)}$ and
$
\lambda_1<\lambda_2<\cdots<\lambda_k
$
are the eigenvalues of the matrix~$\i\xi$ 
with the associated eigenspace decomposition
$
{\C^n=E_1\oplus E_2\oplus\cdots\oplus E_k},
$
then~$A$ is a block matrix with~${A_{ij}:E_j\to E_i}$ and
\begin{equation*}
\i[\xi,A] = \left(\begin{array}{cccc}
0 & (\lambda_1-\lambda_2)A_{12} & \cdots &(\lambda_1-\lambda_k)A_{1k} \\
(\lambda_2-\lambda_1)A_{21} & 0 &  & (\lambda_2-\lambda_k)A_{2k} \\
\vdots &  & \ddots & \vdots \\
(\lambda_k-\lambda_1)A_{k1} & (\lambda_k-\lambda_2)A_{k2} & \cdots  & 0
\end{array}\right).
\end{equation*}
Thus it follows from part~(i) of Lemma~\ref{le:Vweight} 
(with~${\hbar=\tfrac{1}{2}}$) that
\begin{equation*}
w_\mu(A,\xi) = \tfrac{1}{2}\max_{A_{ij}\ne0}(\lambda_i-\lambda_j).
\end{equation*}
If this number is negative for some~$\xi$ 
then the block matrix~${(A_{ij})_{i,j=1,\dots,k}}$ 
is strictly upper triangular and so~${\sigma(A)=\{0\}}$.  
The converse follows by considering the Jordan normal form
and using the invariance of the weights under the action 
of the complexified group (Theorem~\ref{thm:MUMFORD1}).
The moment-weight (in)equality in this example takes the form
\begin{equation}\label{eq:LAmw}
\sup_{\xi\in\csu(n)\setminus\{0\}}
\frac{\min_{A_{ij}\ne0}(\lambda_j-\lambda_i)}{\sqrt{\sum_i\lambda_i^2\dim^c(E_i)}}
= \inf_{g\in\SL(n,\C)}\frac{\abs{[g^{-1}Ag,g^*A^*{g^*}^{-1}]}}{\abs{g^{-1}Ag}^2}, 
\end{equation}
where the~$\lambda_i$ and~$E_i$ are determined 
by~$\xi$ as above. By Corollary~\ref{cor:M} equality 
holds in~\eqref{eq:LAmw} because there are no stable points.

The square of the moment map is the function~${f:\bbP(V)\to\R}$ given by
\begin{equation}\label{eq:LAmusquared}
f(A) = \tfrac{1}{2}\abs{\mu(A)}^2 = \frac{\abs{[A,A^*]}^2}{8\abs{A}^4}
\end{equation}
for~${A\in V\setminus\{0\}}$ and its negative 
gradient flow lines are the solutions of the 
differential equation
\begin{equation}\label{eq:LAgradflow}
\dot A 
=  
- \frac{[[A,A^*],A]}{2\abs{A}^2}
+ \frac{\abs{[A,A^*]}^2}{2\abs{A}^4}A.
\end{equation}
Each solution of~\eqref{eq:LAgradflow}
satisfies~$\trace(A^*\dot A)=0$. 
By Lemma~\ref{le:GRADFLOW} it can be written  
in the form~${A(t)=g(t)^{-1}Ag(t)}$, where~${A:=A(0)}$ 
and the function~${g:\R\to\SL(n,\C)}$ satisfies 
the differential equation
\begin{equation}\label{eq:LAKNflow}
g^{-1}\dot g = \frac{[g^{-1}Ag,g^*A^*(g^*)^{-1}]}{2\abs{g^{-1}Ag}^2},\qquad 
g(0)=\one.
\end{equation}
The curve~${g:\R\to\SL(n,\C)}$ in~\eqref{eq:LAKNflow} 
is a negative gradient flow line of the lifted Kempf--Ness 
function~${\Phi_A:\SL(n,\C)\to\R}$ given by 
\begin{equation}\label{eq:LAKN}
\Phi_A(g) = \tfrac{1}{2}
\log\left(\frac{\abs{g^{-1}Ag}}{\abs{A}}\right)
\end{equation}
for~${g\in\SL(n,\C)}$ (see Lemma~\ref{le:VKN}). 
The quotient~${\rG^c/\rG=\SL(n,\C)/\SU(n)}$
is isometrically isomorphic to the space~$\cP_n$
of positive definite Hermitian matrices with determinant one,
equipped with the Riemannian metric
\begin{equation}\label{eq:Pmetric}
\inner{\Phat_1}{\Phat_2}_p 
:= \tfrac{1}{4}\trace\bigl(P^{-1}\Phat_1P^{-1}\Phat_2\bigr)
\end{equation}
for~$P\in\cP_n$ and~${\Phat_i\in T_p\cP_n}$ 
(so~$\Phat_i$ is Hermitian and~${\trace(P^{-1}\Phat_i)=0}$).
The isometry is given by
$$
\SL(n,\C)/\SU(n)\to\cP_n:g\mapsto P:=gg^*.
$$
Under this isometry the Kempf--Ness function~\eqref{eq:LAKN} 
is given by
\begin{equation}\label{eq:LAKNP}
\Phi_A(P) = \tfrac{1}{4}\log\left(
\frac{\trace(APA^*P^{-1})}{\trace(A^*A)}
\right)
\end{equation}
for~${P\in\cP_n}$ and the differential 
equation~\eqref{eq:LAKNflow} translates into
\begin{equation}\label{eq:LAKNPflow}
\dot P = \frac{APA^* - PA^*P^{-1}AP}{\trace(PA^*P^{-1}A)},\qquad
P(0) = \one.
\end{equation}
If~$A$ is diagonalizable and has nonzero spectrum 
(the polystable case) then the solution of~\eqref{eq:LAKNPflow}
converges to a matrix~${P\in\cP_n}$ as~$t$ tends to 
infinity by Theorem~\ref{thm:KNF}.  
This limit defines a Hermitian inner product
\begin{equation}\label{eq:Pinner}
\inner{z}{z'}_P:=\Re(z^*P^{-1}z')
\end{equation}
on~$\C^n$ with respect to which the matrix~$A$ 
is normal, i.e.\ 
$$
[A,PA^*P^{-1}]=0.
$$
It is also interesting to examine the higher critical points 
of the moment map squared.  These are the solutions of the equation
\begin{equation}\label{eq:LAcrit}
[[A,A^*],A]  = \lambda A,\qquad
\lambda = \frac{\abs{[A,A^*]}^2}{\abs{A}^2}.
\end{equation}
The Kirwan--Ness inequality in Corollary~\ref{cor:KIRNESS1}
asserts that every solution~$A$ of equation~\eqref{eq:LAcrit} 
satisfies the inequality
\begin{equation}\label{eq:LAness}
\frac{\abs{[A,A^*]}^2}{\abs{A}^4}
\le \frac{\abs{[gAg^{-1},(gAg^{-1})^*]}^2}{\abs{gAg^{-1}}^4}
\end{equation}
for all~${g\in\SL(n,\C)}$.

\bigbreak

In this example every unstable~$\rG^c$-orbit contains 
a critical point of the square of the moment map. 
Examples of solutions of~\eqref{eq:LAcrit} 
are matrices of the form
\begin{equation}\label{eq:LAcrit1}
A 
:= 
A_{d_1,\dots,d_m}
:=
\left(\begin{array}{cccc}
A_{d_1} & 0 &\cdots & 0 \\
0 & \ddots & \ddots & \vdots \\
\vdots & \ddots & \ddots & 0  \\
0 & \cdots & 0 & A_{d_m}
\end{array}\right)
\end{equation}
with~${\sum_jd_j=n}$, ${\max_jd_j\ge2}$, ${A_1=0}$, and
\begin{equation}\label{eq:LAcrit2}
A_d := 
\left(\begin{array}{ccccc}
0 & a_1 & 0 & \cdots & 0 \\
0 & 0  & a_2 & \ddots & \vdots \\
\vdots & \ddots   & \ddots & \ddots & 0 \\
\vdots &  & \ddots & \ddots & a_{d-1} \\
0 & \cdots & \cdots & 0 & 0
\end{array}\right),\qquad
a_i := \sqrt{\frac{i(d-i)}{2}},
\end{equation}
for~${d\ge2}$.  This matrix satisfies
\begin{equation*}
[A_d,A_d^*] = 
\left(\begin{array}{ccccc}
\frac{d-1}{2}  & 0 & \cdots & \cdots & 0 \\
0 & \frac{d-3}{2}  & \ddots & & \vdots \\
\vdots & \ddots & \ddots & \ddots & \vdots  \\
\vdots &  & \ddots & \frac{3-d}{2}  & 0 \\
0 & \cdots & \cdots & 0 & \frac{1-d}{2} 
\end{array}\right).
\end{equation*}
Hence~${[[A_d,A_d^*],A_d]=A_d}$ 
and
$$
\Abs{A_d}^2 = \Abs{[A_d,A_d^*]}^2 = \frac{(d-1)d(d+1)}{12}.
$$
Thus~\eqref{eq:LAcrit} holds with~$\lambda=1$ and 
$$
f(A_{d_1,\dots,d_m}) = \frac{3}{2\sum_{j=1}^m(d_j-1)d_j(d_j+1)}.
$$
The matrix~${A=A_{d_1,\dots,d_m}}$ with~${m=n-1}$, ${d_1=2}$, and $d_2=\cdots=d_{n-1}=1$
satisfies~${f(A)=1/4}$ and is an absolute maximum of~$f$. 

We examine the Szekelyhidi criterion in 
Corollary~\ref{cor:SZEK0} for the matrix
\begin{equation*}
A
:=
\left(\begin{array}{cccc}
J_{d_1} & 0 &\cdots & 0 \\
0 & \ddots & \ddots & \vdots \\
\vdots & \ddots & \ddots & 0  \\
0 & \cdots & 0 & J_{d_m}
\end{array}\right),\qquad
J_d:=
\left(\begin{array}{cccc}
0 & 1 &  &  0 \\
\vdots & \ddots & \ddots & \\
\vdots &  & \ddots & 1 \\
0 & \cdots & \cdots & 0
\end{array}\right)
\in\R^{d\times d},
\end{equation*}
in Jordan normal form 
(with $\sum_jd_j=n$ and ${\max_jd_j\ge2}$).  
This matrix is $\mu$-balanced (Definition~\ref{def:balanced}).
Its isotropy subgroup contains a maximal 
torus~$\rT$ of dimension~$m$ consisting of diagonal matrices.
The Lie algebra~${\ct=\Lie(\rT)}$ is the space of all 
diagonal matrices with trace zero whose diagonal 
entries~${\xi_{11},\dots,\xi_{1d_1},\dots,\xi_{m1},\dots,\xi_{md_m}}$
are purely imaginary and satisfy 
the condition~${\i\xi_{ji}-\i\xi_{j,i+1}=\lambda}$ 
for~${j=1,\dots,m}$ and~${i=1,\dots,d_j-1}$ 
and some real number~$\lambda$.  Every matrix 
commuting with~$\ct$ is necessarily diagonal.  
So~${\cg_\rT\cap\ct^\perp}$ generates 
a torus of dimension~${n-m-1}$. 
Let~${\eta\in\cg_\rT\cap\ct^\perp\setminus\{0\}}$ 
and denote its diagonal entries by
${\eta_{11},\dots,\eta_{1d_1},\dots,\eta_{m1},\dots,\eta_{md_m}}$.
Then
$$
\sum_{i=1}^{d_j}\eta_{ji}=0\mbox{ for }j=1,\dots,m,\qquad
\sum_{j=1}^m\sum_{i=1}^{d_j}i\eta_{ji}=0.
$$ 
The weight~$w_\mu(A,\eta)$ is the maximal entry 
of the matrix~${[\i\eta,A]}$ by Lemma~\ref{le:Vweight}.
The nonzero entries are~${\i\eta_{ji}-\i\eta_{j,i+1}}$
with~${1\le j\le m}$ and~${1\le i<d_j}$, 
and at least one them is positive.  Thus~$A$ is conjugate 
to a solution of~\eqref{eq:LAcrit} by Corollary~\ref{cor:SZEK0}. 
This also follows directly from the above discussion.

One can think of a matrix~${A\in\C^{n\times n}}$
as a translation invariant Cauchy--Riemann 
operator~${\bar\p_A=\bar\p+A}$ on~$\C$, 
associated to the Hermitian connection
$$
\nabla = d+\Phi dx+\Psi dy,\qquad
\Phi := \tfrac{1}{2}\bigl(A-A^*),\qquad
\Psi := \tfrac{1}{2\i}\bigl(A+A^*),
$$
so~${A=\Phi+\i\Psi}$. Then~$\SU(n)$ is the group 
of translation invariant unitary gauge transformations 
and it acts on~$\nabla$ by 
conjugation~${u^*\nabla = u^{-1}\circ\nabla\circ u}$.
The complexified action of~$\SL(n,\C)$ is given 
by~${g^*\nabla=d+\tPhi\,dx+\tPsi\,dy}$ with
\begin{equation*}
\begin{split}
\tPhi 
&:= 
\tfrac{1}{2}\Bigl(g^{-1}\Phi g 
+ g^*\Phi(g^*)^{-1} + \i g^{-1}\Psi g 
- \i g^*\Psi(g^*)^{-1}\Bigr),\\
\tPsi 
&:= 
\tfrac{1}{2}\Bigl(g^{-1}\Psi g 
+ g^*\Psi(g^*)^{-1} - \i g^{-1}\Phi g 
+ \i g^*\Phi(g^*)^{-1}\Bigr)
\end{split}
\end{equation*}
for~${g\in\SL(n,\C)}$ and~${\Phi,\Psi\in\csu(n)}$. The curvature of~$\nabla$
is the~$2$-form
$$
F^\nabla=[\Phi,\Psi]dx\wedge dy = -\tfrac{\i}{2}[A,A^*]dx\wedge dy
$$
and hence agrees with the moment map.  Thus~$A$ is diagonalizable
if and only if~$\nabla$ is gauge equivalent to a flat connection by a translation
invariant complex\index{normal matrices|)}  
gauge transformation.
\end{exple}

\bigbreak


\begin{exple}[{\bf Control systems}]\label{ex:ABC}\rm
Consider\index{control systems} the action of ${\rG^c:=\GL(n,\C)}$ on the complex 
vector space~${V:=\C^{n\times n}\times\C^{n\times m}\times\C^{p\times n}}$ by 
$$
g\cdot(A,B,C) := \bigl(gAg^{-1},gB,Cg^{-1}\bigr)
$$
for~${g\in\rG^c}$ and~$({A,B,C)\in V}$.
The action of the compact group~${\rG:=\rU(n)}$
preserves the standard Hermitian inner product on~$V$
and the induced action on~$\bbP(V)$ is Hamiltonian 
and generated by the standard moment map
\begin{equation}\label{eq:CTmu}
\mu(A,B,C) 
= 
-\tfrac{\i}{2}
\frac{[A,A^*]+BB^*-C^*C}{\trace(A^*A+B^*B+CC^*)}
\end{equation}
for~${(A,B,C)\in V\setminus\{0\}}$.
The control system~$(A,B,C)$ 
is $\mu$-stable if and only if it is 
controllable and observable i.e.\ 
\begin{equation}\label{eq:CTcontrobs}
\sum_{i=0}^{n-1}\im A^iB =\C^n,\qquad
\bigcap_{i=0}^{n-1}\ker CA^i=\{0\}
\end{equation}
(see~\cite{HM}).  Assume~${\mu(A,B,C)=0}$
and the isotropy subgroup is discrete, i.e.\ 
\begin{equation}\label{eq:ABCbalanced}
AA^*+BB^* = A^*A+C^*C,
\end{equation}
\begin{equation}\label{eq:CTregular}
[\xi,A]=0,\quad \xi B=0,\quad C\xi=0
\qquad\implies\qquad \xi=0
\end{equation}
for~${\xi\in\cg}$.
Let~${z\in\C^n}$ and~${\lambda\in\C}$ 
such that~${Az=\lambda z}$ and~${Cz=0}$. 
Then
$$
(A-\lambda)(A^*-\ola)z+BB^*z = (A^*-\ola)(A-\lambda)z+C^*Cz=0
$$
by~\eqref{eq:ABCbalanced} and hence~${A^*z=\ola z}$ and~${B^*z=0}$.
Thus the matrix~${\xi:=\i zz^*\in\cg}$
commutes with~$A$ and satisfies~${\xi B=0}$ 
and~${C\xi=0}$, and so~$\xi=0$ by~\eqref{eq:CTregular}.  
Hence~${\ker(A-\lambda)\cap\ker C=\{0\}}$
and~${\im(A-\lambda)+\im B=\C^n}$ for all~${\lambda\in\C}$
and so~$(A,B,C)$ satisfies~\eqref{eq:CTcontrobs}.

Conversely, suppose that~$(A,B,C)$ is controllable 
and observable but not~$\mu$-stable.
Then, by the Hilbert--Mumford criterion 
in Theorem~\ref{thm:HMs}, there exists 
a~${\xi\in\cg\setminus\{0\}}$ 
such that~${w_\mu((A,B,C),\xi)\le0}$. 
Let~${\lambda_1<\cdots<\lambda_k}$ be the eigen\-values of~$\i\xi$,
let~${\C^n=E_1\oplus\cdots\oplus E_k}$ be the corresponding
eigen\-space decomposition, and let~${A_{ij}:E_j\to E_i}$,
${B_i:\C^m\to E_i}$,~${C_j:E_j\to\C^p}$ be the 
corresponding components of the matrices~$A,B,C$. 
Then it follows from Lemma~\ref{le:Vweight} 
that~${A_{ij}=0}$ for~${i>j}$ (so~$A$ is upper triangular),
that~${B_i=0}$ whenever~${\lambda_i>0}$, 
and that~${C_j=0}$ whenever~${\lambda_j<0}$.
Moreover, by controllability and observability
we have~${C_1\ne0}$ and~${B_k\ne 0}$. 
In the case~${k\ge2}$ we 
obtain~${0\le\lambda_1<\lambda_k\le0}$,
and in the case~${k=1}$ we obtain~${\lambda_1=0}$ 
and so~${\xi=0}$, a contradiction in both cases.
This shows that the control system~$(A,B,C)$ 
is $\mu$-stable if and only if it is 
controllable and observable. 

\bigbreak

In control theory the solutions of~\eqref{eq:ABCbalanced} 
are known as {\it balanced control systems}~\cite{HM}. 
The relation between balanced control systems and geometric
invariant theory was noted by Helmke and Moore.
They also showed that the gradient flow of the Kempf--Ness 
function in this setting converges to a balanced control system.
More precisely, the lifted Kempf--Ness 
function~${\Phi_{A,B,C}:\rG^c\to\R}$ is given by
\begin{equation}\label{eq:CTKN}
\Phi_{A,B,C}(g) 
= \tfrac{1}{4}\log\left(
\frac{\abs{g^{-1}Ag}^2+\abs{g^{-1}B}^2+\abs{Cg}^2}
{\abs{A}^2+\abs{B}^2+\abs{C}^2}
\right)
\end{equation}
and its gradient flow takes the form
\begin{equation}\label{eq:CTKNflow}
g^{-1}\dot g = \tfrac{1}{2} 
\frac{[g^{-1}Ag,g^*A^*(g^*)^{-1}]+g^{-1}BB^*(g^*)^{-1}-g^*C^*Cg}
{\abs{g^{-1}Ag}^2+\abs{g^{-1}B}^2+\abs{Cg}^2}
\end{equation}
Under the isometry
$
\rG^c/\rG\to\cP_n:[g]\mapsto P:=gg^*
$
in Example~\ref{ex:LA} the Kempf--Ness function is given by
\begin{equation}\label{eq:CTKNP}
\Phi_{A,B,C}(P) = \tfrac{1}{4}\log\left(
\frac{\trace(APA^*P^{-1}+BB^*P^{-1}+PC^*C)}
{\trace(A^*A+BB^*+C^*C)}
\right)
\end{equation}
for~${P\in\cP_n}$ and its negative 
gradient flow equation has the form
\begin{equation}\label{eq:CTKNPflow}
\dot P =
\frac{APA^* - PA^*P^{-1}AP + BB^* - PC^*CP}
{\trace(APA^*P^{-1}+BB^*P+P^{-1}C^*C)}.
\end{equation}
If~$(A,B,C)$ is controllable and observable,
this flow converges to a positive definite Hermitian 
matrix~${P\in\cP_n}$ by Theorem~\ref{thm:KNF},
and~$(A,B,C)$ is balanced with respect to the Hermitian
structure~\eqref{eq:Pinner} determined by~$P$, i.e.
\begin{equation}\label{eq:CTPbalanced}
APA^*P^{-1}+BB^*P^{-1} = PA^*P^{-1}A+PC^*C.
\end{equation}
The isotropy subgroup of a controllable and observable
system~$(A,B,C)$ is trivial and thus the moduli space
of conjugacy classes of such systems is a K\"ahler manifold.  
The moduli space of controllable pairs~$(A,B)$ fits
into the GIT framework with~${\rG^c=\SL(n,\C)}$ 
and the balancing equation~\eqref{eq:ABCbalanced}
replaced by
\begin{equation}\label{eq:ABbalanced}
[A,A^*]+BB^*-\frac{\trace(BB^*)}{n}\one=0.
\end{equation}
The topology of these and related moduli spaces 
was studied by Uwe Helmke and his colla\-bo\-ra\-tors in the 
1980s (see e.g.~\cite{HELMKE1,HELMKE2,HELMKE3,HH,HM}). 

In~\cite{ADHM} Atiyah--Drinfeld--Hitchin--Manin 
used the solutions of the equations
\begin{equation}\label{eq:ADHM}
[A_1,A_1^*]+[A_2,A_2^*]+BB^*-C^*C=0,\qquad
[A_1,A_2]+BC = 0,
\end{equation}
with~${B,C^*\in\C^{n\times m}}$ and~$A_1,A_2\in\C^{n\times n}$ 
to construct anti-self-dual $\SU(m)$-instantons 
of charge~$n$ on the four-sphere.
This is the {\it ADHM construcion}.
\end{exple}



\begin{exple}[{\bf Partial flag manifolds}]\label{ex:FLAG}\rm
Fix a finite sequence\index{flag manifold}
of positive integers~${n_1,\dots,n_r}$ and 
define~${n:=n_1+\cdots+n_r}$.  A {\bf partial flag of type~${(n_1,\dots,n_r)}$}
is a sequence of subspaces~${\{0\}=F_0\subset F_1\subset\cdots\subset F_r=\C^n}$
such that~${\dim^c(F_i/F_{i-1})=n_i}$ for~${i=1,\dots,r}$.  
The space~$\cF(n_1,\dots,n_r)$ of all partial flags
of type~$(n_1,\dots,n_r)$ can be identified with the quotient space
$$
\cF(n_1,\dots,n_r) \cong \rU(n)/(\rU(n_1)\times\cdots\times\rU(n_r)).
$$
If~${\lambda_1>\lambda_2>\cdots>\lambda_{r-1}>\lambda_r=0}$, 
then the centralizer of the matrix
\begin{equation}\label{eq:xilambda}
\xi
:=
\left(\begin{array}{ccccc}
\bi\lambda_1\one_{n_1} & 0 & \cdots & 0 \\
0 &\bi\lambda_2\one_{n_2} & \ddots & \vdots \\
\vdots & \ddots & \ddots & 0\\
0 & \cdots & 0 &\bi\lambda_r\one_{n_r}
\end{array}\right)
\in\cu(n)
\end{equation}
is the subgroup~${C(\xi)=\rU(n_1)\times\cdots\times\rU(n_r)}$
and so the map
$$
\rU(n)/(\rU(n_1)\times\cdots\times\rU(n_r))\to\cO(\xi):[u]\mapsto u\xi u^*
$$
gives rise to a diffeomorphism from the partial flag manifold~$\cF(n_1,\dots,n_r)$
to the (co)adjoint orbit~${\cO(\xi)=\left\{u\xi u^*\,|\,u\in\rU(n)\right\}\subset\cu(n)}$.
This orbit has the tangent spaces~${T_\eta\cO(\xi)=\im(\ad(\eta))}$ 
and is equipped with a symplectic 
form~${\om_\eta(\etahat_1,\etahat_2)=\tfrac{1}{2}\Re\,\trace(\eta[\zeta_1,\zeta_2])}$
for~${\eta\in\cO(\xi)}$ and~${\etahat_i=[\eta,\zeta_i]}$, ${\zeta_i\in\cu(n)}$. 

The partial flag manifold~$\cF(n_1,\dots,n_r)$ can be described as a symplectic 
quotient as follows.  For~${i=1,\dots,r}$ define~${k_i:=n_1+\cdots+n_i}$ and consider
the action of the group~${\rG^c:=\GL(k_1,\C)\times\cdots\times\GL(k_{r-1},\C)}$
on the complex vector space~${V:=\C^{k_2\times k_1}\times\cdots\times\C^{k_r\times k_{r-1}}}$
by
$$
g\cdot A 
:= \Bigl(
g_2A_1g_1^{-1},g_3A_2g_2^{-1},\dots,
g_{r-1}A_{r-2}g_{r-2}^{-1},A_{r-1}g_{r-1}^{-1}
\Bigr)
$$
for~${g=(g_1,\dots,g_{r-1})\in\rG^c}$ with~${g_i\in\GL(k_i,\C)}$
and~${A=(A_1,\dots,A_{r-1})\in V}$ with~${A_i\in\C^{k_{i+1}\times k_i}}$.
The space~$V$ is equipped with the symplectic form
$$
\om(A,B) = \sum_{i=1}^{r-1}\Im\,\trace(A_i^*B_i)
$$
for~$A=(A_1,\dots,A_{r-1}),\,B=(B_1,\dots,B_{r-1})\in V$ with $A_i,B_i\in\C^{k_{i+1}\times k_i}$,
and the action of~${G:=\rU(k_1)\times\cdots\times\rU(k_{r-1})}$
preserves this symplectic form.  If we identify the Lie 
algebra~${\cg:=\Lie(\rG)=\cu(k_1)\times\cdots\times\cu(k_{r-1})}$ 
with its dual space~$\cg^*$ via the inner 
product~${\inner{\xi}{\eta}=\sum_{i=1}^{r-1}\Re\,\trace(\xi_i^*\eta_i)}$,
then the action of~$\rG$ on $V$ is generated by the moment 
map~${\mu:V\to\cg}$, given by
$$
\mu(A) := \tfrac{\bi}{2}\Bigl(A_1^*A_1,
A_2^*A_2-A_1A_1^*,\dots,
A_{r-1}^*A_{r-1}-A_{r-2}A_{r-2}^*
\Bigr)
$$
for~${A=(A_1,\dots,A_{r-1})\in V}$.  

\bigbreak

Now fix a central element~${\tau:=\tfrac{\bi}{2}(\tau_1\one_{k_1},\dots,\tau_{r-1}\one_{k_{r-1}})\in\cg}$
with~${\tau_i>0}$ for~${i=1,\dots,r-1}$.  Then, with~${A_0:=0}$, we have
$$
\mu^{-1}(\tau) = \left\{A\in V\,\big|\,
A_i^*A_i-A_{i-1}A_{i-1}^*=\tau_i\one_{k_i}
\mbox{ for }i=1,\dots,r-1 \right\}.
$$
The group~$\rG=\rU(k_1)\times\cdots\times\rU(k_{r-1})$
acts freely on~$\mu^{-1}(\tau)$ and the quotient
$$
V\dslash\rG(\tau) = \mu^{-1}(\tau)/\rG
$$
is diffeomorphic to the partial flag manifold~${\cF(n_1,\dots,n_r)}$.
The diffeomorphism sends the equivalence class of an element~${A\in\mu^{-1}(\tau)}$
to the partial flag~${F=(F_i)_{i=1,\dots,r}}$ given by~${F_r=\C^n}$ and
$$
F_i := \im\bigl(A_{r-1}\circ A_{r-2}\circ\cdots\circ A_i:\C^{k_i}\to\C^n\bigr),\qquad
i=1,\dots,r-1.
$$
Moreover, if~${\xi\in\cu(n)}$ is chosen as in~\eqref{eq:xilambda} with
${\lambda_r=0}$ and
$$
\lambda_i := \tau_i+\tau_{i+1}+\cdots+\tau_{r-1},\qquad i=1,\dots,r-1,
$$
then the map
\begin{equation}\label{eq:FLAGdiffeo}
V\dslash\rG(\tau)=\mu^{-1}(\tau)/\rG\to\cO(\xi):
[A]\mapsto \bi A_{r-1}A_{r-1}^*
\end{equation}
is a diffeomorphism.  Indeed, if~${A\in\mu^{-1}(\tau)}$,
then~$\bi\lambda_i$ is an eigenvalue of the skew-adjoint 
matrix~$\bi A_{r-1}A_{r-1}^*$ with the $n_i$-dimensional
eigenspace
$$
E_i := \bigl\{A_{r-1}A_{r-2}\cdots A_iz\,\big|\,z\in\C^{k_i},\,A_{i-1}^*z=0\bigr\}.
$$
Thus the image of the map~\eqref{eq:FLAGdiffeo} is contained in~$\cO(\xi)$.
({\bf Exercise:} Prove that the map~\eqref{eq:FLAGdiffeo} is bijective
and that it is a symplectomorphism.)

The lifted Kempf--Ness function~${\Phi_A:\rG^c\to\R}$
is given by 
\begin{equation}\label{eq:KNFLAG}
\Phi_A(g) 
= 
\tfrac{1}{4}\sum_{i=1}^{r-1}\Bigl(\Abs{g_{i+1}^{-1}A_ig_i}^2 - \Abs{A_i}^2\Bigr) 
- \tfrac{1}{4}\sum_{i=1}^{r-1}\tau_i\log\bigl(\det(g_ig_i^*)\bigr)
\end{equation}
for~${A\in V}$ and~${g\in\rG^c}$.  Here we define~${g_r:=\one_n}$.
The quotient space~$\rG^c/\rG$ can be identified with the space~$\cP$
of $(r-1)$-tuples~${P=(P_1,\dots,P_{r-1})}$ of positive definite
Hermitian matrices~${P_i=P_i^*\in\C^{k_i\times k_i}}$ via~${[g_i]\mapsto P_i:=g_ig_i^*}$.
In this formulation the Kempf-Ness function~${\Phi_A:\cP\to\R}$ is given by 
$$
\Phi_A(P) 
= 
\tfrac{1}{4}\sum_{i=1}^{r-1}
\Bigl(\trace\bigl(A_iP_iA_i^*P_{i+1}^{-1}\bigr) - \trace(A_i^*A_i)
- \tau_i\log\bigl(\det(P_i)\bigr)\Bigr)
$$
and its negative gradien flow with respect 
to the metric~\eqref{eq:Pmetric} has the form
\begin{equation}\label{eq:KNFLAGflow}
\dot P_i = A_{i-1}P_{i-1}A_{i-1}^*
- P_iA_i^*P_{i+1}^{-1}A_i P_i + \tau_iP_i
\end{equation}
for~${i=1,\dots,r-1}$.  Here we use the convention~${P_r:=\one_n}$, ${P_0:=0}$, ${A_0:=0}$.
An element~$A\in V$ is $(\mu-\tau)$-stable if and only~$A_i$ has rank~$k_i$ for each~$i$.
\end{exple}

\bigbreak



\begin{exple}[{\bf Lie algebras}]\label{ex:LIEALGEBRA}\rm
This example is due to Lauret~\cite{LAURET}.
Consider the standard action of the group~$\rG^c=\SL(n,\C)$
on the complex vector space
$
V := \Lambda^2(\C^n)^*\otimes\C^n
$
of all complex bilinear~$2$-forms~${\tau:\C^n\times\C^n\to\C^n}$.
The group action is given by
$$
(g\cdot\tau)(z,z') := g\tau(g^{-1}z,g^{-1}z')
$$
for~${g\in\SL(n,\C)}$, ${\tau\in V}$, and~$z,z'\in\C^n$.
The action of the compact subgroup~${\rG=\SU(n)}$ preserves the 
standard inner product on~$V$, given by 
$$
\inner{\sigma}{\tau} = \sum_{i,j=1}^n \inner{\sigma(e_i,e_j)}{\tau(e_i,e_j)}
$$
for~${\sigma,\tau\in V}$.  Here~${\inner{z}{z'}=\Re\sum_{i=1}^n\oz_iz_i'}$ 
denotes the the standard Hermitian inner product on~$\C^n$ 
and~${e_1,\dots,e_n}$ is the standard basis of~$\C^n$.  
By Lemma~\ref{le:Vmu} the action of~$\SU(n)$ on the projective
space~$\bbP(V)$ with the symplectic form induced 
by the Hermitian structure on the sphere of radius~${r=\sqrt{2\hbar}}$ 
is generated by the moment map~${\mu:\bbP(V)\to\csu(n)^*}$ given by 
\begin{equation}\label{eq:muLIE1}
\begin{split}
\inner{\mu(\tau)}{\xi}
&= 
\frac{\hbar}{\abs{\tau}^2}\inner{\tau}{\bi\xi\cdot\tau} \\
&=
\frac{\hbar}{\abs{\tau}^2}\sum_{i,j=1}^n\inner{\tau(e_i,e_j)}{\i\xi\tau(e_i,e_j)-2\tau(e_i,\i\xi e_j)}
\end{split}
\end{equation}
for~${\tau\in V}$ and~${\xi\in\csu(n)}$.
Now suppose that~$\tau:\C^n\times\C^n\to\C^n$ is a complex Lie bracket
and denote its adjoint representation by~${\ad_\tau:\C^n\to\Der(\C^n,\tau)}$,
so that~${\ad_\tau(z)=\tau(z,\cdot)}$ for~${z\in\C^n}$.   
Then, for all~${\xi\in\csu(n)}$,
\begin{equation}\label{eq:muLIE2}
\begin{split}
\inner{\mu(\tau)}{\xi} 
&\phantom{:}= 
-\frac{\bi\hbar}{\abs{\tau}^2}\trace(M_\tau\xi),\\
M_\tau 
&:= 
\sum_{i=1}^n\Bigl(2\ad_\tau(e_i)^*\ad_\tau(e_i)-\ad_\tau(e_i)\ad_\tau(e_i)^*\Bigr).
\end{split}
\end{equation}
Thus~$\mu(\tau)=0$ if and only if $M_\tau\in\R\one$.
It was shown by Lauret in~\cite[Theorem~4.3]{LAURET} 
that the complexified group orbit~$\rG^c([\tau])\subset\bbP(V)$ 
intersects the zero set of the moment map (i.e.~$\tau$ is polystable)
if and only if the Lie algebra~$(\C^n,\tau)$ is semisimple.
His proof uses Cartan's theorem about the compact real form 
of a semisimple Lie algebra.  This theorem implies that
in the semisimple case there exists a
Hermitian inner product~$\inner{\cdot}{\cdot}'$ on~$\C^n$ 
and a Hermitian orthonormal basis~${e_1',\dots,e_n'}$ of~$\C^n$ 
such that~${\ad_\tau(e_i')^*=-\ad_\tau(e_i')}$ 
and~${\sum_{i=1}^n\ad_\tau(e_i')^2=-c\one}$ 
for some~${c>0}$ and so~${M'_\tau=c\one}$. 
After suitable rescaling the two inner products are related 
by the action of~$\SL(n,\C)$ and so every semisimple
complex Lie bracket~$\tau$ is polystable.

\bigbreak

One can modify this approach following the work of Donaldson in~\cite{DON-lie}
so that Cartan's theorem does not need to be used but instead 
a proof of Cartan's theorem as well as various other standard results
in Lie algebra theory emerge as byproducts of the proof. 
The key idea is to show directly that every simple Lie algebra
is polystable by finding a critical point of the 
lifted Kempf--Ness function~${\Phi_\tau:\SL(n,\C)\to\R}$, given by 
\begin{equation}\label{eq:KNLIE}
\Phi_\tau(g) 
:= \hbar\Bigl(\log\Abs{g^{-1}\cdot\tau}-\log\Abs{\tau}\Bigr)
= \frac{\hbar}{2}\log\left(\frac{\Abs{g^{-1}\cdot\tau}^2}{\Abs{\tau}^2}\right)
\end{equation}
for~${g\in\SL(n,\C)}$ (see Lemma~\ref{le:VKN}). 
The Kempf--Ness function can be replaced by 
the function~$f_\tau$ on the Hadamard 
manifold~$\cP_n\cong\SL(n,\C)/\SU(n)$ of all positive definite
Hermitian matrices of determinant one, defined by
\begin{equation}\label{eq:DONLIE}
f_\tau(gg^*) := \Abs{g^{-1}\cdot\tau}^2
\end{equation}
for~${g\in\SL(n,\C)}$.  The Hadamard manifold~$\cP_n$ can be 
identified with the space~$\cH_n$ of Hermitian inner products on~$\C^n$ 
with a fixed determinant.  Thus the task is to find a Hermitian inner
product on~$\C^n$ which minimizes the norm of the Lie bracket~$\tau$.
This approach was suggested by Cartan~\cite{CARTAN2} 
and carried out by Richardson~\cite{RR} under 
the assumption that the Killing form is nondegenerate.
In~\cite{DON-lie} Donaldson proved that every convex 
function~${f:\cH_n\to\R}$ that is invariant under a 
subgroup of~$\SL(n,\C)$ that acts irreducibly on~$\C^n$ 
has a critical point which is fixed by the subgroup.   
In the case at hand the relevant subgroup
is the group~${\Aut_0(\C^n,\tau)\cap\SL(n,\C)}$
of all Lie algebra automorphisms of~$(\C^n,\tau)$
in the identity component with determinant one.
It acts irreducibly on~$\C^n$ whenever~$(\C^n,\tau)$
is a simple Lie algebra.  Every critical point of~$f_\tau$ 
is then a Hermitian inner product~$h$ on~$\C^n$
for which the space of derivations is invariant 
under the involution~${A\mapsto A^*}$.
For Lie algebras with a trivial center the  existence 
of such an inner product implies that the Killing form of $(\C^n,\tau)$
is nondegenerate, the adjoint repre\-sentation $\ad_\tau:\C^n\to\Der(\C^n,\tau)$ 
is bijective, and~$\Aut_0(\C^n,\tau)$ is actually contained in~$\SL(n,\C)$.
Another byproduct of Donaldson's proof is that 
$$
\rK := \Aut_0(\C^n,\tau)\cap\SU(\C^n,h)
$$
is a maximal compact subgroup of~$\Aut_0(\C^n,\tau)$,
that every compact subgroup of~$\Aut_0(\C^n,\tau)$ 
is conjugate to a subgroup of~$\rK$,
and that~$\Aut_0(\C^n,\tau)$ is the complexification of~$\rK$.
Moreover,~$\rK$ is connected and 
the Lie algebra of~$\rK$ is Cartan's compact real form of $(\C^n,\tau)$.
Once these results have been established for simple Lie algebras,
it is a straight forward matter to deduce that the polystable points
in the present setting are precisely the semisimple Lie algebras.
For more details see~\cite{DON-lie,LAURET}
and also~\cite[\S7.6]{RS-DG}.
\end{exple}

\bigbreak




\appendix


\chapter{Nonpositive sectional curvature}\label{app:NONPOS}

Assume throughout that $M$ is a complete, connected, simply 
connected Riemannian manifold of nonpositive sectional curvature. 
Then Hadamard's theorem asserts that the exponential 
map~${\exp_p:T_pM\to M}$ is a diffeomorphism for every~${p\in M}$
(see e.g.~\cite[Theorem~6.5.2]{RS-DG}). 
The Levi-Civita connection on~$M$ is denoted by~$\nabla$
and the distance function associated to the Riemannian
metric by~${d:M\times M\to[0,\infty)}$.

\begin{lem}\label{le:gamma01}
Let $I\subset\R$ be an interval,
let $\gamma_0,\gamma_1:I\to M$ be smooth curves,
and, for each $t\in I$, denote by $[0,1]\to M:s\mapsto\gamma(s,t)$ 
the unique geodesic with the endpoints 
$$
\gamma(0,t)=\gamma_0(t),\qquad
\gamma(1,t)=\gamma_1(t).
$$ 
Define the function $\rho:I\to\R$ by 
\begin{equation}\label{eq:rho}
\rho(t) := d(\gamma_0(t),\gamma_1(t))
= \int_0^1\abs{\p_s\gamma(s,t)}\,ds.
\end{equation}
If $t\in I$ such that $\rho(t)\ne0$
then $\rho$ is differentiable at $t$ and 
\begin{equation*}
\begin{split}
\dot\rho(t)
&= 
\int_0^1\frac{\inner{\p_s\gamma}{\Nabla{t}\p_s\gamma}}
{\abs{\p_s\gamma}}\,ds  \\
&=
\frac{\inner{\p_s\gamma(1,t)}{\p_t\gamma(1,t)}
- \inner{\p_s\gamma(0,t)}{\p_t\gamma(0,t)}}{\rho(t)}.
\end{split}
\end{equation*}
\end{lem}

\begin{proof}
If $\rho(t)\ne 0$ then $\p_s\gamma(s,t)\ne0$ 
for all $s$ so $\rho$ is differentiable at $t$ and
\begin{equation*}
\begin{split}
\dot\rho(t)
&=
\int_0^1\p_t\abs{\p_s\gamma}\,ds \\
&=
\int_0^1\frac{\inner{\p_s\gamma}{\Nabla{t}\p_s\gamma}}
{\abs{\p_s\gamma}}\,ds \\
&=
\frac{1}{\rho(t)}
\int_0^1\p_s\inner{\p_s\gamma}{\p_t\gamma}\,ds.
\end{split}
\end{equation*}
The last equation follows from the fact that 
$\Nabla{s}\p_s\gamma\equiv0$ and  
${\abs{\p_s\gamma(s,t)}=\rho(t)}$.
Now Lemma~\ref{le:gamma01} follows 
from the fundamental theorem of calculus.
\end{proof}

\begin{lem}\label{le:gradflow}
Let $\Phi:M\to\R$ be a smooth function that is convex 
along geodesics, let $\gamma_0,\gamma_1:\R\to M$ 
be negative gradient flow lines of $\Phi$, and let
$\gamma$ and $\rho$ be as in Lemma~\ref{le:gamma01}. 
Then $\rho$ is nonincreasing and, if 
$$
\rho(t)\ne0,
$$
then 
\begin{equation}\label{eq:rhodot}
\dot\rho(t) = -\frac{1}{\rho(t)}\int_0^1\frac{\p^2}{\p s^2}
(\Phi\circ\gamma)(s,t)\,ds.
\end{equation}
\end{lem}

\begin{proof}
Assume $\rho(t)\ne 0$. Then, by Lemma~\ref{le:gamma01},
$\rho$ is differentiable at $t$ and
\begin{equation*}
\begin{split}
\dot\rho(t)
&=
\frac{1}{\rho(t)}
\Bigl(\inner{\p_s\gamma(1,t)}{-\nabla\Phi(\gamma(1,t))}
- \inner{\p_s\gamma(0,t)}{-\nabla\Phi(\gamma(0,t))}\Bigr) \\
&=
- \frac{1}{\rho(t)}
\Bigl(\p_s(\Phi\circ\gamma)(1,t)
- \p_s(\Phi\circ\gamma)(0,t)\Bigr).
\end{split}
\end{equation*}
This proves~\eqref{eq:rhodot}.
If $\rho(t)=0$ for some $t$ then $\rho(t)=0$ for all $t$. 
If $\rho(t)\ne 0$ for all~$t$ then $\rho$ is nonincreasing
by~\eqref{eq:rhodot}. This proves Lemma~\ref{le:gradflow}.
\end{proof}

\begin{lem}\label{le:distance}
Let $I\subset\R$ be an interval and 
$\gamma_0,\gamma_1:I\to M$ be smooth curves 
such that~$\gamma_0$ is a geodesic.  Define 
the function~${\rho:I\to[0,\infty)}$ by~\eqref{eq:rho}.
If~${\gamma_0(t)=\gamma_1(t)}$ then 
$$
\frac{d}{dt^\pm}\rho(t)
= \lim_{h\to 0\atop h>0}\frac{\rho(t\pm h)}{\pm h} 
= \pm\abs{\dot\gamma_0(t)-\dot\gamma_1(t)}.
$$
If~${\gamma_0(t)\ne\gamma_1(t)}$ 
then
$$
\ddot\rho(t) \ge -\abs{\nabla\dot\gamma_1(t)}.
$$
\end{lem}

\begin{proof}
A theorem in differential geometry (e.g.~\cite[Lemma~4.7.1]{RS-DG}) 
asserts that for every positive constant~${\alpha<1}$ 
there exists a~${\delta>0}$ such that, for all~${v_0,v_1\in T_pM}$,
$$
\abs{v_0},\abs{v_1}<\delta
\quad\implies\quad
\alpha\abs{v_0-v_1}\le d(\exp_p(v_0),\exp_p(v_1)) \le \alpha^{-1}\abs{v_0-v_1}.
$$
Now assume~${\rho(t_0)=0}$, denote~${p_0:=\gamma(t_0)}$, 
and choose a pair of smooth functions~${v_0,v_1:\R\to T_{p_0}M}$
such that
$$
\gamma_0(t_0+t) = \exp_{p_0}(v_0(t)),\qquad
\gamma_1(t_0+t) = \exp_{p_0}(v_1(t))
$$
for all $t$. Then
\begin{equation*}
\begin{split}
\frac{d}{dt^+}\rho(t_0)
&=
\lim_{h\to 0\atop h>0}\frac{\rho(t_0+h)}{h}  \\
&= 
\lim_{h\to 0\atop h>0}\frac{d(\gamma_0(t_0+h),\gamma_1(t_0+h))}{h}  \\
&=
\lim_{h\to 0\atop h>0}\frac{d(\exp_{p_0}(v_0(h)),\exp_{p_0}(v_1(h)))}{h}  \\
&=
\lim_{h\to 0\atop h>0}\frac{\abs{v_0(h)-v_1(h)}}{h}  \\
&=
\left.\frac{d}{dt^+}\right|_{t=0}\abs{v_0(t)-v_1(t)} \\
&=
\abs{\dot v_0(0)-\dot v_1(0)} \\
&=
\abs{\dot\gamma_0(t_0)-\dot\gamma_1(t_0)}.
\end{split}
\end{equation*}
An analogous argument shows that 
$\frac{d}{dt^-}\rho(t_0)=-\abs{\dot\gamma_0(t_0)-\dot\gamma_1(t_0)}$.
This proves the first assertion of Lemma~\ref{le:distance}.

To prove the second assertion, 
define~${\gamma:[0,1]\times I\to M}$ and~${\rho:I\to\R}$ 
as in Lemma~\ref{le:gamma01}. If~${\rho(t)\ne 0}$, 
then Lemma~\ref{le:gamma01} asserts that~$\rho$ 
is differentiable at~$t$ and 
\begin{eqnarray*}
\dot\rho(t)
&=&
\int_0^1\frac{\inner{\p_s\gamma}{\Nabla{t}\p_s\gamma}}
{\abs{\p_s\gamma}}\,ds \\
&=&
\frac{1}{\rho(t)}
\Bigl(\inner{\p_s\gamma(1,t)}{\p_t\gamma(1,t)}
- \inner{\p_s\gamma(0,t)}{\p_t\gamma(0,t)}\Bigr).
\end{eqnarray*}
In particular, by the Cauchy--Schwarz inequality, we have
\begin{equation}\label{eq:distance1}
\dot\rho(t)^2 \le \int_0^1\abs{\Nabla{t}\p_s\gamma}^2\,ds.
\end{equation}
Moreover,~${\frac{d}{dt}(\rho\dot\rho)=\rho\ddot\rho+\dot\rho^2}$ 
and hence
$$
\rho(t)\ddot\rho(t)+\dot\rho(t)^2
= \frac{d}{dt}
\Bigl(\inner{\p_s\gamma(1,t)}{\p_t\gamma(1,t)}
- \inner{\p_s\gamma(0,t)}{\p_t\gamma(0,t)}\Bigr) 
= I+II,
$$
where 
$$
I = \inner{\p_s\gamma(1,t)}{\Nabla{t}\p_t\gamma(1,t)}  
\ge - \rho(t)\abs{\Nabla{t}\dot\gamma_1(t)}
$$
and
\begin{equation*}
\begin{split}
II
&=
\inner{\Nabla{t}\p_s\gamma(1,t)}{\p_t\gamma(1,t)}
- \inner{\Nabla{t}\p_s\gamma(0,t)}{\p_t\gamma(0,t)} \\
&=
\frac{\p}{\p s}\frac{\abs{\p_t\gamma(1,t)}^2}{2}
- \frac{\p}{\p s}\frac{\abs{\p_t\gamma(0,t)}^2}{2} \\
&=
\int_0^1\frac{\p^2}{\p s^2}\frac{\abs{\p_t\gamma}^2}{2}\,ds \\
&=
\int_0^1\Bigl(
\abs{\Nabla{s}\p_t\gamma}^2
+ \inner{\p_t\gamma}{\Nabla{s}\Nabla{s}\p_t\gamma}
\Bigr)\,ds \\
&=
\int_0^1\Bigl(
\abs{\Nabla{t}\p_s\gamma}^2 + \inner{\p_t\gamma}
{\Nabla{s}\Nabla{t}\p_s\gamma-\Nabla{t}\Nabla{s}\p_s\gamma}
\Bigr)\,ds \\
&=
\int_0^1\Bigl(
\abs{\Nabla{t}\p_s\gamma}^2
+ \inner{\p_t\gamma}{R(\p_s\gamma,\p_t\gamma)\p_s\gamma}
\Bigr)\,ds \\
&\ge
\dot\rho(t)^2.
\end{split}
\end{equation*}
Here the last inequality follows from~\eqref{eq:distance1}
and the fact that $M$ has nonpositive sectional curvature.
This proves Lemma~\ref{le:distance}.
\end{proof}

\begin{lem}\label{le:expand}
For all $p\in M$, all $v_0,v_1\in T_pM$, and all $t\ge 1$,
$$
\abs{v_0-v_1}\le d(\exp_p(v_0),\exp_p(v_1)) 
\le \frac{d(\exp_p(tv_0),\exp_p(tv_1))}{t}.
$$
\end{lem}

\begin{proof}
Define the functions~${\gamma_0,\gamma_1:[0,\infty)\to M}$
and~${\rho:[0,\infty)\to[0,\infty)}$ by 
$$
\gamma_0(t):=\exp_p(tv_0),\qquad
\gamma_1(t):=\exp_p(tv_1),\qquad
\rho(t):=d_M(\gamma_0(t),\gamma_1(t)).
$$
By Lemma~\ref{le:distance} the function $\rho$ is convex 
and~${\rho(0)=0}$ and~${\dot\rho(0)=\abs{v_0-v_1}}$. Hence 
$$
d_M(\exp_p(tv_0),\exp_p(tv_1)) 
= \rho(t)\ge t\rho(1) 
= td_M(\exp_p(v_0),\exp_p(v_1))
$$ 
for $t\ge 1$ and 
$$
d_M(\exp_p(v_0),\exp_p(v_1))
= \rho(1) \ge \dot\rho(0) = \abs{v_0-v_1}.
$$
This proves Lemma~\ref{le:expand}.
\end{proof}

\begin{thm}[{\bf Cartan Fixed Point Theorem}]\label{thm:CARTAN}
Let $M$\index{Cartan Fixed Point Theorem} 
be a complete connected simply connected
Riemannian manifold with nonpositive sectional
curvature. Let~$\rG$ be a compact topological group
that acts on~$M$ by isometries.  Then there exists 
an element~${p\in M}$ such that~${gp=p}$ for every~${g\in\rG}$.
\end{thm}

\begin{proof}
See page~\pageref{proof:CARTAN}.
\end{proof}

The proof follows the argument given by Bill Casselmann in~\cite{BC} 
and requires the following two lemmas.  The first lemma asserts that
every manifold of nonpositive sectional curvature is a semi-hyperbolic space 
in the sense of Alexandrov.  

\begin{lem}\label{le:ALEXANDROV}
Let $M$ be a complete connected simply connected
Riemannian manifold with nonpositive sectional
curvature. Let $m\in M$ and $v\in T_mM$ and define
$$
p_0:=\exp_m(-v),\qquad p_1 := \exp_m(v).
$$ 
Then 
$$
2d(m,q)^2 + \frac{d(p_0,p_1)^2}{2} 
\le d(p_0,q)^2+d(p_1,q)^2
$$
for every $q\in M$.
\end{lem}

\begin{proof}
By Hadamard's theorem the exponential map~${\exp_m:T_mM\to M}$ 
is a diffeomorphism (see e.g.~\cite[Theorem~6.5.2]{RS-DG}). 
Hence~${d(p_0,p_1)=2\abs{v}}$. Now let~${q\in M}$.  
Then there exists a unique tangent vector~${w\in T_mM}$ such that
$$
q=\exp_m(w),\qquad d(m,q) = \abs{w}.
$$
Since the exponential map is expanding, 
by Lemma~\ref{le:expand}, we have
$$
d(p_0,q) \ge \abs{w+v},\qquad
d(p_1,q) \ge \abs{w-v}.
$$
Hence 
\begin{eqnarray*}
d(m,q)^2 
&=& 
\abs{w}^2 \\
&=& 
\frac{\abs{w+v}^2+\abs{w-v}^2}{2}-\abs{v}^2 \\
&\le& 
\frac{d(p_0,q)^2+d(p_1,q)^2}{2}-\frac{d(p_0,p_1)^2}{4}. 
\end{eqnarray*}
This proves Lemma~\ref{le:ALEXANDROV}.
\end{proof}

The next lemma is {\bf Serre's Uniqueness Theorem} 
for the {\it circumcentre} of a bounded set in 
a {\it semi-hyperbolic space}. 

\begin{lem}[{\bf Serre}]\label{le:SERRE}
Let $M$\index{Serre's Uniqueness Theorem} 
be a complete connected simply connected
Riemannian manifold with nonpositive sectional
curvature. For $p\in M$ and $r\ge 0$ denote by 
$B(p,r)\subset M$ the closed ball of radius $r$ 
centered at $p$.  Let $\Om\subset M$ be a nonempty 
bounded set and define
$$
r_\Om := \inf \left\{r>0\,|\,\mbox{there exists a }p\in M
\mbox{ such that }\Om\subset B(p,r)\right\}
$$
Then there exists a unique point $p_\Om\in M$ such that
$\Om\subset B(p_\Om,r_\Om)$.
\end{lem}

\begin{proof}
We prove existence.  Choose sequences~${r_i>r_\Om}$ 
and~${p_i\in M}$ such that
$$ 
\Om\subset B(p_i,r_i),\qquad \lim_{i\to\infty} r_i=r_\Om.
$$
Choose~${q\in\Om}$.  Then~${d(q,p_i)\le r_i}$ for every~$i$.
Since the sequence~$r_i$ is bounded and~$M$ is complete,
it follows that~$p_i$ has a convergent subsequence, 
still denoted by~$p_i$.  Its limit~${p_\Om:=\lim_{i\to\infty}p_i}$
satisfies~${\Om\subset B(p_\Om,r_\Om)}$.

We prove uniqueness. Let $p_0,p_1\in M$  
such that 
$$
\Om\subset B(p_0,r_\Om)\cap B(p_1,r_\Om).
$$ 
Since the exponential map $\exp_p:T_pM\to M$ 
is a diffeomorphism there exists a unique 
tangent vector~${v_0\in T_{p_0}M}$ 
such that~${p_1 = \exp_{p_0}(v_0)}$. 
Denote the midpoint between~$p_0$ and~$p_1$ by 
$$
m := \exp_{p_0}\left(\tfrac12 v_0\right).
$$
Then it follows from Lemma~\ref{le:ALEXANDROV} that
\begin{equation*}
\begin{split}
d(m,q)^2 
\le 
\frac{d(p_0,q)^2+d(p_1,q)^2}{2}-\frac{d(p_0,p_1)^2}{4} 
\le 
r_\Om^2 - \frac{d(p_0,p_1)^2}{4}
\end{split}
\end{equation*}
for every $q\in\Om$. Since~${\sup_{q\in\Om}d(m,q)\ge r_\Om}$
by definition of~$r_\Om$, it follows that~${d(p_0,p_1)=0}$ 
and hence~${p_0=p_1}$. This proves Lemma~\ref{le:SERRE}.
\end{proof}

\begin{proof}[Proof of Theorem~\ref{thm:CARTAN}]
\label{proof:CARTAN}
Let $q\in M$ and consider the group orbit 
$$
\Om:=\left\{gq\,|\,g\in\rG\right\}.
$$
Since $\rG$ is compact, this set is bounded.
Let $r_\Om\ge0$ and $p_\Om\in M$ be as in Lemma~\ref{le:SERRE}.
Then 
$$
\Om\subset B(p_\Om,r_\Om).
$$
Since $\rG$ acts on $M$ by isometries, this implies
$$
\Om = g\Om \subset B(gp_\Om,r_\Om)
$$
for all $g\in\rG$.  Hence it follows from
the uniqueness statement in Lemma~\ref{le:SERRE} 
that~${gp_\Om=p_\Om}$ for every~${g\in\rG}$. 
This proves Theorem~\ref{thm:CARTAN}. 
\end{proof}


\chapter{The complexified group}\label{app:Gc}



\begin{defn}
A {\bf complex Lie group}\index{Lie group!complex}
is a Lie group $\rG$ equipped with the structure of a complex manifold
such that the structure maps
$$
\rG\times\rG\to\rG:(g,h)\mapsto gh,\qquad
\rG\to\rG:g\mapsto g^{-1}
$$
are holomorphic. 
\end{defn}

The Lie algebra~${\cg:=\Lie(\rG)}$ of a complex Lie group~$\rG$ 
is equipped with a linear complex structure~${\cg\to\cg:\zeta\mapsto\bi\zeta}$ 
that is preserved by the adjoint action of~$\rG$, so the Lie bracket is complex bilinear. 
Conversely, if the Lie algebra is equipped with a linear 
complex structure~$\bi$ that is preserved by the adjoint action, 
then the formula~${g^{-1}J(g)\ghat:=\bi(g^{-1}\ghat)}$
for~${\ghat\in T_g\rG}$ defines an integrable almost complex 
structure on~$\rG$ with respect to which the structure maps are holomorphic.  
Here integrability follows from the fact that the almost 
complex structure is preserved by the torsion-free 
connection~${g^{-1}\Nabla{t}\ghat=\tfrac{d}{dt}(g^{-1}\ghat)
+ [g^{-1}\dot g,g^{-1}\ghat]}$.

\begin{thm}\label{thm:Gc1}
Let~$\rG$ be a compact Lie group
and let~$\rG^c$ be a complex Lie group
with Lie algebras $\cg:=\Lie(\rG)$ and $\cg^c=\Lie(\rG^c)$. 
Let~${\iota:\rG\to\rG^c}$ be a Lie group homomorphism.
Then the following are equivalent.

\smallskip\noindent{\bf (i)}
For every complex Lie group~$\rH$ and every Lie group
homomorphism ${\rho:\rG\to\rH}$ there exists a unique 
holomorphic homomorphism~${\rho^c:\rG^c\to\rH}$ 
such that~${\rho=\rho^c\circ \iota}$.

\smallskip\noindent{\bf (ii)}
The homomorphism~$\iota$ is injective, its image~$\iota(\rG)$ is a maximal
compact subgroup of~$\rG^c$, the quotient~$\rG^c/\iota(\rG)$
is connected, and the differential~${d\iota(\one):\cg\to\cg^c}$ 
maps~$\cg$ onto a totally real subspace of~$\cg^c$. 
\end{thm}

\begin{proof}
See pages~\pageref{proof:Gc1-1} and~\pageref{proof:Gc1-2}.
\end{proof}

A Lie group homomorphism~${\iota:\rG\to\rG^c}$ 
that satisfies the equivalent conditions of Theorem~\ref{thm:Gc1} 
is called a {\bf complexification} of $\rG$.  
By the universality property in part~(i) of Theorem~\ref{thm:Gc1}, 
the complexification~$(\rG^c,\iota)$ of a compact Lie group~$\rG$ 
is\index{Lie group!complexification} unique up to canonical 
isomorphism.\index{complexified group}
A complex Lie group is called {\bf reductive} if it is the 
complexification\index{Lie group!reductive}
of a compact Lie group.\index{reductive Lie group}

\begin{thm}\label{thm:Gc2}
Every compact Lie group admits a complexification,
unique up to canonical isomorphism.
\end{thm}

\begin{proof}
See page~\pageref{proof:Gc2}.
\end{proof}

\bigbreak

The archetypal example of a complexification is the 
inclusion of the unitary group~$\rU(n)$ into~$\GL(n,\C)$. 
Polar decomposition gives rise to a diffeomorphism 
\begin{equation}\label{eq:GLn}
\phi:\rU(n)\times\cu(n)\to\GL(n,\C),\qquad
\phi(u,\eta):=\exp(\bi\eta)u.
\end{equation}
This example extends to every Lie subgroup of $\rU(n)$.

\begin{thm}\label{thm:Gc3}
Let~${\rG\subset\rU(n)}$ be a Lie subgroup 
with Lie algebra~${\cg\subset\cu(n)}$. 
Then the set
$$
\rG^c :=\left\{\exp(\bi\eta)u\,|\,
u\in\rG,\,\eta\in\cg\right\}\subset\GL(n,\C)
$$
is a complex Lie subgroup of $\GL(n,\C)$ and the inclusion
of $\rG$ into $\rG^c$ satisfies condition~(ii) in Theorem~\ref{thm:Gc1}.
Moreover,~$\rG^c/\rG$ is diffeomorphic to~$\cg$.
\end{thm}

\begin{proof}
The proof has ten steps.

\medskip\noindent{\bf Step~1.}
{\it $\rG^c$ is a closed submanifold of $\GL(n,\C)$.}

\medskip\noindent
This follows from the fact that~\eqref{eq:GLn}
is a diffeomorphism.

\medskip\noindent{\bf Step~2.}
{\it $\one\in\rG^c$ and $T_\one\rG=\cg\oplus\bi\cg=:\cg^c$.}

\medskip\noindent
For~${\xi,\eta\in\cg}$ and~${t\in\R}$ 
define
$$
\gamma(t):=\exp(\bi t\eta)\exp(t\xi)\in\rG^c.
$$ 
Then~${\dot\gamma(0)=\xi+\bi\eta}$. 
Hence
$$
\cg^c\subset T_\one\rG^c
$$
and both spaces have the same dimension. 

\bigbreak

\medskip\noindent{\bf Step~3.}
{\it $T_g\rG^c=g\cg^c$ for every $g\in\rG^c$.}

\medskip\noindent
Both spaces have the same dimension, so it suffices
to prove that~${T_g\rG^c\subset g\cg^c}$.  
Let $\phi$ be the diffeomorphism~\eqref{eq:GLn}. 
Let~$(u,\eta)\in\rG\times\cg$ 
and
$$
g:=\phi(u,\eta)=\exp(\bi\eta)u\in\rG^c.
$$
Then, for every~${\uhat\in T_u\rG}$, we have 
${d\phi(u,\eta)(\uhat,0)=\exp(\bi\eta)u(u^{-1}\uhat)\in g\cg^c}$.
Now let~$\heta\in\cg$.  We must prove that
$d\phi(u,\eta)(0,\heta)\in g\cg^c$.
To see this, define 
$$
\gamma(s,t) := \phi(u,t(\eta+s\heta)) = \exp(\bi t(\eta+s\heta))u
$$
and
$$
\xi(s,t) := \gamma(s,t)^{-1}\p_s\gamma(s,t),\qquad
\eta(s,t) := \gamma(s,t)^{-1}\p_t\gamma(s,t)
$$
for~${s,t\in\R}$.  Then~${\eta(s,t)=u^{-1}\bi(\eta+s\heta)u}$ 
and
$$
\p_t\xi(s,t) = \p_s\eta(s,t) + [\xi(s,t),\eta(s,t)],\qquad
\xi(s,0) = 0
$$ 
for all~${s,t\in\R}$. 
Thus~${\eta(s,t)\in\cg^c}$ and hence~${\xi(s,t)\in\cg^c}$ 
for all~${s,t}$.  In particular, we have
${d\phi(u,\eta)(0,\heta)=\gamma(0,1)\xi(0,1)\in g\cg^c}$ 
and this proves Step~3. 

\medskip\noindent{\bf Step~4.}
{\it Let $g\in\GL(n,\C)$.  Then $g\in\rG^c$ if and only if 
there exists a smooth curve~${\gamma:[0,1]\to\GL(n,\C)}$ 
with~${\gamma(0)\in\rG}$, ${\gamma(1)=g}$, 
and~${\gamma(t)^{-1}\dot\gamma(t)\in\cg^c}$ 
for all~$t$.}

\medskip\noindent
If~$g=\exp(i\eta)u\in\rG^c$, 
then the curve~${\gamma(t):=\exp(\bi t\eta)u}$
satisfies the requirements of Step~4. Conversely, 
let~${\gamma:[0,1]\to\GL(n,\C)}$ be a smooth curve 
with~${\gamma(0)\in\rG}$, ${\gamma(1)=g}$,
and~${\gamma(t)^{-1}\dot\gamma(t)\in\cg^c}$ 
for all~${t\in[0,1]}$. 
Then the set~${I:=\{t\in[0,1]\,|\,\gamma(t)\in\rG^c\}}$
is nonempty because $0\in I$ and is closed because~$\rG^c$ 
is a closed subset of~$\GL(n,\C)$ by Step~1. 
To prove that it is open, 
let~$\eta(t):=\gamma(t)^{-1}\dot\gamma(t)\in\cg^c$
and define the vector field~$v_t$  on~$\C^{n\times n}$ 
by~${v_t(A):=A\eta(t)}$.  By Step~3,~$v_t$ is tangent to $\rG^c$.  
Hence every solution of the differential
equation~${\dot A(t)=A(t)\eta(t)}$ that starts in~$\rG^c$ 
remains in~$\rG^c$ on a sufficiently small time interval.  
Hence~$I$ is open. Thus~${I=[0,1]}$ and so~${g=\gamma(1)\in\rG^c}$. 

\medskip\noindent{\bf Step~5.}
{\it If $g\in\rG^c$ and $\zeta\in\cg^c$ then $g^{-1}\zeta g\in\cg^c$.}

\medskip\noindent
Choose~${\gamma:[0,1]\to\rG^c}$ as in Step~4 with~${\gamma(0)\in\rG}$
and~${\gamma(1)=g}$ and define
$$
\zeta(t) := \gamma(t)^{-1}\zeta\gamma(t),\qquad 
\zeta'(t) := \gamma(t)^{-1}\dot\gamma(t).
$$
Then~$\zeta'(t)\in\cg^c$ for all~$t$ by Step~3 and
$$
\dot\zeta(t)+[\zeta'(t),\zeta(t)]=0,\qquad 
\zeta(0) = \gamma(0)\zeta\gamma(0)^{-1}\in\cg^c.
$$
Hence $\zeta(t)\in\cg^c$ for all~$t$, 
and so~${g^{-1}\zeta g=\zeta(1)\in\cg^c}$. 

\bigbreak

\medskip\noindent{\bf Step~6.}
{\it If $g\in\rG^c$ and $\zeta\in\cg^c$ then $g\zeta g^{-1}\in\cg^c$.}

\medskip\noindent
The linear map $\zeta\mapsto g^{-1}\zeta g$ maps $\cg^c$ to itself,
by Step~5, and it is injective.  Hence the map
$\cg^c\to\cg^c:\zeta\mapsto g^{-1}\zeta g$ is bijective 
and this proves Step~6.

\medskip\noindent{\bf Step~7.}
{\it If $g_0,g_1\in\rG^c$ then $g_0g_1\in\rG^c$.}

\medskip\noindent
Choose two curves $\gamma_0,\gamma_1:[0,1]\to\rG^c$ as in Step~4 
with $\gamma_i(0))\in\rG$ and $\gamma_i(1)=g_i$.
Then the curve $\gamma:=\gamma_0\gamma_1:[0,1]\to\GL(n,\C)$
satisfies
$$
\gamma^{-1}\dot\gamma = \gamma_1^{-1}\dot\gamma_1 
+ \gamma_1^{-1}(\gamma_0^{-1}\dot\gamma_0)\gamma_1,\qquad
\gamma(0)\in\rG.
$$
By Step~5, $\gamma(t)^{-1}\dot\gamma(t)\in\cg^c$
for all $t$ and hence, by Step~4, $g_0g_1=\gamma(1)\in\rG^c$. 

\medskip\noindent{\bf Step~8.}
{\it If $g\in\rG^c$ then $g^{-1}\in\rG^c$.}

\medskip\noindent
Let $\gamma:[0,1]\to\rG^c$ be as in Step~4 
with~${\gamma(0)\in\rG}$ and~${\gamma(1)=g}$,
and define the curve~${\beta:[0,1]\to\GL(n,\C)}$ 
by
$$
\beta(t):=\gamma(t)^{-1}
$$ 
for~${0\le t\le 1}$. Then~${\beta(0)\in\rG}$ and
$$
\beta^{-1}\dot\beta 
= \gamma\frac{d}{dt}\gamma^{-1}
= - \dot\gamma\gamma^{-1} 
= \gamma(-\gamma^{-1}\dot\gamma)\gamma^{-1}.
$$
Hence~$\beta(t)^{-1}\dot\beta(t)\in\cg^c$ for all~$t$ 
by Step~6, and so~$g^{-1}=\beta(1)\in\rG^c$ by Step~4. 

\medskip\noindent{\bf Step~9.}
{\it $\rG^c$ is a complex Lie subgroup of $\GL(n,\C)$.}

\medskip\noindent
By Step~3 $\rG^c$ is a complex submanifold of $\GL(n,\C)$,
and by Steps~7 and~8 it is a subgroup of $\GL(n,\C)$ .

\medskip\noindent{\bf Step~10.}
{\it $\rG$ is a maximal compact subgroup of $\rG^c$ 
and~$\rG^c/\rG$ is diffeomorphic to~$\cg$.}

\medskip\noindent
That~$\rG^c/\rG$ is diffeomorphic to~$\cg$ follows directly from the definition.
Now let~${\rH\subset\rG^c}$ be a subgroup such that~${\rG\subsetneq\rH}$.
Choose an element $h\in\rH\setminus\rG$.  Since~${\rH\subset\rG^c}$,
there is a pair ${(u,\eta)\in\rG\times\cg}$ such that~${h=\exp(\bi\eta)u}$.
Since~${\rG\subset\rH}$ and~$\rH$ is a subgroup of~$\rG^c$
we have~${P:=\exp(\bi\eta)\in\rH}$.  The matrix~$P$ is Hermitian 
and positive definite.  Since~${h\notin\rG}$ we also have~${P\notin\rG}$.
But this implies~${\eta\ne 0}$ and so at least one eigenvalue
of~$P$ is not equal to~$1$.  Hence the sequence
$$
P^k = \exp(\bi k\eta)\in\rH,\qquad k=1,2,3,\dots
$$
has no subsequence that converges to an element 
of~$\GL(n,\C)$.  Thus~$\rH$ is not compact and this 
proves Theorem~\ref{thm:Gc3}.
\end{proof}

\begin{remk}\label{rmk:Gc}\rm
Here is a sketch of a proof of Theorem~\ref{thm:Gc3}
in the intrinsic setting. If~$\rG$ is a Lie group then the formula 
\begin{equation}\label{eq:AG}
A(g)\ghat := g^{-1}\ghat
\end{equation}
for~${g\in\rG}$ and~${\ghat\in T_g\rG}$ defines a 
flat connection~${A\in\Om^1(\rG,\cg)}$ such that

\smallskip\noindent{\bf (a)} 
for every~$g\in\rG$ the linear map~${A(g):T_g\rG\to\cg}$
is bijective,

\smallskip\noindent{\bf (b)} 
for every smooth path~${\zeta:\R\to\cg}$ and every element~${g\in\rG}$
the differential equation
$$
A(\gamma(t))\dot\gamma(t)=\zeta(t),\qquad
\gamma(0)=g, 
$$
has a solution~${\gamma:\R\to\rG}$ (on all of~$\R$),

\smallskip\noindent{\bf (c)} 
the holonomy of~$A$ is trivial, i.e.\ if~${\gamma:[0,1]\to\rG}$
is a smooth curve with~${\gamma(0)=\gamma(1)}$
then every solution~${\zeta:[0,1]\to\cg}$ of
the differential equation
$$
\dot\zeta+[A(\gamma)\dot\gamma,\zeta]=0
$$
satisfies~${\zeta(0)=\zeta(1)}$.

\smallskip\noindent
Conversely, if~$A\in\Om^1(\rG,\cg)$ is a
Lie algebra valued $1$-form on a connected 
manifold~$\rG$ that satisfies~(a), (b), and~(c), 
and~$\one$ is any element of~$\rG$,
then~$\rG$ has the unique structure of a 
Lie group with unit~$\one$ such that~$A$ 
is given by~\eqref{eq:AG}.

Let~$\rG$ be a compact Lie group with 
Lie algebra~$\cg$ and define~${\cg^c:=\cg\oplus\bi\cg}$.
Then there exists a unique flat connection~$B\in\Om^1(\cg,\cg^c)$
such that
\begin{equation}\label{eq:BGC}
[\eta,\etahat]=0 \qquad\implies\qquad B(\eta)\etahat = \bi\etahat
\end{equation}
for all~${\eta,\etahat\in\cg}$. (Define~$B(\eta)\etahat:=\zeta(1)$,
where~${\zeta:[0,1]\to\cg^c}$ is the unique solution 
of~${\dot\zeta+[\bi\eta,\zeta]=\bi\etahat}$ with~${\zeta(0)=0}$.)  
Now define
\begin{equation}\label{eq:AGC}
A(u,\eta)(\uhat,\etahat) := u^{-1}\uhat + u^{-1}(B(\eta)\etahat) u
\end{equation}
for~${u\in\rG}$, ${\uhat\in T_u\rG}$, and~${\eta,\etahat\in\cg}$.
Then~${A\in\Om^1(\rG\times\cg,\cg^c)}$ is a flat connection 
satisfying~(a), (b), (c) and hence gives rise to a 
unique group structure on the manifold
$$
\rG^c:=\rG\times\cg
$$ 
such that the map
$$
\iota:\rG\to\rG^c,
$$ 
defined by~${\iota(u):=(u,0)}$ for~${u\in\rG}$, 
is a group homomorphism. Moreover, the homorphism~${\iota:\rG\to\rG^c}$
satisfies condition~(ii) in Theorem~\ref{thm:Gc1}
and this gives rise to a proof of Theorem~\ref{thm:Gc3} 
in the intrinsic setting. 
\end{remk}

\begin{remk}\label{rmk:CARTAN}\rm
If~$\rG^c$ satisfies condition~(ii) in Theorem~\ref{thm:Gc1},
then the homogeneous space~$\rG^c/\rG$ is simply connected.
The proof requires the following steps.

\smallskip\noindent{\bf Step~1.}
{\it If~${\eta\in\cg}$ and~${\exp(\bi\eta)\in\rG}$ then~${[\xi,\eta]=0}$
for all~${\xi\in\cg}$.}

\smallskip\noindent{\bf Step~2.}
{\it If~${\eta\in\cg}$ and~${\exp(\bi\eta)\in\rG}$ then~${\eta=0}$.}

\smallskip\noindent{\bf Step~3.}
{\it $\rG^c/\rG$ is simply connected.}

\smallskip\noindent
The proof of Step~1 uses the fact that~$\rG^c/\rG$ has 
nonpositive sectional curvature (see Appendix~\ref{app:GcG}).
The proof of Step~2 uses Step~1 and the fact 
that~$\rG$ is a maximal compact subgroup of~$\rG^c$
(which by the Cartan--Iwasawa--Malcev Theorem 
in~\cite[Thm~14.1.3]{HN} implies that the intersection of~$\rG$
with the identity component~$\rG^c_0$ of~$\rG^c$
is a maximal compact subgroup of~$\rG^c_0$).
The proof of Step~3 uses Step~2 and the existence 
of a nonconstant closed geodesic in each nontrivial 
homotopy class. It follows from Step~3 and Hadamard's 
theorem (e.g.~\cite[Theorem~6.5.2]{RS-DG}) that the 
map
$$
\rG\times\cg\to\rG^c:(u,\eta)\mapsto\phi(u,\eta):=\exp(\bi\eta)u
$$
is a diffeomorphism.\index{Cartan Decomposition Theorem} 
This is the {\bf Cartan Decomposition Theorem}.
\end{remk}

\begin{proof}[Proof of Theorem~\ref{thm:Gc1} ``(ii)$\implies$(i)'']
\label{proof:Gc1-1}
If~$\rG^c$ satisfies condition~(ii) in Theorem~\ref{thm:Gc1},
then the homogeneous space~$\rG^c/\rG$ is connected
and simply connected (see Remark~\ref{rmk:CARTAN}). 
Hence the following holds.

\medskip\noindent{\bf (I)}
{\it For every~${g\in\rG^c}$ there exists a smooth 
curve~${\gamma:[0,1]\to\rG^c}$ 
that satisfies~${\gamma(0)\in\rG}$ and~${\gamma(1)=g}$.}

\medskip\noindent{\bf (II)}
{\it Any two paths~${\gamma_0,\gamma_1:[0,1]\to\rG^c}$ 
as in~(I) can be joined by a smooth 
homotopy~${\{\gamma_s\}_{0\le s\le1}}$ 
satisfying~${\gamma_s(0)\in\rG}$ 
and~${\gamma_s(1)=g}$ for all~${s\in[0,1]}$.}

\medskip\noindent
Now let~$\rH$ be a complex Lie group, 
let~${\ch:=\Lie(\rH)}$ be its Lie algebra,
let~${\rho:\rG\to\rH}$ be a Lie group homomorphism,
let~${\Phi:=d\rho(1):\cg\to\ch}$ be the induced Lie algebra 
homomorphism, and denote its complexification
by~${\Phi^c:\cg^c\to\ch}$.
We define the extended map~${\rho^c:\rG^c\to\rH}$ as follows.  
Given an element~${g\in\rG^c}$ choose~$\gamma$ as in~(I),
let~${\beta:[0,1]\to\rH}$ be the unique solution of the 
differential equation 
\begin{equation}\label{eq:rhoc}
\beta^{-1}\dot\beta = \Phi^c(\gamma^{-1}\dot\gamma),\qquad
\beta(0) = \rho(\gamma(0)), 
\end{equation}
and define 
$$
\rho^c(g) := \beta(1).
$$
We prove that $\rho^c$ is well defined, i.e.\
that~$\beta(1)$ does not depend on the choice of the path~$\gamma$. 
By~(II) any two smooth paths~$\gamma_0$ and~$\gamma_1$ 
satisfying~${\gamma_0(0),\gamma_1(0)\in\rG}$
and~${\gamma_0(1)=\gamma_1(1)=g}$ 
can be joined by a smooth 
homotopy~${[0,1]^2\to\rG^c:(s,t)\mapsto\gamma_s(t)=\gamma(s,t)}$
such that~${\gamma_s(0)\in\rG}$ and~${\gamma_s(1)=g}$ 
for all~$s$.  Define~${\beta:[0,1]^2\to\rH}$ by 
$$
\beta^{-1}\p_t\beta = \Phi^c(\gamma^{-1}\p_t\gamma),\qquad
\beta(s,0) = \rho(\gamma(s,0)).
$$
We claim that
\begin{equation}\label{eq:betas}
\beta^{-1}\p_s\beta = \Phi^c(\gamma^{-1}\p_s\gamma).
\end{equation}
To see this, abbreviate 
$$
\zeta_s:=\gamma^{-1}\p_s\gamma,\qquad
\zeta_t:=\gamma^{-1}\p_t\gamma,\qquad
\eta_s:=\beta^{-1}\p_s\beta,\qquad
\eta_t:=\beta^{-1}\p_t\beta.
$$
Then~${\eta_t=\Phi^c(\zeta_t)}$ by definition of~$\beta$ and
$$
\p_t\eta_s = \p_s\eta_t + [\eta_s,\eta_t],\qquad
\p_t\Phi^c(\zeta_s) = \p_s\Phi^c(\zeta_t) + [\Phi^c(\zeta_s),\Phi^c(\zeta_t)].
$$
Moreover, when~${t=0}$ we have
$$
d\rho(\gamma)\p_s\gamma 
= d\rho(\gamma)\gamma\zeta_s 
= \rho(\gamma)\Phi(\zeta_s)
$$
and hence
$$
\eta_s(s,0)
= 
\beta(s,0)^{-1}\p_s\beta(s,0)
= 
\Phi(\gamma(s,0)^{-1}\p_s\gamma(s,0))
= 
\Phi(\zeta_s(s,0)).
$$
Hence both functions $t\mapsto\eta_s(s,t)$ and
$t\mapsto\Phi^c(\zeta_s(s,t))$ satisfy the same 
initial value problem and so they agree. 
This proves~\eqref{eq:betas}. 
Hence
$$
\eta_s(s,1)=\Phi^c(\zeta_s(s,1))=0
$$
and therefore~${\p_s\beta(s,1)=0}$. 
This shows that $\rho^c$ is well defined. 

We prove that, for $g\in\rG^c$ and $\zeta\in\cg^c$, we have 
\begin{equation}\label{eq:rhoc1}
\Phi^c(g^{-1}\zeta g) = \rho^c(g)^{-1}\Phi^c(\zeta)\rho^c(g).
\end{equation}
Choose~$\gamma$ and~$\beta$ as in the definition of~$\rho^c(g)$, 
so that
$$
\beta^{-1}\dot\beta=\Phi^c(\gamma^{-1}\dot\gamma),
$$
and define
$$
\zeta(t) := \gamma(t)^{-1}\zeta\gamma(t),\qquad
\eta(t) := \beta(t)^{-1}\Phi^c(\zeta)\beta(t).
$$
Then~$\eta(t)$ and~$\Phi^c(\zeta(t))$ satsify the same differential equation
$$
\dot\eta(t)+[\beta(t)^{-1}\dot\beta(t),\eta(t)]=0
$$
and the same initial condition
$$
\eta(0)
= \rho(\gamma(0))^{-1}\Phi^c(\zeta)\rho(\gamma(0)) 
= \Phi^c(\zeta(0)).
$$
Hence they have the same endpoints and this proves equation~\eqref{eq:rhoc1}.

\bigbreak

We prove that~$\rho^c$ is a group homomorphism.
Let~${g_0,g_1\in\rG^c}$ and let~$\gamma_i$ and~$\beta_i$ 
be as in the definition of~$\rho^c(g_i)$ for~${i=0,1}$. 
Then~${\rho^c(\gamma_i(t))=\beta_i(t)}$ for~${0\le t\le1}$ and~${i=0,1}$.  
Define~${\gamma:=\gamma_0\gamma_1}$ and~${\beta:=\beta_0\beta_1}$.
Then, by~\eqref{eq:rhoc1}, we have
\begin{equation*}
\begin{split}
\beta^{-1}\dot\beta
&=
\beta_1^{-1}\dot\beta_1 
+ \beta_1^{-1}\beta_0^{-1}\dot\beta_0\beta_1 \\
&= 
\Phi^c(\gamma_1^{-1}\dot\gamma_1)
+ \rho^c(\gamma_1)^{-1}
\Phi^c(\gamma_0^{-1}\dot\gamma_0)
\rho^c(\gamma_1) \\
&= 
\Phi^c(\gamma_1^{-1}\dot\gamma_1
+\gamma_1^{-1}\gamma_0^{-1}\dot\gamma_0\gamma_1) \\
&= 
\Phi^c(\gamma^{-1}\dot\gamma).
\end{split}
\end{equation*}
Hence
$
\rho^c(g_0g_1)
= \beta(1)
= \beta_0(1)\beta_1(1)
= \rho^c(g_0)\rho^c(g_1)
$
and so~$\rho^c$ is a group homomorphism. 

We prove that $\rho^c$ is smooth. 
Consider the commutative diagram
$$
\xymatrix 
@C=50pt 
@R=25pt 
{&\rG^c \ar[dr]^{\rho^c}  \\
\rG\times\cg \ar[ur]^{\phi}\ar[rr] & & \rH}.
$$
Here~$\phi:\rG\times\cg\to\rG^c$ is the diffeomorphism given 
by~$\phi(u,\eta)=\exp(\bi\eta)u$ for~${u\in\rG}$ and~${\eta\in\cg}$
(see Theorem~\ref{thm:Gc3}).
In the intrinsic setting the Cartan Decomposition
Theorem asserts that this map is a diffeomorphism
under the assumption that~$\rG^c$ satisfies condition~(ii) 
in Theorem~\ref{thm:Gc1} (see Remark~\ref{rmk:CARTAN}).
The map~${\rG\times\cg\to\rH}$ is given 
by~${(u,\eta)\mapsto\exp(\bi\Phi^c(\eta))\rho(u)}$
and hence is smooth. That the differential of~$\rho^c$ 
at~$\one$ is given by~$\Phi^c$ follows also from this diagram.
This proves existence, and uniqueness is obvious.
Thus we have proved that~(ii) implies~(i) 
in Theorem~\ref{thm:Gc1}. 
\end{proof}

\begin{proof}[Proof of Theorem~\ref{thm:Gc2}]
\label{proof:Gc2}
By Theorem~\ref{thm:Gc3} 
(respectively Remark~\ref{rmk:Gc} in the intrinsic setting), 
there exists an embedding~${\iota:\rG\to\rG^c}$ into a 
complex Lie group (diffeomorphic to $\rG\times\cg$) 
that satisfies condition~(ii) in Theorem~\ref{thm:Gc1}. 
Since~(ii) implies~(i) in Theorem~\ref{thm:Gc1},
the embedding $\iota:\rG\to\rG^c$ satisfies both~(i) and~(ii) 
in Theorem~\ref{thm:Gc1} and hence is a complexification.
Moreover, any two embeddings of $\rG$ into a complex 
Lie group that satisfy~(i) in Theorem~\ref{thm:Gc1} are 
naturally isomorphic.  This proves Theorem~\ref{thm:Gc2}. 
\end{proof}

\begin{proof}[Proof of Theorem~\ref{thm:Gc1} ``(i)$\implies$(ii)'']
\label{proof:Gc1-2}
Let $\iota:\rG\to\rG^c$ be an embedding into a complex 
Lie group that satisfies~(i).  By Theorem~\ref{thm:Gc2}
there exists an embedding $\widetilde{\iota}:\rG\to\widetilde{\rG}^c$
into a complex Lie group that satisfies both~(i) and~(ii).
Since both embeddings satisfy~(i), there exists a unique
holomorphic Lie group isomorphism $\phi:\rG^c\to\widetilde{\rG}^c$
such that $\phi\circ\iota=\widetilde{\iota}$. 
Since the embedding $\widetilde\iota$ satisfies~(ii),
so does $\iota$. This proves Theorem~\ref{thm:Gc1}.
\end{proof}


\chapter{The homogeneous space $M=\rG^c/\rG$}\label{app:GcG}

Let~${\rG\subset\rU(n)}$\index{G@$\rG^c/\rG$|(} 
be a compact Lie group and denote
by~${\rG^c\subset\GL(n,\C)}$ the complexified group.  
Then the homogeneous space
$$
M := \rG^c/\rG := \left\{\pi(g)\,|\,g\in\rG^c\right\},\qquad
$$
is a connected, simply connected, complete Riemannian manifold 
with nonpositive sectional curvature.   The purpose of the present 
appendix to explain this basic fact.  Denote by~${\pi:\rG^c\to M}$ 
the canonical projection, given by~${\pi(g):=g\rG}$ for~${g\in\rG^c}$.

\begin{thm}\label{thm:GcG}
Choose an invariant inner product on $\cg$ and define 
a Riemannian metric on $M$ by 
\begin{equation}\label{eq:metric}
\inner{v_1}{v_2}_p := \inner{\eta_1}{\eta_2},\qquad
p=\pi(g),\qquad v_i=d\pi(g)g\i\eta_i,
\end{equation}
for $g\in\rG^c$ and $\eta_1,\eta_2\in\cg$.

\smallskip\noindent{\bf (i)} 
Let $g:\R\to\rG^c$ and $\eta:\R\to\cg$ be smooth functions.
Then the covariant derivative of the vector 
field~${X:=d\pi(g)g\i\eta\in\Vect(\gamma)}$ 
along the curve~${\gamma:=\pi\circ g:\R\to M}$ 
is given by 
\begin{equation}\label{eq:LC}
\nabla X = d\pi(g)g\i\Bigl(\dot\eta+[\Re(g^{-1}\dot g),\eta]\Bigr).
\end{equation}

\smallskip\noindent{\bf (ii)}
The geodesics on $M$ 
have the form
$$
\gamma(t) = \pi(g\exp(\i t\eta))
$$
for $g\in\rG^c$ and $\eta\in\cg$. 

\bigbreak

\smallskip\noindent{\bf (iii)}
The Riemann curvature tensor on $\rG^c/\rG$ 
is given by 
$$
R_p(v_1,v_1)v_3 = d\pi(g)g\i [[\eta_1,\eta_2],\eta_3],\qquad
p=\pi(g),\qquad
v_i=d\pi(g)g\i\eta_i.
$$
for $g\in\rG^c$ and $\eta_i\in\cg$. 

\smallskip\noindent{\bf (iv)}
$M$ is a complete, connected, simply connected Riemannian
manifold of non\-positive sectional curvature.
\end{thm}

\begin{proof}
The projection $\pi:\rG^c\to M$ is a principal $\rG$-bundle.
The formula
$$
A_g(\ghat) := \Re(g^{-1}\ghat)
$$
defines a connection $1$-form $A\in\Om^1(\rG^c,\cg)$.
The map 
$$
\rG^c\times\cg\to TM:(g,\eta)\mapsto d\pi(g)g\i\eta.
$$
descends to a vector bundle isomorphism from
the associated bundle $\rG^c\times_\ad\cg$ 
to the tangent bundle of $M$.  Thus $A$ induces a 
connection on $TM$ and this connection is
given by~\eqref{eq:LC}.  Whenever the action of $\rG$
on a vector space preserves the inner product so does
the induced connection.  Hence~\eqref{eq:LC} is a Riemannian
connection on $TM$. We prove that it is torsion free.
Denote by $s$ and $t$ the standard coordinates on $\R^2$.
Choose a smooth function $g:\R^2\to\rG^c$ and denote
$\gamma:=\pi\circ g$. Then
\begin{equation*}
\begin{split}
\Nabla{s}\p_t\gamma
&= 
d\pi(g)g\i\Bigl(
\p_s\Im(g^{-1}\p_tg)+[\Re(g^{-1}\p_sg),\Im(g^{-1}\p_tg)]
\Bigr) \\
&= 
d\pi(g)g\i\Bigl(
\p_t\Im(g^{-1}\p_sg)+[\Re(g^{-1}\p_tg),\Im(g^{-1}\p_sg)]
\Bigr) \\
&=
\Nabla{t}\p_s\gamma
\end{split}
\end{equation*}
Here the second equation follows from the identity
\begin{equation}\label{eq:dsdtg}
\p_s(g^{-1}\p_tg)-\p_t(g^{-1}\p_sg)+[g^{-1}\p_sg,g^{-1}\p_tg]=0.
\end{equation}
This proves~(i).  

We prove part~(ii).
A smooth curve~${\gamma(t)=\pi(g(t))}$ 
is a geodesic if and only if~${\nabla\dot\gamma\equiv0}$. 
By~(i) this is equivalent to the differential equation
$$
\p_t\Im(g^{-1}\dot g) + [\Re(g^{-1}\dot g),\Im(g^{-1}\dot g)] = 0.
$$
A  smooth function~${g:\R\to\rG^c}$ satisfies this equation 
if and only if it has the form~${g(t)=g_0\exp(\i t\eta)u(t)}$ 
for some $g_0\in\rG^c$, $\eta\in\cg$, and $u:\R\to\rG$. 
This proves~(ii).

We prove part~(iii).
Choose maps $g:\R^2\to\rG^c$ and $\eta:\R^2\to\cg$
and define 
\begin{equation}\label{eq:dsdtzeta}
\zeta_s:=g^{-1}\p_sg,\qquad \zeta_t:=g^{-1}\p_tg,\qquad
\p_s\zeta_t-\p_t\zeta_s+[\zeta_s,\zeta_t]=0.
\end{equation}
Here the third equation follows from~\eqref{eq:dsdtg}.  Now define
$$
\gamma :=\pi\circ g,\qquad
Z_s :=\p_s\gamma = d\pi(g)g\zeta_s,\qquad
Z_t :=\p_t\gamma = d\pi(g)g\zeta_s,
$$
and
$$
Y:=d\pi(g)g\i\eta. 
$$
Then, by part~(i), we have
\begin{equation*}
\begin{split}
\Nabla{s}Y &= d\pi(g)g\i\Bigl(\p_s\eta+[\Re(\zeta_s),\eta]\Bigr),\\
\Nabla{t}Y &= d\pi(g)g\i\Bigl(\p_t\eta+[\Re(\zeta_t),\eta]\Bigr).
\end{split}
\end{equation*}
Hence~${R(Z_s,Z_t)Y 
= \Nabla{s}\Nabla{t}Y-\Nabla{t}\Nabla{s}Y 
= d\pi(g)g\i\teta}$, where
\begin{equation*}
\begin{split}
\teta
&=
\p_s\Bigl(\p_t\eta+[\Re(\zeta_t),\eta]\Bigr)
\,+\left[\Re(\zeta_s),\Bigl(\p_t\eta+[\Re(\zeta_t),\eta]\Bigr)\right] \\
&\quad
-\, \p_t\Bigl(\p_s\eta+[\Re(\zeta_s),\eta]\Bigr)
-\left[\Re(\zeta_t),\Bigl(\p_s\eta+[\Re(\zeta_s),\eta]\Bigr)\right] \\
&=
[\Re(\p_s\zeta_t),\eta]
+\left[\Re(\zeta_s),[\Re(\zeta_t),\eta]\right] \\
&\quad
-\,[\Re(\p_t\zeta_s),\eta]
-\left[\Re(\zeta_t),[\Re(\zeta_s),\eta]\right] \\
&= 
\bigl[\Re(\p_s\zeta_t)-\Re(\p_t\zeta_s)+[\Re(\zeta_s),\Re(\zeta_t)],\eta\bigr] \\
&= 
[[\Im(\zeta_s),\Im(\zeta_t)],\eta].
\end{split}
\end{equation*}
Here the last equality follows from~\eqref{eq:dsdtzeta}.
This proves~(iii).  By~(iii), we have
$$
\inner{R(Z_s,Z_t)Z_t}{Z_s}
= 
- \Abs{[\Im(\zeta_s),\Im(\zeta_t)]}^2
\le
0.
$$
This proves Theorem~\ref{thm:GcG}.
\end{proof}

\begin{lem}\label{le:estimate}
If $\xi_0,\xi_1,\eta\in\cg$ satisfy
$$
\exp(-\i\xi_1)\exp(\i\xi_0)\exp(\i\eta) \in\rG,
$$
then $\abs{\xi_0-\xi_1}\le\abs{\eta}$.
\end{lem}

\begin{proof}
Define~${g_0:=\exp(\i\xi_0)}$ and~${g_1:=\exp(\i\xi_1)}$.
Then the unique geodesic in~${M=\rG^c/\rG}$ 
connecting~$\pi(g_0)$ to~$\pi(g_1)$ is given 
by~${\gamma(t):=\pi(g_0\exp(\i t\eta))}$ for~${0\le t\le 1}$. 
This implies~${d(\pi(g_0),\pi(g_1))=\abs{\eta}}$.
Thus by Lemma~\ref{le:expand} 
with~${M=\rG^c/\rG}$, ${p=\pi(\one)}$, 
${v_0=d\pi(\one)\i\xi_0}$, ${v_1=d\pi(\one)\i\xi_1}$, 
${\exp_p(v_0)=\pi(g_0)}$, and ${\exp_p(v_1)=\pi(g_1)}$ 
we have~${\abs{\xi_0-\xi_1}\le\abs{\eta}}$.
This proves Lemma~\ref{le:estimate}.
\end{proof}

\begin{lem}\label{le:KG}
Every compact subgroup of $\rG^c$ is conjugate 
to a subgroup of~$\rG$.
\end{lem}

\begin{proof}
Let $\rK\subset\rG^c$ be a compact subgroup.
Then $\rK$ acts on $\rG^c/\rG$ by isometries 
via~${k\cdot\pi(g):=\pi(kg)}$ for $k\in\rK$ and $g\in\rG^c$. 
By Theorem~\ref{thm:CARTAN} the action of $\rK$ 
on~${\rG^c/\rG}$ has a fixed point $\pi(g)\in\rG^c/\rG$. 
Hence $\pi(kg)=\pi(g)$ and hence $g^{-1}kg\in\rG$ 
for every~${k\in\rK}$. This proves Lemma~\ref{le:KG}.
\end{proof}

\begin{lem}\label{le:Lambda}
Let $\zeta\in\cg^c$.  Then the following are equivalent.

\smallskip\noindent{\bf (i)}
$\zeta$ is semi-simple and has imaginary eigenvalues.

\smallskip\noindent{\bf (ii)}
There is a $g\in\rG^c$ such that $g^{-1}\zeta g\in\cg$.
\end{lem}

\begin{proof}
Assume~(i).  Then the set 
$\rT:=\overline{\{\exp(t\zeta)\,|\,t\in\R\}}\subset\rG^c$
is a (compact) torus.  By Lemma~\ref{le:KG} 
there exists an element $g\in\rG^c$ such that~${g^{-1}\rT g\subset\rG}$.
This implies~${g^{-1}\zeta g=\frac{d}{dt}|_{t=o}g^{-1}\exp(t\zeta)g\in\cg}$.
That~(ii) implies~(i) is obvious.\index{G@$\rG^c/\rG$|)} 
This proves Lemma~\ref{le:Lambda}.
\end{proof}


\chapter{Toral generators}\label{app:Lambda}

This appendix introduces toral generators 
and Mumford's equivalence relation.
Let~${\rG\subset\rU(n)}$ be a compact Lie group
with the complexification 
$$
\rG^c\subset\GL(n,\C)
$$
and denote their Lie algebras by~${\cg:=\Lie(\rG)\subset\cu(n)}$ 
and~${\cg^c:=\cg+\i\cg}$.

\begin{defn}\label{def:toral}
A nonzero\index{toral generator} 
element $\zeta\in\cg^c$ is called a {\bf toral generator}
if it is semi-simple and has purely imaginary eigenvalues. 
This means that the subset
$$
T_\zeta:=\overline{\left\{\exp(t\zeta)\,|\,t\in\R\right\}}
$$
is a torus in~$\rG^c$.  By Lemma~\ref{le:Lambda} 
the set of toral generators is
$$
\sT^c := \ad(\rG^c)(\cg\setminus\{0\}).
$$
Throughout we use the notation
\begin{equation}\label{eq:Lambda}
\begin{split}
\Lambda &:= \left\{\xi\in\cg\setminus\{0\}\,|\,\exp(\xi)=1\right\},\\
\Lambda^c &:= \left\{\zeta\in\cg^c\setminus\{0\}\,|\,\exp(\zeta)=1\right\}.
\end{split}
\end{equation}
Thus
$$
\Lambda\subset\Lambda^c\subset\sT^c.
$$
The elements of~$\Lambda^c$ are in one-to-one correspondence 
with nontrivial one-parameter subgroups~${\C^*\to\rG^c}$.
The set~${\Lambda\cup\{0\}}$ intersects the Lie 
algebra~${\ct\subset\cg}$ of any maximal 
torus~${\rT\subset\rG}$ in a spanning lattice and every 
element of~$\Lambda^c$ is conjugate to an element 
of~${\Lambda\cap\ct}$ (see Lemma~\ref{le:Lambda}).
\end{defn}

\begin{lem}[{\bf Parabolic Subgroups}]\label{le:parabolic1}
For $\zeta\in\sT^c$\index{parabolic subgroup} 
the set
\begin{equation}\label{eq:P}
\rP(\zeta) := \left\{p\in\rG^c\,|\,
\mbox{the limit }\lim_{t\to\infty}
\exp(\i t\zeta)p\exp(-\i t\zeta)\mbox{ exists in }\rG^c\right\}
\end{equation}
is a Lie subgroup of $\rG^c$ with Lie algebra
\begin{equation}\label{eq:p}
\cp(\zeta) := \left\{\rho\in\cg^c\,|\,
\mbox{the limit }
\lim_{t\to\infty}\exp(\i t\zeta)\rho\exp(-\i t\zeta)
\mbox{ exists in }\cg^c\right\}.
\end{equation}
\end{lem}

\begin{proof}
Let $\zeta\in\sT^c$.  Then the matrix $\i\zeta\in\C^{n\times n}$ 
is semi-simple and has real eigenvalues, 
denoted by~${\lambda_1<\lambda_2<\cdots<\lambda_k}$.
Denote  the  eigenspace of~$\lambda_j$ by~$V_j$ so
we have an eigenspace decomposition 
$$
\C^n=V_1\oplus V_2\oplus\cdots\oplus V_k.
$$
Write a matrix $\rho\in\cg^c\subset\cgl(n,\C)$ in the form
$$
\rho = \left(\begin{array}{cccc}
\rho_{11} & \rho_{12} & \cdots & \rho_{1k} \\
\rho_{21} & \rho_{22} & \cdots & \rho_{2k} \\
\vdots   &   \vdots   & \ddots & \vdots \\
\rho_{k1} & \rho_{k2} & \cdots & \rho_{kk}
\end{array}\right),\qquad
\rho_{ij}\in\Hom(V_j,V_i).
$$
Then
$$
\exp(\i t\zeta)\rho\exp(-\i t\zeta) 
= \left(\begin{array}{cccc}
\rho_{11} & e^{(\lambda_1-\lambda_2)t}\rho_{12} 
& \cdots & e^{(\lambda_1-\lambda_k)t}\rho_{1k} \\
 e^{(\lambda_2-\lambda_1)t}\rho_{21} & \rho_{22} 
& \cdots & e^{(\lambda_2-\lambda_k)t}\rho_{2k} \\
  \vdots   &   \vdots   & \ddots & \vdots \\
e^{(\lambda_k-\lambda_1)t}\rho_{k1} 
& e^{(\lambda_k-\lambda_2)t}\rho_{k2} 
& \cdots & \rho_{kk}
\end{array}\right).
$$
Thus $\rho\in\cp(\zeta)$ if and only if $\rho\in\cg^c$
and $\rho_{ij}=0$ for $i>j$.
Likewise, $g\in\rP(\zeta)$ if and only if $g\in\rG^c$
and $g_{ij}=0$ for $i>j$.
Hence $\rP(\zeta)$ is a closed subset of~$\rG^c$.
Since every closed subgroup of a Lie group 
is a Lie subgroup, this proves Lemma~\ref{le:parabolic1}.
\end{proof}

The proof of Lemma~\ref{le:parabolic1} 
shows that $\rP(\zeta)$ is what is called in the theory 
of algebraic groups a {\bf parabolic subgroup}\index{Lie group!parabolic subgroup}
of $\rG^c$ (upper triangular matrices). It also shows that, 
for $\zeta=\xi\in\cg$, its intersection with $\rG$ is the centralizer
$$
\rP(\xi)\cap\rG = \rC(\xi) := \left\{u\in\rG\,|\,u\xi u^{-1}=\xi\right\}.
$$
In this case there is a $\rG$-equivariant 
isomorphism  
\begin{equation}\label{eq:BOREL}
\rG^c/\rP(\xi)\cong\rG/\rC(\xi).
\end{equation}
For a generic element $\xi\in\cg$, 
the identity component $\rB\subset\rP(\xi)$ 
is a {\bf Borel subgroup}\index{Lie group!Borel subgroup}
(a maximal parabolic subgroup)
of~$\rG^c$, the intersection
$$
\rT:=\rB\cap\rG
$$ 
is a maximal torus and is the identity component of the 
centralizer~$C(\xi)$, and equation~\eqref{eq:BOREL} reads 
$$
\rG^c/\rB\cong\rG/\rT.
$$
In general, equation~\eqref{eq:BOREL} can be restated as follows.

\begin{thm}\label{thm:BOREL}
For every~${\xi\in\cg\setminus\{0\}}$ and every~${g\in\rG^c}$
there exists an element~${p\in\rP(\xi)}$ 
such that~${p^{-1}g\in\rG}$.
\end{thm}

\begin{proof}
In Appendix~\ref{app:BOREL} we give a proof 
of Theorem~\ref{thm:BOREL} which does not rely 
on the structure theory of Lie groups.
\end{proof}

\begin{thm}[{\bf Mumford}]\label{thm:MUMFORD}
Define\index{Mumford equivalence relation} 
a relation on $\sT^c$ by
\begin{equation}\label{eq:equiv}
\zeta\sim\zeta'\qquad\stackrel{\text\small{\rm{def}}}{\iff}\qquad
\exists\;p\in\rP(\zeta)\mbox{ such that }p\zeta p^{-1}=\zeta'.
\end{equation}
The formula~\eqref{eq:equiv} defines an equivalence relation on $\sT^c$, 
invariant under conjugation, and every equivalence class contains 
a unique element of $\cg$.
\end{thm}

\begin{proof}
The group $\rG^c$ acts on the space 
$$
\Gamma^c:=\bigsqcup_{\zeta\in\sT^c}\rP(\zeta)
$$
by the diagonal adjoint action.  This determines a groupoid.
The elements of~$\sT^c$ are the objects of the groupoid.
A pair~${(\zeta,p)\in\Gamma^c}$ is a morphism from~$\zeta$
to~$p\zeta p^{-1}$.  The inverse map is given 
by~${(\zeta,p)\mapsto(p\zeta p^{-1},p^{-1})}$
and the composition map sends a composable 
pair of pairs consisting of~${(\zeta,p)}$ and~${(\zeta',p')}$
with~${\zeta'=p\zeta p^{-1}}$ to the pair~${(\zeta,p'p)}$ with
$$
p'p\in\rP(p\zeta p^{-1})p=p\rP(\zeta)=\rP(\zeta).
$$ 
This shows that~\eqref{eq:equiv} is an equivalence relation.

We prove that the equivalence relation~\eqref{eq:equiv}
is invariant under conjugation.  Choose equivalent 
elements~${\zeta,\zeta'\in\sT^c}$ and let~${g\in\rG^c}$.  Then there 
exists an element~${p\in\rP(\zeta)}$ such that~${p\zeta p^{-1}=\zeta'}$.
Hence~${gpg^{-1}\in\rP(g\zeta g^{-1})}$ and
$$
(gpg^{-1})(g\zeta g^{-1})(gpg^{-1})^{-1}=g\zeta'g^{-1},
$$
so~${g\zeta g^{-1}}$ is equivalent to~${g\zeta'g^{-1}}$. 

\bigbreak

Now let $\zeta\in\sT^c$. By Lemma~\ref{le:Lambda}, 
there is a $g\in\rG^c$ such that 
$$
\xi:=g\zeta g^{-1}\in\cg. 
$$
By Theorem~\ref{thm:BOREL}, there is a 
$q\in\rP(\xi)$ such that $u:=q^{-1}g\in\rG$.  Hence
$$
p := u^{-1}g = g^{-1}qg\in\rP(g^{-1}\xi g)=\rP(\zeta)
$$
and
$$
p\zeta p^{-1}=u^{-1}g\zeta g^{-1}u= u^{-1}\xi u\in\cg.
$$
This shows that every equivalence class
in $\sT^c$ contains an element of $\cg$.

We prove uniqueness.
Let $\xi\in\cg\setminus\{0\}$ and $p\in\rP(\xi)$ such that 
${p\xi p^{-1}\in\cg}$. Choose the the eigenvalues 
$\lambda_1<\cdots<\lambda_k$ of $\i\xi$ and the eigenspace 
decomposition $\C^n=V_1\oplus\cdots\oplus V_k$ 
as in the proof of Lemma~\ref{le:parabolic1}. 
Since $\i\xi$ is Hermitian its eigenspaces $V_j$ 
are pairwise orthogonal. Moreover, the subspace 
$V_1\oplus\cdots\oplus V_j$ is invariant under $p$ 
and hence also under $p^{-1}\xi p$ for every~$j$.  
Since $p\xi p^{-1}\in\cg$ is a skew-Hermitian endomorphism 
of $\C^n$ and the complex subspaces $V_1,\dots,V_k$ of $\C^n$ 
are pairwise orthogonal, it follows that
$$
p\xi p^{-1} V_i\subset V_i,\qquad 
i=1,\dots,k.
$$
Hence $p\xi p^{-1}=\xi$. This completes 
the proof of Theorem~\ref{thm:MUMFORD}.
\end{proof}


\chapter{The partial flag manifold $\rG^c/\rP\equiv\rG/\rC$}\label{app:BOREL}

In this appendix we prove\index{flag manifold}
Theorem~\ref{thm:BOREL}.

\begin{lem}\label{le:parabolic2}
Let $N\in\N$. There exist real numbers 
$$
\beta_0(N),\beta_1(N),\dots,\beta_{2N-1}(N)
$$ 
such that $\beta_\nu(N)=0$ when $\nu$ is even 
and, for $k=1,3,5,\dots,4N-1$, 
\begin{equation}\label{eq:betanu}
\sum_{\nu=0}^{2N-1}
\beta_\nu(N) \exp\left({\frac{k\nu\pi\i}{2N}}\right)
=\left\{\begin{array}{rl}
\i,&\mbox{if }0<k<2N,\\
-\i,&\mbox{if }2N<k<4N.
\end{array}\right.
\end{equation}
\end{lem}

\begin{proof}
Define $\lambda:=\exp(\frac{\pi\i}{2N})$ and consider 
the Vandermonde matrix 
$$
\Lambda
:=
\left(\begin{array}{ccccc}
\lambda & \lambda^3 & \lambda^5 
& \dots &\lambda^{2N-1} \\
\lambda^3 & \lambda^9 & \lambda^{15} 
& \dots &\lambda^{6N-3} \\
\lambda^5 & \lambda^{15} & \lambda^{25} 
& \dots &\lambda^{10N-5} \\
\vdots & \vdots & \vdots 
& \ddots & \vdots \\
\lambda^{2N-1} & \lambda^{6N-3} & \lambda^{10N-5} 
& \cdots &\lambda^{(2N-1)^2}
\end{array}\right)
\in\C^{N\times N}.
$$
Its complex determinant is 
$$
\mathrm{det}^c(\Lambda) 
= \lambda^{N(2N-1)}
\prod_{0\le i<j\le N-1}
\left(\lambda^{4j}-\lambda^{4i}\right).
$$
Since $\lambda$ is a primitive $4N$th root of unity,
the numbers $\lambda^{4i}$, $i=0,\dots,N-1$, are 
pairwise distinct.  Hence $\Lambda$ is nonsingular.

Since~$\Lambda$ is nonsingular, there exists a unique 
vector 
$$
z=(z_1,z_3,\dots,z_{2N-1})\in\C^N
$$
such that
\begin{equation}\label{eq:znu1}
\sum_{0<\nu<2N\atop \nu\;odd}
\exp\left({\frac{k\nu\pi\i}{2N}}\right)z_\nu
= \i,\qquad
k= 1,3,\dots,2N-1.
\end{equation}
The numbers $\overline{z}_\nu$ also satisfy 
equation~\eqref{eq:znu1}, because
$$
\sum_{0<\nu<2N\atop \nu\;odd}
\exp\left(\frac{k\nu\pi\i}{2N}\right)\overline{z}_\nu
=
-\sum_{0<\nu<2N\atop \nu\;odd}
\overline{\exp\left(\frac{(2N-k)\nu\pi\i}{2N}\right)z_\nu}
= 
\i
$$
for $k=1,3,\dots,2N-1$. 
Since the solution is unique the $z_\nu$ are real.

Now define the numbers~${\beta_0(N),\beta_1(N),\dots,\beta_{2N-1}(N)}$ by
$$
\beta_\nu(N)
:= \left\{\begin{array}{ll}
z_\nu,&\mbox{for }\nu=1,3,\dots,2N-1,\\
0,&\mbox{for }\nu\mbox{ even}.
\end{array}\right.
$$
These numbers satisfy~\eqref{eq:betanu} 
for~${k=1,3,\dots,2N-1}$ by~\eqref{eq:znu1}. 
Moreover, we have~${\exp(k\pi\i)=-1}$ 
for every odd integer~$k$,
equation~\eqref{eq:betanu} also holds 
for~${k=2N+1,2N+3,\dots,4N-1}$.
This proves Lemma~\ref{le:parabolic2}.
\end{proof}


\begin{lem}\label{le:parabolic3}
Let $m\in\N$ and $N:=2^m$.  
There exist real numbers 
$$
\alpha_0(N),\alpha_1(N),\dots,\alpha_{2N-1}(N)
$$ 
such that $\alpha_\nu(N)=0$ when $\nu$ is even 
and, for every $k\in\{0,1,\dots,2N-1\}$, 
\begin{equation}\label{eq:alphanu}
\sum_{\nu=0}^{2N-1}\alpha_\nu(N) \exp\left({\frac{k\nu\pi\i}{N}}\right)
=\left\{\begin{array}{rl}
\i,&\mbox{if }1\le k\le N-1,\\
-\i,&\mbox{if }N+1\le k\le2N-1,\\
0,&\mbox{if }k=0\mbox{ or }k=N.
\end{array}\right.
\end{equation}
\end{lem}

\begin{proof}
The proof is by induction on $m$. For $m=1$ and $N=2^m=2$ 
choose 
$$
\alpha_1(2):=\tfrac{1}{2},\qquad
\alpha_3(2):=-\tfrac{1}{2}.  
$$
Then 
$$
\sum_{\nu=0}^{3}\alpha_\nu(2) \exp\left({\frac{k\nu\pi\i}{2}}\right)
=
\frac{\i^k-(-\i)^k}{2}  
=
\left\{\begin{array}{rl}
\i,&\mbox{for }k=1,\\
-\i,&\mbox{for }k=3,\\
0,&\mbox{for }k=0,2.
\end{array}\right.
$$
Now let $m\in\N$ and define $N:=2^m$.  Assume, by induction, 
that the numbers $\alpha_\nu(N)$, $\nu=0,1,\dots,2N-1$, 
have been found such that~\eqref{eq:alphanu} holds for 
$k=0,1,\dots,2N-1$.  Let $\beta_\nu(N)$, $\nu=0,1,\dots,2N-1$, 
be the constants of Lemma~\ref{le:parabolic2}.  
Define 
$$
\alpha_{2N+\nu}(N) := \alpha_\nu(N),\qquad
\beta_{2N+\nu}(N) := -\beta_\nu(N),
$$
for $\nu=0,1,2,\dots,2N-1$ and 
\begin{equation}\label{eq:alphanu3}
\alpha_\nu(2N) 
:= \frac{\alpha_\nu(N)  + \beta_\nu(N)}{2},\qquad
\nu=0,1,2,\dots,4N-1.
\end{equation}
Then 
$$
\sum_{\nu=0}^{4N-1}
\alpha_\nu(2N)\exp\left(\frac{k\nu\pi\i}{2N}\right) 
= A_k + B_k,
$$
where
\begin{equation*}
\begin{split}
A_k 
&:=
\frac{1}{2}\sum_{\nu=0}^{4N-1}
\alpha_\nu(N)\exp\left(\frac{k\nu\pi\i}{2N}\right), \\
B_k
&:=
\frac{1}{2}
\sum_{\nu=0}^{4N-1}\beta_\nu(N)
\exp\left(\frac{k\nu\pi\i}{2N}\right).
\end{split}
\end{equation*}
Since $\alpha_{2N+\nu}(N)=\alpha_\nu(N)$, we have 
\begin{equation*}
\begin{split}
A_k
&=
\frac{1}{2}\sum_{\nu=0}^{4N-1}
\alpha_\nu(N)
\exp\left(\frac{k\nu\pi\i}{2N}\right) \\
&=
\frac{1+\exp(k\pi\i)}{2}\sum_{\nu=0}^{2N-1}
\alpha_\nu(N)\exp\left(\frac{k\nu\pi\i}{2N}\right) \\
&=
\frac{1+(-1)^k}{2}\sum_{\nu=0}^{2N-1}
\alpha_\nu(N)\exp\left(\frac{k\nu\pi\i}{2N}\right).
\end{split}
\end{equation*}
If $k$ is odd the right hand side vanishes.
If $k$ is even, then by the induction hypothesis,
\begin{equation*}
\begin{split}
A_k
&=
\sum_{\nu=0}^{2N-1}
\alpha_\nu(N)\exp\left(\frac{(k/2)\nu\pi\i}{N}\right) \\
&=
\left\{\begin{array}{rl}
\i,&\mbox{for } k=2,4,\dots,2N-2,\\
-\i,&\mbox{for } k=2N+2,\dots,4N-2,\\
0,&\mbox{for } k =0,2N.
\end{array}\right.
\end{split}
\end{equation*}
Since $\beta_{2N+\nu}(N)=-\beta_\nu(N)$, we have
\begin{equation*}
\begin{split}
B_k
&=
\frac{1}{2}\sum_{\nu=0}^{2N-1}
\beta_\nu(N)
\left(
\exp\left(\frac{k\nu\pi\i}{2N}\right) 
- \exp\left(\frac{k(2N+\nu)\pi\i}{2N}\right)
\right) \\
&=
\frac{1-\exp(k\pi\i)}{2}
\sum_{\nu=0}^{2N-1}
\beta_\nu(N)
\exp\left(\frac{k\nu\pi\i}{2N}\right) \\
&=
\frac{1-(-1)^k}{2}
\sum_{\nu=0}^{2N-1}
\beta_\nu(N)
\exp\left(\frac{k\nu\pi\i}{2N}\right).
\end{split}
\end{equation*}
If $k$ is even the right hand side vanishes.
If $k$ is odd, then by Lemma~\ref{le:parabolic2},
$$
B_k
=
\sum_{\nu=0}^{2N-1}
\beta_\nu(N)\exp\left(\frac{k\nu\pi\i}{2N}\right) 
= \left\{\begin{array}{rl}
\i,&\mbox{if } k=1,3,\dots,2N-1,\\
-\i,&\mbox{if } k =2N+1,\dots,4N-1.
\end{array}\right.
$$
Combining the formulas for $A_k$ and $B_k$ we find
$$
A_k+B_k =
\left\{\begin{array}{rl}
\i,&\mbox{for } k=1,2,3,\dots,2N-1,\\
-\i,&\mbox{for } k=2N+1,2N+2,\dots,4N-1,\\
0,&\mbox{for } k=0,2N,
\end{array}\right.
$$
and this proves Lemma~\ref{le:parabolic3}.
\end{proof}

\begin{lem}\label{le:parabolic4}
Let $\xi,\eta\in\cg\subset\csu(n)$ and assume $\exp(\xi)=\one$.
Then there exists an element $\zeta\in\cp(\xi)$ 
such that $\zeta-\i\eta\in\cg$.
\end{lem}

\begin{proof}
Let~${\lambda_1<\cdots<\lambda_k}$ 
be the eigenvalues of the Hermitian matrix $\i\xi$ and 
denote the corresponding eigenspace decomposition 
by $\C^n=V_1\oplus\cdots\oplus V_k$. Then 
$$
\lambda_i-\lambda_j = 2\pi m_{ij},\qquad m_{ij}\in\Z,
$$
with $m_{ij}>0$ for $i>j$ and $m_{ij}<0$ for $i<j$.  
Choose $m\in\N$ such that
$$
N:=2^m > m_{k1} 
= \frac{\lambda_k-\lambda_1}{2\pi}.
$$
Choose $\alpha_0,\dots,\alpha_{2N-1}\in\R$ 
as in Lemma~\ref{le:parabolic3}.  Let
$$
\eta 
= 
\left(\begin{array}{cccc}
\eta_{11} & \eta_{12} & \cdots & \eta_{1k} \\
\eta_{21} & \eta_{22} & \cdots & \eta_{2k} \\
\vdots   &   \vdots   & \ddots & \vdots \\
\eta_{k1} & \eta_{k2} & \cdots & \eta_{kk}
\end{array}\right)
\in\cg,\qquad
\eta_{ij}\in\Hom(V_j,V_i).
$$
Define 
$$
\zeta := 
\i\eta - \sum_{\nu=0}^{2N-1}\alpha_\nu
\exp\left(-\frac{\nu}{2N}\xi\right)\eta
\exp\left(\frac{\nu}{2N}\xi\right)
\in\cg^c.
$$
Then, for $i>j$, we have
\begin{equation*}
\begin{split}
\zeta_{ij} 
&=
\i\eta_{ij} - \sum_{\nu=0}^{2N-1}
\alpha_\nu \exp\left(\frac{\nu}{2N}\i(\lambda_i-\lambda_j)\right)\eta_{ij} \\
&=
\left(\i - \sum_{\nu=0}^{2N-1}
\alpha_\nu \exp\left(\frac{m_{ij}\nu\pi\i}{N}\right)\right)\eta_{ij} \\
&=
0.
\end{split}
\end{equation*}
Here the last equality uses Lemma~\ref{le:parabolic3}
and the fact that ${1\le m_{ij}\le N-1}$ for $i>j$.
Since $\zeta_{ij}=0$ for $i>j$ it follows from the 
proof of Lemma~\ref{le:parabolic1} that~${\zeta\in\cp(\xi)}$. 
Moreover, by construction we have~${\bi\eta - \zeta\in\cg}$.
This proves Lemma~\ref{le:parabolic4}.
\end{proof}

\begin{proof}[Proof of Theorem~\ref{thm:BOREL}]
Assume first that $\xi\in\Lambda$ so that $\exp(\xi)=\one$.
Define
$$
A := \left\{g\in\rG^c\,|\,
\exists\; p\in\rP(\xi)\mbox{ such that }
p^{-1}g\in\rG\right\}.
$$
We will prove by an open and closed argument that $A=\rG^c$.

We prove that $A$ is a closed subset of $\rG^c$ (in the relative topology).
Let $g_i\in A$ be a sequence which converges to an element $g\in\rG^c$.  
Then there exists a sequence~${p_i\in\rP(\xi)}$ such that 
$u_i := p_i^{-1}g_i\in\rG$.  Since $\rG$ is compact there exists 
a subsequence (still denoted by $u_i$) which converges to an 
element $u\in\rG$. Since $\rP(\xi)$ is a closed subset 
of $\rG^c$, we have
$$
p := gu^{-1} = \lim_{i\to\infty}g_iu_i^{-1} = \lim_{i\to\infty}p_i\in\rP(\xi).
$$
Hence $p^{-1}g=u\in\rG$ and so $g\in A$.  
Thus $A$ is a closed subset of~$\rG^c$.

We prove that the function $f:\rP(\xi)\times\rG\to\rG^c$, defined by 
$$
f(p,u) := pu
$$
for $p\in\rP(\xi)$ and $u\in\rG$, is a submersion.
Let $p\in\rP(\xi)$ and $u\in\rG$ and denote
$$
g:=f(p,u)=pu.
$$
Let $\ghat\in T_g\rG^c$ and denote 
\begin{equation}\label{eq:ghatzeta}
\tzeta := p^{-1}\ghat u^{-1} 
= u(g^{-1}\ghat)u^{-1} \in\cg^c.
\end{equation}
Let $\eta\in\cg$ be the imginary part of $\tzeta$ 
so that $\tzeta-\i\eta\in\cg$. By Lemma~\ref{le:parabolic4}, 
there exists an element $\zeta\in\cp(\xi)$ such that
$\zeta - \i\eta\in\cg$ and hence
$\tzeta - \zeta \in\cg$. Define
$$
\phat := p\zeta,\qquad 
\uhat := \left(\tzeta-\zeta\right)u.
$$
Then $\phat\in T_p\rP$, $\uhat\in T_u\rG$, 
and 
$$
df(p,u)(\phat,\uhat)
= \phat u + p\uhat \\
= p\tzeta u \\
= \ghat.
$$
Here the last equation follows from~\eqref{eq:ghatzeta}.
Thus we have proved that the differential 
$df(p,u):T_p\rP\times T_u\rG\to T_{pu}\rG^c$ 
is surjective for every $p\in\rP(\xi)$ and every $u\in\rG$.
Hence $f$ is a submersion as claimed.

We prove that $A=\rG^c$.
The set $A$ contains $\rG$ by definition.
Moreover, we have proved that it is closed and that
it is the image of a submersion and hence is open.
Since $\rG^c$ is homeomorphic to $\rG\times\cg$
and $A$ contains $\rG\cong\rG\times\{0\}$, it follows that
$A$ intersects each connected component of $\rG^c$ 
in a nonempty open and closed set. Hence $A=\rG^c$. 
This proves Theorem~\ref{thm:BOREL} for $\xi\in\Lambda$.

Now let $\xi\in\cg\setminus\{0\}$.  Choose sequences $\xi_i\in\Lambda$ 
and $s_i\in\R$ such that $s_i\xi_i$ converges to $\xi$.
By the first part of the proof there exist 
sequences $p_i\in\rP(\xi_i)$ and $u_i\in\rG$ 
such that $u_ip_i = g$
for every $i$.  Passing to a subsequence if necessary 
we may assume that $u_i$ converges to $u\in\rG$. 
Hence $p_i=u_i^{-1}g$ converges to $p=u^{-1}g\in\rG^c$.
Examining the eigenspace decompositions of $\xi$ 
and $\xi_i$ we find that $p\in\rP(\xi)$.
This proves Theorem~\ref{thm:BOREL}.
\end{proof}


\backmatter



\newpage
\addcontentsline{toc}{chapter}{Index}

\input{momentweight-book.ind}



\end{document}